%% file: ms.tex
\documentclass{memo-l}

\usepackage{amssymb}
\usepackage{appendix}
\usepackage{enumitem}
\usepackage[latin1]{inputenc} 
\usepackage{hyperref}
\usepackage{setspace}

\newtheorem{theorem}{Theorem}[chapter]
\newtheorem{lemma}[theorem]{Lemma}
\newtheorem{corollary}[theorem]{Corollary}
\newtheorem{proposition}[theorem]{Proposition}

\newtheorem{notation}[theorem]{Notation}
\newtheorem{question}[theorem]{Question}
\newtheorem{conjecture}[theorem]{Conjecture}
\newtheorem{problem}[theorem]{Problem}

\newtheorem{claim}[theorem]{Claim}
\newtheorem{fact}[theorem]{Fact}

\theoremstyle{definition}
\newtheorem{example}[theorem]{Example}
\newtheorem{definition}[theorem]{Definition}

\numberwithin{section}{chapter}
\numberwithin{equation}{chapter}

\newcommand\N{\mathbb{N}}
\newcommand\F{\mathbb{F}}

\newcommand\e{\epsilon}

\newcommand\age{\text{age}}
\newcommand\bip{\mathrm{Bip}}
\newcommand\trip{\mathrm{Trip}}

\newcommand{\abar}{\bar{a}}
\newcommand{\bbar}{\bar{b}}
\newcommand{\cbar}{\bar{c}}

\newcommand{\ebar}{\bar{e}}
\newcommand{\fbar}{\bar{f}}

\newcommand{\ubar}{\bar{u}}
\newcommand{\vbar}{\bar{v}}
\newcommand{\wbar}{\bar{w}}
\newcommand{\xbar}{\bar{x}}
\newcommand{\ybar}{\bar{y}}
\newcommand{\zbar}{\bar{z}}

\newcommand{\Ubar}{\overline{U}}
\newcommand{\Hbar}{\overline{H}}
\newcommand{\HP}{\textnormal{HP}}
\newcommand{\HOP}{\mathrm{HOP}}
\newcommand{\FOP}{\mathrm{FOP}}

\newcommand{\VC}{\textnormal{VC}}

\newcommand{\GS}{\textnormal{GS}}

\newcommand{\disc}{\textnormal{disc}}

\newcommand{\cycle}{\textnormal{cycle}}
\newcommand{\oct}{\textnormal{oct}}
\newcommand{\dev}{\textnormal{dev}}
\newcommand{\vdisc}{\textnormal{vdisc}}
\newcommand{\msd}{\textnormal{msd}}

\newcommand{\calL}{\mathcal{L}}
\newcommand{\calF}{\mathcal{F}}
\newcommand{\calH}{\mathcal{H}}
\newcommand{\calC}{\mathcal{C}}
\newcommand{\calM}{\mathcal{M}}

\newcommand{\calA}{\mathcal{A}}
\newcommand{\calB}{\mathcal{B}}

\newcommand{\calP}{\mathcal{P}}

\newcommand{\calG}{\mathcal{G}}

\newcommand{\calI}{\mathcal{I}}
\newcommand{\calJ}{\mathcal{J}}
\newcommand{\calQ}{\mathcal{Q}}
\newcommand{\calR}{\mathcal{R}}
\newcommand{\calS}{\mathcal{S}}

\newcommand{\calU}{\mathcal{U}}
\newcommand{\calV}{\mathcal{V}}
\newcommand{\calW}{\mathcal{W}}
\newcommand{\calX}{\mathcal{X}}

\newcommand{\rk}{{\mathrm{rk}}}
\newcommand{\Mix}{{\mathrm{Mix}}}
\newcommand{\Irr}{{\mathrm{Irr}}}

\def\Forb{\operatorname{Forb}}

\def\HOP{\operatorname{HOP}}
\def\NHOP{\operatorname{NHOP}}
\def\FOP{\operatorname{FOP}}
\def\NFOP{\operatorname{NFOP}}
\def\IP{\operatorname{IP}}
\def\NIP{\operatorname{NIP}}
\def\SNIP{\operatorname{WNIP}}
\def\triads{\operatorname{Triads}}
\def\graph{\operatorname{Graph}}
\def\SNIP{\operatorname{SNIP}}
\def\SVC{\operatorname{SVC}}

\newenvironment{proofof}[1]{\indent{\scshape Proof of #1}:~~}{\qed}

\makeindex

\begin{document}

\frontmatter

\title{Irregular Triads in $3$-Uniform Hypergraphs}


\author{C. Terry}
\address{Department of Mathematics, Statistics, and Computer Science, University of Illinois Chicago, Chicago IL 60607, USA}
\curraddr{}
\email{caterry@uic.edu}
\thanks{}

\author{J. Wolf}
\address{Department of Pure Mathematics and Mathematical Statistics, Centre for Mathematical Sciences, Wilberforce Road, Cambridge CB3 0WB, UK}
\curraddr{}
\email{julia.wolf@dpmms.cam.ac.uk}
\thanks{}

\date{11 May 2022}




\begin{abstract}
Over the past several years, numerous authors have explored model theoretically motivated combinatorial conditions that ensure that a graph has an efficient regular decomposition in the sense of Szemer\'edi. In this paper we set out a research program that explores a corresponding set of questions for 3-uniform hypergraphs, a setting in which useful notions of regularity are significantly more intricate. 

The main results in this paper concern certain combinatorial properties which arose as natural higher-order generalizations of the order property in parallel work of the authors in the arithmetic setting. Interpreted in the context of 3-uniform hypergraphs, these are tightly connected to the nature of irregular triads. Specifically, we show that a hereditary property of 3-uniform hypergraphs admits regular decompositions with so-called ``linear error" if and only if it does not have the \emph{functional order property} ($\FOP_2$). 

Along the way, we show that a hereditary property of $3$-uniform hypergraphs is homogeneous (i.e. all regular triads have density near $0$ or near $1$) if and only it has bounded $\VC_2$-dimension. This provides a quantitative version of a recent result of Chernikov and Towsner. 

We also address several questions arising from prior work on tame regularity in hypergraphs. In particular, we characterize the hereditary properties of $3$-uniform hypergraphs admitting the type of regular partitions appearing in work of Fox et al. as those that have bounded slicewise $\VC$-dimension.  This is again analogous to a recent non-quantitative result of Chernikov and Towsner.
\end{abstract}

\maketitle

\tableofcontents


\mainmatter
\include{chapter1.tex}
\include{chapter2.tex}
\include{chapter3.tex}
\include{chapter4.tex}
\include{chapter5.tex}
\include{chapter6.tex}
\include{chapter7.tex}

\appendix
\include{appendices.tex}

\backmatter

\providecommand{\bysame}{\leavevmode\hbox to3em{\hrulefill}\thinspace}

\printindex

\end{document}

%% file: chapter1.tex
\chapter{Introduction}

Szemer\'{e}di's regularity lemma \cite{Szemeredi.1975} is a theorem with a wide range of applications in graph theory.  Roughly speaking, it says that any large graph can be partitioned into a small number of pieces, so that in between most pairs of pieces the graph behaves like a random graph.  More precisely, Szemer\'{e}di's regularity lemma states that for all $\e>0$ there is some $M=M(\e)$ so that for every sufficiently large graph $G=(V,E)$, there is some $m\leq M$ and an equipartition of $V=V_1\cup \ldots \cup V_m$ such that all but at most $\e m^2$ pairs $(V_i,V_j)$ are \emph{$\e$-regular} (see Definition \ref{def:quasi}).  In order to utilize such a regular decomposition, one needs a so-called \emph{counting lemma}, which says that the number of triangles appearing in a triple of $\e$-regular pairs is about what one would expect in a $3$-partite random graph of the same density.  An important early application of the counting lemma was the so-called \emph{triangle removal lemma} \cite{Frankl.1986}, which says that a graph with a small number of triangles can be made triangle free by removing a small number of edges.  We refer the reader to the excellent surveys \cite{Komlos.1996, Conlon.2013} for more details.

Despite its myriad applications, the regularity lemma is limited by the potential existence of irregular pairs, as well as the enormous size of the bound $M(\e)$, which appears as a tower-type function in Szemer\'{e}di's proof.  Investigation into these limitations showed them to be necessary.  Indeed, the folklore example of the half-graph showed the necessity of the irregular pairs, while a construction of Gowers \cite{Gowers.1997} proved the necessity of tower-type bounds.

On the other hand, under certain assumptions, strenghtened versions of the regularity lemma can be proved.  For instance, several authors \cite{Alon.2007, Lovasz.2010} showed that graphs of bounded $\VC$-dimension allow for regular decompositions with very good bounds.  These results yield a structural dichotomy which is best stated in terms of hereditary graph properties (see Theorem \ref{thm:vc} below).  We require a few definitions to state this result.  Throughout, we will use the following notation.  Given $x\neq y$, we will write $xy$ as shorthand for the two-element set $\{x,y\}$.  Similarly, for three distinct elements $x,y,z$, we write $xyz$ as shorthand for the three-element set $\{x,y,z\}$.   For sets $V_1,V_2$, we let $K_2[V_1,V_2]=\{xy: x\in V_1, y\in V_2, x\neq y\}$.

A \emph{hereditary graph property} is a class $\calG$ of finite graphs which is closed under isomorphism and induced subgraphs.  The results in this paper will be interesting mainly in the case of \emph{non-trivial} hereditary  properties, meaning those containing arbitrarily large graphs.  The $\VC$-dimension of a graph $G=(V,E)$ is the largest $k$ for which there exist $a_1,\ldots, a_k\in V$, and $b_S\in V$ for each $S\subseteq [k]$, such that $a_ib_S\in E$ if and only if $i\in S$.   We then say that a hereditary graph property $\calG$ \emph{has $\VC$-dimension at most $k$} if every $G\in \calG$ has $\VC(G)\leq k$, and $\calG$ has \emph{bounded $\VC$-dimension} (or is NIP\footnote{In model theory, a witness of unbounded $\VC$-dimension is known as the ``independence property", explaining the acronym.}) if it has $\VC$-dimension at most $k$ for some $k$.  Otherwise, we say $\calG$ \emph{has unbounded $\VC$-dimension} (or IP). Theorem \ref{thm:vc} sums up the work of Alon, Fischer and Newman \cite{Alon.2007} (part (1)) and Gowers \cite{Gowers.1997} (part (2)), where the bounds quoted are improvements due to Fox  et al. \cite{Fox.2017bfo} (part (1)) and Fox and Lov\'{a}sz \cite{Fox.2014} (part (2)).  The existence of this dichotomy was first pointed out by Alon, Fox, and Zhao in \cite{Alon.2018is}.

\begin{theorem}\label{thm:vc}
Suppose $\calG$ is a hereditary graph property.  Then one of the following holds.
\begin{enumerate}[label=\normalfont(\arabic*)]
\item $\calG$ has bounded $\VC$-dimension.  In this case, if $k=\max\{\VC(G):G\in \calG\}$, then for all $\e>0$, any $G\in \calG$ has a regular partition with $O(\e^{-k})$ parts, and where moreover, all regular pairs have density in $[0,\e)\cup (1-\e,1]$.
\item $\calG$ has unbounded $\VC$-dimension. In this case, there exists $\e>0$ and arbitrarily large $G\in \calG$, such that any $\e$-regular partition of $G$ has at least $\mathrm{Tw}(2+\e^{-2}/16)$ parts.   
\end{enumerate}
\end{theorem}

A similar dichotomy exists with regards to irregular pairs.  A graph $G=(V,E)$ is said to have the \emph{$k$-order property} if there exist $a_1,\ldots, a_k\in V$ and $b_1,\ldots, b_k\in V$ so that $a_ib_j\in E$ if and only if $i\leq j$.  If $G$ does not have the $k$-order property, then we say $G$ is \emph{$k$-stable}.  A hereditary graph property $\calG$ is called $k$-stable if every $G\in \calG$ is $k$-stable, and  $\calG$ is \emph{stable} if it is $k$-stable for some $k$.  Otherwise, $\calG$ is called \emph{unstable}.  

While it was known since at least the 1990s \cite{Komlos.1996} that a sufficiently large instance of the order property  forces the existence of irregular pairs (part (2) below), it was not until 2011 that Malliaris and Shelah \cite{Malliaris.2014} showed that the converse (part (1) below) also holds.  The proof of this direction is highly non-trivial, and uses finitized versions of model theoretic tools which were originally developed to study \emph{stable theories}. 

\begin{theorem}\label{thm:stable}
Suppose $\calG$ is a hereditary graph property.  Then one of the following holds.
\begin{enumerate}[label=\normalfont(\arabic*)]
\item $\calG$ is stable.  In this case, the following holds, with $k$ minimal such that $\calG$ is $k$-stable.  For all  $\e>0$, every sufficiently large $G\in \calG$ has an equipartition into $O(\e^{-2^{2^k}})$ parts, such that each pair is $\e$-regular with density in $[0,\e)\cup (1-\e,1]$.
\item $\calG$ is unstable. In this case, for all $\e>0$ and $M$ there are arbitrarily large $G\in \calG$ such that every partition of $G$ into at most $M$ parts contains $\Omega(M)$ irregular pairs.
\end{enumerate}
\end{theorem}

Model theoretic stability theory was first developed by Shelah \cite{Shelah.1990o5n} to answer an old open question about counting non-isomorphic models of a given first-order theory.  Stable theories and the tools which have resulted from their study are now foundational parts of modern model theory. The paper of Malliaris and Shelah \cite{Malliaris.2014} was the first to show that these ideas have deep connections to finite combinatorics, and has been extremely influential, setting off an entire area of related research which examines connections between infinitary model theory and finite combinatorics. For example, subsequent to \cite{Malliaris.2014} appearing on the arXiv, several papers have explored how various model theoretic notions yield different kinds of ``tame" regularity lemmas \cite{Chernikov.20156sq, Chernikov.2016zb, Chernikov.2020,Ackerman.2017}.  The model theory community also became aware of related theorems (in the realm of bounded $\VC$-dimension) proved earlier by researchers outside of model theory (e.g. \cite{Alon.2007, Lovasz.2010, Alon.20058t, Basu.2009,Fox.2017bfo}).  As a result of this work, tame regularity lemmas are now a widely appreciated point of connection between model theory and combinatorics.  Work in this area has also demonstrated explicit connections between  regularity lemmas and model theoretic tools, such as compact domination, indiscernible sequences, and definability of types (see \cite{Malliaris.2014, Pillay.2016, Pillay.2020}).  

The work described above shows that in the setting of graphs, looking for dichotomies in the behavior of hereditary properties with respect to regular partitions leads directly to notions of primary importance in model theory.  Specifically, the work of Malliaris and Shelah showed that dichotomies related to irregular pairs lead directly to stability theory, a core part of modern model theory.  The aim of this paper is to investigate analogous results in the setting of $3$-uniform hypergraphs.

Analogues of Szemer\'{e}di's regularity lemma for hypergraphs were considered desirable for many years, with the earliest examples appearing in the late 1980s \cite{Haviland.1989} and early 1990s \cite{Chung.1990}.   This early version of hypergraph regularity involved partitioning the vertex set of a hypergraph so that most triples of vertex parts satisfy a certain quasirandomness condition with respect to the hyperedges.  This type of regularity is now often referred to as \emph{weak hypergraph regularity}. 

In the early 2000s, more sophisticated notions of hypergraph regularity were conceived, which went hand in hand with powerful counting lemmas and increasingly general removal lemmas  \cite{Frankl.2002, Nagle.2003, Rodl.2005, Gowers.20063gk,Gowers.2007, Rodl.2004, Nagle.2006}.  For a detailed history, we refer the reader to \cite{Nagle.2013}.  As in the case of graphs, these counting lemmas have found robust applications, including the first fully finitary proof of the multidimensional Szemer\'{e}di's theorem \cite{Gowers.2007} and other, more classical applications in extremal combinatorics \cite{Balogh.2012}.  We will refer to this type of regularity as \emph{strong hypergraph regularity}.  In the setting of a $3$-uniform hypergraph $H=(V,E)$, a strong regular decomposition consists of a partition $\calP_{vert}=\{V_1,\ldots, V_t\}$ of the vertex set $V$, along with, for each $ij\in {[t]\choose 2}$, a partition of the pairs $K_2[V_i,V_j]=\bigcup_{\alpha\leq \ell}P_{ij}^{\alpha}$, so that certain quasirandomness properties hold.  The essential ``unit'' in such a regular decomposition is a so-called \emph{triad}, which is a $3$-partite graph of the form $(V_i\cup V_j\cup V_k, P_{ij}^{\alpha}\cup P_{ik}^{\beta}\cup P_{jk}^{\gamma})$.  The strong regularity lemma then roughly says that every large $3$-uniform hypergraph has an equitable decomposition so that almost all triples of vertices are in a regular triad (for details, see Section \ref{subsec:quasi1}).

Several results mentioned above  \cite{Ackerman.2017, Chernikov.2016zb,Fox.2017bfo} give improved \emph{weak} regular partitions for hypergraphs under various tameness assumptions.  The only previous treatment of strong hypergraph regularity and model theoretic tameness occurs in work of Chernkiov and Towsner \cite{Chernikov.2020}, where higher arity analogues of Theorem \ref{thm:vc} are proved in the context of strong regularity, using infinitary techniques.  The results in \cite{Chernikov.2020} are in terms of a higher arity analogue of $\VC$-dimension called $\VC_k$-dimension, first defined by Shelah \cite{Shelah1, Shelah2}, and subsequently developed in \cite{Chernikov.2019}, with applications to algebra appearing in \cite{ChernikovHempel, Hempel.2014}, and to combinatorics in \cite{Terry.2018}.

Despite the extensive work on model theoretic regularity lemmas of the last ten years, the following basic problem remained open.

\begin{problem}\label{prob:main}
Prove an analogue of Theorem \ref{thm:stable} in the realm of strong regularity for hypergraphs. 
\end{problem}

This is an interesting problem from a model theoretic perspective, as it is asking one to develop a higher arity analogue of stability, a problem which has been discussed in the model theory community for several years (see for example \cite{Shelah1,Takeuchi.2018}). 

\begin{question}\label{prob:main2}
What higher arity notion generalizes stability, in the same way that $\VC_k$-dimension generalizes $\VC$-dimension?
\end{question}

A convincing answer to Question \ref{prob:main2} has also remained open.   Given the important applications of strong hypergraph regularity in combinatorics, as well as its connections to arithmetic structure \cite{Gowers.2007}, one might reasonably expect an answer to Problem \ref{prob:main} to lead to the ``right'' answer to Question \ref{prob:main2}, as well as important new tools and ideas in model theory.

The focus of this paper is Problem \ref{prob:main} and Question \ref{prob:main2} in the setting of $3$-uniform hypergraphs, where our contributions are two-fold.  First, we prove the first analogue of Theorem \ref{thm:stable} in the realm of strong hypergraph regularity in terms of a new higher arity generalization of stability (see Theorem \ref{thm:fopintro} below).  Our work gives substantial evidence that this particular analogue of Theorem \ref{thm:stable} is robust, and the generalization of stability we define is the ``right'' answer to Question \ref{prob:main2}.  

Our second contribution is more complicated.  We show that there are several distinct ways to formulate generalizations of Theorem \ref{thm:stable} in the setting of strong regularity for $3$-uniform hypergraphs, suggesting that Problem \ref{prob:main} does not have a single answer.  We explore these distinct generalizations, related notions of higher arity stability, and specific examples.  The sum of these results is the first significant contribution to understanding the full landscape of Problem \ref{prob:main}.  Our work teases apart several potential generalizations of stability via their connections to Problem \ref{prob:main}, and the examples we give rule out certain natural ideas for how to generalize stability.  Indeed, working out these subtle distinctions and uncovering the relevant examples were significant hurdles to correctly addressing Question \ref{prob:main2}, and we view this as a major contribution of this paper.  

We will now give a more detailed overview of our results (for full details, see Section \ref{subsec:results}).  We begin by taking a slightly different perspective on Theorem \ref{thm:stable}.  A naive attempt to generalize Theorem \ref{thm:stable} to the setting of $3$-uniform hypergraphs would be to insist that there are \emph{no} irregular triads.  However, we will show that doing so reduces to a question about \emph{weak} regularity, rather than strong regularity (see Theorem \ref{thm:vdisceq}).  Therefore, we start by giving a more flexible formulation of Theorem \ref{thm:stable}, which will illuminate the proper generalization to higher arities.  In particular, the following is an exercise which gives an equivalent formulation of part (1) of Theorem \ref{thm:stable}.

\begin{fact}\label{fact:folk}
The following are equivalent for a hereditary graph property $\calG$.
\begin{enumerate}[label=\normalfont(\roman*)]
\item  For all $\e>0$, there is $M$ so that for all $G=(V,E)\in \calG$ with $|V|\geq M$, there is some $t\leq M$ and an equipartition $V=V_0\cup V_1\cup \ldots \cup V_t$ such that $|V_0|\leq \e |V|$ and for all $1\leq i\neq j\leq t$, $(V_i,V_j)$ is  $\e$-regular.
\item For all $\e>0$, there is $M$ so that for all $G=(V,E)\in \calG$ with $|V|\geq M$, there is some $t\leq M$ and an equipartition $V=V_1\cup \ldots \cup V_t$ such that  for all $ij\in {[t]\choose 2}$,  $(V_i,V_j)$ is  $\e$-regular.
\end{enumerate}
\end{fact}

That (ii) implies (i) is trivial.  For the reverse direction, one chooses the parameters carefully in (i), then evenly distributes the error set $V_0$ among the other parts in the partition to obtain a partition as in (ii).  Fact \ref{fact:folk} shows that having no irregular pairs in a regular partition of a graph can be reinterpreted as saying that there is a single, small set in the vertex partition which is used in every irregular pair.  Thinking of sets of vertices as ``lower order'' objects in a graph, this is informally saying one can use a lower order set to control the irregular pairs.  It turns out that this simple reinterpretation is much more conducive to higher arity generalization. We thus reinterpret Problem \ref{prob:main} as follows.

\begin{problem}\label{prob:restatement}
Characterize hereditary properties of $3$-uniform hypergraphs admitting regular partitions where the irregular triads are constrained by ``lower order" sets.  
\end{problem}

Problem \ref{prob:main} is made complicated by the fact that there are several distinct ways to constrain irregular triads using lower order sets.  We define three natural ways of constraining the error triads in a regular decomposition, which we will refer to as zero, binary, and linear $\disc_{2,3}$-error.  The use of the term ``linear'' stems form the authors' work in the arithmetic context, where this type of error set arises from the purely linear part of a quadratic factor (see \cite[Corollary 3.30]{Terry.2021a}).  One of the main theorems of the present paper will characterize this kind of error via a new generalization of the order property (see  Theorem \ref{thm:fopintro}).  

We will now give rough definitions of these notions, and refer the reader to Definition \ref{def:error} for details.  Given a $3$-uniform hypergraph $H=(V,E)$ and a decomposition $\calP$ with $\calP_{vert}=\{V_1, \ldots, V_t\}$, we say that 
\begin{enumerate}
\item $\calP$ has \emph{zero $\disc_{2,3}$-error} if there are no irregular triads in $\calP$ with respect to $H$;
\item $\calP$ has \emph{binary $\disc_{2,3}$-error} if there is a small set $\Sigma\subseteq {[t]\choose 2}$ such that every irregular triad from $\calP$ is contained in some $V_iV_jV_k$ with one of the pairs $ij$, $jk$, or $ik$ in $\Sigma$;
\item $\calP$ has \emph{linear $\disc_{2,3}$-error} if there is a small set $\Sigma\subseteq {[t]\choose 3}$ such that every irregular triad from $\calP$ is contained in some $V_iV_jV_k$ with $ijk\in \Sigma$.  
\end{enumerate}
The above notions describe the behavior of the irregular triads in a fixed decomposition $\calP$ of a hypergraph $H$.  This hypergraph may have several distinct regular decompositions, and it is possible that some of these may have  zero/binary/linear $\disc_{2,3}$-error while others do not.  We will be interested in the case where at least one decomposition exhibits some special behavior.  For this reason, we say that the hypergraph $H$ \emph{admits zero/binary/linear $\disc_{2,3}$-error} if there exists a decomposition $\calP$ which has zero/binary/linear $\disc_{2,3}$-error with respect to $H$.   Extending this definition to the setting of hereditary classes, we say a hereditary property of $3$-uniform hypergraphs $\calH$ \emph{admits zero/binary/linear $\disc_{2,3}$-error} if all the hypergraphs in $\calH$ have regular decompositions with zero/binary/linear $\disc_{2,3}$-error, respectively.   

Our main theorem provides a complete characterization of the properties which admit linear error in terms of a new generalization of the order property, which we call the \emph{functional order property} ($\FOP_2$, see Definition \ref{def:fop}).  We state this result below (see also Theorem \ref{thm:FOP}) and refer the reader to Section \ref{subsec:results} for the precise definitions.

\begin{theorem}\label{thm:fopintro}
Suppose $\calH$ is a hereditary property of $3$-uniform hypergraphs. Then the following are equivalent.
\begin{enumerate}
\item There is some $k$ so that no $H\in \calH$ has the $k$-functional order property.
\item  For all $\e_1>0$, $\e_2:\mathbb{N}\rightarrow (0,1]$, and $t_0,\ell_0\geq 1$, there are $T,L,N$ such that the following hold.  For all $H=(V,E)\in \calH$ with $|V|\geq N$, there are $t_0\leq t\leq T$, $\ell_0\leq \ell \leq L$, and $\calP$ a $(t,\ell,\e_1,\e_2(\ell))$ decomposition of $V$ which is $(\e_1,\e_2(\ell))$-regular and which has linear $\disc_{2,3}$-error with respect to $H$.
\end{enumerate}
\end{theorem}

This is the first analogue of Theorem \ref{thm:stable} to be proved in the setting of strong hypergraph regularity.  Our companion results in the arithmetic setting \cite{Terry.2021a} show a correspondence between forbidding the functional order property in groups, and quadratic arithmetic regularity lemmas where the ``error part" is supported on a small number of \emph{linear} atoms.  These results show that this particular generalization of Theorem \ref{thm:stable} has strong connections to algebraic structure.  Our proof techniques in this paper will also show that the relationship between the functional order property and stability is analogous to that between $\VC_2$-dimension and $\VC$-dimension (see Section \ref{subsec:fopchar}).

Our other main results focus on zero and binary $\disc_{2,3}$-error, and their relationships to various examples and possible generalizations of stability.  We show (Theorem \ref{thm:vdischbark}) that a hereditary property which admits zero error must be close to slicewise stable (Definition \ref{def:weaklystable}). We conjecture that the converse holds as well (Conjecture \ref{conj:weaktriads}).  We also show that the notion of zero error reduces to a problem about weak regularity (Theorem \ref{thm:vdisceq}), and cannot be characterized by stability, as defined in \cite{Ackerman.2017, Chernikov.2016zb}.   We show that  being close to slicewise stable  implies that a property admits binary error (Proposition \ref{prop:wsbinary}), but that the converse to this statement is false (Proposition \ref{prop:slvdiscbin}). 

We give two distinct examples of $3$-uniform hypergraphs which do not admit regular decompositions with binary error. The first arises from a combinatorial definition called the \emph{hyperplane order property} ($\HOP_2$), and the other from an example of Green and Sanders \cite{Green.2015qy4} in elementary abelian $p$-groups. Based on these examples, we prove a general sufficient condition for when a property cannot admit binary error (Chapter \ref{sec:binary}).  

The hyperplane order property has a strong geometric appeal as a potential generalization of the order property.  The existence of the second example in elementary abelian $p$-groups, however, shows that $\HOP_2$ is unlikely to be a fundamental definition generalizing the order property.  We consider this observation crucial for the future investigation of this sphere of problems.

Our results imply that the notions of zero, binary, and linear error are pairwise distinct.  Further, we show that binary error cannot be characterized by stability or slicewise stability, and cross-cuts the theories of bounded $\VC$-dimension.  

The main results of this paper also require several auxiliary results, some of which are of independent interest.  These include finitary versions of the results of Chernikov and Towsner \cite{Chernikov.2020} for $3$-uniform hypergraphs (see Theorems \ref{thm:vc2finite} and \ref{thm:vdischom}), as well as a strong form of the stable graph regularity lemma (see Theorem \ref{thm:functionstablereg}).  

\section*{Acknowledgements} The first author would like to thank M. Malliaris for conversations about stable graphs which informed several proofs in this paper.  The first author would also like to thank G. Conant for many helpful discussions around the paper's main results. The first author was partially supported by NSF grants DMS-2115518 and DMS-2239737 and by a Sloane Fellowship. 

Both authors are profoundly grateful to the referees for their very careful reading of the initial manuscript and their numerous helpful comments and corrections.

The authors' collaboration has been supported by several travel grants over the years, including by the London Mathematical Society, the Simons Foundation, and the Association for Women in Mathematics. 

%% file: chapter2.tex
\chapter{Formal statements of problems and results}\label{sec:mainresults}

In this chapter we state our results in full detail (see Section \ref{subsec:results}).  This requires the definitions in Section \ref{subsec:quasi1}.  The reader familiar with these may jump ahead to Section \ref{subsec:results}.  We begin here with a brief overview of notational conventions appearing in this section.  A more exhaustive account of our notation appears in Chapter \ref{sec:prelims}.  

\label{notation1} Given $x\neq y$, we will write $xy$ as shorthand for the two-element set $\{x,y\}$.  Similarly, for three distinct elements $x,y,z$, we write $xyz$ as shorthand for the three-element set $\{x,y,z\}$.   For sets $V_1,V_2$, we let $K_2[V_1,V_2]=\{xy: x\in V_1, y\in V_2, x\neq y\}$.

A \emph{$3$-uniform hypergraph} is a pair $H=(V,F)$ where $V$ is a vertex set, and $F\subseteq {V\choose 3}$ is an edge set (here ${V\choose 3}=\{Y\subseteq V: |Y|=3\}$). Throughout, we will refer to $3$-uniform hypergraphs as simply \emph{$3$-graphs}.  Given a $3$-graph $H$, an \emph{induced sub-$3$-graph} is a $3$-graph $(V',F')$ where $V'\subseteq V$ and $F'=F\cap {V'\choose 3}$.  A \emph{hereditary $3$-graph property} is a class of finite $3$-graphs closed under induced sub-$3$-graphs and isomorphism\footnote{In the introduction, we referred to these as hereditary properties of $3$-uniform hypergraphs.}.  

In the following sections, we will frequently include translations to model theoretic language to aid the model theoretic reader.  Any such undefined terminology is not necessary for the rest of the paper.  To facilitate these explanations, we will use the fact that every hereditary $3$-graph property $\calH$ is the class of finite models of a universal theory $T_{\calH}$ in the language $\calL=\{R(x,y,z)\}$.  In particular, $T_{\calH}$ contains the axioms stating that $R(x,y,z)$ is symmetric and irreflexive, as well as the sentence stating ``there is no substructure isomorphic to $S$," for every finite $3$-graph $S$ which is not in $\calH$.  This means that, given a $3$-graph $H$ specified using the combinatorial notation, $H=(V,F)$, we have 
$$
F=\{xyz\in {V\choose 3}: H\models R(x,y,z)\}.
$$
We also note that given a $3$-graph $H=(V,E)$, an induced sub-$3$-graph of $H$ is simply a substructure of $H$ in the language $\calL$.

\section{Quasirandomness and regularity in $3$-graphs}\label{subsec:quasi1}

In this section we introduce two notions of quasirandomness for $3$-graphs. We largely follow the notation of \cite{Frankl.2002,Nagle.2013}, with some modifications.   We begin by recalling the definition of quasirandomness for graphs (recall that given sets $V_1,V_2$, we let $K_2[V_1,V_2]=\{xy: x\in V_1, y\in V_2, x\neq y\}$).

\begin{definition}[$\disc_2$]\label{def:disc2}
Suppose $\e>0$ and $G=(U\cup V, E)$ is a bipartite graph with $|E|=d|U||V|$.  We say $G$ \emph{has $\disc_2(\e)$} if for all $U'\subseteq U$ and $V'\subseteq V$, 
$$
\big||E\cap K_2[U',V']|-d|U'||V'|\big| \leq \e |U||V|.
$$
We say that $G$ \emph{has $\disc_2(\e;d')$} if it has $\disc_2(\e)$ and $||E|-d'|U||V||\leq \e |U||V|$. (Here ``disc'' is short for discrepancy.)
\end{definition}

The above notion of $\disc_2$-quasirandomness is essentially equivalent to regularity and several other notions of quasirandomness in graphs.  We will present some of these in more detail in Section \ref{subsec:quasi2}.   Turning to $3$-graphs, the most naive generalization of $\disc_2$ yields Definition \ref{def:vdisc} below.  Given sets $V_1,V_2,V_3$, we let 
$$
K_3[V_1,V_2,V_3]=\{uvw: u\in V_1, v\in V_2, w\in V_3, u\neq v, u\neq w, v\neq w\}.
$$ 

\begin{definition}[$\vdisc_3$]\label{def:vdisc}
Suppose $\e>0$ and $H=(U\cup V\cup W,F)$ is a $3$-partite, $3$-graph with $|F|=d|U||V||W|$.  We say $H$ has $\vdisc_3(\e)$ if for all $U'\subseteq U$, $V'\subseteq V$, and $W'\subseteq W$,
$$
\Big||E\cap K_3[U',V',W']|-d|U'||V'||W'|\Big|\leq \e |U||V||W|.
$$
\end{definition}

The ``v'' in vdisc is meant to emphasize that this definition considers quasirandomness with respect to \emph{vertex} subsets only.  This notion of quasirandomness was introduced in the early 1990s \cite{Chung.1991} as a generalization of quasirandomness for graphs.  There is a corresponding regularity lemma, first proved by Chung \cite{Chung.1991} and stated as Theorem \ref{thm:vreg} below.

\begin{definition}[$\vdisc_3$-regular triples and partitions]\label{def:vdisc3reg}
Suppose  $H=(V,F)$ is a $3$-graph and $\calP=\{V_1,\ldots, V_t\}$ is a partition of $V$.
\begin{enumerate}
\item Given $V_iV_jV_k\in {\calP\choose 3}$, we say $V_iV_jV_k$ \emph{has $\vdisc_3(\e)$ with respect to $H$} if the $3$-partite $3$-graph $(V_i\cup V_j\cup V_k, F\cap K_3[V_i,V_j,V_k])$ has $\vdisc_3(\e)$.
\item We say $\calP$ is a \emph{$\vdisc_3(\e)$-regular partition for $H$} if 
\begin{enumerate}
\item $\calP$ is an equipartition, and 
\item for all but at most $\e{t\choose 3}$ many $ijk\in {[t]\choose 3}$, the triple $V_iV_jV_k$ has $\vdisc_3(\e)$ with respect to $H$.  
\end{enumerate}
\end{enumerate}
\end{definition}

\begin{theorem}[$\vdisc_3$-regularity lemma for 3-graphs, \cite{Chung.1991}]\label{thm:vreg}
For all integers $t_0\geq 1$ and $\e>0$ there exists $T=T(t_0,\e)$ and $N=N(t_0,\e)$ such that the following holds.  If $H=(V,F)$ is a $3$-graph with $|V|\geq N$, there is some $t$ with $t_0\leq t\leq T$ and a $\vdisc_3(\e)$-regular partition for $H$ with $t$ parts. 
\end{theorem}

It was well known that $\vdisc_3$-quasirandomness does not support a general counting lemma. It was shown later \cite{Kohayakawa.2010} that it does suffice for counting copies of \emph{linear} hypergraphs (see also \cite{Person.2009} for an application of this counting lemma).  A general counting lemma was the motivation for a more sophisticated notion of quasirandomness \cite{Frankl.2002,Gowers.20063gk}, which considers not only sets of vertices but also sets of pairs of vertices.  Given a graph $G=(V,E)$, we let $K_3^{(2)}(G)$ denote the set of triples from $V$ forming a triangle in $G$. In other words,\label{k3}
 $$
 K_3^{(2)}(G)=\left\{xyz\in {V\choose 3}: xy, yz, xz\in E\right\}.
 $$
 Given a $3$-graph $H=(V,F)$ on the same vertex set, we say that $G$ \emph{underlies $H$} if $F\subseteq K_3^{(2)}(G)$, i.e. if all the edges of $H$ sit atop a triangle from $G$.  A subgraph $G'\subseteq G$ is a graph $(V',E')$ where $V'\subseteq V$ and $E'\subseteq E$.  
 
 Our next definition describes when a $3$-graph $H$ is quasirandom relative to an underlying quasirandom graph $G$.  Informally, such a pair $(H,G)$ will be considered  quasirandom if the graph $G$ is quasirandom, and further, the ternary edges $E$ of $H$ are ``uniformly distributed" on the triangles of $G$, in the sense that restricting $E$ to the triangles of any \emph{subgraph} of $G$ does not drastically change the relative edge density. We now give the formal definition.
 
\begin{definition}[$\disc_{2,3}$]\label{def:regtrip}
Let $\e_1,\e_2>0$, and suppose $H=(V,F)$ is a $3$-graph. Assume $G=(V_1\cup V_2\cup V_3,E)$ is a $3$-partite graph underlying $H$, and suppose $d_3$ is defined by $|F|=d_3|K^{(2)}_3(G)|$.  We say that \emph{$(H,G)$ has $\disc_{2,3}(\e_1,\e_2)$} if there is $d_2\in (0,1)$ such that for each $1\leq i<j\leq 3$, the graph $G[V_i,V_j]$ has $\disc_2(\e_2;d_2)$, and for every subgraph $G'\subseteq G$,
\[||F\cap K_3^{(2)}(G')|-d_3|K^{(2)}_3(G')||\leq \e_1 d_2^3 |V_1||V_2||V_3|.\]
\end{definition}

To aid the reader's understanding, we give here an example of a pair $(H,G)$ which is \emph{not} quasirandom in the sense of Definition \ref{def:regtrip}.  Let $G=(V_1\cup V_2\cup V_3,E)$ be any tripartite graph where $G[V_1,V_2]$, $G[V_1,V_3]$, and $G[V_2,V_3]$ are quasirandom as graphs.  Choose an arbitrary partition $A\cup B$ of $R\cap K_2[V_2,V_3]$ so that $A$ and $B$ have roughly equal size, and let $H$ be the $3$-graph whose edges are the triangles of $G$ containing a pair of vertices belonging to $A$, i.e. $H=\{xyz\in K_3^{(2)}(G): yz\in A\}$.  Then there are two naturally defined subgraphs of $G$, $G_A$ and $G_B$, obtained from $G$ by deleting the edges in $B$ and $A$, respectively. The relative density of $H$ on $G_A$ is $1$, while the relative density of $H$ on $G_B$ is zero.  Because of this, the pair $(H,G)$ will not satisfy Definition \ref{def:regtrip} for sufficiently small $\e_1,\e_2$. 

Our notational convention is to use single subscripts of $2$ or $3$ for notions of strictly binary or strictly ternary quasirandomness (e.g. $\disc_2$ and $\vdisc_3$), and the subscript 2,3 as above to emphasize that $\disc_{2,3}$ quasirandomness has both a binary and ternary quasirandomness requirement.  Definition \ref{def:regtrip} has several equivalent formulations, and a more complete account of these will appear in Section \ref{subsec:quasi2}.  We now state the definition which will serve as a higher order analogue of an equitable  partition of a vertex set.  In contrast to the partition appearing in Theorem \ref{thm:vreg}, Definition \ref{def:decomp} will involve not only a partition of the vertex set but also a partition of the set of pairs of vertices. 

\begin{definition}[$(t,\ell,\e_1,\e_2)$-decomposition]\label{def:decomp}
Let $V$ be a vertex set, $t,\ell \in \mathbb{N}^{>0}$, and $\e_1,\e_2>0$.  A \emph{$(t,\ell,\e_1,\e_2)$-decomposition} $\calP$ for $V$ consists of a partition $V=V_1\cup \ldots \cup V_t$ and for each $1\leq i\neq j\leq t$, a collection $P^1_{ij},\ldots, P^\ell_{ij}$ of disjoint subsets of ${V\choose 2}$ such that the following hold. 
\begin{enumerate}
\item $|V_1|\leq \ldots \leq |V_t|\leq |V_1|+1$,
\item For each $1\leq i\neq j\leq t$, $K_2[V_i,V_j]=\bigcup_{\alpha\leq \ell} P^\alpha_{ij}$,
\item For all but $\e_1 |V|^3$ triples $xyz\in {V\choose 3}$, the following holds.  There is some $ijk\in {[t]\choose 3}$ and $1\leq \alpha,\beta,\gamma\leq \ell$ such that $xy\in P^\alpha_{ij}$, $xz\in P^\beta_{ik}$, $yz\in P^\gamma_{jk}$, and such that each of $P^\alpha_{ij}, P^\beta_{ik}, P^\gamma_{jk}$ have $\disc_{2}(\e_2,1/\ell)$.  

\end{enumerate}
\end{definition}

\label{triadnotation}In the notation of Definition \ref{def:decomp}, a \emph{triad of $\calP$} is a $3$-partite graph of the form
$$
G^{\alpha,\beta,\gamma}_{ijk}:=(V_i\cup V_j\cup V_k, P_{ij}^\alpha\cup P_{ik}^\beta\cup P_{jk}^\gamma)
$$
for some $ijk\in {[t]\choose 3}$ and $\alpha,\beta,\gamma\leq \ell$.

To ease notation, we define $\triads(\calP):=\{G_{ijk}^{\alpha,\beta,\gamma}: ijk\in {[t]\choose 3}, \alpha,\beta,\gamma \in [\ell]\}$.  When $H=(V,F)$ is a $3$-graph, then for each triad $G_{ijk}^{\alpha,\beta,\gamma}\in \triads(\calP)$, we write
$$
H_{ijk}^{\alpha,\beta,\gamma}:=(V_i\cup V_j\cup V_k, F\cap K_3^{(2)}(G_{ijk}^{\alpha,\beta,\gamma})).
$$

Note that by definition, $H_{ijk}^{\alpha,\beta,\gamma}$ is a $3$-partite $3$-graph, and $G_{ijk}^{\alpha,\beta,\gamma}$ underlies $H_{ijk}^{\alpha,\beta,\gamma}$.  We next define regular triads and decompositions, in analogy with Definition \ref{def:decomp}.  Informally speaking, given a hypergraph $H=(V,F)$ and a decomposition $\calP$ of $V$, we will call $\calP$ regular with respect to $H$ if for almost all triples of vertices $x,y,z$ from $V$, there is a triad $G$ from $\calP$ so that $x,y,z$ forms a triangle in $G$, and so that $G$ is regular (in the sense of Definition \ref{def:regtrip})) with respect to the restriction of $H$ to the triangles of $G$. We now give the precise definitions.

\begin{definition}[$\disc_{2,3}$-regular triads and decompositions]\label{def:regdec}
Suppose $H=(V,F)$ is a $3$-graph and $\calP$ is a $(t,\ell,\e_1,\e_2)$-decomposition for $V$.  
\begin{enumerate}
\item We say $G_{ijk}^{\alpha,\beta,\gamma}\in \triads(\calP)$  \emph{has $\disc_{2,3}(\e_1,\e_2)$ with respect to} the 3-graph $H$ if  $(H_{ijk}^{\alpha,\beta,\gamma}, G^{\alpha,\beta,\gamma}_{ijk})$ has $\disc_{2,3}(\e_1,\e_2)$.  

\item We say that $\calP$ is \emph{$(\e_1,\e_2)$-regular} for $H$ if for all but $\e_1 n^3$ many triples $xyz\in {V\choose 3}$, there is some $G_{ijk}^{\alpha,\beta,\gamma}\in \triads(\calP)$, so that $xyz\in K_3^{(2)}(G_{ijk}^{\alpha,\beta,\gamma})$ and which has  $\disc_{2,3}(\e_1,\e_2)$ with respect to $H$.  
\end{enumerate}
\end{definition}

Definition \ref{def:regdec}, along with its equivalent formulations (see Section \ref{subsec:quasi2}), are now considered standard notions of regular triads and decompositions for $3$-graphs (see  \cite{Gowers.20063gk, Gowers.2007, Rodl.2005,Nagle.2013}).  We now state the regularity lemma for $\disc_{2,3}$-quasirandomness, which is taken from \cite{Frankl.2002}.  We note that in this regularity lemma, the $\e_2$ parameter is allowed to depend on the parameter $\ell$.

\begin{theorem}[$\disc_{2,3}$-regularity lemma for 3-graphs, \cite{Frankl.2002,Gowers.20063gk}]\label{thm:reg2} For all $\e_1>0$, every function $\e_2:\mathbb{N}\rightarrow (0,1]$, and every $t_0\in \mathbb{N}$, there exist positive integers $T_0$ and $n_0$ such that for any $3$-graph $H=(V,F)$ on $n\geq n_0$ vertices, there exists a $(\e_1,\e_2(\ell))$-regular $(t,\ell, \e_1,\e_2(\ell))$-decomposition $\calP$ for $H$ with $t_0\leq t \leq T_0$ and $\ell\leq T_0$.
\end{theorem}

We end this section by setting up some terminology around irregular triads in the sense of Definition \ref{def:decomp}.  

\begin{definition}[Irregular triads]\label{def:regtriads2}
Suppose $H=(V,E)$ is a $3$-graph and $\calP$ is a $(t,\ell,\e_1,\e_2)$-decomposition for $V$.
\begin{enumerate}
\item We say $G_{ijk}^{\alpha,\beta,\gamma}\in \triads(\calP)$ is \emph{$\disc_{2,3}$-regular with respect to $(H,\calP)$} if it has $\disc_{2,3}(\e_1,\e_2(\ell))$ with respect to $H$. Otherwise, we say $G_{ijk}^{\alpha,\beta,\gamma}$ is \emph{$\disc_{2,3}$-irregular}.
\item We say $G_{ijk}^{\alpha,\beta,\gamma}\in \triads(\calP)$ is \emph{$\disc_2$-regular with respect to $(H,\calP)$} if each of $P_{ij}^\alpha, P_{ik}^\beta, P_{jk}^\gamma$ have $\disc_2(\e_2(\ell);1/\ell)$.  Otherwise, we say $G_{ijk}^{\alpha,\beta,\gamma}$ is \emph{$\disc_{2}$-irregular}.
\item We say $G_{ijk}^{\alpha,\beta,\gamma}\in \triads(\calP)$ is \emph{$\disc_{3}$-irregular with respect to $(H,\calP)$} if it is $\disc_2$-regular but not $\disc_{2,3}$-regular with respect to $(H,\calP)$.
\end{enumerate}
\end{definition}

Note that in the notation of Definition \ref{def:regtriads2}, $G_{ijk}^{\alpha,\beta,\gamma}$ is $\disc_{2,3}$-irregular with respect to $(H,\calP)$ if and only if it is $\disc_2$-irregular with respect to $(H,\calP)$ or $\disc_3$-irregular with respect to $(H,\calP)$.

There are several other notions of quasirandomness for hypergraphs which we do not consider (see \cite{Lenz.2012, Towsner.2016}).  Our main focus in this paper is to consider $\disc_{2,3}$-regularity, due to the important applications of its corresponding counting lemma.  As our results will show, one of our main problems related to $\disc_{2,3}$ (namely that centering around ``zero $\disc_{2,3}$-error") will in fact collapse to a problem about $\vdisc_3$ (see Theorem \ref{thm:vdisceq}).  For this reason, we will also be considering $\vdisc_3$-regularity. 

\section{Statements of problems}\label{subsec:results}

In this section we formally state the problems of interest in this paper.  We begin with two basic definitions  concerning $\vdisc_3$-regularity.  Roughly speaking, we will say that a property admits \emph{zero $\vdisc_3$-error} if the $3$-graphs in $\calH$ admit $\vdisc_3$-regular partitions with \emph{no} irregular triples.  On the other hand, we will say $\calH$ admits \emph{binary $\vdisc_3$-error} if the $3$-graphs in $\calH$ admit $\vdisc_3$-regular partitions where the irregular triples can be ``controlled" by a small number of pairs of vertex parts. 

\begin{definition}[zero and binary $\vdisc_3$-error]\label{def:vdiscerror}
Let $\calH$ be a hereditary $3$-graph property.  
\begin{enumerate}
\item We say $\calH$ \emph{admits zero $\vdisc_3$-error}, if for all $\e>0$ and $t_0\geq 1$, there is $T=T(\e,t_0)$ and $N=N(\e,t_0)$ such that for all $H=(V,F)\in \calH$ with $|V|\geq N$, there is an equipartition $\calP$ of $V$ with some $t_0\leq t\leq T$ many parts, such that every triple $XYZ\in {\calP\choose 3}$ has $\vdisc_3(\e)$ with respect to $H$. 

\item We say $\calH$ \emph{admits binary $\vdisc_3$-error}, if for all $\e>0$ and $t_0\geq 1$, there is $T=T(\e,t_0)$ and $N=N(\e,t_0)$ such that for all $H=(V,F)\in \calH$ with $|V|\geq N$, there is an equipartition $\calP$ of $V$ with some $t_0\leq t\leq T$ many parts, along with a set $\Sigma \in {\calP\choose 2}$ of size at most $\e t^2$, such that  every triple $XYZ\in {\calP\choose 3}$ with $XY, YZ, XZ\notin \Sigma$ has $\vdisc_3(\e)$ with respect to $H$. 
 \end{enumerate}
\end{definition}

When $\calH$ fails (1) (respectively, (2)), we say it \emph{requires non-zero (respectively, non-binary) $\vdisc_3$-error}.   Clearly, if $\calH$ satisfies (1), it also satisfies (2).  

We next define three related notions for $\disc_{2,3}$-regularity.  We give an informal account of these definitions ahead of the precise statements below.  In close analogy to Definition \ref{def:vdiscerror}(1), we will say that a property $\calH$ admits \emph{zero $\disc_{2,3}$-error} if the $3$-graphs in $\calH$ admit $\disc_{2,3}$-regular partitions with \emph{no} irregular triads.  Similarly to Definition \ref{def:vdiscerror}(2), we say that $\calH$ admits \emph{binary $\disc_{2,3}$-error} if the $3$-graphs in $\calH$ admit $\disc_{2,3}$-regular partitions where the irregular triads can be ``controlled" by a small number of pairs of vertex parts.  We will define one additional notion which (for reasons that will be apparent shortly) has no counterpart in Definition \ref{def:vdiscerror}.  In particular, we will say that $\calH$ admits \emph{linear $\disc_{2,3}$-error} if the $3$-graphs in $\calH$ admit $\disc_{2,3}$-regular partitions where the irregular triads can be ``controlled" by a small number of triples of vertex parts.

\begin{definition}[zero, binary, and linear $\disc_{2,3}$-error]\label{def:error}
Suppose $\calH$ is a hereditary $3$-graph property.  
\begin{enumerate}
\item We say $\calH$  \emph{admits zero $\disc_{2,3}$-error} if for all $t_0,\ell_0\geq 1$, $\e_1>0$, all $\e_2:\mathbb{N}\rightarrow (0,1]$, there are $N,T,L$ such that for every $H=(V,F)\in \calH$ with $|V|\geq N$, there exists $t_0\leq t\leq T$, $\ell_0\leq \ell\leq L$, and a $(t,\ell, \e_1,\e_2(\ell))$-decomposition $\calP$ of $V$, such that every triad in $\triads(\calP)$ has $\disc_{2,3}(\e_1,\e_2(\ell))$ with respect to $H$. 
\item  We say $\calH$  \emph{admits binary $\disc_{2,3}$-error} if for all $t_0,\ell_0\geq 1$, $\e_1>0$, all $\e_2:\mathbb{N}\rightarrow (0,1]$, there are $N,T,L$ such that for every $H=(V,F)\in \calH$ with $|V|\geq N$, there exists $t_0\leq t\leq T$, $\ell_0\leq \ell\leq L$, and a $(t,\ell, \e_1,\e_2(\ell))$-decomposition $\calP$ of $V$ such that the following holds.  There is a set $\Gamma\subseteq{[t]\choose 2}$ with $|\Gamma|\leq \e_1 t^2$, such that for every $G_{ijk}^{\alpha,\beta,\gamma}\in \triads(\calP)$ with $ij, ik, jk\notin \Gamma$, $G_{ijk}^{\alpha,\beta,\gamma}$ has $\disc_{2,3}(\e_1,\e_2(\ell))$ with respect to $H$. 
\item  We say $\calH$   \emph{admits linear\footnote{The use of the term ``linear'' stems form the authors' work in the arithmetic context, where this type of error set arises from the purely linear part of a quadratic factor, see \cite[Corollary 3.30]{Terry.2021a}.} $\disc_{2,3}$-error} if for all $t_0,\ell_0\geq 1$, $\e_1>0$, all $\e_2:\mathbb{N}\rightarrow (0,1]$, there are $N,T,L$ such that for every $H=(V,F)\in \calH$ with $|V|\geq N$, there exists $t_0\leq t\leq T$, $\ell_0\leq \ell\leq L$, and a $(t,\ell, \e_1,\e_2(\ell))$-decomposition $\calP$ of $V$ such that the following holds. There is a set $\Sigma\subseteq{[t]\choose 3}$ with $|\Sigma|\leq \e_1 t^3$, such that for every $G_{ijk}^{\alpha,\beta,\gamma}\in \triads(\calP)$ with $ijk\notin \Sigma$, $G_{ijk}^{\alpha,\beta,\gamma}$ has $\disc_{2,3}(\e_1,\e_2(\ell))$ with respect to $H$. 
\end{enumerate}
\end{definition} 

If $\calH$ fails (1) (respectively (2), (3)), we say $\calH$ \emph{requires non-zero (respectively non-binary, non-linear) $\disc_{2,3}$-error}.   By definition, if $\calH$ admits zero $\disc_{2,3}$-error, it also admits binary $\disc_{2,3}$-error and linear $\disc_{2,3}$-error.  It is also not difficult to see that if $\calH$ admits binary $\disc_{2,3}$-error then it admits  linear $\disc_{2,3}$-error (one just takes $\Sigma$ to be the set of all triples with at least one pair in $\Gamma$).  Thus $(1)\Rightarrow (2)\Rightarrow (3)$.  The reader may be wondering why the possibility of a linear error appears in Definition \ref{def:error} but not in Definition \ref{def:vdiscerror}.  This is because Theorem \ref{thm:vreg} implies that \emph{every} hereditary $3$-graph property admits linear $\vdisc_3$-error.

One can restate parts (1) and (2) of Definition \ref{def:error} in the rough form ``sufficiently large elements in $\calH$ have regular decompositions with no irregular triads of a certain type.''    It is clear that (1) is already in this form.  Indeed, informally speaking, (1) says that sufficiently large elements in $\calH$ have $\disc_{2,3}$-regular decompositions with \emph{no} $\disc_{2,3}$-irregular triads.  It turns out that (2) is equivalent to saying that all sufficiently large elements in $\calH$ have $\disc_{2,3}$-regular decompositions with no $\disc_3$-irregular triads (see Definition \ref{def:regtriads2}).  More specifically, we will prove the following in Chapter \ref{sec:binary}.

\begin{theorem}\label{thm:disc3}
$\calH$ admits binary $\disc_{2,3}$-error if and only if for all $t_0,\ell_0\geq 1$, $\e_1>0$, all $\e_2:\mathbb{N}\rightarrow (0,1]$, there are $n,T,L$ such that for every $H=(V,F)\in \calH$ with $|V|\geq N$, there are $t_0\leq t\leq T$, $\ell_0\leq \ell\leq L$ and an $(\e_1,\e_2(\ell))$-regular $(t,\ell, \e_1,\e_2(\ell))$-decomposition $\calP$ for $H$, such that no triad of $\calP$ is $\disc_{3}$-irregular with respect to $H$.
\end{theorem}

On the other hand, we do not know of a way to restate (3) using a statement of the form ``there are no irregular triads of a certain type."  It is for this reason that we transition away from thinking about ``no irregular triads", and instead aim to \emph{control} irregular triads using lower order sets. 

Finally, the reader may be wondering why none of these definitions address $\disc_2$-irregular triads.  It turns out that these are simply not very interesting, since any $\disc_{2,3}$-regular partition $\calP$ can be modified so that it has no $\disc_2$-irregular triads.  This is implicit in the work of Frankl and R\"{o}dl \cite{Frankl.2002}, but we will also include a proof for completeness.  Specifically, we prove the following in Section \ref{subsec:intersecting}.

\begin{proposition}\label{prop:binaryboring}
For every hereditary $3$-graph property $\calH$ the following holds. For all $\ell_0,t_0\geq 1$, $\e_1>0$, all $\e_2:\mathbb{N}\rightarrow (0,1]$, there are $n_0,T_0,L_0$ such that for all $n\geq n_0$, every $G\in \calH$ with at least $n_0$ vertices has an $(\e_1,\e_2(\ell))$-regular $(t,\ell, \e_1,\e_2(\ell))$-decomposition for some $t_0\leq t\leq T_0$, and $\ell_0\leq \ell\leq L_0$, such that every triad of $\calP$ is $\disc_{2}$-regular.\end{proposition}

Thus we are left with three potentially interesting notions with regards to $\disc_{2,3}$-irregular triads: zero, binary, and linear $\disc_{2,3}$-error. In summary then, the following are the problems concerning irregular triples/triads in $3$-graphs which we will consider in this paper.

\begin{problem}\label{prob:triads}
Suppose $\calH$ is a hereditary $3$-graph property.
\begin{enumerate}[label=\normalfont(\arabic*)]
\item Characterize when $\calH$ admits zero $\vdisc_3$-error. 
\item Characterize when $\calH$ admits binary $\vdisc_3$-error. 
\item Characterize when $\calH$ admits zero $\disc_{2,3}$-error. 
\item Characterize when $\calH$ admits binary $\disc_{2,3}$-error. 
\item Characterize when $\calH$ admits linear $\disc_{2,3}$-error. 
\end{enumerate}
\end{problem}

It will turn out that in order to address these problems, we must also consider the following general question, which is also motivated by Theorem \ref{thm:vc}: which hereditary $3$-graph properties have regular decompositions where all regular triples/triads have density near $0$ or $1$?  We will informally refer to triads/triples with density near $0$ or $1$ as \emph{homogeneous}.  Homogeneity is stronger than regularity, in the sense that a homogeneous triple $(V_i,V_j,V_k)$ will be $\vdisc_3$-regular, and a homogeneous triad $(V_i\cup V_j\cup V_k, P_{ij}^\alpha\cup P_{jk}^\beta\cup P_{ik}^\gamma)$ will be $\disc_{2,3}$-regular (see Propositions \ref{prop:homimpliesrandomv} and \ref{prop:homimpliesrandome}). A decomposition $\calP$ of $V$ will be said to be homogeneous with respect to a 3-graph $H=(V,F)$ if almost all triples of vertices lie in a homogeneous triad.

\begin{definition}
Suppose $H=(V,F)$ is a $3$-graph with $|V|=n$ and $\mu>0$.
\begin{enumerate}
\item  Suppose $\calP$ is a partition of $V$.  We say $\calP$ is  \emph{$\mu$-homogeneous with respect to $H$}, if at least $(1-\mu){n\choose 3}$ triples $xyz\in {V\choose 3}$ satisfy the following.  There is some $V_iV_jV_k\in {\calP\choose 3}$ such that $xyz\in K_3[V_i,V_j,V_k]$ and 
$$
|F\cap K_3[V_i,V_j,V_k]|/|K_3[V_i,V_j,V_k]|\in [0,\mu)\cup (1-\mu,1].
$$
\item Suppose $t,\ell\geq 1$ and $\calP$ is a $(t,\ell)$-decomposition of $V$.  We say that $\calP$ is \emph{$\mu$-homogeneous with respect to $H$} if at least $(1-\mu){n\choose 3}$ triples $xyz\in {V\choose 3}$ satisfy the following.  There is some triad $G_{ijk}^{\alpha,\beta,\gamma}\in \triads(\calP)$ such that $xyz\in K_3^{(2)}(G_{ijk}^{\alpha,\beta,\gamma})$ and 
$$
|F\cap K_3^{(2)}(G_{ijk}^{\alpha,\beta,\gamma})|/|K_3^{(2)}(G_{ijk}^{\alpha,\beta,\gamma})|\in [0,\mu)\cup (1-\mu,1].
$$
\end{enumerate}
\end{definition}

\begin{definition}[$\vdisc_3$- and $\disc_{2,3}$-homogeneity]\label{def:hom}
Suppose $\calH$ is a hereditary $3$-graph property.  
\begin{enumerate}
\item We say $\calH$ is \emph{$\vdisc_3$-homogeneous} if for all $\e>0$ and $t_0$, there are $T=T(\e,t_0)$ and $N=N(\e,t_0)$ such that for all $H\in \calH$ on at least $N$ vertices, there is a partition $\calP$ of $V$ which is $\vdisc_3(\e)$-regular and $\e$-homogeneous with respect to $H$.
\item We say $\calH$ is \emph{$\disc_{2,3}$-homogeneous} if for all $t_0,\ell_0\geq 1$, $\e_1>0$ and $\e_2:\mathbb{N}\rightarrow (0,1]$, there are $T, L, N$ such that for $H\in \calH$ on at least $N$ vertices, there are $t_0\leq t\leq T$, $\ell_0\leq \ell \leq L$, and a $(t,\ell,\e_1,\e_2(\ell))$-decomposition $\calP$ of $V$ which is $(\e_1,\e_2(\ell))$-regular and $\e_1$-homogeneous with respect to $H$.
\end{enumerate}
\end{definition}

We can now state concrete versions of our homogeneity problems.
\begin{problem}\label{prob:hom}
Suppose $\calH$ is a hereditary $3$-graph property.
\begin{enumerate}[label=\normalfont(\arabic*)]
\item Characterize when $\calH$ is $\vdisc_3$-homogeneous.
\item Characterize when $\calH$ is $\disc_{2,3}$-homogeneous.
\end{enumerate}
\end{problem}

\section{Statements of results}\label{subsec:stateresults}
In this section we state the main results of the paper. We will give full answers to Problem \ref{prob:hom} and Problem \ref{prob:triads} (5), and partial results  for Problem \ref{prob:triads} (1)-(4).  

To state some of our results, we require a few more definitions.  First, we can recharacterize stability and $\VC$-dimension in terms of the bipartite counterpart of a graph.  
\begin{definition}[$\bip(G)$]
Given a graph $G=(V,E)$, let $\bip(G)$ denote the bipartite graph $(U\cup W, E')$, where $U=\{u_v:v\in V\}$ and $W=\{w_v:v\in V\}$ are disjoint copies of $V$, and $E'=\{u_vw_{v'}: vv'\in E\}$.  
\end{definition}

Note that the $\VC$-dimension of a finite graph $G$ is then the largest $k$ such that $\bip(G)$ contains an induced copy of the \emph{$k$th power set graph},
$$
U(k)=(\{a_i: i\in [k]\}\cup \{b_S:S\subseteq [k]\},\{a_ib_S: i\in S\}).
$$ 
Similarly, a graph $G$ has the $k$-order property if and only if $\bip(G)$ contains an induced copy of the \emph{half-graph of height $k$}, 
\[H(k)=(\{a_i : i\in [k]\}\cup \{b_j : j\in [k]\} , \{a_ib_j: i\leq j\}).\]  

We will also use a natural equivalence relation on hereditary properties.  Given two $3$-graphs $H=(V,F)$ and $H'=(V,F')$ on the same vertex set, we say that $H$ and $H'$ are \emph{$\delta$-close} if $|F\Delta F'|\leq \delta |V|^3$.

\begin{definition}[Closeness of hereditary properties]\label{def:close}
Suppose $\calH$ and $\calH'$ are hereditary $3$-graph properties.  We say \emph{$\calH$ is close to $\calH'$} if for all $\delta>0$, there is $N=N(\delta)$ such that if $n\geq N$, then every $G\in \calH$ on at least $n$ vertices is $\delta$-close to an element of $\calH'$ with the same vertex set.  

We say $\calH$ and $\calH'$ are \emph{close}, denoted $\calH\sim \calH'$, if $\calH$ is close to $\calH'$, and $\calH'$ is close to $\calH$.
\end{definition}

This notion appears in several other places in the literature.  For example it is referred to as indistinguishability in \cite{Alon.2000}.   Clearly $\sim$ is an equivalence relation.  It turns out that the behavior of a hereditary property with respect to Problems \ref{prob:hom} and \ref{prob:triads} depends only on its $\sim$-equivalence class.

\begin{proposition}\label{prop:simclasses2}
Suppose $\calH$ and $\calH'$ are close.  Then the following are equivalent.
\begin{enumerate}[label=\normalfont(\arabic*)]
\item  $\calH$ is $\vdisc_3$-homogeneous (is $\disc_{2,3}$-homogeneous, admits zero/binary\\ $\vdisc_3$-error, admits zero/binary/linear $\disc_{2,3}$-error). 
\item   $\calH'$ is  $\vdisc_3$-homogeneous (is $\disc_{2,3}$-homogeneous, admits zero/binary $\vdisc_3$-error, admits zero/binary/linear $\disc_{2,3}$-error). 
\end{enumerate}
\end{proposition}

The proof of Proposition \ref{prop:simclasses2} is straightforward from standard arguments combined with the induced removal lemma of R\"{o}dl and Schacht \cite{Rodl.2009}, and appears in Appendix \ref{app:simclasses}.

To answer Problem \ref{prob:hom} (2), we characterize $\disc_{2,3}$-homogeneity in terms of $\VC_2$-dimension, which we define below.  We first need some notation for the tripartite counterpart of a 3-graph.

\begin{definition}[$\trip(H)$]
Given a $3$-graph $H=(V,F)$, let $\trip(H)$ denote the $3$-partite $3$-graph $(X\cup Y\cup Z, E')$, where $X=\{x_v:v\in V\}$, $Y=\{y_v: v\in V\}$, and $Z=\{z_v: v\in V\}$ are disjoint copies of $V$, and $F'=\{x_vy_{v'}w_{v''}: vv'v''\in F\}$.  

For a hereditary $3$-graph property $\calH$, we define $\trip(\calH)=\{\trip(H): H\in \calH\}$. 
\end{definition}

It is straightforward to check that $\trip(\calH)$ is also a hereditary $3$-graph property. 

\begin{definition}[$V(k)$]\label{def:vk}
Given $k\geq 1$, define $V(k)$ to be the $3$-partite $3$-graph 
\[V(k)=(\{b_1,\ldots, b_k\}\cup \{c_1,\ldots, c_k\}\cup\{a_S: S\subseteq [k]^2\}, \{a_Sb_vc_w: (v,w)\in S\}).\]
\end{definition}

In other words, $V(k)$ encodes the power set of the Cartesian product $[k]\times [k]$.

\begin{definition}[$\VC_2$-dimension of a 3-graph]\label{def:vc2dim}
For a finite $3$-graph $H$, the \emph{$\VC_2$-dimension of $H$}, denoted by $\VC_2(H)$, is the largest $k$ such that $\trip(H)$ contains an induced copy of $V(k)$.    
\end{definition}

In other words, given a finite $3$-graph $H=(V,F)$, $\VC_2(H)$ is the largest $k$ for which there exists $b_1,\ldots, b_k,c_1,\ldots, c_k\in V$, and $a_S\in V$ for each $S\subseteq[k]^2$, so that $a_Sb_ic_j\in F$ if and only if $(i,j)\in S$.  

Given a hereditary $3$-graph property $\calH$, we say it has \emph{bounded $\VC_2$-dimension}\footnote{For the model theorist, $\calH$ has bounded $\VC_2$-dimension if and only if the edge relation is an $\NIP_2$-formula in every model of $T_{\calH}$.} (is $\NIP_2$) if there is some $k\geq 1$ such that $\VC_2(H)\leq k$ for all $H\in \calH$.  In this case, we define the \emph{$\VC_2$-dimension of $\calH$}, $\VC_2(\calH)$, to be $\max\{\VC_2(H): H\in \calH\}$.  Otherwise, we say that $\calH$ has \emph{unbounded $\VC_2$-dimension} (has $\IP_2$), and write $\VC_2(\calH)=\infty$.   This notion has several equivalent formulations, for example the following (see \cite{Chernikov.2019}).

\begin{fact}\label{fact:vc2universal}
$\calH$ has unbounded $\VC_2$-dimension if and only if $\trip(\calH)$ contains every finite $3$-partite $3$-graph. 
\end{fact}

This is analogous to the graph setting, where a hereditary graph property $\calH$ has unbounded $\VC$-dimension if and only if $\bip(\calH)$ contains every finite bipartite graph (this is folklore, and an easy exercise).  In the case of graphs, with a bit more work on top of Theorem \ref{thm:vc}, one can show that a hereditary graph property has bounded $\VC$-dimension if and only if it has regular decompositions where all \emph{regular} pairs are homogeneous (note this still allows for the possibility of irregular pairs).  We prove the analogous result for hereditary $3$-graph properties and $\VC_2$-dimension.  In particular, we first show the following result for $3$-graphs of uniformly bounded $\VC_2$-dimension. This result appears in Chapter \ref{sec:dischom}, with the proof beginning on page \pageref{proof:vc2finite}.

\begin{theorem}\label{thm:vc2finite}
For all $k\geq 1$, $\e_1>0$, $\e_2:\mathbb{N}\rightarrow (0,1]$, and $\ell_0,t_0\geq 1$, there are $N,T,L\geq 1$ such that the following hold.  Suppose $H=(V,F)$ is a $3$-graph with $|V|\geq N$ and $\VC_2(H)<k$.  Then there exist $t_0\leq t\leq T$, $\ell_0\leq \ell\leq L$, and a $(t,\ell,\e_1,\e_2(\ell))$-decomposition of $V$ which is $(\e_1,\e_2(\ell))$-regular and $\e_1$-homogeneous with respect to $H$.
\end{theorem}

 A closely related result was recently proved by Chernikov and Towsner \cite{Chernikov.2020}.  Specifically, they showed that if an $r$-uniform hypergraph equipped with a Keisler measure has bounded $\VC_{r-1}$-dimension, then it has a decomposition where all regular triads density near $0$ or $1$ with respect to said measure \cite{Chernikov.2020}.  They employ infinitary techniques, and their results are non-quantitative.  Our proof of Theorem \ref{thm:vc2finite}, on the other hand, is based on Theorem \ref{thm:reg2} and thus yields the same bounds as can be achieved there (see \cite{Moshkovitz.2019}). Subsequent to the results of Chapter \ref{sec:vdischom}, the first author \cite{Terry.2022} produced an additional combinatorial argument to show that $L$ in Theorem \ref{thm:vc2finite} can be taken to be of the form $\e_1^{-O_k(1)}$.  

An immediate corollary of Theorem \ref{thm:vc2finite} is that any $\NIP_2$ hereditary $3$-graph property is  $\disc_{2,3}$-homogeneous.   We will show that this is a characterization (see page \pageref{proof:dischom} for the proof).

\begin{theorem}\label{thm:dischom}
Suppose $\calH$ is a hereditary $3$-graph property.  Then $\calH$ is $\disc_{2,3}$-homogeneous if and only if $\calH$ has bounded $\VC_2$-dimension. 
\end{theorem}

Combining \cite{Terry.2018} with the lower bound construction of Moshkovitz and Shapira \cite{Moshkovitz.2019}, we see that bounded $\VC_2$-dimension is also characterized in terms of the bounds in the regularity lemma.  

We now turn to $\vdisc_3$-homogeneity and Problem \ref{prob:hom} (1).  We begin by stating prior results. Recall that the $\VC$-dimension of a graph $G$ is the largest $k$ such that $U(k)$ appears as an induced subgraph of the bipartite graph $\bip(G)$. A simple way to extend this notion to the setting of $3$-graphs is via the $\VC$-dimension of an auxiliary bipartite graph associated to a $3$-graph.

\begin{definition}[Graph associated to a 3-graph]\label{def:assgraph}
Suppose $H=(V,F)$ is a $3$-graph.  The \emph{graph associated to $H$} is $\graph(H):=(V\cup {V\choose 2}, E)$, where $E=\{ae: \{a\}\cup e\in F\}$.
\end{definition}

The \emph{$\VC$-dimension of $H$}, $\VC(H)$, is then defined as the $\VC$-dimension of $\graph(H)$.  We say a hereditary $3$-graph property has \emph{unbounded $\VC$-dimension} (or $\IP$) if for all $k$, there is $H\in \calH$ with $\VC(H)\geq k$.  Otherwise we say $\calH$ has \emph{bounded $\VC$-dimension}\footnote{For the model theorist, $\calH$ has bounded $\VC$-dimension if and only if the edge relation is an NIP formula for every partition of the variables, in every model of $T_{\calH}$.} (or $\NIP$), and let $\VC(\calH)=\max\{\VC(H):H\in \calH\}$.   In \cite{Fox.2017bfo}, Fox et al. showed that $3$-graph properties of bounded $\VC$-dimension are $\vdisc_3$-homogeneous.  In fact, they showed something much stronger, namely that any $k$-graph with $\VC$-dimension at most $d$ has an equitable $\e$-homogeneous partition with $O_{d,k}(\e^{-d})$ many parts (see \cite{Fox.2017bfo} for details).  Similar results with a weaker bound were also obtained by Chernikov and Starchenko in \cite{Chernikov.2016zb}.  Thus, having bounded $\VC$-dimension is sufficient for $\vdisc_3$-homogeneity.  However, it turns out that it is not necessary. 

\begin{proposition}\label{prop:weakhomex}
There is a hereditary $3$-graph property $\calH$ which is $\vdisc_3$-homogeneous but which has unbounded $\VC$-dimension.
\end{proposition}

The example used to prove Proposition \ref{prop:weakhomex} is extremely simple.  Concretely, we define the infinite $3$-graph
$$
\calW_1=(\{b_i,c_i: i\in \mathbb{N}\}\cup \{a_S: S\in \calP(\mathbb{N})\}, \{a_Sb_jc_k: j=k\in S\}),
$$
and then let $\calH_{\calW_1}$ contain all finite $3$-graphs isomorphic to a finite induced sub-$3$-graph of $\calW_1$.  The key feature of this example is that, while the edge relation has unbounded $\VC$-dimension in $\calW_1$, for any finite induced sub-$3$-graph $H$ of $\calW_1$, the copies of $U(k)$ in $\graph(H)$ are extremely sparse.  This contrasts with the following example of a 3-graph, in which copies of $U(k)$ occur more densely.

\begin{definition}[$\Ubar(k)$]\label{def:ubark}
Define $\Ubar(k)$ to be the $3$-partite $3$-graph 
\[\Ubar(k)=(\{a_1,\ldots, a_k\}\cup \{b_S: S\subseteq [k]\}\cup \{c_1,\ldots, c_k\},\{a_ib_Sc_j: j\in S\}).\]  
\end{definition}

One can think of $\Ubar(k)$ as arising from adjoining $k$ new vertices to $U(k)$ to create a $3$-graph.  We will show that $\vdisc_3$-homogeneity is characterized by whether or not there is a bound on the size of $k$ for which $\Ubar(k)\in \trip(\calH)$ (see page \pageref{proof:vdiscubark} for the proof).

\begin{theorem}\label{thm:vdiscubark}
$\calH$ is $\vdisc_3$-homogeneous if and only if there is $k\geq 1$ such that $\Ubar(k)\notin \trip(\calH)$. 
\end{theorem}

Note that if a hereditary $3$-graph property $\calH$ has $\VC$-dimension less than $k$, then $\Ubar(k)\notin \trip(\calH)$.  However, the converse of this is false, as witnessed by $\calH_{\calW_1}$.  

There are several ways a 3-graph $H=(V,F)$ can have large $\VC$-dimension, based on how elements of $V\cup {V\choose 2}$ are chosen to build a copy of $U(k)$ in $\graph(H)$.   We now single out a particular way that will be important for stating our results on $\vdisc_3$-homogeneity. 
 
 \begin{definition}[$U^*(k)$]\label{def:ustar}
 Given $k\geq 1$, define $U^*(k)$ to be the $3$-partite $3$-graph 
 \[U^*(k)=(\{a\}\cup \{b_S: S\subseteq [k]\}\cup \{c_1,\ldots, c_k\},\{ab_Sc_j: j\in S\}).\] 
 \end{definition}
 
 Note that $U^*(k)$ arises by adding a single new vertex to $U(k)$ to make a $3$-graph.  It is easy to see that if $U^*(k)\in \trip(\calH)$, then $\calH$ has $\VC$-dimension at least $k$.  On the other hand, it is not hard to show that even though the example $\calH_{\calW_1}$ above has unbounded $\VC$-dimension, $U^*(k)\notin \trip(\calH_{\calW_1})$ for any $k>1$. Thus, not every $3$-graph $H$ with large $\VC$-dimension contains a copy of  $U^*(k)$ in $\trip(H)$.  

\begin{definition}[Slicewise $\VC$-dimension]\label{def:ustark}
Suppose $H$ is a finite $3$-graph. The \emph{slicewise $\VC$-dimension of $H$}\footnote{In previous versions of this paper, this was referred to as \emph{weak $\VC$-dimension}. We have here adopted terminology introduced in \cite{Chernikov.2020} to avoid excessive use of the adjective ``weak".}
 is defined to be 
 \[\SVC(H):=\max\{k: \text{$U^*(k)$ is an induced sub-$3$-graph of $\trip(H)$}\}.\] 
\end{definition}
In other words, the slicewise $\VC$-dimension of $H=(V,F)$ is the largest $k$ for which there exist $a, c_1,\ldots, c_k\in V$, and $b_S\in V$ for each $S\subseteq[k]$, so that $ac_ib_S\in F$ if and only if $i\in S$.  Given $a\in V$, define the \emph{slice graph at $a$} to be $H_a:=(V,N(a))$, where $N(a)=\{bc\in {V\choose 2}: abc\in F\}$. Then the slicewise $\VC$-dimension of $H$ is also the maximum $\VC$-dimension of the graphs $H_a$, as $a$ ranges over all elements in $V$.  

We say that a hereditary $3$-graph property $\calH$ has \emph{bounded $\SVC$-dimension}\footnote{For the model theorist, $\calH$ has bounded $\SVC$-dimension if and only if for every $M\models T_{\calH}$, and every $a\in M$, $R(a,x,y)$ is an NIP formula in $\calM$.}  (or is SNIP) if there is some $k$ so that $\SVC(H)\leq k$ for all $H\in \calH$.  In this case, we define the \emph{$\SVC$-dimension of $\calH$} to be  $\SVC(\calH)=\max\{\SVC(H): H\in \calH\}$.  Otherwise, we say $\calH$ has \emph{unbounded slicewise $\VC$-dimension} (or has SIP), and write $\SVC(\calH)=\infty$. 

It turns out that a hereditary $3$-graph property is close to a SNIP property if and only if $\Ubar(k)\notin \trip(\calH)$ for some $k\geq 1$ (see Section \ref{subsec:closenip}).

\begin{theorem}\label{thm:simclasswnip}
Suppose $\calH$ is a hereditary $3$-graph property.  Then $\calH$ is close to some $\calH'$ with  $\SVC(\calH')<\infty$ if and only if there is $k$ such that $\Ubar(k)\notin \trip(\calH)$.
\end{theorem}

This provides us with another way of identifying the $\sim$-classes which are $\vdisc_3$-homogeneous.

\begin{theorem}\label{thm:vdischom} 
The following are equivalent.
\begin{enumerate}[label=\normalfont(\arabic*)]
\item $\calH$ is $\vdisc_3$-homogeneous.
\item $\trip(\calH)$ omits $\Ubar(k)$ for some $k$.
\item There is a property $\calH'$ with $\SVC(\calH')<\infty$ such that $\calH\sim \calH'$. 
\end{enumerate}
\end{theorem}

An interesting open question is to determine the optimal bounds corresponding to regularity lemmas for $\vdisc_3$-homogeneous properties.  Specifically, Fox et al. \cite{Fox.2017bfo} obtained polynomial bounds under the assumption of bounded $\VC$-dimension.  Can one still obtain such bounds if one only assumes bounded $\SVC$-dimension? 

In \cite{Chernikov.2020}, Chernikov and Towsner produced results similar to Theorem \ref{thm:vdischom} via a different proof which uses infinitary techniques, as well as related results for $k$-uniform hypergraphs satisfying various generalizations of the notion of having bounded slicewise $\VC$-dimension.

We now turn to Problem \ref{prob:triads}, which asks about irregular triples/triads. We begin with Problem \ref{prob:triads} (1) concerning zero $\vdisc_3$-error.  Previous work has shown that this problem is related to stability in $3$-graphs, as defined below. 

\begin{definition}[Stability of a 3-graph property]\label{def:stab3}
Given a $3$-graph $H=(V,F)$, we say $H$ is \emph{$k$-stable} if $\graph(H)$ is $k$-stable. 

We say a hereditary $3$-graph property $\calH$ is \emph{stable}\footnote{For the model theorist, $\calH$ is stable if and only if the edge relation is stable, for any bipartition of the variables, in every model of $T_{\calH}$.} if there is $k$ such that every $H\in \calH$ is $k$-stable.
\end{definition}

It was shown in \cite{Ackerman.2017} that $k$-stable $3$-graphs have $\vdisc_3$-regular decompositions with no irregular triples. Similar  results were also obtained in \cite{Chernikov.2016zb} without equitability conditions on the partition.  The results in \cite{Ackerman.2017} directly imply the following.

\begin{theorem}[Stable 3-graphs admit zero $\vdisc_3$-error, \cite{Ackerman.2017}]\label{thm:AFP}
Suppose $\calH$ is a hereditary $3$-graph property.  If $\calH$ is stable, then $\calH$ admits zero $\vdisc_3$-error.
\end{theorem}

However, it turns out that stability is too strong to give a characterization, as we shall show by way of the following proposition.

\begin{proposition}\label{prop:fullweakex}
There is a hereditary $3$-graph property which admits zero $\vdisc_3$-error but which is not stable.
\end{proposition}

The example used to prove Proposition \ref{prop:fullweakex} is analogous to $\calH_{\calW_1}$.  Concretely, we define the infinite $3$-graph
$$
\calW_2=(\{a_i,b_i,c_i: i\in \mathbb{N}\}, \{a_ib_jc_k: i\leq j\text{ and }j=k\}),
$$
and then let $\calH_{\calW_2}$ be the class of all finite $3$-graphs isomorphic to an induced sub-$3$-graph of $\calW_2$.    The important quality of this example is that while the edge relation is unstable in $\calW_2$, given any finite induced sub-$3$-graph $H$ of $\calW_2$, the copies of $H(k)$ in $\graph(H)$ are extremely sparse.  The following is a simple example of a $3$-graph where half-graphs occur more densely.

\begin{definition}[$\Hbar(k)$]\label{def:hbark}
Given $k\geq 1$, define $\Hbar(k)$ to be the $3$-graph 
\[\Hbar(k)=(\{a_i,b_i,c_i: i\in [k]\},\{a_ib_jc_k: j\leq k\}).\]
\end{definition}

One can think of $\Hbar(k)$ as constructed by adjoining $k$ new vertices to $H(k)$ to obtain a $3$-graph. We prove that the presence of these $3$-graphs is sufficient to imply that a property requires non-zero $\vdisc_3$-error. 

\begin{theorem}\label{thm:vdischbark}
If $\Hbar(k)\in \trip(\calH)$  for every $k$, then $\calH$ requires non-zero $\vdisc_3$-error.
\end{theorem}

We conjecture that the converse is true as well. 

\begin{conjecture}\label{conj:weaktriads}
If there is some $k\geq 1$ such that $\Hbar(k)\notin \trip(\calH)$, then $\calH$ admits zero $\vdisc_3$-error.
\end{conjecture}

We will also give an equivalent formulation of the conjecture by characterizing the $\sim$-classes $[\calH]$ with the property that $\Hbar(k)\notin \trip(\calH)$ for some $k$.  This characterization centers on the $3$-graph obtained by adjoining a single vertex to $H(k)$ to obtain a $3$-graph.  

\begin{definition}[$H^*(k)$]\label{def:hstar}
Given $k\geq 1$, define $H^*(k)$ to be the $3$-graph 
\[H^*(k)=(\{a\}\cup \{b_i,c_i: i\in [k]\},\{ab_ic_j: i\leq j\}).\]
\end{definition}

Note that if $H^*(k)$ appears as an induced sub-$3$-graph of $G$, then $G$ is not $k$-stable. However, not every $G$ with the $k$-order property has an induced copy of $H^*(k)$ in $\trip(G)$ (the $3$-graph $\calW_2$ used to prove  Proposition \ref{prop:fullweakex} is such an example).  We now define the notion of \emph{slicewise stability}, in analogy to slicewise $\VC$-dimension.

\begin{definition}[Slicewise stability]\label{def:weaklystable}
A $3$-graph $H=(V,F)$ is \emph{slicewise $k$-stable}\footnote{In previous versions of this paper, this notion was referred to as \emph{weakly stable}. Here we have adopted terminology introduced in \cite{Chernikov.2020} to avoid excessive use of the adjective ``weak".} if and only if $H^*(k)$ is not an induced sub-$3$-graph of $ \trip(H)$.

Note that a $3$-graph $H=(V,F)$ is slicewise $k$-stable if and only if the slicegraphs $H_a$ are stable for all $a\in V$.  A hereditary $3$-graph property $\calH$ is \emph{slicewise stable}\footnote{For the model theorist, $\calH$ is slicewise stable if and only if for all $\calM\models T_{\calH}$ and $a\in M$, $R(a,x,y)$ is a stable formula in $\calM$.} if there is some $k$ such that every $H\in \calH$ is slicewise $k$-stable.
\end{definition}

In other words, a $3$-graph $H=(V,F)$ is slicewise $k$-stable if for all $a\in V$, the graph $H_a:=(V,N(a))$ is $k$-stable, where $N(a)=\{bc: abc\in F\}$.   It turns out that slicewise stability characterizes the $\sim$-class $[\calH]$ such that $\trip(\calH)$ omits $\Hbar(k)$ for some $k$ (see Section \ref{subsec:closenip}).

\begin{theorem}\label{thm:simclassws} Suppose $\calH$ is a hereditary $3$-graph property.  The following are equivalent.
\begin{enumerate}[label=\normalfont(\arabic*)]
\item There is some $k\geq 1$ such that $\Hbar(k)\notin \trip(\calH)$. 
\item There is some slicewise stable hereditary $3$-graph property $\calH'$ such that $\calH\sim \calH'$. 
\end{enumerate}
\end{theorem}

Thus the following (a priori weaker) conjecture is in fact equivalent to Conjecture \ref{conj:weaktriads}.

\begin{conjecture}\label{conj:weaktriads2}
If $\calH$ is slicewise stable, then $\calH$ admits zero $\vdisc_3$-error.
\end{conjecture}

We will prove something weaker than Conjecture \ref{conj:weaktriads} in Section \ref{subsec:wsbinary}, namely that a property which is close to a slicewise stable property admits binary $\vdisc_3$-error.

\begin{proposition}\label{prop:wsbinary}
Suppose $\calH$ is close to some slicewise stable $\calH'$.  Then $\calH$ admits binary $\vdisc_3$-error.
\end{proposition}

We also show in Proposition \ref{prop:slvdiscbin} that the converse to Proposition \ref{prop:wsbinary} is false.

One tool we will use in the proof of Proposition \ref{prop:wsbinary} is the following set of conditions that are equivalent to having zero and binary $\vdisc_3$-error, respectively.  

\begin{proposition}\label{prop:wnipreduction}
Suppose $\calH$ is a hereditary $3$-graph property.  
\begin{enumerate}[label=\normalfont(\arabic*)]
\item $\calH$ admits zero $\vdisc_3$-error if and only if $\calH$ is close to a $\SNIP$ property and admits zero $\disc_{2,3}$-error.
\item $\calH$ admits binary $\vdisc_3$-error if and only if $\calH$ is close to a $\SNIP$ property and admits binary $\disc_{2,3}$-error.
\end{enumerate}
\end{proposition}

The proof of Proposition \ref{prop:wnipreduction} (which begins on page \pageref{proof:wnipreduction}) uses Theorem \ref{thm:vdischom} alongside some standard techniques from extremal combinatorics.  To prove Proposition \ref{prop:wsbinary}, we use Proposition \ref{prop:wnipreduction}, and ideas from the proof of the stable graph regularity lemma of \cite{Malliaris.2014}, as well as the strong graph regularity lemma (see e.g. \cite{Rodl.2007,  Lovasz.2007,  Alon.2000}).  We believe this argument provides useful techniques for showing that a given hereditary 3-graph property admits binary error.

Note that the combination of Proposition \ref{prop:wnipreduction} and Theorem \ref{thm:vdischom} immediately implies that any property admitting zero or binary $\vdisc_3$-error is also $\vdisc_3$-homogeneous.

\begin{corollary}\label{cor:wnipcor}
If $\calH$ admits zero or binary $\vdisc_3$-error, then $\calH$ is $\vdisc_3$-homogeneous.
\end{corollary}

Proposition \ref{prop:wnipreduction} will also play a role in showing that Problems \ref{prob:triads} (1) and (2) are in fact the same question.  In particular, a property admits zero $\disc_{2,3}$-error if and only if it admits zero $\vdisc_3$-error.

\begin{theorem}\label{thm:vdisceq}
$\calH$ admits zero $\disc_{2,3}$-error if and only if $\calH$ admits zero $\vdisc_3$-error.
\end{theorem}

In addition to Proposition \ref{prop:wnipreduction}, the proof of Theorem \ref{thm:vdisceq} will use Theorems \ref{thm:vdischom} and \ref{thm:vdischbark} (see page \pageref{thm:vdischom}). 

We now turn to part (4) of Problem \ref{prob:triads}, which asks for a characterization of the hereditary $3$-graph properties admitting linear $\disc_{2,3}$-error.  The motivation for this question, as well as the definition of linear error, stems from \cite{Terry.2021a}, where such decompositions occur naturally in the ternary sum graphs associated to certain quadratically structured subsets of $\mathbb{F}_p^n$. As in that paper, we show that such decompositions are related to the following combinatorial notion, first defined by the authors here and in \cite[Definition 1.10]{Terry.2021a}.

\begin{definition}[$F(\ell)$]\label{def:fell}
Given $\ell\geq 1$, let $F(\ell)$ be the $3$-graph 
\[F(\ell)=(\{a^f_i,b_i, c_i: i\in [\ell], f:[\ell]^2\rightarrow [\ell]\}, \{a^f_ib_jc_k: k\leq f(i,j)\}).\]
\end{definition}

\begin{definition}[Ternary functional order property ($\FOP_2$)]\label{def:fop}
Given $\ell\geq 1$ and a $3$-graph $H$, we say $H$ \emph{has the $\ell$-functional order property of dimension 3} (or $\ell$-$\FOP_2$) if there is an induced copy of $F(\ell)$ in $\trip(H)$.
\end{definition}
In other words, $H=(V,E)$ has $\ell$-$\FOP_2$ if and only if there exist $b_1,\ldots, b_{\ell}$, $c_1,\ldots, c_{\ell}\in V$ and for every function $f:[\ell]^2\rightarrow [\ell]$, vertices $a_1^f,\ldots, a_{\ell}^f$, such that $a^f_ib_jc_k\in E$ if and only if $k\leq f(i,j)$.  We then say that a hereditary $3$-graph property $\calH$ \emph{has $\FOP_2$} if for all $\ell\geq 1$, there is $H\in \calH$ with $\ell$-$\FOP_2$.  Otherwise, we say $\calH$ is $\NFOP_2$. 

It is not difficult to see (see \cite[Lemma 5.11]{Terry.2021a}) that if $\calH$ has $\FOP_2$ then it has unbounded $\VC$-dimension.  It is further clear from Fact \ref{fact:vc2universal}  that if $\calH$ has unbounded $\VC_2$-dimension, then it has $\FOP_2$.  We will show this implication is strict in Appendix \ref{app:fopgen}.  For further properties of general $\NFOP_2$ formulas, including a proof that $\NFOP_2$ is closed under finite boolean combinations, we refer the reader to Appendix \ref{app:fopgen}.

In \cite{Terry.2021a} we show that given a set $A\subseteq \F_p^n$, if the $3$-graph $\Gamma_A=(\F_p^n, \{xyz: x+y+z\in A\})$ does not have $\ell$-$\FOP_2$, then there is a decomposition of $\F_p^n$ defined by a high rank quadratic factor, so that on almost all parts in the decomposition, the set $A$ has density near $0$ or $1$.  Moreover, the parts of this decomposition on which the density of $A$ is bounded away from 0 and 1 are contained in a small number of linearly defined pieces \cite[Theorem 1.12]{Terry.2021a}, hence the term ``linear error''.  This structure theorem can be translated directly into the statement that there exists a regular decomposition for $\Gamma_A$ with linear error, in the sense of Definition \ref{def:error} (see \cite[Corollary 3.30]{Terry.2021a}).  Here, we prove an analogous result for $3$-graphs, namely that if $H$ is a large $3$-graph that does not have $\ell$-$\FOP_2$, then it has a $\disc_{2,3}$-regular decomposition with linear error.

\begin{theorem}\label{thm:FOPfinite}
For all $k\geq 1$, $\e_1>0$, $\e_2:\mathbb{N}\rightarrow (0,1]$, and $t_0,\ell_0\geq 1$, there are $T,L,N$ such that the following hold.  For all $H=(V,E)$ with $|V|\geq N$ such that $H$ does not have $k$-$\FOP_2$, there are $t_0\leq t\leq T$, $\ell_0\leq \ell \leq L$, and $\calP$ a $(t,\ell,\e_1,\e_2(\ell))$ decomposition of $V$ which is $(\e_1,\e_2(\ell))$-regular and which has linear $\disc_{2,3}$-error with respect to $H$.
\end{theorem}

This result will be proved from page \pageref{proof:FOPfinite} onwards.  Theorem \ref{thm:FOPfinite} implies, by definition, that any $\NFOP_2$ hereditary $3$-graph property admits linear $\disc_{2,3}$-error.  We will show that $\NFOP_2$ hereditary $3$-graph properties are in fact characterized by the admission of  linear $\disc_{2,3}$-error.

\begin{theorem}\label{thm:FOP}
Suppose $\calH$ is a hereditary $3$-graph property.  Then $\calH$ is $\NFOP_2$ if and only if $\calH$ admits linear $\disc_{2,3}$-error.
\end{theorem}

The proof of Theorem \ref{thm:FOP} begins on page \pageref{proof:FOP}.  We finally arrive at parts (2) and (3) of Problem \ref{prob:triads}, which ask about $3$-graph properties admitting binary $\vdisc_3$- and binary $\disc_{2,3}$-error, respectively.  For this problem, we will give two examples of properties which require binary $\disc_{2,3}$-error and binary $\vdisc_3$-error.  The first example is the following.

\begin{definition}[$\HP(k)$]\label{def:pk}
Given $k\geq 1$, define $\HP(k)$ to be the $3$-partite $3$-graph 
\[\HP(k)=(\{a_i,b_i,c_i: i\in [k]\},\{a_ub_vc_w: u+v+w\geq k+2\}).\]
\end{definition}

The $3$-graph $\HP(k)$ is a natural higher-order analogue of the half graph.  Indeed, observe that $H(k)$ is isomorphic to the graph with vertex set $\{a_i,b_i: i\in [k]\}$ and edge set $\{a_ub_v: u+v\geq k+1\}$.  Based on this analogy, it is natural to make the following definition.

 \begin{definition}[Ternary hyperplane order property ($\HOP_2$)]\label{def:hop3}
We say a $3$-graph $H$ has the \emph{$k$-hyperplane order property of dimension 3} (or $k$-$\HOP_2$) if $\HP(k)$ is an induced sub-$3$-graph of $\trip(H)$. 
 
A hereditary $3$-graph property $\calH$ has $k$-$\HOP_2$ if there is $H\in \calH$ which has $k$-$\HOP_2$.  We say that $\calH$ has \emph{$\HOP_2$} if it has $k$-$\HOP_2$ for all $k\geq 1$.  Otherwise, we say $\calH$ has $\HOP_2$.
 \end{definition}

Besides the obvious geometric reasons for viewing $\HOP_2$ as a ternary analogue of a half-graph, there are algebraic reasons why Definition \ref{def:pk} might be considered a natural ternary analogue of the order property. Specifically, the parts of the partition approximating an $\NFOP_2$-subset of $\F_p^n$ established in \cite{Terry.2021a} do not have $2$-$\HOP_2$, just as the parts of the partition approximating a stable subset of $\F_p^n$ in \cite{Terry.2019} do not have the $2$-order property. We refer the reader to \cite[Section 5]{Terry.2021a} for an extensive discussion of this. In the hypergraph setting, one compelling piece of evidence is the fact that if a hereditary $3$-graph property has $\HOP_2$, then it requires non-binary $\disc_{2,3}$-error. 

\begin{theorem}\label{thm:ternarytriads}
Suppose $\calH$ is a hereditary $3$-graph property with $\HOP_2$.  Then $\calH$ requires non-binary $\disc_{2,3}$-error. 
\end{theorem}

However, the authors' work on analogous questions in the arithmetic setting \cite{Terry.2021a} identified a second, distinct example whose presence requires non-binary $\disc_{2,3}$-error.  The example is defined as a disjoint union of cosets of a nested chain of subgroups in an elementary abelian $p$-group.

\begin{definition}[Green-Sanders examples]\label{def:gsexamples}
Suppose $p\geq 3$.  Given $n\geq 1$, define
\[ A(p,n):=\{\xbar\in \mathbb{F}_p^n: \exists \; i \in [n] \text{ such that }x_i=1\text{ and for all }1\leq j<i, x_j=0\}.\]
The \emph{Green-Sanders $p$-example of dimension $n$} is the pair $(\mathbb{F}_p^n, A(p,n))$.
 \end{definition}

In the case when $p=3$, Green and Sanders \cite{Green.2015qy4} showed that every subspace of $\mathbb{F}_3^n$ has a non-trivial Fourier coefficient with respect to $A(3,n)$, and in \cite{Terry.2021a}, the authors showed that for any prime $p\geq 3$, the sets $A(p,n)\subseteq \mathbb{F}_p^n$ cannot be written as a union of quadratic atoms with ``no error'' (see \cite[Proposition 5.21]{Terry.2021a}).  We refer the reader to our forthcoming note \cite{Terry.2021c} for further examples with similar properties.  The sets $A(p,n)$ of Definition \ref{def:gsexamples} naturally give rise to $3$-partite $3$-graphs as follows.

\begin{definition}[$\GS_{p}(n)$]\label{def:gspi}
We define $\GS_{p}(n)$ to be the $3$-graph 
\[\GS_{p}(n)=(\{a_g, b_g,c_g: g\in \mathbb{F}_p^n\},\{a_gb_{g'}c_{g''}: g+g'+g''\in A(p,n)\}).\] 
\end{definition}

We prove that these examples also require non-binary $\disc_{2,3}$-error.

\begin{theorem}\label{thm:ternarytriads1}
If a hereditary $3$-graph property $\calH$ has $\GS_p(n)\subseteq \trip(\calH)$ for every $n$, then $\calH$ requires non-binary $\disc_{2,3}$-error. 
\end{theorem}

The authors show in \cite[Appendix A.3]{Terry.2021a} that $\GS_3(n)$ does not have $4$-$\HOP_2$, and thus, $\NHOP_2$ cannot characterize the properties which admit binary $\disc_{2,3}$-error.  Instead, it seems that the hypotheses of Theorems \ref{thm:ternarytriads} and \ref{thm:ternarytriads1} should be distinct instantiations of another, more general tameness condition (see Conjecture \ref{conj:x}).

Let $\calH_{\HP}$ be the hereditary $3$-graph obtained by closing $\{\HP(k):k\geq 1\}$ under isomorphism and induced sub-3-graphs, and let $\calH_{\GS_p}$ be the hereditary $3$-graph property obtained by closing $\{\GS_p(n):n\geq 1\}$ under isomorphism and induced sub-$3$-graphs. Theorems \ref{thm:ternarytriads} and \ref{thm:ternarytriads1} tell us that if $\calH_{\HP}$ or $\calH_{\GS_p}$ are contained in $\trip(\calH)$, then $\calH$ requires binary $\disc_{2,3}$-error.   Interestingly, both these examples have bounded (in fact, very small) $\VC$-dimension.

\begin{proposition}\label{prop:finitevcnotdisc3tame}
$\calH_{\HP}$ has $\VC$-dimension $1$ and for all $p\geq 3$, $\calH_{\GS_p}$ has $\VC$-dimension at most $3$.
\end{proposition}

We will give proofs of these bounds in Chapter \ref{sec:binary}.   Proposition \ref{prop:finitevcnotdisc3tame} also plays a role in our proofs of Theorems \ref{thm:ternarytriads} and \ref{thm:ternarytriads1}.  In particular, since properties with finite $\VC$-dimension also have finite $\SVC$-dimension, both $\calH_{\GS_p}$ and $\calH_{\HP}$ are SNIP.  Thus by Proposition \ref{prop:wnipreduction}, in order to show that they require non-binary $\disc_{2,3}$-error, it suffices to show that both require non-binary $\vdisc_3$-error.  We will do this using a general sufficient condition satisfied by both examples (see Section \ref{subsec:special}).   

Proposition \ref{prop:finitevcnotdisc3tame} shows that a hereditary $3$-graph property can have bounded $\VC$-dimension while also requiring non-binary $\disc_{2,3}$-error.  On the other hand, we will show that there are properties admitting binary $\disc_{2,3}$-error which not only have  unbounded $\VC$-dimension, but which contain $\Ubar(k)$ for all $k$.

\begin{proposition}\label{prop:disc3vsVC}
There is a hereditary $3$-graph property $\calH$ which admits binary $\disc_{2,3}$-error, and which contains $\Ubar(k)$ for all $k$. 
\end{proposition}

Thus, whether or not a property admits binary $\disc_{2,3}$-error has no general relationship to whether or not it has bounded $\VC$-dimension. This might seem surprising at first, but it reflects the fact that that hypergraph regularity is really about ternary structure, while $\VC$-dimension is closely related to binary structure.  Proposition \ref{prop:disc3vsVC} also shows that a property can admit binary $\disc_{2,3}$-error while not admitting binary $\vdisc_3$-error.  Indeed, by Theorem \ref{thm:vdischom}, if a property contains $\Ubar(k)$ for all $k$, it is not $\vdisc_3$-homogeneous, and therefore cannot admit binary $\vdisc_3$-error.  Combining this with Proposition \ref{prop:wnipreduction}, we see that admitting binary $\vdisc_3$-error is strictly stronger than admitting binary $\disc_{2,3}$-error.

 In contrast to this, there is a clear implication between stability and binary $\disc_{2,3}$-error.  In particular, if $\calH$ is stable, then by Theorem \ref{thm:AFP},  $\calH$ admits zero $\vdisc_3$-error, so by Proposition \ref{prop:wnipreduction}, it admits zero $\disc_{2,3}$-error.   If Conjecture \ref{conj:weaktriads2} were true, then a similar argument would show that any slicewise stable $\calH$ admits zero $\disc_{2,3}$-error, and hence binary $\disc_{2,3}$-error.   We show directly in this paper (Proposition \ref{prop:wsbinary}) that any slicewise stable property admits binary $\disc_{2,3}$-error (in fact, binary $\vdisc_3$-error). 
 
 Finally, our results imply that admitting any good control over the $\disc_{2,3}$-irregular triads implies $\disc_{2,3}$-homogeneity (see Chapter \ref{sec:binary}).
 
  \begin{corollary}\label{cor:disc23hom}
 If $\calH$ admits zero, binary, or linear $\disc_{2,3}$-error, then $\calH$ is $\disc_{2,3}$-homogeneous.
 \end{corollary}
 
 We end the introduction with a discussion of some open problems. As mentioned above, we conjecture that $\calH_{\HP}$ and $\calH_{\GS_p}$ are two manifestations of an underlying combinatorial definition, which will characterize the properties admitting binary $\disc_{2,3}$-error. 
\begin{conjecture}\label{conj:x}
There exists a definition $\mathrm{XOP}_2$ (for ``Unknown Order Property") which is defined by an infinite scheme of existential sentences in the language $\calL=\{R(x,y,z)\}$ such that a hereditary $3$-graph property admits binary $\disc_{2,3}$-error if and only if it does not satisfy $\mathrm{XOP}_2$.  
\end{conjecture}

We also conjecture a connection between this property and quadratic structure in the arithmetic setting.  For a discussion of this, we direct the reader to \cite[Section 1.3]{Terry.2021a}.

In addition to the interesting problem above for $3$-graphs, several open questions arise form this work with regard to $k$-graphs for $k>3$. For example, the results of Fox et al. \cite{Fox.2017bfo} give strong regularity lemmas for $k$-graphs with bounded $\VC$-dimension.  Our work suggests that an analogue of slicewise $\VC$-dimension should characterize hereditary $k$-graph properties with such partitions.  Similarly, while Theorem \ref{thm:AFP} was proven in the setting of $k$-graphs for $k\geq 2$, it remains open to give a characterization of which hereditary $k$-graph properties admit regular partitions as described there.  It is fairly clear that bounded $\VC_k$-dimension (see \cite{Chernikov.2019} for the definition) will correspond to the appropriate higher order analogue of $\disc_{2,3}$-homogeneity.  Indeed, all the same tools used in Chapter \ref{sec:dischom} have their analogues for $k$-graphs.  One direction of this has been proved, non-quantitatively, by Chernikov and Towsner \cite{Chernikov.2020}.

It is also straightforward to define a higher-dimensional analogue of the ternary functional order property.

\begin{definition}[$\FOP_k$]\label{def:FOPk}
Given a $(k+1)$-graph $H=(V,E)$, we say $H$ \emph{has $\ell$-$\FOP_k$} if there are $\{x^u_j: j\in [\ell], u\in [k]\}\subseteq V$ and for every function $f:[\ell]^{k}\rightarrow [\ell]$, there are $z_{1}^f,\ldots, z_{\ell}^f\in V$ such that $z_j^fx_{i_1}^1\ldots x_{i_k}^k\in E$ if and only if $j< f(i_1,\ldots, i_{k})$. 

We say a hereditary $(k+1)$-graph property $\calH$ \emph{has $\FOP_k$} if for all $\ell\geq 1$, there is $H\in \calH$ with $\ell$-$\FOP_k$.  Otherwise, we say $\calH$ is $\NFOP_k$. 
\end{definition}

We conjecture that $\NFOP_k$ should have a characterization analogous to Theorem \ref{thm:FOP}.  There are also natural $k$-dimensional generalizations of the $\HP(n)$ and $\GS_p(n)$.  

\begin{definition}\label{def:hopk}
Given $k,n\geq 1$, define
$$
\HP_k(n)=(\{a_i^j: i\in [n], j\in [k]\}, \{a_{i_1}^1\ldots a_{i_k}^k: i_1<i_2+\ldots +i_k\}.
$$
Let $\calH_{\HP_k}$ be the closure of $\{\HP_k(n):n\geq 1\}$ under isomorphism and induced sub-$k$-graph.
\end{definition}

\begin{definition}[$\GS_p^{(k)}(n)$]\label{def:gskprop}
Suppose $\ell,k\geq 1$ and $p\geq 3$ is a prime.  Define 
$$
\GS_p^{(k)}(n)=(\{a^i_g: i\in [k], g\in \mathbb{F}_p^n\},\{a^1_{g_1}\ldots a^k_{g_k}:g_1+\ldots +g_k\in A(p, n)\}).
$$  
Let $\calH_{\GS_p^{(k)}}$ be the closure of $\{\GS_p^{(k)}(n):n\geq 1\}$ under isomorphism and induced sub-$k$-graph. 
\end{definition}

We conjecture that both $\calH_{\HP_k}$ and  $\calH_{\GS^{(k)}_p}$ fail to admit the appropriate generalization of binary $\disc_{2,3}$-error.\label{higherconj}  Clearly this also leads to a higher order version of Conjecture \ref{conj:x}.  

\section{Outline of paper}\label{subsec:outline}

We now give a brief outline of what we do in each section.  In Chapter \ref{sec:prelims}, we set up notation, and give further background on regularity lemmas for $3$-graphs, counting lemmas, and other tools we use throughout the paper.  In Chapter \ref{sec:dischom} we prove Theorem \ref{thm:dischom}, the equivalence between $\disc_{2,3}$-homogeneity and finite $\VC_2$-dimension.  In Chapter \ref{sec:vdischom}, we prove our results about $\vdisc_3$-homogeneity, slicewise $\VC$-dimension, and slicewise stability, including Theorems \ref{thm:vdiscubark}, \ref{thm:vdischom}, \ref{thm:simclasswnip}, \ref{thm:vdischbark}, \ref{thm:simclassws} and Propositions \ref{prop:weakhomex}, \ref{prop:fullweakex},\ref{prop:wsbinary}, \ref{prop:wnipreduction}, and Corollary \ref{cor:wnipcor}.  Chapter \ref{sec:vdischom} also contains strong and generalized versions of the regularity lemma for stable graphs (see Theorems \ref{thm:functionstablereg} and \ref{thm:goodstrong}).  In Chapter \ref{sec:fop} we prove our results on linear error and $\FOP_2$, including Theorem \ref{thm:FOP}.  In this chapter there are also several general results about $\FOP_2$-formulas and  recharacterizations of both $\FOP_2$ and $\IP_2$ hereditary $3$-graph properties. Chapter \ref{sec:fop} also contains a strong removal lemma for stable graphs, with an accompanying structure theorem for graphs which contain few half-graphs (see Theorem \ref{thm:stableremoval} and Lemma \ref{lem:closetostableinWweak}).  In Chapter \ref{sec:binary}, we prove our main results regarding examples of properties requiring binary error, including Theorems \ref{thm:ternarytriads} and \ref{thm:ternarytriads1}.  It is in this chapter that we prove Theorem \ref{thm:disc3}, which recharacterizes binary $\disc_{2,3}$-error in terms of no $\disc_3$-irregular triads.

\section{Summary of problems and updates}\label{subsec:problems} To conclude this chapter, we list below, for the convenience of the reader, the problems and conjectures posed in this paper, alongside pointers to work towards them that has appeared since a version of this manuscript was first posted on the arXiv.

\begin{enumerate}
\item Characterize zero $\vdisc_3$-error or, equivalently by Theorem \ref{thm:vdisceq}, zero $\disc_{2,3}$-error (see Problem \ref{prob:triads} and Conjecture \ref{conj:weaktriads}). Results in this paper related to this problem are Proposition \ref{prop:fullweakex}, Theorem \ref{thm:vdischbark}, Proposition \ref{prop:wsbinary}, Proposition \ref{prop:wnipreduction} (a), and Proposition \ref{prop:disc3vsVC}.  A conjectured solution is given in Conjecture \ref{conj:weaktriads} (and equivalently Conjecture \ref{conj:weaktriads2}). This conjecture was recently shown to be false by Chernikov and Towsner \cite[Theorem G]{Chernikov.2024}.
\item Characterize binary $\vdisc_3$-error (see Problem \ref{prob:triads}). Results in this paper related to this problem are Proposition \ref{prop:wsbinary}, Proposition \ref{prop:wnipreduction} (b), and Theorem \ref{prop:slvdiscbin}.
\item Characterize binary $\disc_{2,3}$-error (see Problem \ref{prob:triads} and Conjecture \ref{conj:x}). Results in this paper related to this problem are Proposition \ref{prop:wnipreduction} (b), Theorem \ref{thm:ternarytriads}, Theorem \ref{thm:ternarytriads1}, and Proposition \ref{prop:disc3vsVC}.
\item We conjecture that there is a characterisation of $\NFOP_k$ for $k\geq 3$ analogous to Theorem \ref{thm:FOP} (see the remark following Definition \ref{def:FOPk}). In particular, we conjecture that there are higher order analogues of Theorems \ref{thm:FOPfinite}, \ref{thm:ternarytriads} and \ref{thm:ternarytriads1}.
\item Following Theorem \ref{thm:vdischom}, we asked whether a quantitatively improved version of Theorem \ref{thm:vdischom} holds.  A substantial quantitative improvement was recently obtained by the first author in \cite{Terry.2024}, which shows that a double exponential bound can be obtained, and moreover, that the bound cannot be improved past a single exponential. The correct form of the bound remains open.
\end{enumerate}

%% file: chapter3.tex
\chapter{Preliminaries}\label{sec:prelims}

This chapter contains the necessary background on quasirandomness and regularity. We begin with some notation.  Much of the notation below is new, and some is repetition of notation fixed in earlier in the paper.  In this case we will endeavor to provide references to where related notation has already appeared.  

Given parameters $a$ and $b$, we will write $a\ll b$ to mean that $a$ sufficiently small compared to $b$.  Similarly, we write $b\gg a$ to mean $b$ is sufficiently large compared to $a$.  Given $\e>0$ and $d,r\in \mathbb{R}$, we will write $d=r\pm \e$ to mean $d\in (r-\e, r+\e)$.  Given a positive integer $n$, let $[n]:=\{1,\ldots, n\}$.  Given a set $X$ and $k\geq 1$, define
$$
{X\choose k}:=\{Y\subseteq X: |Y|=k\}.
$$
Some of the following notation will already have appeared in Chapter \ref{sec:mainresults}, but is repeated here for the convenience of the reader.  For sets $X,Y,Z$, we set
\begin{align*}
K_2[X,Y]&=\{\{x,y\}: x\in X, y\in Y, x\neq y\}\text{ and }\\
K_3[X,Y,Z]&=\{\{x,y,z\}: x\in X, y\in Y, z\in Z, x\neq y, y\neq z, x\neq z\}.
\end{align*}
When $x,y\in X$ are distinct, we will write $xy$ to denote the set $\{x,y\}$.  Similarly, if $x,y,z\in X$ are pairwise distinct, we write $xyz$ to denote the set $\{x,y,z\}$. 

A \emph{$k$-uniform hypergraph} is a pair $(V,E)$ where $V$ is a set of vertices and $E\subseteq {V\choose k}$ is a set of edges.  When $k=2$, a $k$-uniform hypergraph is a \emph{graph}.  To ease notation, throughout the paper, we will call a $k$-uniform hypergraph a \emph{$k$-graph}.  Given a $k$-graph $G$, $V(G)$ denotes the vertex set of $G$, $E(G)$ denotes its edge set, and we set $v(G):=|V(G)|$ and $e(G):=|E(G)|$.  We say a $k$-graph $G'=(V',E')$ is a \emph{sub-$k$-graph of $G$} if $V'\subseteq V$ and $E'\subseteq E$.  In this case we write $G'\subseteq G$.   Note that a sub-$k$-graph is not necessarily induced.  We say $G'$ is an \emph{induced sub-$k$-graph of $G$} if $V'\subseteq V$ and $E'=E\cap {V'\choose k}$.  We will occasionally call a sub-$k$-graph a \emph{non-induced sub-$k$-graph} to emphasize what is meant.  When $V''\subseteq V$, we define $G[V'']=(V'',E\cap {V''\choose k})$.  For a $k$-graph $G''$, an \emph{induced copy of $G''$ in $G$} is an induced sub-$k$-graph of $G$ isomorphic to $G''$.  A  \emph{non-induced copy of $G''$ in $G$} is a non-induced sub-$k$-graph of $G$ isomorphic to $G''$.   

Suppose $G=(V,E)$ is a $k$-graph.  If $E={V\choose k}$, then $G=(V,E)$ is called \emph{complete}.  For any $2\leq k\leq r$, we say $G=(V,E)$ is \emph{$r$-partite} if there is a partition $V=V_1\cup \ldots \cup V_r$ such that for all $e\in E$, and $i\in [r]$, $|e\cap V_i|\leq 1$.  In this case we write $G=(V_1\cup \ldots \cup V_r, E)$ to denote that $G$ is $r$-partite, with choice of partition $V=V_1\cup \ldots \cup V_r$.  When $r=2$, we say that $G$ is \emph{bipartite}. 

Suppose now that $G=(V,E)$ is a graph.  For $v\in V$, the \emph{neighborhood of $v$} is 
$$
N_G(v)=\{u\in V :uv\in E\}.
$$
Given sets $U,W\subseteq V$, set $E_G[U,W]:=E\cap K_2[U,W]$, and $e_G(U,W):=|E_G[U,W]|$.  The  \emph{density} between $U$ and $W$ is the number $d_G(U,W)\in [0,1]$ such that 
$$
e_G(U,W)=d_G(U,W)|U||W|.
$$
When $G$ is clear from context, we will occasionally drop the subscripts above.  We define $G[U,W]$ to be the graph $(U\cup W,E_G[U,W])$.  Similarly, if $U,W,Z\subseteq V$, 
$$
G[U,W,Z]:=(U\cup W\cup Z,E_G[U,W]\cup E_G[U,Z]\cup E_G[W,Z]).
$$

Suppose that $H=(V,F)$ is a $3$-graph.  For sets $X,Y,Z\subseteq V$, we define 
$$
E_H[X,Y,Z]=\{xyz\in F: x\in X, y\in Y, z\in Z\}.
$$
Given  $v\in V$, $N_H(v)=\{bc\in {V\choose 2}:vbc\in F\}$ is the \emph{neighborhood of $v$}, and given $uv\in {V\choose 2}$, $N_H(uv)=\{b\in V: uvb\in F\}$ is the \emph{neighborhood of $uv$}.

Let $G=(V,E)$ be a graph. For $p\geq 2$, define
\[K_p^{(2)}(G)=\left\{e\in {V\choose p}: e'\in E \text{ for all }e'\in {e\choose 2}\right\}.\]
In other words, $K_p^{(2)}(G)$ is the collection of complete $p$-element subgraphs of $G$ (recall this was previously defined for $p=3$ on  page \pageref{k3}).   We recall from Section \ref{subsec:quasi1} that given a $3$-graph $H=(V,F)$ on the same vertex set as $G$, then we say $G$ \emph{underlies} $H$ if $F\subseteq K_3^{(2)}(G)$.  In other words, $G$ underlies $H$ if every edge of $H$ sits atop a triangle of $G$.  Similarly, we let
$$
K_p^{(3)}(H)=\left\{e\in {V\choose p}: f\in F \text{ for all }f\in {e\choose 3}\right\},
$$
i.e. $K_p^{(3)}(H)$ is the collection of complete $p$-element sub-$3$-graphs of $H$.

If $\calH$ is a class of finite $k$-graphs, we say $\calH$ is \emph{closed under induced sub-$k$-graphs} if for all $H=(V,F)\in \calH$ and $V'\subseteq V$, $H[V']:=(V',H\cap {V\choose k})\in \calH$.

\begin{definition}
Given $k\geq 2$, a \emph{hereditary $k$-graph property} is a class of finite $k$-graphs which is closed under isomorphism and induced sub-$k$-graphs.
\end{definition}

When $H,H'$ are $k$-graphs, we say that $H$ \emph{omits $H'$} or is \emph{$H'$-free} if no induced sub-$k$-graph of $H$ is isomorphic to $H'$.  Given a class of $k$-graphs $\calF$, let $\Forb(\calF)$ denote the class of finite $k$-graphs which omit every element of $\calF$.  Every class of the form $\Forb(\calF)$ is a hereditary $k$-graph property.  Conversely, it is well known that for every hereditary $k$-graph property $\calH$, there is a class $\calF_\calH$ of finite $k$-graphs such that $\calH=\Forb(\calF_\calH)$.

Recall that given a graph $G=(V,E)$, let $\bip(G)$ denote the bipartite graph $(U\cup W, E')$, where $U=\{u_v:v\in V\}$, $W=\{w_v:v\in V\}$, and $E'=\{u_vw_{v'}: vv'\in E\}$. To head off possible misinterpretation, we emphasize to the reader that $\bip(G)$ denotes the ``bipartite version" of an arbitrary graph $G$, not the bipartition of a bipartite graph $G$.  Similarly, for a $3$-graph $H=(V,E)$, let $\trip(H)$ denote the $3$-partite $3$-graph $(X\cup Y\cup Z, F')$, where $X=\{x_v:v\in V\}$, $Y=\{y_v: v\in V\}$, $Z=\{z_v: v\in V\}$, and $F'=\{x_vy_{v'}w_{v''}: vv'v''\in F\}$.

Given a hereditary graph property $\calG$, let $\bip(\calG)=\{\bip(G): G\in \calG\}$, and given a hereditary $3$-graph property $\calH$, we let $\trip(\calH)=\{\trip(H): H\in \calH\}$.   Given $n\geq 1$, and a hereditary $k$-graph property, $\calH$, set $\calH_n=\{G\in \calH: V(G)=[n]\}$.

For a $k$-graph $H$,  $\age(H)$ denotes the class of finite $k$-graphs isomorphic to an induced sub-$k$-graph of $H$.  Note that $\age(H)$ is always a hereditary $k$-graph property. 

Given a vertex set $V$, and $t,\ell\geq 1$, a $(t,\ell)$-decomposition $\calP$ of $V$ consists of a partition $V=V_1\cup \ldots \cup V_t$ and for each $ij\in {[t]\choose 2}$, a partition $K_2[V_i,V_j]=\bigcup_{\alpha\in [\ell]}P_{ij}^{\alpha}$.  Observe that Definition \ref{def:decomp} is a special kind of $(t,\ell)$-decomposition.   Much of the notation introduced below is related to that introduced on page \pageref{triadnotation}, although here the context is more general (namely that of an arbitrary $(t,\ell)$-decomposition).  We begin by setting notation for the two components of a $(t,\ell)$-decomposition.

\begin{notation}[$\calP_{vert}$ and $\calP_{edge}$]\label{not:pedge}\label{not:pedgeref}
Suppose that $\calP$ is a $(t,\ell)$-decomposition consisting of the vertex partition $V=V_1\cup \ldots \cup V_t$ and for each $ij\in {[t]\choose 2}$, the partition $K_2[V_i,V_j]=\bigcup_{\alpha\in [\ell]}P_{ij}^{\alpha}$.  We then define 
\begin{align*}
\calP_{vert}=\left\{V_1,\ldots, V_t\right\}\text{ and }\calP_{edge}=\left\{P_{ij}^{\alpha}: ij\in {[t]\choose 2}, \alpha\in [\ell]\right\}.
\end{align*}
\end{notation}
Continuing with the notation above, we will write $G_{ijk}^{\alpha,\beta,\gamma}$ to denote the $3$-partite graph
$$
G_{ijk}^{\alpha,\beta,\gamma}=(V_i\cup V_j\cup V_k, P_{ij}^{\alpha}\cup P_{ik}^{\beta}\cup P_{jk}^{\gamma}),
$$
called a \emph{triad} of $\calP$.  When $V$ is the vertex set of a $3$-graph $H=(V,F)$, and no confusion might otherwise arise, we will write $H_{ijk}^{\alpha,\beta,\gamma}$ for the $3$-partite $3$-graph
$$
H_{ijk}^{\alpha,\beta,\gamma}=(V_i\cup V_j\cup V_k, F\cap K_3^{(2)}(G_{ijk}^{\alpha,\beta,\gamma})).
$$
We say a triad $G_{ijk}^{\alpha,\beta,\gamma}$ satisfies $\disc_{2,3}(\e_1,\e_2)$ with respect to $H$ if $(H_{ijk}^{\alpha,\beta,\gamma}, G_{ijk}^{\alpha,\beta,\gamma})$ satisfies $\disc_{2,3}(\e_1,\e_2)$, and we say that $G_{ijk}^{\alpha,\beta,\gamma}$ is \emph{$\e$-homogeneous with respect to $H$} if 
$$
\frac{|F\cap K_3^{(2)}(G_{ijk}^{\alpha,\beta,\gamma})|}{|K_3^{(2)}(G_{ijk}^{\alpha,\beta,\gamma})|}\in [0,\e)\cup (1-\e,1].
$$ 

Similarly, given a partition $\calQ=\{V_1,\ldots, V_t\}$ of $V$, we say that a triple $V_iV_jV_k$ satisfies $\vdisc_3(\e)$ with respect to $H$ if $H[V_i,V_j,V_k]$ satisfies $\vdisc_3(\e)$, and we say that $V_iV_jV_k$ is \emph{$\e$-homogeneous with respect to $H$} if 
$$
\frac{|F\cap K_3[V_i,V_j,V_k]|}{|K_3[V_i,V_j,V_k]|}\in [0,\e)\cup (1-\e,1].
$$
Note that the two definitions above are related to Definition \ref{def:hom}.

\section{Notions of quasirandomness for graphs and hypergraphs}\label{subsec:quasi2}

In this section we give a more complete account of quasirandomness for graphs and $3$-graphs.  We begin with quasirandomness for graphs, of which there are several equivalent formulations.  One of them, $\disc_2$, was already defined in the introduction (see Definition \ref{def:disc2}), and we now define three more which turn out to be equivalent.  Recall that $C_4$ is the cycle on $4$-vertices.

\begin{definition}\label{def:quasi}
Let $\e,d_2>0$, and let $G=(V\cup U,E)$ be a bipartite graph with $|E|=d|U||V|$.
\begin{enumerate}
 \item We say that $G$ has $\cycle_2(\e;d_2)$ if $d\in (d_2-\e, d_2+\e)$ and $G$ contains at most $(d_2^4+\e)|U|^2|V|^2$ many non-induced copies of $C_4$.  
 \item We say that $G$ has $\dev_2(\e; d_2)$ if $d\in (d_2-\e,d_2+\e)$ and the following holds, where $g(u,v)=1-d_2$ if $uv\in E$ and $g(u,v)=-d_2$ if $uv\notin E$.
 $$
 \sum_{u_0,u_1\in U}\sum_{v_0,v_1\in V}\prod_{i\in \{0,1\}}\prod_{j\in \{0,1\}}g(u_i,v_j)\leq \e |U|^2|V|^2.
 $$
\item We say that $G$ is \emph{$\e$-regular} if for all $U'\subseteq U$ with $|U'|\geq \e|U|$ and $W'\subseteq W$ with $|W'|\geq \e |W|$, 
$$
|e_G(U',W')-d|U'||W'||\leq \e |U||W|.
$$
 \end{enumerate}
\end{definition}

In the notation above, we say simply that $G$ has $\disc_2(\e)$ (respectively $\cycle_2(\e)$, $\dev_2(\e)$) if it has $\disc_2(\e;d_2)$ (respectively $\cycle_2(\e;d_2)$, $\dev_2(\e;d_2)$) for some $d_2$.  

Given a graph $G=(V,E)$ and pair of subsets $U,V\subseteq G$, we say $(U,V)$ is an \emph{$\e$-regular pair} (respectively that the pair \emph{has $\disc_2(\e)$, $\cycle_2(\e)$, $\dev_2(\e)$}) if  $G[U,V]$ is $\e$-regular (respectively has $\disc_2(\e)$,  $\cycle_2(\e)$, $\dev_2(\e)$).  In the late 1980s, Chung and Graham \cite{Chung.1988} showed that these notions are equivalent in the sense summarized in Proposition \ref{prop:translate1}.

\begin{proposition}\label{prop:translate1}
For all $\e, d_2>0$, there exists $\delta_2>0$ and $N\geq 1$ such that if $G=(U\cup V,E)$ is a bipartite graph with $|U|=|V|=n\geq N$ then the following hold.
\begin{enumerate}[label=\normalfont(\arabic*)]
\item If $G$ has $\disc_2(\delta_2;d_2)$, then $G$ has $\cycle_2(\e;d_2)$.
\item If $G$ has $\cycle_2(\delta_2;d_2)$, then $G$ has $\dev_2(\e;d_2)$.
\item If $G$ has $\dev_2(\delta;d_2)$, then $G$ is $\e$-regular.
\item If $G$ is $\delta_2$-regular with density $d_2\pm \e$, then $G$ has $\cycle_2(\e,d_2)$.
\end{enumerate}
\end{proposition}

These notions of quasirandomness capture the behavior of a random graph.  For example, for every $\eta>0$, a sufficiently large random bipartite graph will have $\disc_2(\eta;1/2)$ with high probability.   This can be used to prove the following fact.

\begin{fact}\label{fact:rg}
For all $\eta>0$, there is $n_0$ such that for all $n\geq n_0$, almost every bipartite graph $G=(U\cup W, E)$ with $|U|=|W|=n$ has $\disc_2(\eta;1/2)$. 
\end{fact}

Quasirandomness is useful because it allows us to count copies of small fixed subgraphs. For example, quasirandom graphs have about the same number of triangles as one would expect in a random graph, see for instance \cite[Theorem 18]{Rodl.2010}.

\begin{proposition}[Counting Lemma]\label{prop:counting}
For every $t\geq 2$ and $\gamma>0$ there exists $\e>0$ such that the following holds. Let $G=(V_1\cup \ldots \cup V_t,E)$ be a $t$-partite graph such that for each $1\leq i\neq j\leq t$, $G[V_i,V_j]$ has $\disc_2(\e;d_{ij})$. Then 
$$
\Big| |K_t^{(2)}(G)|- \prod_{1\leq i<j\leq t}d_{ij}\prod_{i=1}^t|V_i|\Big|\leq \gamma \prod_{i=1}^t|V_i|.
$$
\end{proposition}

It will be convenient to have a slightly simpler version of this proposition packaged in the following way.

\begin{corollary}\label{cor:counting}
For all $t\geq 2$ and $\e, r>0$, there is $\mu=\mu(k,\e,r)>0$ such that the following hold.  Suppose $G=(U,E)$ is a graph and $U=U_1\cup \ldots \cup U_t$ is a partition so that for each $1\leq i<j\leq t$, $G[U_i,U_j]$ has $\disc_2(\mu;r)$.  

Then the number of $\ubar\in \prod_{i\in [t]}U_i$ with $u_iu_j\in E$ for each $ij\in {[t]\choose 2}$ is 
\[(1\pm \e) r^{{k\choose 2}} \prod_{i\in [t]}|U_i|.\] 
\end{corollary}

It is an exercise to deduce Corollary \ref{cor:counting} from Proposition \ref{prop:counting} (see Appendix \ref{app:slicing}).  One important application of counting lemmas are so-called removal lemmas, which say that if a graph $G$ has very few copies of a small graph $H$, then $G$ can be made $H$-free by changing a small number of edges.  We will use an induced version of this result, due to Alon, Fischer, Krivelevich, and Szegedy \cite{Alon.2000}.  Recall that two graphs $G$ and $G'$ on the same vertex set are \emph{$\delta$-close} if $|E(G)\Delta E(G')|\leq \delta |V(G)|^2$.

\begin{theorem}[Induced graph removal]\label{lem:indremgraph}
For any graph $H$ and $\e>0$, there are $N$ and $\delta$ such that if $G$ is a graph on $n\geq N$ vertices, containing at most $\delta n^{v(H)}$ copies of $H$, then $G$ is $\e$-close to some $G'$ on the same vertex set, which contains no induced copies of $H$.
\end{theorem}

We now turn to notions of quasirandomness for $3$-graphs.  One notion of $3$-graph quasirandomness already defined in the introduction is $\vdisc_3$, and this can be seen as a natural generalization of $\disc_2$.  As was shown by Kohayakawa, R\"{o}dl, and Skokan,  this notion of quasirandomness does not correspond to a general counting lemma, but rather to one which counts copies of \emph{linear hypergraphs} \cite{Kohayakawa.2010}.  

It is also natural, and for certain applications necessary, to consider notions which take into account the quasirandomness of both the pairs and triples in a $3$-graph.  The next definition contains three such notions, which are equivalent to one another (but not to $\vdisc_3$, see \cite{Nagle.2013}).  We note that one of these, $\disc_{2,3}$ was defined in the introduction in Definition \ref{def:regtrip}.

\begin{definition}
Let $\e>0$, $n\geq 1$, and suppose $H=(V,F)$ is a $3$-graph. Assume $G=(V_1\cup V_2\cup V_3,E)$ is a $3$-partite graph underlying $H$, and assume  $d_3$ is defined by $|E(H)|=d_3|K^{(2)}_3(G)|$. 
\begin{enumerate}
\item We say that \emph{$(H,G)$ has $\disc_{2,3}(\e)$} if for each $1\leq i<j\leq 3$, the graph $G[V_i,V_j]$ has $\disc_2(\e)$, and for every (not necessarily induced) subgraph $G'\subseteq G$,
\[\big||F\cap K_3^{(2)}(G')|-d_3|K^{(2)}_3(G')|\big|\leq \e |V_1||V_2||V_3|.\]
\item We say that \emph{$(H,G)$ has $\dev_{2,3}(\e)$} if for each $1\leq i<j\leq 3$, the graph $G[V_i,V_j]$ has $\dev_2(\e)$, and 
$$
\sum_{u_0,u_1}\sum_{v_0,v_1}\sum_{w_0,w_1}\prod_{\e \in \{0,1\}^3}f_{(H,G)}(u_{\e_1},v_{\e_2}, w_{\e_3})\leq \e |V_1|^2|V_2|^2|V_3|^2,
$$
where $f_{(H,G)}(u,v,w)=1_{K_3^{(2)}(G)}(u,v,w)(1_{F}(u,v,w)-d_3)$.
\item We say that \emph{$(H,G)$ has $\oct_{2,3}(\e)$} if for each $1\leq i<j\leq 3$, the graph $G[V_i,V_j]$ has $\cycle_2(\e;d_2)$, and $H$ contains at most $(d_3^8d_2^{12}+\e)|V_1|^2|V_2|^2|V_3|^2$ copies of $K_{2,2,2}^{(3)}$.
\end{enumerate}
\end{definition}

Like quasirandomness in graphs, these notions reflect the behavior of random $3$-graphs.  For example, for all $\eta>0$, a sufficiently large $3$-partite $3$-graph has $\disc_{3}(\eta)$ with high probability.  This can be used to prove the following well known fact which we will use later on (specifically in the proof of Theorem \ref{thm:dischom} on page \pageref{fact:rhgref}).

\begin{fact}\label{fact:rhg}
For all $\eta>0$, there is $N\geq 1$ such that for all $n\geq N$ there exists a $3$-partite $3$-graph $H_n=(U\cup W\cup Z, F)$ such that $(H_n, K_3[U,W,Z])$ has $\disc_{3}(\eta)$ and density $1/2\pm \eta$.
\end{fact}

As mentioned in the introduction,  to obtain a general counting lemma, one requires a more subtle definition, which  maintains control over the relative quasirandomness of the edge decomposition.  

\begin{definition}\label{def:regtrip2}
Let $\e_1,\e_2>0$, and suppose $H=(V,F)$ is a $3$-graph. Assume $G=(V_1\cup V_2\cup V_3,E)$ is a $3$-partite graph underlying $H$, and suppose $d_3$ is defined by $|F|=d_3|K^{(2)}_3(G)|$.
\begin{enumerate}
\item We say that \emph{$(H,G)$ has $\disc_{2,3}(\e_1,\e_2)$} if there is $d_2\in (0,1]$ such that for each $1\leq i<j\leq 3$, the graph $G[V_i,V_j]$ has $\disc_2(\e_2;d_2)$, and for every subgraph $G'\subseteq G$ (not necessarily induced),
\[\big||F\cap K_3^{(2)}(G')|-d_3|K^{(2)}_3(G')|\big|\leq \e_1(d_2)^3 |V_1||V_2||V_3|.\]
\item We say that \emph{$(H,G)$ has $\dev_{2,3}(\e_1,\e_2)$} if there is $d_2\in (0,1]$ such that  for each $1\leq i<j\leq 3$, the graph $G[V_i,V_j]$ has $\dev_2(\e_2;d_2)$, and 
$$
\sum_{u_0,u_1}\sum_{v_0,v_1}\sum_{w_0,w_1}\prod_{\e \in \{0,1\}^3}f_{(H,G)}(u_{\e_1},v_{\e_2}, w_{\e_3})\leq \e_1 d_2^{12}|V_1|^2|V_2|^2|V_3|^2.
$$
where $f_{(H,G)}(u,v,w)=1_{K_3^{(2)}(G)}(u,v,w)(1_{F}(u,v,w)-d_3)$.
\item  We say that \emph{$(H,G)$ has $\oct_{2,3}(\e_1,\e_2)$} if there is $d_2\in (0,1]$ such that for each $1\le i\neq j\leq 3$, the graph $G[V_i,V_j]$ has $\cycle_2(\e_2;d_2)$, and the number of copies of $K_{2,2,2}^{(3)}$ in $H$ is at most $d_3^8d_2^{12} |V_1|^2|V_2|^2|V_3|^2+\e_1d_2^{12} |V_1|^2|V_2|^2|V_3|^2$.
\end{enumerate}
\end{definition}

In applications, $\e_2$ will usually be much smaller than $\e_1$.  The history of these notions and their corresponding regularity and counting lemmas is explained in full in \cite{Nagle.2013}, for example.  Concerning the results that we shall use, a regularity and counting lemma for $\dev_{2,3}$ and $\oct_{2,3}$  was first proved by Gowers in  \cite{Gowers.20063gk}, while R\"{o}dl and Frankl proved a regularity lemma for $\disc_{2,3}$ in addition to a counting lemma for a stronger version of $\disc_{2,3}$ in \cite{Frankl.2002}.   However, it has been shown that the notions in Definition \ref{def:regtrip} and Definition \ref{def:regtrip2} are all equivalent in the following sense. 

\begin{proposition}\label{prop:qrsame}
For all $d_3,\e_3>0$ there exists $\delta_3>0$ such that for all $d_2,\e_2>0$ there exists $\delta_2>0$ and $n_0\geq 1$ such that the following hold.  Suppose $H=(V,F)$ is a $3$-graph. Assume $G=(V_1\cup V_2\cup V_3,E)$ is a $3$-partite graph underlying $H$ with $|V_1|=|V_2|=|V_3|=n$.  Suppose the density of each $G[V_i,V_j]$ is in $(d_2-\delta_2,d_2+\delta_2)$ and $e(H)=d_3|K_3^{(2)}(G)|$.
\begin{enumerate}[label=\normalfont(\arabic*)]
\item {\rm(Nagle et al. \cite{Nagle.2013})} If $(H,G)$ has $\disc_{2,3}(\delta_3,\delta_2)$ then it also has $\oct_{2,3}(\e_3,\e_2)$.
\item {\rm(Dementieva et al. \cite{Dementieva.2002})} If $(H,G)$ has $\oct_{2,3}(\delta_3,\delta_2)$ then it also has $\disc_{2,3}(\e_3,\e_2)$.
\item {\rm(Nagle et al. \cite{Nagle.2013})} If $(H,G)$ has $\oct_{2,3}(\delta_3,\delta_2)$ then it also has $\dev_{2,3}(\e_3,\e_2)$.
\item {\rm(Gowers \cite{Gowers.20063gk})} If $(H,G)$ has $\dev_{2,3}(\delta_3,\delta_2)$ then it also has $\oct_{2,3}(\e_3,\e_2)$.
\end{enumerate}
\end{proposition}

Proposition \ref{prop:qrsame} shows that Definition \ref{def:error} could be equivalently formulated with any of these three notions of quasirandomenss. 

In order to apply a hypergraph regularity lemma, one needs a corresponding counting lemma. Theorem \ref{thm:countingHG} was first proven by Gowers in \cite{Gowers.20063gk} for $\dev_{2,3}$, but we state it in terms of $\disc_{2,3}$ (via Proposition \ref{prop:qrsame}).  

\begin{theorem}[Counting lemma for regular 3-graphs]\label{thm:countingHG}
For all $t\in \mathbb{N}$ and $\xi$, $d_3>0$, there exists $\delta_3>0$ such that for every $d_2>0$, there exists $\delta_2>0$ and $n_0\in \N$ such that the following holds.

Let $G=(V,E)$ be a $t$-partite graph with vertex partition $V=V_1\cup \ldots \cup V_t$, where $|V_1|=\ldots =|V_t|=n\geq n_0$, and let $H=(V,F)$ be a $3$-graph on vertex set $V$ such that $F\subseteq K_3^{(2)}(G)$. Suppose that for each $ij\in {[t]\choose 2}$, $G[V_i,V_j]$ has $\disc_2(\delta_2;d_2)$.  For each $ijk\in {[t]\choose 3}$, let $G^{ijk}=G[V_i,V_j,V_k]$, $H^{ijk}=(V_i\cup V_j\cup V_k, K_3^{(2)}(G^{ijk})\cap F)$, and let $d_{ijk}$ be such that $e(H^{ijk})=d_{ijk}|K^{(2)}_3(G^{ijk})|$.  Suppose that for each $ijk\in {[t]\choose 3}$, $d_{ijk}\geq d_3$, and $(H^{ijk},G^{ijk})$ satisfies $\disc_{2,3}(\delta_3,\delta_2)$.  Then 
\[\left| |K^{(3)}_t(H)|-d_2^{t\choose 2} d n^t\right|\leq \xi d_2^{t\choose 2}d n^t\]
where $d=\prod_{ijk\in {[t]\choose 3}}d_{ijk}$.
\end{theorem}

An induced analogue of the above was proven by R\"{o}dl and Schacht in \cite{Rodl.2009}.

\begin{theorem}[Induced Counting Lemma]\label{thm:countinginducedHG}
For all $t\in \mathbb{N}$ and $\xi$, $d_3>0$, there exists $\delta_3>0$ such that for every $d_2>0$, there exists $\delta_2>0$ and $n_0\in \N$ such that the following holds.

Let $H_C=([t],C)$ be a $t$-partite, $3$-graph.  Let $G=(V,E)$ be a $t$-partite graph with vertex partition $V=V_1\cup \ldots \cup V_t$, where $|V_1|=\ldots =|V_t|=n\geq n_0$, and let $H=(V,F)$ be a $3$-graph on vertex set $V$ such that $F\subseteq K_3^{(2)}(G)$. Suppose that for each $ij\in {[t]\choose 2}$, $G[V_i,V_j]$ has $\disc_2(\delta_2;d_2)$. For each $ijk\in {[t]\choose 3}$, let $G^{ijk}=G[V_i,V_j,V_k]$, $H^{ijk}=(V_i\cup V_j\cup V_k, K_3^{(2)}(G^{ijk})\cap F)$, and let $d_{ijk}$ be such that $e(H^{ijk})=d_{ijk}|K^{(2)}_3(G^{ijk})|$.

Suppose that for each $ijk\in C$, $d_{ijk}\geq d_3$ and for each $ijk\in {[t]\choose 3}\setminus C$, $d_{ijk}\leq 1-d_3$.  Suppose further that $(H^{ijk},G^{ijk})$ satisfies $\disc_{2,3}(\delta_3,\delta_2)$ for each $ijk\in {[t]\choose 3}$.  Let $m_C$ be the number of tuples $(v_1,\ldots, v_t)\in \prod_{i=1}^tV_i$ such that $v_iv_jv_k\in F$ if and only if $ijk\in C$.  Then  
\[\left| m_C-d_2^{t\choose 2} d n^t\right|\leq \xi d_2^{t\choose 2}d n^t\]
where $d=\prod_{ijk\in {[t]\choose 3}}d_{ijk}$.
\end{theorem}
As an application of Theorem \ref{thm:countinginducedHG}, R\"{o}dl and Schacht \cite{Rodl.2009} proved the following induced removal lemma.

\begin{theorem}[Induced removal for $3$-graphs]\label{thm:indrem}
For all $\eta>0$ and finite $3$-graphs $G$, there are $c>0$, $C>0$, and $n_0$ such that the following holds.  Suppose $H$ is a $3$-graph on $n\geq n_0$ vertices, and $H$ contains at most $c n^{v(G)}$-many induced copies of $G$. Then $H$ is $\eta$-close to some $H'$ containing no induced copies of $G$.
\end{theorem}

In fact, R\"{o}dl and Schacht's results hold for $k$-graphs for any $k\geq 3$, but we will only use the case $k=3$.  We will also use the following, even stronger result for hereditary properties, also from \cite{Rodl.2009}.

\begin{theorem}[Induced removal for hereditary properties]\label{thm:indremprop}
Suppose $\calH$ is a hereditary $3$-graph property.  Then for all $\eta>0$ there are $c>0$, $C>0$, and $n_0$ such that the following holds.  Suppose $H$ is a $3$-graph on $n\geq n_0$ vertices, and for every $G\in \calF_\calH$ on at most $C$ vertices, $H$ contains at most $c n^{v(G)}$ copies of $G$. Then $H$ is $\eta$-close to some $H'\in \calH_n$.
\end{theorem}

This theorem can be used to give another characterization of $\sim$-classes (see Definition \ref{def:close}).

\begin{definition}
Given a hereditary $3$-graph property $\calH$, let $\Delta(\calH)$ denote all finite $3$-graphs $G$ such that for all $\e>0$, there is $N=N(G,\calH,\e)$ so that if $n\geq N$ and $H\in \calH_n$, then $H$ contains at most $\e n^{v(G)}$-many induced sub-$3$-graphs isomorphic to $G$.  
\end{definition}

In other words, $\Delta(\calH)$ is the set of finite $3$-graphs $G$ so that large elements of $\calH$ cannot contain too many copies of $G$.   Theorem \ref{thm:indremprop} implies that this definition characterizes $\sim$-classes.  

\begin{proposition}\label{prop:simclasses1}
$\calH\sim \calH'$ if and only if $\Delta(\calH)=\Delta(\calH')$.
\end{proposition}

For a proof, see \cite{Rodl.2009}.

\section{Lemmas for slicing, intersecting, and combining decompositions}\label{subsec:intersecting}

This section contains several facts which will be useful for building regular decompositions. While most of these are fairly standard, we include several proofs in Appendix \ref{app:slicing} for the sake of completeness.  First, Lemma \ref{lem:sl} tells us that large sub-pairs in quasirandom graphs are still quasirandom.

\begin{proposition}[Sub-pairs lemma]\label{lem:sl} Suppose $G=(A\cup B, E)$ is a bipartite graph and $|E|=d|A||B|$.  Suppose $A'\subseteq A$ and $B'\subseteq B$ satisfy $ |A'| \geq \gamma |A|$ and $|B'| \geq \gamma|B|$ for some $\gamma \geq \e$, and $G$ satisfies $\disc_2(\epsilon;d)$.  Then $G':=(A'\cup B', G[A',B'])$ satisfies $\disc_2(\epsilon';d')$ where $\epsilon'=2\gamma^{-2}\e$ and $d'\in(d -\gamma^{-2}\e,d + \gamma^{-2}\e)$.
\end{proposition}

We refer the reader to \cite{Komlos.1996} for a proof.  Lemma \ref{lem:subpairs} is a similar result for triples, and says that large sub-triples of $\disc_{2,3}$-regular triads are still somewhat regular. 

\begin{lemma}[Sub-triples lemma]\label{lem:subpairs}
For all $k\geq 1$ and $\delta_3>0$, there exists $\delta_3'>0$ such that for all $d_2,d_3\in (0,1]$ and $0<\delta_2\leq d^3_2/8$, there are $\delta_2'>0$ and $m_0\geq 1$ so that the following holds.  Suppose $V=V_1\cup V_2\cup V_3$ where $|V_1|, |V_2|, |V_3|\geq m_0$, and $G=(V_1\cup V_2\cup V_3, E)$ is a $3$-partite graph such that  for each $1\leq i< j\leq 3$, $G[V_i,V_j]$ has $\disc_2(\delta'_2;d_2)$.  Assume $H=(V,F)$ is a $3$-graph, and $(H',G)$ has $\disc_{2,3}(\delta'_2,\delta'_3)$ with density $d_3$, where $H'=(V, F\cap K_3^{(2)}(G))$.

Suppose that for each $i\in [3]$, $V_i'\subseteq V_i$ satisfies $|V_i'|\geq |V_i|/k$.  Then  $(H'',G')$ has $\disc_{2,3}(\delta_3, \delta_2)$ with density $d_3'=d_3\pm \delta_3$, where $G'=G[V_1'\cup V'_2\cup V'_3]$, and $H''=H'[V_1'\cup V_2'\cup V_3']$.
\end{lemma}

The proof consists of standard arguments and appears in Appendix \ref{app:slicing}.  The next couple of propositions tell us that homogeneity (i.e. density near zero or one) implies quasirandomness, both in the $\disc_{2,3}$ and $\vdisc_3$-sense, thereby generalizing the following easy fact about graphs.

\begin{fact}\label{fact:homimpliesrandombinary}
Suppose $H=(U\cup V, E)$ is a bipartite graph. For all $\e>0$, if $|E\cap K_2[U,V]|\leq \e|K_2[U,V]|$, then $H$ satisfies $\disc_2(\e)$.
\end{fact}

\begin{proposition}\label{prop:homimpliesrandomv}
If $H=(V_1\cup V_2\cup V_3, F)$ is a $3$-partite $3$-graph and $|F\cap K_3[V_1,V_2,V_3]|\leq \e |V_1||V_2||V_3|$, then $H$ has $\vdisc_3(\e)$.
\end{proposition}
\begin{proof}
Suppose for each $i\in [3]$, $V_i'\subseteq V_i$.   Let $d$ be such that 
\[|F\cap K_3[V_1,V_2,V_3]|=d|K_3[V_1,V_2,V_3]|.\]  
Since $d\leq \e$,
$$
|F\cap K_3[V_1',V_2',V_3']|\leq \e|V_1||V_2||V_3| \leq \e|V_1||V_2||V_3|+d|V_1'||V_2'||V_3'|.
$$
Clearly $|F\cap K_3[V_1',V_2',V_3']|\geq 0>d|V_1'||V_2'||V_3'|-\e|V_1||V_2||V_3|$.  Thus
$$
\big||F\cap K_3[V_1',V_2',V_3']|-d|V_1'||V_2'||V_3'|\big|\leq \e|V_1||V_2||V_3|,
$$
as desired.
\end{proof}

\begin{proposition}\label{prop:homimpliesrandome}
For all $0<\e<1/2$ and $d_2>0$, there is $\delta_2>0$ such that for all $0<\delta\leq \delta_2$, there is $N$ such that the following holds.  Suppose $H=(V_1\cup V_2\cup V_3, F)$ is a $3$-partite $3$-graph on $n\geq N$ vertices, and for each $i,j\in [3]$, $||V_i|-|V_j||\leq \delta|V_i|$.  Suppose $G=(V_1\cup V_2\cup V_3, E)$ is a $3$-partite graph, where for each $1\leq i< j\leq 3$, $G[V_i,V_j]$ has $\disc_2(\delta;d_2)$, and assume
$$
|F\cap K_3^{(2)}(G))|\leq \e|K_3^{(2)}(G)|.
$$
Then $(H',G)$ has $\disc_{2,3}(\delta,\sqrt{\e})$, where $H'=(V_1\cup V_2\cup V_3, F\cap K_3^{(2)}(G))$.
\end{proposition}
\begin{proof}
Suppose $0<\e<1/2$ and $0<d_2$.  Let $\delta_1$ and $N$ be from Corollary \ref{cor:counting} for $\gamma=1/32$, and set $\delta=\min\{\e^2/32, \delta_2/32\}$.  Suppose $0<\delta\leq \delta_2$,  $H=(V_1\cup V_2\cup V_3, F)$ is a $3$-partite $3$-graph on $n\geq N$ vertices such that for each $i,j\in [3]$, $||V_i|-|V_j||\leq \delta|V_i|$.  Suppose $G=(V_1\cup V_2\cup V_3, E)$ is a $3$-partite graph, where for each $1\leq i< j\leq 3$, $G[V_i,V_j]$ has $\disc_2(\delta;d_2)$, and assume
$$
|F\cap K_3^{(2)}(G))|\leq \e|K_3^{(2)}(G)|.
$$
By Corollary \ref{cor:counting}, $|K_3^{(2)}(G)|=d^3_2(1\pm  \gamma)|V_1||V_2||V_3|$. Suppose $G'\subseteq G$ is a subgraph, and let $d$ be such that $|F\cap K_3^{(2)}(G))|=d|K_3^{(2)}(G)|$.  Then 
$$
|F\cap K_3^{(2)}(G')|\leq \e|K_3^{(2)}(G)| \leq  \e d^3_2(1+\gamma)^3|V_1||V_2||V_3|\leq \sqrt{\e}d_2^3|V_1||V_2||V_3|,
$$
where the last inequality is by our choice of $\gamma$ and since $\e<1/2$.  Using that $\e|K_3^{(2)}(G)| \leq \sqrt{\e}d_2^3|V_1||V_2||V_3|$, $d\leq \e$, and $K_3^{(2)}(G')\subseteq K_2^{(2)}(G)$, we have that
\begin{align*}
|F\cap K_3^{(2)}(G')|\geq 0\geq d|K_3^{(2)}(G')|-\e |K_3^{(2)}(G)|\geq d|K_3^{(2)}(G')|-\sqrt{\e}d^3_2|V_1||V_2||V_3|.
\end{align*}
Thus $||F\cap K_3^{(2)}(G')|- d|K_3^{(2)}(G')||\leq \sqrt{\e}d^3_2|V_1||V_2||V_3|$.  This concludes the proof that $(H,G)$ has $\disc_{2,3}(\sqrt{\e},\delta)$. 
\end{proof}

Our next set of lemmas deals with combining and partitioning quasirandom graphs.  First, we observe that adding or subtracting quasirandom graphs results in graphs which are still somewhat quasirandom.

\begin{fact}\label{fact:adding}
Suppose $E_1$ and $E_2$ are disjoint subsets of $K_2[U,V]$.
\begin{enumerate}[label=\normalfont(\arabic*)]
\item If $(U\cup V, E_1)$ has $\disc_2(\e_1;d_1)$, $(U\cup V, E_2)$ has $\disc_2(\e_2;d_2)$, then $(U\cup V, E_1\cup E_2)$ has $\disc_2(\e_1+\e_2; d_2+d_1)$.  
\item If $(U\cup V, E_1\cup E_2)$ has $\disc_2(\e;d)$ and $(U\cup V, E_1)$ has $\disc_2(\e_1;d_1)$, then $(U\cup V, E_2)$ has $\disc_2(\e+\e_1; d-d_1)$.
\end{enumerate}
\end{fact}

For a proof, see Appendix \ref{app:slicing}.  Our next lemma is due to Frankl and R\"{o}dl \cite[Lemma 3.8]{Frankl.2002}, and allows us to partition quasirandom graphs into several equally sized, quasirandom subgraphs.

\begin{lemma}\label{lem:3.8}
For all $\e>0$, $\rho\geq 2\e$, $0<r<\rho/2$, and $\delta>0$, there is $m_0=m_0(\e, \rho, \delta)$ such that the following holds.  Suppose $|U|=|V|=m\geq m_0$, and $G=(U\cup V, E)$ is a bipartite graph  satisfying $\disc_2(\e)$ with density $\rho$.  Then if $\ell=[1/r]$, and $\e\geq 10(1/\ell m)^{1/5}$, then there is a partition $E=E_0\cup E_1\cup \ldots \cup E_\ell$ such that 
 \begin{enumerate}[label=\normalfont(\roman*)]
 \item for each $1\leq i\leq \ell $, $(U\cup V, E_i)$ has $\disc_2(\e)$ with density $\rho r(1\pm \delta)$, and
 \item $|E_0|\leq \rho r(1+\delta)m^2$.
 \end{enumerate}
 Further, if $1/r\in \mathbb{Z}$, then $E_0=\emptyset$.
\end{lemma} 

Lemma \ref{lem:3.8} will come in handy frequently as we build decompositions.  An example of such an application is Lemma \ref{lem:disc2} below, which shows that given a decomposition $\calP'$ of a vertex set $V$, we can change a small number of elements from $\calP'_{edge}$ to obtain a decomposition $\calP$ of $V$ which has no $\disc_2$-irregular elements in $\calP_{edge}$ (recall Notation \ref{not:pedge} on page \pageref{not:pedgeref}). This shows that irregularity involving only the $\disc_2$-irregularity of $\calP_{edge}$ does not produce a meaningful notion of ternary irregularity.

\begin{lemma}\label{lem:disc2}
For all $\e_1>0$ and $\e_2:\mathbb{N}\rightarrow (0,1]$, there are $\e_1'>0$ and $\e_2':\mathbb{N}\rightarrow (0,1]$ such that for all $\ell, t\geq 1$, there is $N$ such that the following holds.  Suppose $V$ is a set of size $n\geq N$, and $\calP'$ is an $(t,\ell, \e_1',\e_2'(\ell))$-decomposition of $V$.  Then there is a $(t,\ell, \e_1,\e_2(\ell))$-decomposition $\calP$ of $V$ such that
\begin{enumerate}[label=\normalfont(\arabic*)]
\item  $\calP_{vert}=\calP'_{vert}$;
\item every $P\in \calP_{edge}$ has $\disc_2(\e_2(\ell);1/\ell)$;
\item if $P'\in \calP_{edge}'$ satisfied $\disc_2(\e'_2(\ell);1/\ell)$, then $P'\in \calP_{edge}$. 
\end{enumerate}
\end{lemma}
\begin{proof}
Fix $\e_1>0$ and $\e_2:\mathbb{N}\rightarrow (0,1]$.   Without loss of generality, we may assume $\e_2$ is non-increasing (e.g. we can replace $\e_2$ with any increasing function bounded above by the original $\e_2$).  Set $\e_1'=\e_1^2/4$, and for all $x\in \mathbb{N}$, set $\e_2'(x)=\e_2( x^2)/x^3$.  Given $\ell,t\geq 1$, let $N=\max\{tm_0(s\e_2'(\ell),1-s/\ell ,s\e_2'(\ell)): 1\leq s\leq (1-\sqrt{\e_1'})\ell\}$, where $m_0(\cdot,\cdot,\cdot)$ is from Lemma \ref{lem:3.8}.

Suppose $V$ is a set of size $n\geq N$, and $\calP'$ is a $(t,\ell,\e_1',\e'_2(\ell))$-decomposition of $V$, consisting of $\{V_i:i\in [t]\}$ and $\{P_{ij}^{\alpha}:ij\in {[t]\choose 2}, \alpha\leq \ell\}$.  Set $m=n/t$, and define $\Delta=\{P\in \triads(\calP'): P\text{ has }\disc_2(\e'_2(\ell);1/\ell)\}$.  Then let
$$
\Omega=\{xyz\in {V\choose 3}:\text{ for some }P\in \triads(\calP')\setminus \Delta, P\cap\{xy, yz,xz\}\neq \emptyset\}.
$$
In other words, $\Omega$ is the set of triples which intersect a $\disc_2$-irregular element from $\calP'_{edge}$.  Then define
$$
\Sigma=\Big\{ij\in {[t]\choose 2}:|K_2[V_i,V_j]\setminus (\bigcup \Delta)|\geq \sqrt{\e_1'}n^2/t^2\Big\}.
$$
Note that if $ij\in \Sigma$, then for all $k\in [t]\setminus \{i,j\}$, $K_3[V_i,V_j,V_k]$ contains at least $\sqrt{\e_1'}n^3/t^3$ many triples from $\Omega$, and thus, at least $(t-2)\sqrt{\e_1'}n^3/t^3$ triples in $\Omega$ intersect $K_2[V_i,V_j]$.  Since by assumption, $|\Omega|\leq \e_1' n^3$, this implies $|\Sigma|\leq 3\sqrt{\e_1'}t^2$.  We now define a new partition of $K_2[V_i,V_j]$, for each $ij\in {[t]\choose 2}$.  We do this in cases, depending on whether $ij\in \Sigma$ or $ij\notin \Sigma$. 

Suppose first $ij\in \Sigma$.  In this case, apply Lemma \ref{lem:3.8} to $K_2[V_i,V_j]=E_1$ with $\e=\delta=\e'_2(\ell)$ and $p=1/\ell$ to obtain a partition $K_2[V_i,V_j]=Q_{ij}^1\cup \ldots \cup Q_{ij}^{\ell}$ such that for each $1\leq \alpha\leq \ell$, $Q_{ij}^\alpha$ has  has $\disc_2(2\e'_2(\ell);1/\ell)$.  By definition of $\e_2'(\ell)$, this shows each such $Q_{ij}^{\alpha}$ has $\disc_2(\e_2(\ell);1/\ell)$.  

Suppose now $ij\notin \Sigma$.  Let $P_{ij}^{\alpha_1},\ldots, P_{ij}^{\alpha_s}$ enumerate the elements of $\Delta_{ij}$ and set $E_0=K_2(V_i,V_j)\setminus \bigcup_{u=1}^sP_{ij}^{\alpha_u}$.  If $s=\ell$, then let $Q_{ij}^{\alpha}=P_{ij}^{\alpha}$ for each $\alpha\leq \ell$.  Assume on the other hand that $s<\ell$.  Then since $ij\notin \Sigma$, $|E_0|\leq \sqrt{\e_1'}m^2$, so we must have that $s\geq (1-\sqrt{\e_1'})\ell$.  By Fact \ref{fact:adding}, $E_1\cup \ldots \cup E_s$ has $\disc_2(s\e'_2(\ell); s/\ell)$, so its compliment in $K_2[V_i,V_j]$, $E_0$, has $\disc_2(s\e'_2(\ell); 1-s/\ell)$.  Setting $\rho_0=|E_0|/m^2$, this implies $\rho_0\geq \frac{1}{\ell}(1-\e_2'(\ell))$, and consequently, by definition of $\e_2'(\ell)$, $\rho_0\geq 2s\e'_2(\ell)$.  Let $u=2\ell(\ell-s)$ and $p=1/u$.  Then
$$
p=\frac{1}{2\ell(\ell-s)}<\frac{1}{2(\ell-s)}<\frac{\ell-s}{2\ell}\leq \rho_0/2,
$$
where the second inequality is because $s\geq (1-\sqrt{\e_1'})\ell$.   Consequently, we may apply Lemma \ref{lem:3.8} to $E_0$ with this $u$ and $p$, and $\e=\delta=s\e'_2(\ell)$ to obtain a partition $E_0=E_1'\cup \ldots \cup E_u'$ where for each $i\in [u]$, $E_i'$ has of density $\frac{1}{u}(1\pm s\e_2'(\ell))$ and satisfies $\disc_2(s\e'_2(\ell))$.  Let $I_1\cup \ldots \cup I_{\ell-s}$ be an equipartition of $[u]$ and for each $i\in [\ell-s]$, set $E_i''=\bigcup_{j\in I_i}E_j'$.  By Fact \ref{fact:adding}, each $E_i''$ has $\disc_2(2\ell s\e_2'(\ell), \frac{2(\ell-s)}{u}(1\pm 2\ell s\e'_2(\ell)))$.  By our choice of $\e_2'(\ell)$ and the definition of $u$, this means each $E_i''$ has $\disc_2(\e_2(\ell);1/\ell)$.  Now define $\calP$ by setting $\calP_{vert}=\{V_1,\ldots, V_t\}$ and $\calP_{edge}=\{P_{ij}^{\alpha_1},\ldots, P_{ij}^{\alpha_s}, E_1'',\ldots, E_{\ell-s}''\}$.   

Clearly each element of $\calP_{edge}$ satisfies $\disc_2(\e_2(\ell);1/\ell)$.  Moreover, by construction, if $P_{ij}^\alpha\in \calP'_{edge}$ had $\disc_2(\e'_2(\ell);1/\ell)$ in $\calP'$, then $\calP_{ij}^{\alpha}\in \calP_{edge}$, as desired. 
\end{proof}

 As a corollary, we now prove Proposition \ref{prop:binaryboring}, which tells us that every sufficiently large $3$-graph has regular decompositions with no $\disc_2$-irregular triads. 
\vspace{2mm}

\begin{proofof}{Proposition \ref{prop:binaryboring}}
Fix a hereditary $3$-graph property $\calH$.  Fix $\e_1>0$ and $\e_2:\mathbb{N}\rightarrow (0,1]$. Let $\e_1'>0$ and $\e_2':\mathbb{N}\rightarrow (0,1]$ be as in Lemma \ref{lem:disc2} for $\e_1$ and $\e_2$.  Given $t_0,\ell_0\geq 1$, choose $n_0, L_0,T_0$ as in Theorem \ref{thm:reg2} for $\e_1'$, $\e'_2$, $t_0,\ell_0$, and let $n_1$ be as in Lemma \ref{lem:disc2} for $\e_1,\e_2,T_0,L_0$.  Set $N=\{n_0,n_1\}$.

Suppose $H\in \calH$ is a $3$-graph on $n\geq N$ vertices.  By Theorem \ref{thm:reg2}, there are $t_0\leq t\leq T_0$, $\ell_0\leq \ell\leq L_0$, and a $(t,\ell,\e_1',\e'_2(\ell))$-decomposition $\calP'$ of $V(H)$ which is  $(\e_1',\e'_2(\ell))$-regular with respect to $H$. By Lemma \ref{lem:disc2}, there is a $(t,\ell,\e_1,\e_2(\ell))$-decomposition $\calP$ of $V$ which has  no $\disc_2$-irregular triads, and so that the elements from $\calP'_{edge}$ which had $\disc_2(\e_1';\e_2'(\ell))$ in $\calP'$ remain unchanged in $\calP$.  Consequently, every triad of $\calP'$ which had $\disc_{2,3}(\e'_1,\e_2'(\ell))$ with respect to $(H,\calP')$ is unchanged in $\calP$.  Since $\e_1'\leq \e_1$ and $\e_2'(\ell)\leq \e_2(\ell)$, this implies that $\calP$ is $(\e_1,\e_2(\ell))$-regular with respect to $H$. 
\end{proofof}

\vspace{2mm}

Our next series of lemmas will be useful for situations where we intersect two different decompositions  of a vertex set $V$.  First, we make a definition that describes when one decomposition approximately refines another.

\begin{definition}\label{def:refinement}
Suppose $V$ is a set and $\calP$, $\calQ$ are decompositions of $V$. Then $\calP$ is an \emph{$(\e_1,\e_2)$-approximate refinement of $\calQ$} if there exists $\Sigma\subseteq {\calP_{vert}\choose 2}$ such that the following hold.
\begin{enumerate}
\item For all $V_i\in \calP_{vert}$ there is $W_j\in \calQ_{vert}$ such that $|V_i\setminus W_j|<\e_1|V_i|$, 
\item For all  $P\in \calP_{edge}$ with $V_iV_j\notin \Sigma$, there $Q\in \calQ_{edge}$ with  $|P\setminus Q|\leq \e_1 |P|$, and such that $P\setminus Q$ has $\disc_2(\e_2)$.
\item $|\Sigma|\leq \e_1 {\calP_{vert}\choose 2}$.
\end{enumerate}
\end{definition}

Suppose $\calP$ and $\calQ$ are decompositions of the same vertex set $V$, say $\calP_{vert}=\{V_1,\ldots, V_t\}$, $\calQ_{vert}=\{U_1,\ldots, U_s\}$, $\calP_{edge}=\{P_{ij}^{\alpha}: ij\in {[t]\choose 2}, \alpha\in [\ell]\}$,  and $\calQ_{edge}=\{Q_{ij}^{\beta}: ij\in {[s]\choose 2}, \beta\in [r]\}$.  We then let $\calP\wedge \calQ$\footnote{This usage of wedge notation is unrelated to that appearing in Section \ref{subsec:strongstable}, page \pageref{not:wedge}.} denote the \emph{common refinement of $\calP$ and $\calQ$}, meaning the decomposition with
$$
(\calP\wedge \calQ)_{vert}=\{V_i\cap U_j: i\in [t],j\in [s]\},
$$
 and with $(\calP\wedge \calQ)_{edge}=\{P_{ij}^{\alpha}\cap Q_{i'j'}^{\beta}: ij\in {[t]\choose 2}, i'j'\in {[s]\choose 2}, \alpha\in [\ell], \beta\in [r]\}$.  Given such decompositions $\calP$ and $\calQ$, we can always find an approximate refinement of $\calP\wedge \calQ$ which is also a $(t',\ell',\e_1,\e_2(\ell'))$-decomposition.

\begin{lemma}\label{lem:refinement}
For all $\e_1>0$ and $\e_2:\mathbb{N}\rightarrow (0,1]$, there are functions $f, g:\mathbb{N}^4\rightarrow (0,1]$  such that for all $\ell,t,s,r\geq 1$, there is $N\geq 1$ such that the following holds.  

Suppose $V$ is a set of size $n\geq N$, $\calP$ is an $(\ell,t)$-decomposition of $V$, and $\calQ$ is an $(r,s)$-decomposition of $V$.  Then there are $\ell'\leq g(r,s,\ell,t)$ and $t'\leq f(r,s,\ell,t)$, and a $(t',\ell',\e_1,\e_2(\ell'))$-decomposition $\calR$ of $V$ such that  every $R_{ij}^{\alpha}\in \calR_{edge}$ has $\disc_2(\e_2(\ell');1/\ell')$, and so that $\calR$ is an $(\e_1,\e_2(\ell'))$-approximate refinement of $\calP\wedge \calQ$. 
\end{lemma}

The proof of Lemma \ref{lem:refinement} appears in Appendix \ref{app:slicing}.  Finally, it will be useful to know that approximate refinements of homogeneous decompositions are still fairly homogeneous.

\begin{lemma}\label{lem:intersecting}
For all $\e_1>0$ and $\e_2:\mathbb{N}\rightarrow (0,1]$, there are $\e_1'>0$ and $\e_2':\mathbb{N}\rightarrow (0,1]$ such that for all $\ell, t, \ell',t'\geq 1$, there is $N$ such that the following holds.  Suppose $H=(V,E)$ is a $3$-graph on $n\geq N$ vertices, and $\calQ$ is a $(t,\ell,\e'_1,\e'_2(\ell))$-decomposition of $V$ which is $(\e'_1,\e'_2(\ell))$-regular and $\e'_1$-homogeneous with respect to $H$.  

Suppose $\calP$ is a $(t',\ell',\e_1',\e_2'(\ell'))$-decomposition of $V$ which is a $(\e'_1,\e'_2(\ell))$-approx\-imate refinement of $\calQ$, witnessed by the set $\Sigma\subseteq {[t']\choose 2}$.  Then $\calP$ is a $(t',\ell',\e_1,\e_2(\ell'))$-decomposition of $V$ which is $(\e_1,\e_2(\ell'))$-regular and $\e_1$-homogeneous with respect to $H$.
\end{lemma}

For a proof, see Appendix \ref{app:slicing}.  

\section{Ramsey facts}\label{subsec:ramseyfacts}

It will be convenient at various points to know some simple Ramsey facts about hereditary properties and the combinatorial definitions from the introduction.  For example, if $G$ is a $3$-graph and $\Hbar(k)$ is an induced sub-$3$-graph of $\trip(G)$, then by definition there are $a_1,\ldots, a_k\in V(G)$ and $b_1,\ldots, b_k,d\in V(G)$ and so that $a_ib_jd\in E(G)$ if and only if $i\leq j$. However, it is possible that the sets $\{a_1,\ldots, a_k\}$ and $\{b_1,\ldots, b_k\}$ are not disjoint, so we do not necessarily know the size of the set $\{a_1,\ldots, a_k,b_1,\ldots, b_k\}$. This can be inconvenient at times, and is the motivation for the following definition.

\begin{definition}
Suppose $H=(V_1\cup V_2\cup V_3, F)$ is a $3$-partite $3$-graph.  A \emph{clean copy of $H$} is a $3$-graph $G=(V,E)$ such that there exists a bijective map $f:V_1\cup V_2\cup V_3\rightarrow V$ so that for all $uvw\in K_3[V_1,V_2,V_3]$, $uvw\in F$ if and only if $f(u)f(v)f(w)\in E$.  
\end{definition}

Using standard Ramsey arguments, we will show that we can usually assume that there are clean copies of the configurations we are interested in.  For example, the following is a stronger version of Fact \ref{fact:vc2universal}, the universality property corresponding to $\IP_2$.

\begin{fact}\label{fact:vc2universal2}
If $\calH$ has $\IP_2$, then $\calH$ contains a clean copy of every $3$-partite $3$-graph. 
\end{fact}

For a proof, see Appendix \ref{app:ramsey}. Similarly, we will use the following facts.

\begin{lemma}\label{lem:clean}
Suppose $\calH$ is a hereditary $3$-graph property.
\begin{enumerate}[label=\normalfont(\arabic*)]
\item Suppose $H^*(k)\in \trip(\calH)$ for all $k\geq 1$.  Then $\calH$ contains a clean copy of $H^*(k)$ for all $k\geq 1$.
\item If $\calH$ has $\FOP_2$, then $\calH$ contains a clean copy of $F(k)$ for all $k\geq 1$.
\item Suppose $\Hbar(k)\in \trip(\calH)$ for all $k\geq 1$.  Then $\calH$ contains a clean copy of $\Hbar(k)$ for all $k\geq 1$.
\item Suppose $U^*(k)\in \trip(\calH)$ for all $k\geq 1$.  Then $\calH$ contains a clean copy of $U^*(k)$ for all $k\geq 1$.
\end{enumerate}
\end{lemma}

For proofs of these facts, see Appendix \ref{app:ramsey}.  Before being able to state the condition we will use for $\Ubar(k)$, we need to make one definition.

\begin{definition}
Suppose $n\geq 1$ and $G=(U\cup W, E)$ is a bipartite graph.  Set
$$
n\otimes G:=(\{U\cup W\cup \{c_1,\ldots, c_n\}, \{c_iuw: uw\in E\}).
$$
\end{definition}

In other words, $n\otimes G$ is the $3$-graph obtained by adjoining $n$ new vertices to the graph $G$.  For example, note that $\Hbar(k)=k\otimes H(k)$ and $\Ubar(k)=k\otimes U(k)$.  We will show that arbitrarily large $\Ubar(k)$ come with a universality property for $3$-graphs of the form $n\otimes G$.

\begin{lemma}\label{lem:ubarkuniversal}
Suppose $\calH$ is a hereditary $3$-graph property and $\Ubar(k)\in \trip(\calH)$ for all $k\geq 1$.  Then for every bipartite graph $G$ and $n\geq 1$, $\calH$ contains a clean copy of $n\otimes G$.
\end{lemma}

Again, we refer the reader to Appendix \ref{app:ramsey} for a proof. 

%% file: chapter4.tex
\chapter{$\VC_2$-dimension and $\disc_{2,3}$-homogeneity}\label{sec:dischom}

The goal of this chapter is to show that a hereditary $3$-graph property is $\disc_{2,3}$-homogeneous if and only if it has finite $\VC_2$-dimension.  

As a first step, we show that a sufficiently $\disc_{2,3}$-regular triad in a large $3$-graph of bounded $\VC_2$-dimension must have density near $0$ or $1$ (see Proposition \ref{prop:trivdensity}).  The crucial ingredient for this is the following counting lemma for induced $3$-partite $3$-graphs.  It follows easily from Theorem \ref{thm:countinginducedHG}.   For completeness, we include a proof in Appendix \ref{app:slicing}.

\begin{theorem}[Counting lemma for induced $3$-partite $3$-graphs]\label{thm:counting}
For all $\xi$, $d_3>0$ and $t\geq 1$, there exists $\delta_3>0$ such that for every $d_2>0$, there exists $\delta_2>0$ and $N\geq1$ such that the following holds.

Suppose $t_1,t_2,t_3\geq 1$, $t_1+t_2+t_3=t$, and $F=(U\cup W\cup Z,R_F)$ is a $3$-partite $3$-graph with $U=\{u_1,\ldots, u_{t_1}\}$, $W=\{w_1,\ldots, w_{t_2}\}$, and $Z=\{z_1,\ldots, z_{t_3}\}$.  Suppose $V$ is a set of size $n\geq N$, $H=(V,R)$ is a $3$-graph, and  $G=(V,E)$ is a $t$-partite graph, with vertex partition 
$$
V=A_1\cup \ldots \cup A_{t_1}\cup B_1\cup \ldots \cup B_{t_2}\cup C_1\cup \ldots \cup C_{t_3},
$$
 where each part has size $\frac{n}{t}(1\pm \delta_2)$.    For each $(i,j,s)\in [t_1]\times [t_2]\times [t_3]$, set $G^{ijs}=G[A_i,B_j,C_s]$ and $H^{ijs}=(A_i\cup B_j\cup C_s, K_3^{(2)}(G^{ijs})\cap R)$, and let $d_{ijs}$ be such that $e(H^{ijs})=d_{ijs}|K^{(2)}_3(G^{ijs})|$.  Assume that for each $(i,j,s)\in [t_1]\times [t_2]\times [t_3]$, each of $G[A_i,B_j]$, $G[A_i,C_s]$, and $G[B_j,C_s]$ have $\disc_2(\delta_2;d_2)$, $(H^{ijs},G^{ijs})$ satisfies $\disc_{2,3}(\delta_3,\delta_2)$, and 
\[
\begin{cases} d_{ijs}\geq d_3&\text{ if }u_iv_jz_s\in R_F,\text{ and }\\
d_{ijs}\leq (1-d_3)&\text{ if }\text{ if }u_iv_jz_s\notin R_F.
\end{cases}\]
Then the number of tuples $\abar\bbar\cbar\in \prod_{i=1}^{t_1}A_i\times \prod_{i=1}^{t_2}B_i\times \prod_{i=1}^{t_3}C_i$ such that $a_ib_jc_s\in E(H)$ if and only if $u_iw_jz_s\in R_F$ is 
$$
d_2^{t\choose 2}d(1\pm \xi)\prod_{i=1}^{t_1}|A_i| \prod_{i=1}^{t_2}|B_i| \prod_{i=1}^{t_3}|C_i|,
$$
where $d=(\prod_{u_iv_jw_s\in R_F}d_{ijs})(\prod_{u_iv_jw_s\in R_F}(1-d_{ijs}))$.
\end{theorem}

We now prove that $\disc_{2,3}$-regular triads in $3$-graphs of bounded $\VC_2$-dimension have density near $0$ or near $1$.  We do this by showing that a large $\disc_{2,3}$-regular triad with intermediate density can be broken up into sub-triads, all of which are regular and of intermediate density (this uses Lemma \ref{lem:subpairs}).  We then apply Theorem \ref{thm:counting} to the resulting configuration to find a copy of $V(k)$, and thus a contradiction.

\begin{proposition}[$\VC_2<\infty$ implies regular triples have trivial density]\label{prop:trivdensity}
For all $k\geq 1$ and $\e\in (0,1/2)$, there exists a $\delta_3>0$ such that for every $d_2>0$ there is $\delta_2>0$ and $M_0$ such that the following holds.  Suppose $V=V_1\sqcup V_2\sqcup V_3$, where $V$ is a set of size $n\geq M_0$, and $|V_i|=\frac{n}{3}(1\pm \delta_2)$ for each $i\in [3]$.  Suppose $G=(V_1\cup V_2\cup V_3, E)$ is a $3$-partite graph, and $H=(V,R)$ is a $3$-graph with $\VC_2(H)<k$.   Assume that for each $1\leq i< j\leq 3$, $G[V_i,V_j]$ has $\disc_2(\delta_2;d_2)$.  Let $H'=(V, R\cap K_3^{(2)}(G))$, and assume $(H',G)$ has $\disc_{2,3}(\delta_2,\delta_3)$. Then
$$
\frac{|R\cap K_3^{(2)}(G)|}{|K_3^{(2)}(G)|}\in [0,\e)\cup (1-\e,1].
$$
\end{proposition}

\begin{proof}. Fix $k\geq 1$ and $\e\in (0,1/2)$.  First, apply Theorem \ref{thm:counting} to obtain $\delta'''_3$ for $t=2k+2^k$, $\xi=1/2$, and $d_3=\e/2$.   Fix any $d_2>0$.  Apply Theorem \ref{thm:counting} to obtain $\delta'''_2$ and $n_0$ for $t=2^{k^2}+k^2$, $\xi$, $d_3$, $\delta_3'''$, and $d_2$.  Let  $\delta_3''=\min\{\delta_3''', \e/4\}$ and $\delta_2''=\min\{\delta_2''',d_2^3/8\}$, and apply Lemma \ref{lem:subpairs} to $d_2, d_3, \delta_3'', \delta_2''$ to obtain $\delta'_2, \delta'_3$ and $m_0$.  Finally, set $\delta_2=\min\{\delta'_2/(3k2^{1+k^2}), \delta_2''/(3k2^{1+k^2})\}$, $\delta_3=\delta'_3$, and choose $M_0\geq \delta_2^{-1}\max\{m_0,n_0\}2^k$ sufficiently large so that $d_2^{t\choose 2}(\e/2)^{t\choose 3}(M_0)^t/2>0$.

Now  suppose $V=V_1\cup V_2\cup V_3$ where $|V|=n\geq M_0$, and for each $i\in [3]$, $|V_i|=\frac{n}{3}(1\pm \delta_2)$.   Assume $G=(V_1\cup V_2\cup V_3, E)$ a $3$-partite graph, and $H=(V,R)$ is a $3$-graph with $\VC_2(H)< k$.  Set $H'=(V, R\cap K_3^{(2)}(G))$.  Assume that for each $1\leq i<j\leq 3$, $G[V_i,V_j]$ has $\disc_2(\delta_2;d_2)$, and assume $(H',G)$ has $\disc_{2,3}(\delta_2,\delta_3)$.  Suppose towards a contradiction that 
$$
d:=\frac{|R\cap K_3^{(2)}(G)|}{|K_3^{(2)}(G)|}\in [\e,1-\e].
$$

For each $i\in \{1,2,3\}$, let $V_{i,1},\ldots, V_{i,2^k}$ be an equipartition of $V_i$.  Note each $V_{i,j}$ has size $\frac{|V_i|}{2^k}\pm 1= \frac{n}{3(2^k)}(1\pm \delta_22^{k+1})=\frac{n}{3(2^k)}(1\pm \delta_2')\geq \max\{m_0,n_0\}$.

For each $(u,v,w)\in [k]\times [k]\times [2^k]$, set $G_{uvw}=G[V_{1,u},V_{2,v},V_{3,w}]$ and 
$$
H_{uvw}=(V_{1,u}\cup V_{2,v}\cup V_{3,w}, R\cap K_3^{(2)}(G_{uvw})).
$$
By Lemma \ref{lem:sl}, each of the graphs $G[V_{1,u},V_{2,v}]$, $G[V_{1,u},V_{3,w}]$, and $G[V_{2,v},V_{3,w}]$ has $\disc_2(\delta_2'';d_2)$.  By Lemma \ref{lem:subpairs}, $(H_{uvw},G_{uvw})$ has $\disc_{2,3}(\delta''_2,\delta''_3)$ with density $d_{uvw}\in [\e-\delta_3'',1-\e+\delta_3'']$.  Since $\delta_3''\leq \e/4$, this implies $d_{uvw}\in (3\e/4, 1-3\e/4)$.  

Since $\delta_2''\leq \delta_2'''$ and $\delta_3''\leq \delta_3'''$, Theorem \ref{thm:counting} implies that there are at least $d_2^{t\choose 2}(\e/2)^{t\choose 3}m^t/2$ tuples $\vbar\wbar\zbar \in (\prod_{j=1}^{2^k}V_{1,j}) \times (\prod_{j=1}^{k}V_{2,j}) \times (\prod_{j=1}^{k}V_{3,j})$ such that $v_{i}w_{j}z_{s}\in R$ if and only if $v_iu_jw_s\in E(V(k))$.   By choice of $m$, $d_2^{t\choose 2}(\e/2)^{t\choose 3}m^t/2>0$, so $\trip(H)$ contains an induced copy of $V(k)$, contradicting that $\VC_2(H)< k$.
\end{proof}

\vspace{2mm}

We are now in a position to prove Theorem \ref{thm:vc2finite}, which shows that $3$-graphs with bounded $\VC_2$-dimension have $\disc_{2,3}$-homogeneous decompositions.

\vspace{2mm}

\begin{proofof}{Theorem \ref{thm:vc2finite}}\label{proof:vc2finite}
Fix $k\geq 1$.   Fix $\ell_0,t_0\geq 1$, $\e_1>0$ and $\e_2:\mathbb{N}\rightarrow (0,1]$.  Without loss of generality, let us assume $\e_2$ is non-increasing (e.g. by replacing $\e_2$ with a non-increasing function bounded above by the original $\e_2$).  Apply Proposition \ref{prop:trivdensity} to $\e=\e_1$ to obtain $\delta_3>0$.  For each $\ell\in \mathbb{N}$, let $\delta_2(\ell)$ and $M_0(\ell)$ be obtained by applying from Proposition \ref{prop:trivdensity} to $\e$, $\delta_3$, and $d_2=1/\ell$.  Let $\e_2''(x)$ be as in Corollary \ref{cor:counting} for $t=3$, $\e_1^2$ and density $1/x$.  Let $\delta_2'(x)$ be from Lemma \ref{lem:sl} for $\e_2''(x)$ and with density $1/x$.  Then define $\e_2':\mathbb{N}\rightarrow (0,1]$ by setting $\e_2'(\ell)=\min\{\e_2(\ell),\e_2''(\ell), \delta_2(\ell),\delta'_2(\ell)\}$ for each $\ell \in \mathbb{N}$, and let $\e_1'=\min\{\e_1,\delta_3\}$.  Choose $T=T(\e_1',\e_2',t_0)$ and $N=N(\e_1',\e_2',t_0)$  from Theorem \ref{thm:reg2}, and define $m_0=\max\{M_0(T), N\}$.  

Suppose $H=(V,E)\in \calH$ and $|V|\geq m_0$.  By Theorem \ref{thm:reg2}, there are $t_0\leq t\leq T$ and $\ell_0\leq \ell\leq L$ and a $(t,\ell,\e_1',\e_2'(\ell))$-decomposition, $\calP$, of $V$ which is $(\e_1',\e_2'(\ell))$-regular  with respect to $H$.   We show every $\disc_{2,3}$-regular triad of $\calP$ is $\e_1$-homogeneous with respect to $H$.  Suppose $G_{ijk}^{\alpha,\beta,\gamma}\in \triads(\calP)$, $H_{ijs}^{\alpha,\beta,\tau}:=(V_i\cup V_j\cup V_s, E\cap K_3^{(2)}(G_{ijk}^{\alpha,\beta,\tau}))$, and assume  $(H_{ijk}^{\alpha,\beta,\tau}, G_{ijk}^{\alpha,\beta,\tau}))$ has $\disc_{2,3}(\delta_3, \delta_2(\ell))$.  By Proposition \ref{prop:trivdensity}, the density of edges in this triad is in $[0,\e_1)\cup (1-\e_1,1]$.  Thus $\calP$ is $\e_1$-homogeneous.  Since $\e_1'<\e_1$ and $\e'_2(\ell)<\e_2(\ell)$, $\calP$ is also $(\e_1,\e_2(\ell))$-regular with respect to $H$. 
\end{proofof}

\vspace{2mm}

Equipped with Theorem \ref{thm:vc2finite}, we now prove Theorem \ref{thm:dischom}, which says that a hereditary $3$-graph property is $\disc_{2,3}$-homogeneous if and only if it has bounded $\VC_2$-dimension. 

\vspace{2mm}

\begin{proofof}{Theorem \ref{thm:dischom}}
\label{proof:dischom} Suppose first that $\VC_2(\calH)<\infty$, say $\VC_2(\calH)<k$.  Fix $\ell_0,t_0\geq 1$, $\e_1>0$ and  $\e_2:\mathbb{N}\rightarrow (0,1]$.  As in the proof of Theorem \ref{thm:vc2finite}, we may assume, without loss of generality, that $\e_2$ is non-increasing.  By Theorem \ref{thm:vc2finite}, there are $T,L,N$ so that for all $H=(V,E)\in \calH$ with $|V|\geq m_0$, there are $t_0\leq t\leq T$ and $\ell_0\leq \ell\leq L$ and a $(t,\ell,\e_1,\e_2(\ell))$-decomposition, $\calP$, of $V$ which is $(\e_1,\e_2(\ell))$-regular and $\e_1$-homogeneous with respect to $H$.   This shows $\calH$ is $\disc_{2,3}$-homogeneous.

For the converse, suppose $\VC_2(\calH)=\infty$. Fix $\ell_0,t_0\geq 1$ and $\e_1\in (0,8^{-4})$. Define $\e_2:\mathbb{N}\rightarrow (0,1]$ as follows.  Given $\ell\geq 1$, let $\gamma_\ell=\e_1^3/4\ell^3$.  Apply Corollary \ref{cor:counting} to $\gamma_\ell$ to obtain $\e_\ell>0$.  Let $\e_2(\ell)$ be as in Lemma \ref{lem:sl} for $d=1/\ell$ and $\delta=\e_{\ell}\e_1^3/\ell^3$.  

Fix any $T\geq t_0$ and $L\geq \ell_0$.  Set $\e=\e_2(L)/T^3L^3$, let $n_0$ be as in Fact \ref{fact:rhg} for $\e$.   We show that for all $n\geq n_0$ there is $H=(V,E)\in \calH$ with $|V|\geq n_0$ such that for every $\ell_0\leq \ell\leq L$ and $t_0\leq t\leq T$, and every  $(t,\ell,\e_1,\e_2(\ell))$-decomposition $\calP$ of $V$, if $\calP$ is $(\e_1,\e_2(\ell))$-regular with respect to $H$, then at least $\frac{1}{9}(1-\e_1^{1/5}){n\choose 3}$ triples $xyz\in {V\choose 3}$ are in a triad of $\calP$ with density in $(1/8,7/8)$ (so in particular $\calP$ cannot be $\e_1$-homogeneous with respect to $H$).
 
Fix $n\geq n_0$.  Let $H_n=(U\cup W\cup Z,E)$ be as in Fact \ref{fact:rhg} \label{fact:rhgref} for $\e$, where $U=\{u_1,\ldots, u_n\}$, $W=\{w_1,\ldots, w_n\}$, and $Z=\{z_1,\ldots, z_n\}$. By assumption, $(H_n, K_3[U,W,Z])$ has $\disc_3(\e)$ and density $1/2\pm \e$.  By Fact \ref{fact:vc2universal2}, there is a clean copy of $H_n$ in $\calH$.  In other words, there is $H=(A\sqcup B\sqcup C, R)\in\calH$ where $A=\{a_1,\ldots, a_n\}$, $B=\{b_1,\ldots, b_n\}$, and $C=\{c_1,\ldots, c_n\}$ such that $a_ib_jc_k\in R$ if and only if $u_iw_jz_k\in E$.  Let $H'=(A\cup B\cup C, R\setminus ({A\choose 2}\cup {B\choose 2}\cup {C\choose 2}))$, and note $H'$ satisfies $\disc_3(\e)$.

Now suppose $t_0\leq t\leq T$, $\ell_0\leq \ell\leq L$ and $\calP$ is a $(t,\ell,\e_1,\e_2(\ell))$-decomposition for $V=A\cup B\cup C$.  Let $V=V_1\cup \ldots \cup V_t$ be the corresponding vertex partition, and for each $1\leq i\leq t$, set $A_i=A\cap V_i$, $B_i=B\cap V_i$ and $C_i=C\cap V_i$.  Let $\Sigma$ be the set of $i_1i_2i_3\in {[t]\choose 3}$ such that for some $\sigma\in S_3$, $\min\{|A_{i_{\sigma(1)}}|, |B_{i_{\sigma(2)}}|, |C_{i_{\sigma(3)}}|\}\geq \sqrt{\e_1}n/t$.

Suppose $i_1i_2i_3\in \Sigma$, and $\alpha,\beta,\gamma\leq \ell$ are such that $(H',G')$ has $\disc_{2,3}(\e_1,\e_2(\ell))$ where $G'=G_{i_1i_2i_3}^{\alpha,\beta,\gamma}$ and $H'=H_{i_1i_2i_3}^{\alpha,\beta,\gamma}$.   Let $\sigma\in S_3$ be such that \[\min\{|A_{i_{\sigma(1)}}|, |B_{i_{\sigma(2)}}|, |C_{i_{\sigma(3)}}|\}\geq \sqrt{\e_1}n/t.\]  Set $G''=G'[A_{i_{\sigma(1)}}\cup B_{i_{\sigma(2)}}\cup C_{i_{\sigma(3)}}]$.  By Lemma \ref{lem:sl}, Corollary \ref{cor:counting}, and because each of $A_{i_{\sigma(1)}}$, $B_{i_{\sigma(2)}}$, $C_{i_{\sigma(3)}}$ are large,
\begin{align}\label{align:tri}
|K_3^{(2)}(G'')|=\frac{1}{\ell}(1\pm \e_1^2)|A_{i_{\sigma(1)}}||B_{i_{\sigma(2)}}||C_{i_{\sigma(3)}}|.
\end{align}
Since $H'$ has $\disc_3(\e)$, 
\begin{align}\label{dens}
\big||R\cap K_3(G'')|-\frac{1}{2}|K_3(G'')|\big|\leq \e n^3.
\end{align}
Since $(H',G')$ has $\disc_{2,3}(\e_1,\e_2(\ell))$, if $d$ is such that $|R\cap K_3^{(2)}(G')|=d|K_3^{(2)}(G')|$, then 
$$
\Big| |R\cap K_3^{(2)}(G'')|-d|K_3^{(2)}(G'')|\Big|\leq \frac{\e_1}{\ell^3}(n/t)^3.
$$
By the triangle inequality, (\ref{align:tri}), and the lower bounds on the sizes of $A_{i_{\sigma(1)}}$, $B_{i_{\sigma(2)}}$, and $C_{i_{\sigma(3)}}$, we can conclude
\begin{align*}
\Big| d|K_3^{(2)}(G'')|-\frac{1}{2}|K_3^{(2)}(G'')|\Big|&\leq \frac{\e_1}{\ell^3}(n/t)^3+\e n^3\\
&\leq 2\e_1^{1/4}|A_{i_{\sigma(1)}}||B_{i_{\sigma(2)}}||C_{i_{\sigma(3)}}|/\ell^3\\
&\leq 3\e_1^{1/4}|K_3^{(2)}(G'')|.
\end{align*}
Thus, $|d-\frac{1}{2}|\leq 3\e_1^{1/4}$, so $d\in (1/8,7/8)$ since $\e_1<8^{-4}$. 

Call a set $A_i, B_i, C_i$ \emph{trivial} if it has size at most $\e_1^{1/4} n/t$.  For each $i\in [t]$, at most two of $A_i, B_i, C_i$ can be trivial.  Consequently, the number of vertices in a trivial set is at most $3\e_1^{1/4} n$.  Let $\Gamma$ be the set of $xyz\in K_3[A,B,C]$ such that none of $x$, $y$ or $z$ are in a trivial set.  Then  $|\Gamma|\geq n^3-3\e_1^{1/4}n^3$.  Since $\calP$ is $(\e_1,\e_2(\ell))$-regular with respect to $H$, if $\Gamma'$ is the set of $xyz\in {V\choose 3}$  in a $\disc_{2,3}$-regular triad of $\calP$, then $|\Gamma'|\geq {3n\choose 3}-\e_1 (3n)^3\geq 9n^3(1-\e_1)$.  Therefore $|\Gamma|+|\Gamma'|>{3n\choose 3}+n^3(1-3\e_1^{1/4}-9\e_1)$, so
$$
|\Gamma\cap \Gamma'|\geq n^3(1-3\e_1^{1/4}-9\e_1)\geq n^3(1-4\e^{1/4})\geq {|V|\choose 3}\frac{1}{9}(1-4\e^{1/4}).
$$
Since every $xyz\in \Gamma\cap \Gamma'$ is in a $\disc_{2,3}$-regular triad $G_{ijk}^{\alpha,\beta,\tau}$ with density in $(1/8,7/8)$, this shows $\calP$ is not $\e_1$-homogeneous, and thus $\calH$ is not $\disc_{2,3}$-homogeneous.
\end{proofof}

\vspace{2mm}

The proof of Theorem \ref{thm:dischom} shows that properties with unbounded $\VC_2$-dim\-ension fail $\disc_{2,3}$-homogeneity in a rather strong way. In particular, we have shown that there exist $3$-graphs in which a $\frac{1}{9}$-proportion of triads are always non-homogeneous, and which are guaranteed to exist in any property of  unbounded $\VC_2$-dimension. 

%% file: chapter5.tex
\chapter{Slicewise NIP, slicewise stability, and $\vdisc_3$-homogeneity}\label{sec:vdischom}

This chapter contains the results on slicewise stable and SNIP properties.  We will first show in Section \ref{subsec:closenip} that the $\sim$-classes of slicewise stable properties are exactly  those $\calH$ for which $\Hbar(k)\notin \trip(\calH)$ for some $k$ (recall Definition \ref{def:close}). The analogous results characterizing SNIP in terms of $\Ubar(k)$ are proved similarly and appear  in the appendix.  In Section \ref{subsec:vdiscwnip}, we prove Theorem \ref{thm:vdischom}, characterizing the $\vdisc_3$-homogeneous properties as those close to a SNIP property.  In Section \ref{subsec:nonzerror} we show that a property admits zero $\disc_{2,3}$-error if and only if it admits zero $\vdisc_3$-error.  In the same section we show that any property containing $\Hbar(k)$ for all $k\geq 1$ requires non-zero $\disc_{2,3}$-error.  In Section \ref{subsec:strongstable}, we prove a strong version of the stable regularity lemma.  This will serve as a warm-up for the proof of Proposition \ref{prop:wsbinary}, and will also play a key role in the proof of Theorem \ref{thm:FOPfinite}.  Finally, in Section \ref{subsec:wsbinary}, we will prove Proposition \ref{prop:wsbinary}, which says that any property which is close to a slicewise stable property admits binary $\disc_{2,3}$-error.

\section{The $\sim$-classes of slicewise stable and slicewise NIP properties}\label{subsec:closenip}

In this section, we prove Theorem \ref{thm:simclassws}, which says that a $3$-graph property $\calH$ is close to a slicewise stable property if and only if for some $k\geq 1$, $\Hbar(k)\notin \trip(\calH)$.  The proof of Theorem \ref{thm:simclasswnip} (the analogous result for SNIP) is extremely similar, and is thus relegated to Appendix \ref{app:simclasses}.  We note that Theorems \ref{thm:simclassws} and Theorem \ref{thm:simclasswnip} are almost surely particular instances of a more general result showing that a property $\calH$ is close to $\calF$-free if and only if $\calH$ omits certain blow-ups of the elements in $\calF$.  As we were not able to find this exact type of general result proved in the literature, we include here the two special cases we require (namely Theorems \ref{thm:simclassws} and Theorem \ref{thm:simclasswnip}).

Given a $3$-partite $3$-graph $F$ on $n$ vertices, let $\calB(F)$ be the set of clean copies of $F$ with vertex set $[n]$.  Note that $\calB(F)$ is finite.

\vspace{2mm}

\proofof{Theorem \ref{thm:simclassws}.} Suppose $\Hbar(k)\in \trip(\calH)$ for all $k$.  Fix $\calH'$ a slicewise stable hereditary $3$-graph property.  We show that $\calH$ is far from $\calH'$.  

By assumption, there is an integer $k\geq 1$ such that $H^*(k)\notin \trip(\calH')$. Set $\e=(2k)^{-2k-1}/|\calB(H^*(k))|$, and consider any $n\gg \e^{-1}k$.  By Lemma \ref{lem:clean}, there is $G\in \calH$ a clean copy of $\Hbar(k)$.  In other words, $V(G)=U\sqcup V\sqcup W$, for some $U=\{u_1,\ldots, u_n\}$, $V=\{v_1,\ldots, v_n\}$, $W=\{w_1,\ldots, w_n\}$, and $u_iv_jw_m\in E(G)$ if and only if $j\leq m$.  For each $1\leq i\leq k$, let $V_i=\{v_{n i/k},\ldots, v_{n(i+1)/k}\}$ and $W_i=\{w_{n i/k},\ldots, w_{n(i+1)/k}\}$.  Let $m=|U\cup V\cup W|$, and note $m\geq n$.  We construct a large set $S$ of induced subgraphs of $G$ as follows. 
\begin{itemize}
\item Choose $u\in U$. There are $n$ choices.
\item For each $i\in [k]$, choose $v_i'\in V_i\setminus \{u\}$. There are at least $(\frac{n}{k}-1)^k\geq (n/2k)^k$ choices.
\item For each $1\leq i\leq k$, choose $w_i\in W_i\setminus \{u,v_1',\ldots, v_k'\}$.  There are at least $(\frac{n}{k}-k-1)^k\geq (\frac{n}{2k})^k$ ways to do this.  
\item Put $G[\{u\}\cup \{v_i': i\in [k]\}\cup \{w_i: i\in [k]\}]$ in $S$.
\end{itemize}
Clearly every element of $S$ is isomorphic to an element of $\calB(H^*(k))$ and
$$
|S|\geq n (n/2k)^{2k}>\e |\calB(H^*(k))|n^{2k+1},
$$
where the inequality is by definition of $\e$ and because $n$ is large.  By the pigeonhole principle, there is some $G_n\in \calB(H^*(k))$, such that $S$ contains at least $\e n^{2k+1}$ elements isomorphic to $G_n$.  Consequently, there is some $G'\in \calB(H^*(k))$ so that $\calB(H^*(k))=G_n$ for arbitrarily large $n$.  By definition, this implies $G'\in \Delta(\calH)$.  Since $H^*(k)\cong \trip(G')$ and $H^*(k)\notin \trip(\calH')$, $G'\notin \Delta(\calH')$ (in fact $G'\notin \calH$).  By Proposition \ref{prop:simclasses1}, $\calH\nsim \calH'$.

Suppose on the other hand that $\calH\nsim \calH'$ for every slicewise stable $\calH'$.  We will use the fact that, by Lemma \ref{lem:clean}, a property $\calH'$ is slicewise stable if and only if there exists $k\geq 1$ so that $\calH'$ contains no clean copy of $H^*(k)$.  

Fix $k\geq 1$.   By assumption, $\calH$ is not close to any slicewise stable property, so there is some $\delta>0$ such that for arbitrarily large $n$, there is $G(n)\in \calH_n$ and $\Gamma(n)\in \calB(H^*(k))$, such that $G(n)$ is not $\delta/|\calB(H^*(k))|$-close to $\Gamma(n)$-free. For each $\Gamma\in \calB(H^*(k))$, let $c_k(\Gamma), N_k(\Gamma)$ be as in Theorem \ref{thm:indrem} for $\Gamma$ and $\delta/|\calB(H^*(k))|$, then set $c_k=\min\{c_k(\Gamma): \Gamma\in \calB(H^*(k))\}$ and $N_k=\max\{N_k(\Gamma): \Gamma\in \calB(H^*(k))\}$.  Choose any $n\gg k, N_k$.  By Theorem \ref{thm:indrem}, $G(n)$ contains at least $c_k n^{2k+1}$ induced copies of some $\Gamma(n)$.  Say $V(\Gamma(n))=\{a,b_1,\ldots,b_k,c_1,\ldots, c_k\}$ and $ab_ic_j\in E(\Gamma(n))$ if and only if $i\leq j$.  Let $S$ be the set of induced sub-$3$-graphs isomorphic to $\Gamma(n)$ in $G(n)$.  For all $H\in S$, let $a^H, b_1^H,\ldots, b_k^H, c_1^H,\ldots, c_k^H\in V(H)$ be such that $a^Hb_i^Hc_j^H\in E(G(n))$ if and only if $i\leq j$. By assumption $|S|\geq c_k n^{2k+1}$.  Consequently, there is  $\ebar\fbar\in [n]^{2k}$ and $X\subseteq S$ with $|X|\geq \e n$, such that for all $H\in X$, $(b_1^H,\ldots, b_k^H, c_1^H,\ldots, c_k^H)=\ebar\fbar$.  Let $d_1,\ldots, d_k$ be any $k$ distinct elements in $X'$.  Then $d_ue_if_j\in E$ if and only if $i\leq j$, and consequently, $\trip(G)$ contains an induced copy of $\Hbar(k)$, so $\Hbar(k)\in \trip(\calH)$. 
\qed

\vspace{2mm}

\section{Slicewise $\VC$-dimension and $\vdisc_3$-homogeneity}\label{subsec:vdiscwnip}
In this section we prove Theorem \ref{thm:vdischom}, characterizing the properties which are $\vdisc_3$-homogeneous.  We begin by showing that if $H$ is a large $3$-graph such that $\trip(H)$ omits $\Ubar(k)$ for some $k\geq 1$, then there are restrictions on the behavior of $\disc_{2,3}$-regular triads in any sufficiently regular decomposition of $H$.

\begin{proposition}\label{prop:wnipdens}
For all $k\geq 1$ and $\e\in (0,1/2)$ there are $\e_1>0$ and $\e_2:\mathbb{N}\rightarrow (0,1]$ such that for all $T,L\geq 1$, there is $M$ such that the following holds.  

Suppose $H=(V,E)$ is a $3$-graph on at least $M$ vertices such that $\Ubar(k)$ is not an induced sub-$3$-graph of $\trip(H)$.  Suppose $ t\leq T$, $\ell\leq L$, and $\calP$ is a $(t,\ell,\e_1,\e_2(\ell))$-decomposition of $V$ which is  $(\e_1,\e_2(\ell))$-regular with respect to $H$.  Suppose $ijs\in {[t]\choose 3}$ is such that at least $(1-\sqrt{\e_1})|V_i||V_j||V_s|$ many $xyz\in K_3[V_i,V_j,V_s]$ are in a $\disc_{2,3}$-regular triad of $\calP$.  Then $d_{ijs}\in [0,\e)\cup (1-\e,1]$ where $d_{ijs}$ is such that $|E\cap K_3[V_i,V_j,V_s]|=d_{ijs}|V_i||V_j||V_s|$.
\end{proposition}
\begin{proof}
Fix $k\geq 1$ and $\e\in (0,1/2)$.  Choose $\delta_3^{\ref{prop:trivdensity}}$ by applying Proposition \ref{prop:trivdensity} to $k$ and $\e$.  Then for each $\ell\geq 1$, choose $\delta^{\ref{prop:trivdensity}}_2(\ell)$, and $M_\ell$ by applying Proposition \ref{prop:trivdensity} to $k$, $\e$, $\delta_3^{\ref{prop:trivdensity}}$, and $d_2=1/\ell$.  For each $\ell$, let $\e_\ell$ be obtained by applying Corollary \ref{cor:counting} to $\e^2/\ell^3$ and density $1/\ell$.  Apply Theorem \ref{thm:counting} to $t=2k+2^k$, $\xi=1/2$, and $d_3=1/2$ to obtain $\delta_3^{\ref{thm:counting}}$.  For each $\ell\geq 1$, let $\delta_2^{\ref{thm:counting}}(\ell)$ and $m_\ell$ be obtained by applying Theorem \ref{thm:counting} to $t=2k+2^k$, $\xi=1/2$, $d_3=1/2$, and $d_2=1/\ell$.  Now let $\delta_3^{\ref{lem:subpairs}}$ be obtained by applying Lemma \ref{lem:subpairs} to $t=2k+2^k$, then for each $\ell\geq 1$, let $\delta_2^{\ref{lem:subpairs}}(\ell)$ be obtained by applying Lemma \ref{lem:subpairs} to $t=2k+2^k$, $\delta_2=\delta_2'(\ell)$, and $\delta_3=\delta_3^{\ref{thm:counting}}$.  Now set $\e_1=\min\{\delta^{\ref{lem:subpairs}}_3,\delta^{\ref{prop:trivdensity}}_3, 3^{-16}, \e^{-10}/10\}$,  $\e_2(\ell)=\min\{\delta_2^{\ref{lem:subpairs}}(\ell),\delta_2^{\ref{prop:trivdensity}}(\ell),\e_\ell\}$, and let $M=TM_Tm_T$.

Suppose $H$ is a $3$-graph on at least $M$ vertices such that $\Ubar(k)$ is not an induced sub-$3$-graph of $\trip(H)$.  By Fact \ref{fact:vc2universal}, $\VC_2(H)\leq k$. Suppose $ t\leq T$, $\ell\leq L$, and $\calP$ is a $(t,\ell,\e_1,\e_2(\ell))$-decomposition of $V$ which is $(\e_1,\e_2(\ell))$-regular with respect to $H$.  For each triad $G_{ijs}^{\alpha,\beta,\tau}$ of $\calP$, let $d_{ijs}^{\alpha,\beta,\tau}$ be such that $|E\cap K_3^{(2)}(G_{ijs}^{\alpha,\beta,\tau})|=d_{ijs}^{\alpha,\beta,\tau}|K_3^{(2)}(G_{ijs}^{\alpha,\beta,\tau})|$.   For each $ijs\in {[t]\choose 3}$, let $I_{ijs}$ be the set of $(\alpha,\beta,\tau)\in [\ell]^3$ such that $(H_{ijk}^{\alpha,\beta,\gamma}, G_{ijs}^{\alpha,\beta,\tau})$ satisfies $\disc_{2,3}(\e_2(\ell),\e_1)$.  Let 
$$
\Sigma=\Big\{ijs\in {[t]\choose 3}: \Big|\bigcup_{(\alpha,\beta,\gamma)\in I_{ijs}}K_3^{(2)}(G_{ijs}^{\alpha,\beta,\gamma})\Big|\geq (1-\sqrt{\e_1})|K_3^{(2)}(G_{ijs})|\Big\}.
$$
Note that by our choice of parameters and by Proposition \ref{prop:trivdensity}, for every $(\alpha,\beta,\gamma)\in I_{ijs}$,  there is $\delta_{ijs}^{\alpha,\beta,\tau}\in \{0,1\}$ such that $|\delta_{ijs}^{\alpha,\beta,\tau}-d_{ijs}^{\alpha,\beta,\tau}|\leq \e$.

\begin{claim}\label{claim:noswitch}
For all $ijs\in \Sigma$, if $\tau\neq \tau'$ and $(\alpha,\beta,\tau),(\alpha,\beta,\tau')\in I_{ijs}$, then $\delta_{ijs}^{\alpha,\beta,\tau}=\delta_{ijs}^{\alpha,\beta,\tau'}$.  
\end{claim}
\begin{proofof}{Claim \ref{claim:noswitch}}
Suppose towards a contradiction there are $(\alpha,\beta,\tau)$ and $(\alpha,\beta,\tau')$ in $I_{ijs}$ with $\tau\neq \tau'$ such that $\delta_{ijs}^{\alpha,\beta,\tau}=0$ while $\delta_{ijs}^{\alpha,\beta,\tau'}=1$.  Let $m=n/t2^k$, and choose equipartitions $V_i=V_{i,1}\cup \ldots \cup V_{i,2^k}$, $V_j=V_{j,1}\cup \ldots \cup V_{j,2^k}$, and $V_s=V_{s,1}\cup \ldots \cup V_{s,2^k}$.  Let $\{S_1,\ldots, S_{2^k}\}$ enumerate the elements of $\calP([k])$.  We define a $(2k+2^k)$-partite graph $G$ with vertex set $V'=\bigcup_{u=1}^{2^k}V_{i,u}\cup \bigcup_{u=1}^k V_{j,u}\cup \bigcup_{u=1}^k V_{s,u}$ as follows.  For each $1\leq a,b\leq k$, and $1\leq u\leq 2^k$, set $G[V_{i,u},V_{j,a}]=P_{ij}^\alpha$,  $G[V_{i,u},V_{s,b}]=P_{is}^{\beta}$, and
\[
G[V_{j,a},V_{s,b}]=\begin{cases}P_{js}^{\tau'}&\text{ if }a\in S_u\\
P_{js}^\tau&\text{ if }a\notin S_u.
\end{cases}\]
Now define $H'$ to be the $3$-partite $3$-graph with vertex set $V'$ and whose edge set is defined as follows.
$$
E(H):=\bigcup_{1\leq a,b\leq k, 1\leq u\leq 2^k} E\cap K_3^{(2)}(G[V_{i,u},V_{j,a},V_{s,b}])).  
$$
By Lemma \ref{lem:subpairs}, and our choice of parameters, we have that for all $1\leq a,b\leq k$ and $1\leq u\leq 2^k$, $(H', G[V_{i,u},V_{j,a},V_{s,b}])$ has $\disc_{2,3}(\delta^{\ref{thm:counting}}_2(\ell),\delta^{\ref{thm:counting}}_3)$ and by assumption, the density of edges in $H'$ on $K_3(G[V_{i,u},V_{j,a},V_{s,b}])$ is at least $\e$ if $a\in S_u$ and at most $1-\e$ if $a\notin S_u$.  But now Theorem \ref{thm:counting} implies $\trip(H)$ contains an induced copy of $\Ubar(k)$, a contradiction.
\end{proofof}
\vspace{2mm}

Note that the same argument shows that if $ijs\in \Sigma$ and $(\alpha,\beta',\tau),(\alpha,\beta',\tau)\in I_{ijs}$ for some $\beta\neq \beta'$ (respectively $(\alpha',\beta,\tau),(\alpha',\beta,\tau)\in I_{ijs}$ for some  $\alpha\neq \alpha'$), then $\delta_{ijs}^{\alpha,\beta,\tau}=\delta_{ijs}^{\alpha,\beta',\tau}$ (respectively $\delta_{ijs}^{\alpha,\beta,\tau}=\delta_{ijs}^{\alpha',\beta,\tau}$). 

Now fix $ijs\in \Sigma$.  Let $J_{ijs}$ be the set of $\alpha\in [\ell]$ such that there are at least $(1-\e_1^{1/4})\ell^2$ many pairs $(\beta,\tau)\in [\ell]^2$, with $(\alpha,\beta,\tau)\in I_{ijs}$.  By an easy counting argument, $|J_{ijs}|\geq (1-\e_1^{1/4})\ell$.  Given $\alpha\in J_{ijs}$, let $J_{ijs}^\alpha$ be the set of $\beta\in [\ell]$ such that there are at least $(1-\e_1^{1/16})\ell$ many $\tau\in [\ell]$ such that $(\alpha,\beta,\tau)\in I_{ijs}$.  By an easy counting, $|J_{ijs}^\alpha|\geq (1-\e_1^{1/16})\ell$.   Given $\beta\in J_{ijs}^\alpha$, let $J_{ijs}^{\alpha,\beta}$ be the set of $\tau\in [\ell]$ such that $(\alpha,\beta,\tau)\in I_{ijs}$.  By definition, $|J_{ijs}^{\alpha,\beta}|\geq (1-\e^{1/16})\ell$.  

We show that there is a $\delta_{ijs}\in \{0,1\}$ so that for all $(\alpha,\beta,\gamma)$ with $\alpha\in J_{ijs}$, $\beta\in J_{ijs}^{\alpha}$, and $\gamma\in J_{ijs}^{\alpha,\beta}$, $\delta_{ijs}^{\alpha,\beta,\gamma}=\delta_{ijs}$.   Suppose $\alpha\in J_{ijs}$.  By Claim \ref{claim:noswitch}, for all $\beta\in J_{ij}^\alpha$, there is $\delta_{ijs}^{\alpha,\beta}\in \{0,1\}$ such that for all $\tau\in J_{ijs}^{\alpha,\beta}$, $\delta_{ijs}^{\alpha,\beta,\tau}=\delta_{ijs}^{\alpha,\beta}$.  Suppose $\beta\neq \beta'\in J_{ijs}^{\alpha}$.  By definition of $J_{ijs}^{\alpha}$, $| J_{ijs}^{\alpha,\beta}\cap J_{ijs}^{\alpha,\beta'}|\geq (1-2\e^{1/16})\ell$.  Thus there is some $\tau \in J_{ijs}^{\alpha,\beta}\cap J_{ijs}^{\alpha,\beta'}$.  Claim \ref{claim:noswitch} applied to $(\alpha,\beta,\tau)$ and $(\alpha,\beta',\tau)$ implies that $\delta_{ijs}^{\alpha,\beta}=\delta_{ijs}^{\alpha,\beta'}$.  Thus there is $\delta_{ijs}^\alpha$, such that for all $\beta\in J_{ijs}^\alpha$, $\delta_{ijs}^{\alpha,\beta}=\delta_{ijs}^{\alpha}$.  Now suppose $\alpha \neq \alpha' \in J_{ijs}$.  By definition of $J_{ijs}$, $|J_{ijs}^\alpha\cap J_{ijs}^{\alpha'}|\geq (1-2\e^{1/4})\ell$.  Thus there is some $\beta\in J_{ijs}^\alpha\cap J_{ijs}^{\alpha'}$.  By definition, $|J_{ijs}^{\alpha,\beta}\cap J_{ijs}^{\alpha',\beta}|\geq (1-2\e^{1/16})\ell$, so there is some $\tau \in J_{ijs}^{\alpha,\beta}\cap J_{ijs}^{\alpha',\beta}$.  Now Claim \ref{claim:noswitch} applied to $(\alpha,\beta,\tau)$ and $(\alpha',\beta,\tau)$ implies that $\delta_{ijs}^{\alpha,\beta}=\delta_{ijs}^{\alpha',\beta}$, and thus $\delta_{ijs}^\alpha=\delta_{ijs}^{\alpha'}$.  We have shown there is $\delta_{ijs}$ such that for all $\alpha\in J_{ijs}$, $\delta_{ijs}^\alpha=\delta_{ijs}$.  In other words, for all $\alpha\in J_{ijs}$, $\beta\in J_{ijs}^{\alpha}$, and $\gamma\in J_{ijs}^{\alpha,\beta}$, $\delta_{ijs}^{\alpha,\beta,\gamma}=\delta_{ijs}$.  

We now show that the triple $(V_i,V_j,V_s)$ has density near $\delta_{ijs}$ in $H$.  For simplicity, let us assume $\delta_{ijs}=1$ (the case $\delta_{ijs}=0$ is similar).  First, given $\alpha\in J_{ijs}$ and writing $G_{ijs}^\alpha=(V_i\cup V_j\cup V_s, P_{is}^\alpha\cup K_2[V_i,V_j]\cup K_2[V_j,V_s])$, we have that
\[|E\cap K_{3}^{(2)}(G_{ijs}^\alpha)| \geq  \sum_{\beta\in J_{ijs}^\alpha}\sum_{\tau\in J_{ijs}^{\alpha,\beta}}(1-\e)|K_3^{(2)}(G_{ijs}^{\alpha,\beta,\tau})|,\]
which, by Corollary \ref{cor:counting} and the bounds on $|J^{\alpha}_{ijs}|, |J_{ijs}^{\alpha,\beta}|$, is at least
\begin{align*}
&(1-\e^{1/16})(1-\e^{1/4})\ell^2\frac{1}{\ell^3}(1-\ell^2\e^2/\ell^2)|V_i||V_j||V_s|\\
&\geq (1-\e^{1/16})(1-\e^{1/4})(1-\e^2)\frac{1}{\ell}|V_i||V_j||V_s|,
\end{align*}
Thus
\begin{align*}
|E\cap K_{3}[V_i,V_j,V_s]|&\geq |J_{ijs}|(1-\e^{1/16})(1-\e^{1/4})(1-\e^2)\frac{1}{\ell}|V_i||V_j||V_s|\\
&\geq (1-\e^{1/16})(1-\e^{1/4})^2(1-\e^2)|V_i||V_j||V_s|\\
&\geq (1-\e)|V_i||V_j||V_s|,
\end{align*}
where the last inequality is by choice of $\e$.  This shows that for every $ijs\in \Sigma$, the triple $(V_i,V_j,V_s)$ is $\e$-homogeneous with respect to $H$, as desired.
\end{proof}

We are now in a position to prove Theorem \ref{thm:vdiscubark}, which states that a hereditary $3$-graph property $\calH$ is $\vdisc_3$-homogeneous if and only if there is $k\geq 1$ such that $\Ubar(k)\notin \trip(\calH)$.

\vspace{2mm}

\proofof{Theorem \ref{thm:vdiscubark}}
\label{proof:vdiscubark} Suppose first that there is $k\geq 1$ such that $\Ubar(k)\notin \trip(\calH)$.  Fix $0<\e<1/100$ and $t_0\geq 1$.   Let $\delta_3'$ and $\delta_3''$ be obtained by applying Proposition \ref{prop:trivdensity} and Proposition \ref{prop:wnipdens} respectively to $\e$ and $k$, and set $\delta_3=\min\{\delta_3',\delta_3''\}$.  Given $\ell$, let $\delta'_2(\ell)$, $m'_\ell$ be obtained by applying Proposition \ref{prop:trivdensity}  to $\e$, $k$, $\delta'_3$ and $d_2=1/\ell$.  Given $\ell$, let  $\delta_2''(\ell)$ and $m_{\ell}''$ be obtained by applying Proposition \ref{prop:wnipdens} to $\e$, $k$, and $\delta_3''$ with $d_2=1/\ell$.

Set $\delta_2(\ell)=\min\{\delta'_2(\ell),\delta_2''(\ell)\}$, and $\e_1=\min\{\delta_3,\e^2\}$, and define $\e_2:\mathbb{N}\rightarrow (0,1]$ so that $\e_2(\ell)=\min\{\delta_2(\ell),\e/(3\ell)^3\}$.   Let $T=T(\e_1,\e_2,t_0,1)$, $L=L(\e_1,\e_2,t_0,1)$ and $N_1=N_1(\e_1,\e_2,t_0,1)$ be as in Theorem \ref{thm:reg2}.  Let $N_2$ be sufficiently large so that ${N_2\choose 3}\geq N_2^3/10$.  Now let $M=\max\{N_1,N_2,m'_TT, m''_TT\}$.  

Now suppose $H=(V,E)\in \calH$ with $|V|=n\geq M$.  By Theorem \ref{thm:reg2}, there is $t_0\leq t\leq T$ and $1\leq \ell\leq L$ and a $(t,\ell,\e_1,\e_2(\ell))$-decomposition $\calP$ of $V$ which is $(\e_1,\e_2(\ell))$-regular with respect to $H$.   Let $\Gamma$ be the set of $xyz\in {V\choose 3}$, such that $xyz$ is not in a $\disc_{2,3}$-regular triad of $\calP$.  Let $\Sigma$ be the set of $ijs\in {[t]\choose 3}$ such that $|\Gamma\cap K_3[V_i,V_j,V_s]|\geq \sqrt{\e_1}|V_i||V_j||V_s|$.  Clearly $|\Sigma|\leq \sqrt{\e_1} t^3$ because $\calP$ is $(\e_1,\e_2(\ell))$-regular with respect to $H$.  By Proposition \ref{prop:wnipdens}, for all $ijs\in {[t]\choose 3}\setminus \Sigma$, the triple $(V_i,V_j,V_k)$ is $\e$-homogeneous with respect to $H$.  Since $\sqrt{\e_1}\leq \e$, this shows that $\calP_{vert}$ is an $\e$-homogeneous partition with respect to $\calH$.  We have now shown that $\calH$ is $\vdisc_3$-homogeneous.

For the converse, suppose $\Ubar(k)\in \trip(\calH)$ for all $k\geq 1$.  Fix $\e<1/8$ and $t_0\geq 1$, and set $M_0=\max\{\e^{-1},t_0\}$.   By Fact \ref{fact:rg}, for all $T\geq t_0$ there is $n_0$ such that for all $n\geq n_0$, there is a bipartite graph $G=(U\cup V, E)$ with $|U|=|V|=n$, and which satisfies $\disc_2(\e^2/T^2;1/2)$.  

Suppose $T\geq t_0$ and $n\geq n_0M_0$.  Fix a bipartite graph $G=(U\cup V, E)$ with $U=\{u_1,\ldots, u_n\}$, $W=\{w_1,\ldots, w_n\}$ which satisfies $\disc_2(\e^2/T;1/2)$.   By Lemma \ref{lem:ubarkuniversal}, $\calH$ contains a clean copy of $n\otimes G$.  Thus there exists some $H=(V,R)\in \calH$, and sets $B=\{b_1,\ldots, b_n\}, A=\{a_1,\ldots, a_n\}, C=\{c_1,\ldots, c_n\}$ so that $V=A\sqcup B\sqcup C$, and such that $a_sb_ic_j\in R$ if and only if $u_iw_j\in E$.  We show there is no equipartition of $V$ into at least $t_0$ and at most $T$ parts which is $\e$-homogeneous with respect to $H$.  

Suppose towards a contradiction there is an integer $t_0\leq t\leq T$, an equipartition $\{V_1,\ldots, V_t\}$ of $V$,  and a set $\Sigma\subseteq {[t]\choose 3}$, such that $|\Sigma|\leq \e t^3$ and such that for all $ijk\in {[t]\choose 3}\setminus \Sigma$, $|R\cap K_3[V_i,V_j,V_k]|/|K_3[V_i,V_j,V_k]|\in [0,\e)\cup (1-\e,1]$.  Given $ijs\in {[t]\choose 3}$, set 
$$
d_{ijs}:=\frac{|E\cap K_3(V_i,V_j,V_s)|}{|V_i||V_j||V_s|}.
$$
For each $i\in [t]$, let $A_i=A\cap V_i$, $B_i=B\cap V_i$, and $C_i=C\cap V_i$.  Let $\Sigma'$ be the set of $i_1i_2i_3\in {[t]\choose 3}$ such that for some $\sigma \in S_3$, $\min\{|A_{i_{\sigma(1)}}|,|B_{i_{\sigma(2)}}|,|C_{i_{\sigma(3)}}|\}\geq \e^{1/4}n/t$.

We show that for all $i_1i_2i_3\in \Sigma'$, $d_{i_1i_2i_3} \in (\e,1-\e)$.  Fix $i_1i_2i_3\in \Sigma'$ and let $\sigma \in S_3$ be such that $\min\{|A_{i_{\sigma(1)}}|,|B_{i_{\sigma(2)}}|,|C_{i_{\sigma(3)}}|\}\geq \e^{1/4}n/t$. Let $U'=\{u_j: b_j\in B_{i_{\sigma(2)}}\}$ and $W'=\{w_j: c_j\in C_{i_{\sigma(3)}}\}$.  Then since $G$ has $\disc_2(\e^2/T^2; 1/2)$, the following holds.
\begin{align}\label{rg}
\big||E\cap K_2(U',V')|-|U'||V'|/2\big|\leq \frac{\e^2}{T^2} n^2.
\end{align}
Note that because $i_{\sigma(2)}\neq i_{\sigma(3)}$, $B_{i_{\sigma(2)}}$, $C_{i_{\sigma(3)}}$, and $A_{\sigma(1)}$  are disjoint, thus
\[|R\cap K_3[A_{i_{\sigma(1)}},B_{i_{\sigma(2)}}, C_{i_{\sigma(3)}}]|=|A_{\sigma(1)}||E\cap K_2[U',V']\]
and
\[K_3[A_{i_{\sigma(1)}},B_{i_{\sigma(2)}}, C_{i_{\sigma(3)}}]\setminus R|=|A_{\sigma(1)}||K_2[U',V']\setminus E|.\]
Combining this with (\ref{rg}), we have the following.
$$
\frac{|A_{\sigma(1)}|}{(n/t)^3}\Big(\frac{|U'||V'|}{2}-\frac{\e^2n^2}{T^2}\Big)\leq d_{i_1i_2i_3}\leq\frac{|A_{\sigma(1)}|}{(n/t)^3}\Big(\frac{|U'||V'|}{2}+\frac{\e^2n^2}{T^2}\Big).
$$
By assumption, we have $\min\{|A_{i_{\sigma(1)}}|,|U'|,|W'|\}\geq \e^{1/4}n/t$, whilst at the same time $n/t\geq \max\{|A_{i_{\sigma(1)}}|,|U'|,|W'|\}$.  Consequently,
$$
\e^{3/4}/2-\e^{5/4}t^2/T^2\leq d_{ijs}\leq \frac{1}{2}+\frac{\e^2t^2}{T^2}.
$$
Since $t\leq T$ and $\e<1/8$, this implies $d_{ijs}\in (\e,1-\e)$.  

We have now shown that for all $ijs\in \Sigma'$, $d_{ijs}\in (\e,1-\e)$, while by assumption, for all $ijs\notin \Sigma$, $d_{ijs}\in [0,\e)\cup (1-\e,1]$.  Consequently, $\Sigma'\subseteq \Sigma$, and therefore, $|\Sigma'|\leq \e t^3$.  

Given $i\in [t]$, we say $A_i$, $B_i$, or $C_i$ is \emph{trivial} if it has size at most $\e^{1/4}|V_i|$.  Note that for each $i\in [t]$, at most two of $A_i$, $B_i$, or $C_i$ can be trivial.   Consequently, if $V_0$ is the union of all trivial $A_i,B_i,C_i$, then $|V_0|\leq 2\e^{1/4}|V|$.   Note that if $ijs\in {[t]\choose 3}\setminus \Sigma'$, then $K_3[V_i,V_j,V_s]\subseteq K_3[V_0,V,V]$, and consequently $|\bigcup_{ijs\in {[t]\choose 3}\setminus \Sigma'}K_3[V_i,V_j,V_s]|\leq 2\sqrt{\e}|V|^3$.  But now we have that
\[|\bigcup_{ijs\in {[t]\choose 3}\setminus \Sigma'}K_3[V_i,V_j,V_s]|+ |\bigcup_{ijs\in \Sigma'}K_3[V_i,V_j,V_s]|\leq 2\sqrt{\e}|V|^3+\e t^3(|V|/t)^3 \leq 3\sqrt{\e}|V|^3,\]
Clearly this is a contradiction, as the left hand side is bounded below by ${|V|\choose 3}$, and our assumptions on $\e$ and $|V|$ imply that ${|V|\choose 3}>3\sqrt{\e}|V|^3$.  This finishes our proof that there is no equipartition of $V$ into $t_0\leq t\leq T$ parts which is $\e$-homogeneous with respect to $H$.  This shows that $\calH$ is not $\vdisc_3$-homogeneous.
\qed

\vspace{2mm}

Theorem \ref{thm:vdischom} follows immediately from Theorem \ref{thm:vdiscubark} and Theorem \ref{thm:simclasswnip}.   We end this section by proving Proposition \ref{prop:wnipreduction}, which gives an equivalent characterization of properties admitting zero/binary $\vdisc_{3}$-error.  For this, it will be useful to have a uniform way of talking about decompositions with zero and binary $\vdisc_3$-error.

\begin{definition}
Given $G=(V,E)$ a $3$-graph, an \emph{$\e$-homogeneous partition of $G$ with zero (respectively, binary) error} is an equipartition $V=V_1\cup \ldots \cup V_t$ along with a set $\Sigma \subseteq {[t]\choose 2}$ such that the following hold.
\begin{enumerate}
\item $|\Sigma|=0$ (respectively, $|\Sigma|\leq \e t^2$).
\item For all $ijk\in {[t]\choose 3}$ with $ij,jk,ik\notin \Sigma$, $\frac{|E\cap K_3[V_i,V_j,V_k]|}{|K_3[V_i,V_j,V_k]|}\in [0,\e)\cup (1-\e,1]$.
\end{enumerate}
\end{definition}

Observe that a hereditary $3$-graph property $\calH$ is $\vdisc_3$-homogeneous with zero (respectively, binary) error if and only if for all $\e>0$ and $t_0\geq 1$ there are $T,N\geq 1$ such that for every $G=(V,E)\in \calH$ with $|V|\geq N$, there is an equipartition of $V$ into at most $T$ parts which is a $\e$-homogeneous partition of $G$ with zero (respectively, binary error).

We now prove Proposition \ref{prop:wnipreduction}, which shows that a hereditary $3$-graph property $\calH$ admits zero (respectively binary) $\vdisc_3$-error if and only if $\calH$ is close to a SNIP property and admits zero (respectively binary) $\disc_{2,3}$-error.

\vspace{2mm}

\begin{proofof}{Proposition \ref{prop:wnipreduction}}
\label{proof:wnipreduction}Suppose first that $\calH$ is a hereditary $3$-graph property which is SNIP and admits zero (respectively, binary) $\disc_{2,3}$-error.  Say $\SVC(\calH)<k$.  By Theorem \ref{thm:vdischom}, $\calH$ is $\vdisc_3$-homogeneous.   Fix $\e>0$ and $t_0\geq 1$.  Without loss of generality, assume $\e<1/2$.  Choose $\e_1>0$ and $\e_2:\mathbb{N}\rightarrow (0,1]$ as in Proposition \ref{prop:wnipdens} for $k$ and $\e/2$. Let $M_1, T, L$ be as in the definition of $\calH$ is admitting zero (respectively binary) $\disc_{2,3}$-error for $\e_1,\e_2, t_0$ and $\ell_0=1$.  Let $M_2=M(k,\e,T,L)$ from Proposition \ref{prop:wnipdens}.  

Now suppose $G=(V,E)\in \calH$ satisfies $|V|\geq \max\{M_1,M_2\}$.  By assumption, there are $t_0\leq t\leq T$, $1\leq \ell\leq L$, and $\calP$ a $(t,\ell,\e_1,\e_2(\ell))$-decomposition of $V$ which is $(\e_1,\e_2(\ell))$-regular with respect to $H$, and $\Sigma\subseteq {[t]\choose 2}$ of size $0$ (respectively, of size at most $\e_1 t^2$), such that for every $P_{ijs}^{\alpha,\beta,\gamma}\in \triads(\calP)$ which is not $\disc_{2,3}$-regular with respect to $H$, one of $ij, is, jk\in \Sigma$.  By Proposition \ref{prop:wnipdens}, for every $ijs\in {[t]\choose 3}$ with $ij, js, is\notin \Sigma$, $|E\cap K_3[V_i,V_j,V_s]|/|K_3[V_i,V_j,V_s]|\in [0,\e)\cup (1-\e,1]$.   Thus $V_1\cup \ldots \cup V_t$ is an $\e$-homogeneous partition of $G$ with zero (respectively, binary) error.   This shows that if $\calH$ is SNIP and admits zero (respectively, binary) $\disc_{2,3}$-error, then $\calH$ admits zero (respectively, binary) $\vdisc_3$-error.  It follows immediately from Proposition \ref{prop:simclasses2} that the same holds for $\calH'$ close to a SNIP property.  

Suppose conversely that $\calH$ admits zero (respectively, binary) $\vdisc_{3}$-error.  By Theorem \ref{thm:vdischom}, $\calH$ is close to a SNIP property.  Fix $\e_1>0$, $\e_2:\mathbb{N}\rightarrow (0,1]$, and $t_0, \ell_0\geq 1$.  Given $x\in \mathbb{N}$, let $\delta_2(x)$ be the $\delta_2$ from Proposition \ref{prop:homimpliesrandome} for $\e_1/2\ell_0$ and $d_2=1/x$.  Let $\e_2''$ be as in Corollary \ref{cor:counting} for density $1/\ell_0$ and $\min\{\e_2(\ell_0)/\ell_0^3, \delta_2(\ell_0)\}$.  Let $\e_2'=\min\{\e_2'', \e_2(\ell_0)/\ell_0, \delta_2(\ell_0)\}$.  Let $T$ and $N_1$ be as in the definition of $\calH$ admitting zero (respectively, binary) $\vdisc_{3}$-error for $\e=\e^2_1/2\ell^3_0$ and $t_0$.  Let $N_2$ be the $m(\e'_2,\e_2',\e_2', \ell_0)$ from Lemma \ref{lem:3.8}, and set $N=\max\{N_1,N_2T\}$.

Suppose $H=(V,E)\in \calH$ and $|V|\geq N$.  By assumption, there exists $t_0\leq t\leq T$ and $\calP=\{V_1,\ldots, V_t\}$ an $(\e^2_1/2\ell^3_0)$-homogeneous partition for $H$ with zero (respectively, binary) error.  For each $ij\in {[t]\choose 2}$, by Lemma \ref{lem:3.8}, there exists a partition $K_2[V_i,V_j]=\bigcup_{\alpha\leq \ell_0}P_{ij}^{\alpha}$ such that each $P_{ij}^{\alpha}$ has $\disc_2(\e_2'(\ell_0);1/\ell_0)$.   Suppose $ijk\in {[t]\choose 3}$ has the property that none of $ij, ik,jk$ are $\Sigma$.  Then by assumption, there is $u\in \{0,1\}$ such that $|E^u\cap K_3[V_i,V_j,V_k]|\geq (1-\e^2_1/2\ell^3_0)|V_i||V_j||V_k|$.  By Corollary \ref{cor:counting}, for any $\alpha,\beta,\gamma\leq \ell_0$, 
$$
|K_3^{(2)}(G_{ijk}^{\alpha,\beta,\gamma})|=\frac{1}{\ell_0^3}(1\pm \e_2(\ell_0)/\ell^3_0)|V_i||V_j||V_k|,
$$
 so 
 $$
 |E^u\cap K_3^{(2)}(G_{ijk}^{\alpha,\beta,\gamma})|\geq |K_3^{(2)}(G_{ijk}^{\alpha,\beta,\gamma})|-(\e^2_1/2\ell^3_0)|V_i||V_j||V_k|\geq (1-\e^2_1)|K_3^{(2)}(G_{ijk}^{\alpha,\beta,\gamma})|.
 $$
 By Proposition \ref{prop:homimpliesrandome}, $(H_{ijk}^{\alpha,\beta,\gamma},G_{ijk}^{\alpha,\beta,\gamma})$ satisfies $\disc_{2,3}(\e_1,\e_2(\ell_0))$, as desired.  This shows that if $\calH$ admits zero (respectively, binary) $\vdisc_3$-error, then $\calH$ is close to a SNIP property and admits zero (respectively, binary) $\vdisc_{2,3}$-error.
\end{proofof}

\section{Slicewise stability and non-zero error}\label{subsec:nonzerror}
In this section, we prove results concerning slicewise stability and the notions of zero $\disc_{2,3}$-error and zero $\vdisc_3$-error.  We begin by showing that a property may be unstable but still admit zero $\disc_{2,3}$-error.  

\begin{example}\label{ex:stable}
Let $H$ be the $3$-graph with vertex set $\{a_i,b_i,c_i:i\in \mathbb{N}\}$ and edge set $\{a_ib_jc_k:i\leq j= k\}$.  Let  $\calH=\age(H)$. Then $R(x;y,z)$ is unstable in $\calH$, but $\calH$ admits zero $\disc_{2,3}$-error. 
\end{example}
\begin{proof}
Fix $\ell_0,t_0\geq 1$, $\e_1>0$ and $\e_2:\mathbb{N}\rightarrow (0,1]$.  Let $\e_2':\mathbb{N}\rightarrow (0,1]$ be such that for all $\ell\geq 1$, $\e'_2(\ell)\leq \min\{\e_2(\ell),\e^4_1/\ell^4,\e_\ell\}$, where $\e_{\ell}$ comes from applying Corollary \ref{cor:counting} to $\gamma=\e_1^4/10\ell^3$ and $t=3$.   Let $N_1,T,L$ be as in Proposition \ref{prop:binaryboring} for $\e_1^4$, $\e_2'$ and $t_0,\ell_0$.  Let $N_2$ be sufficiently large so that $N_2^2<\e_1^4 N_2^3/2T^3L^3$.  Let $N=\max\{N_1,N_2\}$.  

Suppose $H=(V,E)\in \calH$, where $|V|=n\geq N$.  Then by definition of $\calH$, we can write  $V(H)=V_A\cup V_B\cup V_C$ where $V_A=\{a_{i_1},\ldots, a_{i_s}\}$, $V_B=\{b_{j_1},\ldots, b_{j_t}\}$ and $V_C= \{c_{r_1},\dots, c_{r_m}\}$ for some $s+t+m=n$, and $E(G)=\{(a_{i_u},b_{j_v},c_{r_k}): i_u\leq j_v=r_k\}$.  Note $|E(H)|\leq s\max\{t,m\}\leq n^2/4$, where the second inequality is by the AM-GM inequality.  Since $n\geq N_1$, there exists $t_0\leq t\leq T$ and $\ell_0\leq \ell\leq L$, and a $(t,\ell,\e^4_1,\e'_2(\ell))$-decomposition $\calP$ of $V$ which is $(\e^4_1,\e'_2(\ell))$-regular with respect to $H$, and which has no $\disc_2$-irregular triads. Fix a triad $G_{ijk}^{\alpha,\beta,\tau}$ of $\calP$.  By assumption, each of $P_{ij}^\alpha,P_{ik}^\beta,P_{jk}^\tau$ has $\disc_2(\e'_2(\ell);1/\ell)$.  Combining Corollary \ref{cor:counting} with the bound on $E(H)$ implies the following.
\begin{align*}
\frac{|E\cap K_3^{(2)}(G_{ijk}^{\alpha,\beta,\gamma})|}{|K_3^{(2)}(G_{ijk}^{\alpha,\beta,\gamma})|}\leq \frac{n^2}{4|K_3^{(2)}(G_{ijk}^{\alpha,\beta,\gamma})|}\leq \frac{\e^4_1|V_i||V_j||V_k|}{T^3L^3}\leq \e^2_1,
\end{align*}
where the second inequality is by our choice of $N$ and Corollary \ref{cor:counting}.  By Proposition \ref{prop:homimpliesrandome}, $(H_{ijk}^{\alpha,\beta,\gamma}, G_{ijk}^{\alpha,\beta,\gamma})$ satisfies $\disc_{2,3}(\e_1,\e_2(\ell))$.  This shows $\calH$ admits zero $\disc_{2,3}$-error.
\end{proof}

 Our next goal is to show that if a hereditary $3$-graph property $\calH$ has $\Hbar(k)\in \trip(\calH)$ for all $k$, then $\calH$ requires non-zero $\disc_{2,3}$-error.  The main reason for this is the following theorem, which shows that any equipartition of the vertex set of $\Hbar(k)$ has a triple with intermediate density.

We would like to acknowledge here a (no longer available) blog post written by Aleksandar Makelov, which presented a proof that $H(k)$ requires irregular pairs in the graph regularity lemma.  A proof of this fact does not appear in the published literature (to the authors' knowledge), and we thus found this blog post a very useful starting point as we developed the following proof, as well the other proofs in this paper showing that certain examples require certain kinds of irregular triples/triads.

\begin{theorem}\label{thm:hbarksimple}
For all $0<\e_1\leq 2^{-18}$, and all $T,N\geq 1$, there is $n\geq N$ such that for every $t\leq T$ and equipartition $\calP=\{V_1,\ldots, V_t\}$ of $V(\Hbar(n))$, the following holds.  There is a triple $ijk\in {[t]\choose 3}$, and sets $A_{i}\subseteq V_i$, $B_{j1},B_{j0}\subseteq V_j$, and $C_{k1},C_{k0}\subseteq V_k$, each of size at least $\e^{1/9}(3n/t)$, such that $K_3[A_{i},B_{j1},C_{k1}]\subseteq E$ and  $K_3[A_{i},B_{j0},C_{k0}]\cap E=\emptyset$.
\end{theorem}
\begin{proof}
Fix $0<\e_1<2^{-18}$ and $T,N\geq 1$.   Let  $n_0=\e_1^{-100}T^2$ and assume $n\gg n_0$ and $t\leq T$.  Suppose $\calP=\{V_1,\ldots, V_t\}$ is an equipartition of $V(\Hbar(n))$.  Say $V=V_1\cup \ldots \cup V_t$ is the partition from $\calP$.   Fix a labeling $V(\Hbar(k))=A\cup B\cup C$ where $A=\{a_i: i\in [n]\}$, $B=\{b_i: i\in [n]\}$, $C=\{c_i: i\in [n]\}$ are such that $b_ic_j\in E(\Hbar(k))$ if and only if $i\leq j$.  For each $1\leq i\leq t$, define $A_i=A\cap V_i$, $B_i=B\cap V_i$, and $C_i=C\cap V_i$, and set 
\[\calA=\{i\in [T]:|A_i|\geq \e_1^{1/9}n/t\},\text{ }\calB=\{i\in [T]:|B_i|\geq \e_1^{1/9}n/t\},\text{ and }\]
\[\calC=\{i\in [T]:|C_i|\geq \e_1^{1/9}n/t\}.\]
Set $A'=\bigcup_{i\in \calA}A_i$, $B'=\bigcup_{i\in \calB}B_i$, and $C'=\bigcup_{i\in \calC}C_i$.  Observe that
$$
|A\setminus A'|\leq \sum_{i\notin \calA}\e_1^{1/9}n/t\leq \e_1^{1/9}n.
$$
Similarly, $|B\setminus B'|\leq \e_1^{1/9}n$ and  $|C\setminus C'|\leq \e_1^{1/9}n$.  Given $i\in [t]$, and $\mu>0$, let
\[B_i(\mu)=\{j\in [n]:b_j\in B_i, \min\{|\{b_1,\ldots, b_j\}\cap B_i|,|\{b_{j+1},\ldots, b_n\}\cap B_i|\}\geq \mu|B_i|\}
\]
and 
\[C_i(\mu)=\{j\in [n]:c_j\in C_i, \min\{|\{c_1,\ldots, c_j\}\cap C_i|,|\{c_{j+1},\ldots, c_n\}\cap C_i|\}\geq \mu|C_i|\}.\]
Observe that for any $\mu>0$ and $i\in [t]$, set
$$
(B_i\setminus B_i(\mu))\subseteq (B_i\cap \{b_1,\ldots, b_{\min (B_i(\mu))-1}\})\cup (B_i\cap \{b_{\max (B_i(\mu))+1}, \ldots, n\}),
$$
and consequently, $|B_i(\mu)|\geq (1-2\mu)|B_i|$. A similar argument shows  $|C_i(\mu)|\geq (1-2\mu)|C_i|$.

Suppose first there exist $i\neq j$ such that $i\in \calB$, $j\in \calC$, and $B_i(\e_1^{1/9})\cap C_j(\e_1^{1/8})\neq \emptyset$.  Say $w\in B_i(\e_1^{1/9})\cap C_j(\e_1^{1/9})$. Set  $C_{j0}=C_j\cap \{c_{w+1},\ldots, c_n\}$, $C_{j1}=C_j\cap \{c_1,\ldots, c_w\}$, $B_{i0}=B_i\cap \{b_{w+1},\ldots, b_n\}$, and $B_{i1}=B_i\cap \{b_1,\ldots, b_w\}$.  Fix any $k\in \calA\cap ([t]\setminus \{i,j\})$.  Then $K_3[A_k,C_{j1}, B_{i1}]\subseteq E$ while $K_3[A_k,C_{j0}, B_{i0}]\cap E=\emptyset$, so we are done.

This leaves us with the case where, for all $i\in \calB$ and $j\in \calC$, if $i$ and $j$ are distinct, then $B_i(\e^{1/9})\cap C_j(\e^{1/9})= \emptyset$. In this case we have that
\begin{align*}
n&\geq \Big|\bigcup_{i\in \calB}B_i(\e_1^{1/9})\cup \bigcup_{j\in \calC}C_j(\e_1^{1/9})\Big|\\
&\geq \sum_{i\in \calB}|B_i(\e_1^{1/9})|+\sum_{j\in \calC}|C_j(\e_1^{1/9})|-|\sum_{s\in \calB\cap \calC} B_s(\e_1^{1/9})\cap C_s(\e_1^{1/9})|,\\
\end{align*}
which is at least
\begin{align*}
(1-2\e_1^{1/9})&\Big(\sum_{i\in \calB}|B_i|+\sum_{j\in \calC}|C_j|\Big)-|\sum_{s\in \calB\cap \calC} B_s(\e^{1/9})\cap C_s(\e^{1/9})|\\
&=(1-2\e_1^{1/9})(|B'|+|C'|)-|\sum_{s\in \calB\cap \calC} B_s(\e^{1/9})\cap C_s(\e_1^{1/9})|\\
&\geq (1-2\e_1^{1/9})2n(1-\e_1^{1/9})-|\sum_{s\in \calB\cap \calC} B_s(\e_1^{1/9})\cap C_s(\e_1^{1/9})|.
\end{align*}
This implies $|\sum_{s\in \calB\cap \calC} B_s(\e_1^{1/9})\cap C_s(\e_1^{1/9})|\geq (1-2\e_1^{1/9})2n(1-\e_1^{1/9})- n\geq n(1-\e_1^{1/18})$.   Since $n$ is sufficiently large, this implies there exist $i\neq j$, with both in $i,j\in \calB\cap \calC$, and such that $|B_{i}(\e_1^{1/9})\cap C_{i}(\e_1^{1/9})|\neq \emptyset$ and $|B_{j}(\e_1^{1/9})\cap C_{j}(\e_1^{1/9})|\neq \emptyset$.

Fix any $u\in B_{i}(\e^{1/9})\cap C_{i}(\e_1^{1/9})$ and $v\in B_{j}(\e^{1/9})\cap C_{j}(\e_1^{1/9})$.  Since $i\neq j$, we know $u\neq v$.  Without loss of generality, let us assume $u<v$ (the other case is similar).  Define $B_{i1}=B_{i}\cap \{b_1,\ldots, b_u\}$, $C_{i1}=C_{i}\cap \{c_1,\ldots, c_u\}$, $B_{j0}=B_{j}\cap \{b_{v+1},\ldots, b_n\}$, and $C_{j0}=C_{j}\cap \{c_{v+1},\ldots, c_n\}$.  Choose any $k\in \calA\setminus \{i,j\}$.  Note that $K_3[A_k,B_{i1}, C_{j1}]\subseteq E$, while $K_3[A_k, B_{i0}, C_{j0}]\cap E=\emptyset$, so we are done.
\end{proof}

As a corollary, we show that if $\trip(\calH)$ contains arbitrarily large $\Hbar(k)$, then $\calH$ requires non-zero $\disc_{2,3}$-error.  

\begin{corollary}\label{cor:dischbark}
Suppose $\calH$ is a hereditary property of $3$-graphs such that $\Hbar(k)\in \trip(\calH)$ for all $k\geq 1$.  Then $\calH$ requires non-zero $\disc_{2,3}$-error.  
\end{corollary}
\begin{proof}
Suppose towards a contradiction there is some $\calH$ with $\Hbar(k)\in \trip(\calH)$ for all $k\geq 1$, such that $\calH$ admits zero $\disc_{2,3}$-error.  Let $\e=2^{-28}$ (anything small enough for Theorem \ref{thm:hbarksimple} will do).  Let $\e_1$ be as in Proposition \ref{prop:trivdensity} for $k$ and $\e^2$.  Given $\ell\geq 1$, define $\e_2'(\ell)$ as in Corollary \ref{cor:counting} for $t=3$, density $1/\ell$ and $\e^2$.  Then define $\e_2''(\ell)$ be as in Lemma \ref{lem:sl} for density $1/\ell$, $\gamma=\e^{1/9}$, and $\e_2(\ell)$, and let $\e'_1$, $\e_2'''(\ell)$ be as in Lemma \ref{lem:subpairs} for $\e_1$, $\gamma=\e^{1/9}$, $3$, density $1/\ell$, $\e_2'(\ell)$.  Set $\e_2(\ell)=\min\{\e_2'(\ell),\e_2''(\ell),\e_2'''(\ell)\}$.  

By assumption there are $T, L, N\geq 1$ such that if $G\in \calH$ has at least $N$ vertices, then there exists $1\leq \ell\leq L$, $1\leq t\leq T$, and a $(t,\ell,\e_1,\e_2(\ell))$-decomposition of $V(G)$ which is  $(\e_1,\e_2(\ell))$-regular with respect to $G$, with zero $\disc_{2,3}$-error.  

Fix $n\geq N$.  By Lemma \ref{lem:clean}, there is $H=(V,E)\in \calH$ with $|V|=3n$ and such that $\trip(G)\cong \Hbar(n)$.  This means $V=A\sqcup B\sqcup C$ where $A=\{a_i: i\in [n]\}$, $B=\{b_i:i\in [n]\}$ and $C=\{c_i:i\in [n]\}$, and $a_ib_jc_k\in E$ if and only if $i\leq j$.  

By assumption, there are $1\leq \ell\leq L$, $1\leq t\leq T$, and a $(t,\ell,\e'_1,\e_2(\ell))$-decomposition $\calP$ of $V$ which is $(\e'_1,\e_2(\ell))$-regular with respect to $H$ with zero $\disc_{2,3}$-error.  By Theorem \ref{thm:hbarksimple}, there are $ijk\in {[k]\choose 3}$ and sets $A_{i}\subseteq V_i\cap A$, $B_{j1},B_{j0}\subseteq V_j\cap B$, and $C_{k1},C_{k0}\subseteq V_k\cap C$, each of size at least $\e^{1/9}(n/t)$, such that $K_3[A_{i},B_{j1},C_{k1}]\subseteq E$ and  $K_3[A_{i},B_{j0},C_{k0}]\cap E=\emptyset$.    Let $A'=A_i$, $B'=B_{j1}\cup B_{j0}$, $C'=C_{k1}\cup C_{k0}$, and $V'=A'\cup B'\cup C'$.  Then define $H''=(V', E\cap K_3[A',B',C'])$.  By Proposition \ref{lem:subpairs}, for any $1\leq \alpha,\beta,\gamma\leq \ell$, $(H'', G_{ijk}^{\alpha,\beta,\gamma}[V'])$ satisfies $\disc_{2,3}(\e_1,\e_2'(\ell))$.   But we observe that for any $1\leq \alpha,\beta,\gamma\leq \ell$, 
\[\min\{|E\cap K_3(G_{ijk}^{\alpha,\beta,\gamma}[V'])|, |K_3(G_{ijk}^{\alpha,\beta,\gamma}[V'])\setminus E|\}\]
is at least the minimum of 
\[|E\cap K_3(G_{ijk}^{\alpha,\beta,\gamma})\cap K_3[A_i\cup B_{j0}\cup C_{k0}]|\]
and
\[|(K_3(G_{ijk}^{\alpha,\beta,\gamma})\cap K_3[A_i\cup B_{j1}\cup C_{k1}])\setminus E|,\]
and thus at least 
\[\min\{|K_3(G_{ijk}^{\alpha,\beta,\gamma})\cap K_3[A_i\cup B_{j0}\cup C_{k0}]|,  |K_3(G_{ijk}^{\alpha,\beta,\gamma})\cap K_3[A_i\cup B_{j1}\cup C_{k1}]|\}.\]
By Lemma \ref{lem:sl}, for each $u\in \{0,1\}$, each of the graphs $P_{ij}^{\alpha}[A_i,B_{ju}]$, $P_{ik}^{\beta}[A_i,C_{ku}]$ and $P_{jk}^{\gamma}[B_{ju},C_{ku}]$ has $\disc_2(\e_2'(\ell); 1/\ell)$.  Consequently, by Corollary \ref{cor:counting}, for each $u\in \{0,1\}$,
\[|K_3(G_{ijk}^{\alpha,\beta,\gamma})\cap K_2[A_i\cup B_{ju}\cup C_{ku}]|\]
equals
$$
\frac{1}{\ell^3}(1\pm \e^2)|A_i||B_{ju}||C_{ku}|\geq \frac{1}{\ell^3}(1\pm \e^2)\e^{1/3}|V_i||V_j||V_k|.
$$
Combining this with the above, we have that
\begin{align*}
\min\{|E\cap K_3(G_{ijk}^{\alpha,\beta,\gamma}[V'])|, |K_3(G_{ijk}^{\alpha,\beta,\gamma}[V'])\setminus E|\}&\geq \e^{1/3}\frac{1}{\ell^3}(1\pm \e^2)\e^{1/3}|A'||B'||C'|\\
&\geq \e^{2/2}|K_2^{(2)}(G_{ijk}^{\alpha,\beta,\gamma}[V'])|,
\end{align*}
where the last inequality is by Corollary \ref{cor:counting}.  This shows that  $(H'', G_{ijk}^{\alpha,\beta,\gamma}[V'])$ satisfies $\disc_{2,3}(\e_1,\e_2(\ell))$, and has density in $(\e^2,1-\e^2)$.  By our choice of $\e_1$ and $\e_2$, this implies $H''$ has $\VC_2$-dimension at least $k$.  By definition $H''$ is an induced sub-$3$-graph of $\Hbar(n)$, but it is not difficult to show that $\VC_2(\Hbar(n))<3$. Thus we have arrived at a contradiction.
\end{proof}

We can now give a very short proof of Theorem \ref{thm:vdisceq}, which shows that a hereditary $3$-graph property admits zero  $\disc_{2,3}$-error if and only if it admits zero $\vdisc_3$-error.

\vspace{2mm}

\noindent{\bf Proof of Theorem \ref{thm:vdisceq}.}
\label{proof:vdisceq} Suppose first that $\calH$ admits zero $\disc_{2,3}$-error. By Corollary \ref{cor:dischbark}, there is some $k$ such that $\Hbar(k)\notin \trip(\calH)$.  Clearly this implies $\Ubar(k)\notin \trip(\calH)$, so by Theorem \ref{thm:simclasswnip}, $\calH$ is close to a SNIP property and $\vdisc_3$-homogeneous.  By Proposition \ref{prop:wnipreduction}, $\calH$ admits zero $\vdisc_3$-error.  

Suppose on the other hand that $\calH$ admits zero $\vdisc_3$-error.  By Theorem \ref{thm:vdischom}, $\calH$ admits zero $\disc_{2,3}$-error.
\qed

\vspace{2mm}

Combining Theorem \ref{thm:vdisceq} with Corollary \ref{cor:dischbark} yields Theorem \ref{thm:vdischbark}, and combining Theorem \ref{thm:vdisceq} with Example \ref{ex:stable} yields Proposition \ref{prop:fullweakex}.

\section{A generalized strong stable regularity lemma}\label{subsec:strongstable}

The goal of the final section in this chapter will be to prove Proposition \ref{prop:wsbinary}, which shows that a slicewise stable hereditary $3$-graph property admits binary $\vdisc_3$-error and binary $\disc_{2,3}$-error.  The proof will be partially based on the proof of the following strong regularity for stable graphs.

\begin{theorem}\label{thm:functionstablereg}
For all $k\geq 1$, $\e>0$ and $f:\mathbb{N}\rightarrow (0,1]$ non-increasing, there is $M=M(\e,f,d)$ such that the following holds.  

Suppose $G=(V,E)$ is a $k$-stable graph.  Then there is $m\leq M$ and a partition $W=W_0\cup \ldots\cup W_m$ such that $|W_0|\leq \e |W|$, and such that for each $1\leq i,j\leq m$, $(W_i,W_j)$ satisfies $\disc_2(f(m))$ and $|W_i|=|W_j|$.
\end{theorem}

Note that the degree of quasirandomness in the pairs avoiding $W_0$ is allowed to depend on the number of parts in the partition.  Theorem \ref{thm:functionstablereg} is similar to the general strong regularity lemma (see \cite{Rodl.2007,  Lovasz.2007,  Alon.2000, Alon.2008}), however the conclusions here are strictly stronger than can be obtained in general.  Indeed Theorem \ref{thm:functionstablereg} implies the existence of a regular partition with no irregular pairs, so $H(n)$ could not have such a partition. 

We will prove Theorem \ref{thm:functionstablereg} in this section.  This serves as a warm-up for the proof of Proposition \ref{prop:wsbinary} in the next section, which is very similar.  This will also lead naturally to the proof of a result needed in Chapter \ref{sec:fop}, namely Theorem \ref{thm:goodstrong}.

The proof of Theorem \ref{thm:functionstablereg} uses ideas from  \cite{Malliaris.2014}, where the regularity lemma for stable graphs was first proved.  We begin with some terminology from that paper.  Given a graph $G=(V,E)$, a subset $X\subseteq V$, and a vertex $v\in V$, we say that $X$ is \emph{$\e$-good with respect to $v$} if there is some $u=u(v)\in \{0,1\}$ such that 
$$
\Big|\frac{|N(v)\cap X|}{|X|} - u\Big|\leq \e.
$$
The set $X$ is called \emph{$\e$-good} if it is $\e$-good with respect to every $v\in V$.  In other words, $X$ is $\e$-good if every vertex is connected to all but at most $\e|X|$ many elements of $X$, or is non-connected to all but at most $\e|X|$ many elements of $X$.  In \cite{Malliaris.2014}, Malliaris and Shelah showed that in a stable graph, any large set of vertices can by partitioned into $\e$-good subsets, all of which have about the same size.

\begin{theorem}[Theorem 5.18 in \cite{Malliaris.2014}]\label{thm:MSgood}
For all $k\geq 1$ and $\e>0$, there exist $K=K(\e,k)$ and $M=M(k,\e)$ such that the following holds.  Suppose that $G=(V,E)$ is a $k$-stable graph and $W\subseteq V$ has size at least $M$.  Then there is an equipartition of $W$ into at most $K$ subsets, all of which are $\e$-good.
\end{theorem}

In fact, in \cite{Malliaris.2014} a stronger version of Theorem \ref{thm:MSgood} is proved, in which the notion of goodness is upgraded to ``excellence".  This result is iterated to create a partition into excellent sets, which is then turned into an equipartition of excellent sets using a random sampling argument. Finally, a symmetry lemma is invoked, which says that any pair of excellent sets has density near $0$ or $1$, and thus regularity for all pairs is achieved. This is the stable regularity lemma of Malliaris and Shelah, which inspired the set of questions investigated in this paper.

Here, on the other hand, we will iterate Theorem \ref{thm:MSgood} to create a partition into good sets. We will then argue that all pairs of this partition are regular using a version of the symmetry lemma for good sets (Lemma \ref{lem:twosticks} below). We will make no effort to obtain an equipartition, although this could be done using \cite[Claim 5.13(1)]{Malliaris.2014}.

The symmetry lemma for excellent sets used in \cite{Malliaris.2014} follows essentially from the definition of excellence, while the symmetry lemma we use for good sets requires a proof.  However, the combinatorics behind Lemma \ref{lem:twosticks} does already appear in \cite{Malliaris.2014}, in particular inside the proof of Claim 5.4 there.

\begin{lemma}[Symmetry lemma for good sets]\label{lem:twosticks}
For all $0<\e<1/4$ there is $n$ such that the following holds.  Suppose $G=(U\cup W, E)$ is a bipartite graph satisfying $|U|,|W|\geq n$, and assume $U'\subseteq U$, $W'\subseteq W$ satisfy $|U'|\geq (1-\e)|U|$ and $|W'|\geq (1-\e)|W|$.  Suppose that for all $u\in U'$, 
$$
\max\{|N(u)\cap W|, |W\setminus N(u)|\}\geq (1-\e)|W|,
$$ 
and for all $w\in W'$, 
$$
\max\{|N(w)\cap U|, |U\setminus N(w)|\}\geq (1-\e)|U|.
$$ 
Then $|E|/|U||W|\in [0,2\e^{1/2})\cup (1-2\e^{1/2},1]$. 
\end{lemma}
\begin{proof}
For each $\alpha\in \{0,1\}$, define
\[U_\alpha=\{u\in U: |N^\alpha(u)\cap W|\geq (1-\mu)|W|\}\]
and
\[W_\alpha=\{w\in W: |N^\alpha(w)\cap U|\geq (1-\mu)|U|\}.\]
It suffices to show there is $\alpha\in \{0,1\}$ such that $|U_\alpha|\geq (1-\sqrt{\e})|U|$, since then 
$$
|E^{\alpha}|\geq (1-\e)|W||U_{\alpha}|\geq (1-\e)(1-\sqrt{\e})|W||U|\geq (1-2\sqrt{\e})|U||W|.
$$

Suppose towards a contradiction that $\max\{|U_0|, |U_1|\}<(1- \sqrt{\e})|U|$.  Note that for each $\alpha\in \{0,1\}$, $|E^{\alpha}\cap K_2[U_{\alpha},W]|\geq (1-\e)|W||U_{\alpha}|$,  so by standard arguments, there is $W_\alpha'\subseteq W$ of size at least $(1-\sqrt{\e})|W|$ such that for all $w\in W_{\alpha}$, $|N^{\alpha}(w)\cap U_{\alpha}|\geq (1-\sqrt{\e})|U_{\alpha}|$. Then $W''=W_0\cap W_1$ has size at least $(1-2\sqrt{\e})|W|$, and for all $w\in W''$,
\begin{align*}
\min\{|N^1(w)\cap U|, |N^0(w)\cap U|\}&\geq \min\{|N^1(w)\cap U_1|, |N^0(w)\cap U_0|\}\\
&\geq \min\{(1-\sqrt{\e})|U_0|, (1-\sqrt{\e})|U_1|\}\\
&>\e|U|.
\end{align*}
This implies that $W''\cap W'=\emptyset$.  However, this is not possible since $|W'|\geq (1-\e)|W|$, $|W''|\geq (1-2\sqrt{\e})|W|$, and $\e<1/4$.  
\end{proof}

The phenomenon described in Lemma \ref{lem:twosticks} shows up in many places in model theory (see e.g. \cite[Lemma 2.8]{Pillay.1996}) and graph theory (see e.g. \cite[Lemma 27]{Balogh.2000}).  An immediate corollary of Lemma \ref{lem:twosticks} is that between any two good sets in a graph, the density of edges is always close to $0$ or close to $1$.

\begin{corollary}[A pair of good sets is homogeneous]\label{cor:goodpair}
For all $0<\e<1/4$ there is $n\geq 1$ such that the following holds.  Suppose $G=(V,E)$ is a graph and $U,W\subseteq V$ each have size at least $n$. If $U$ and $W$ are both $\e$-good, then $|E\cap K_2[U,W]|/|K_2[U,W]|\in [0,2\e^{1/2})\cup (1-2\e^{1/2},1]$. 
\end{corollary}

One can now quickly sketch a proof of a stable regularity lemma as follows.  Suppose $G=(V,E)$ is a large $k$-stable graph.  By Theorem \ref{thm:MSgood}, there is an equipartition $\calP$ of $V$ into $\e$-good sets. By Corollary \ref{cor:goodpair}, every pair of parts from $\calP$ has density in $[0,2\e^{1/2})\cup (1-2\e^{1/2},1]$.  By Fact \ref{fact:homimpliesrandombinary}, every pair of parts from $\calP$ satisfies $\disc_2(2\e^{1/2})$.  

We will use a similar strategy to prove Theorem \ref{thm:functionstablereg}.  The first step is to prove a strong version of Theorem \ref{thm:MSgood}.  We will work with an equivalent characterization of stability based on the Shelah $2$-rank (see \cite{Shelah.1990o5n}).  To state this equivalence we need some more notation.

Given an integer $d\geq 1$, let $2^d=\{0,1\}^d$.  In other words, $2^d$ consists of sequences of $0$s and $1$s of length $d$.  By convention, we let $2^0=\{<>\}$, where $<>$ denotes the so-called empty string.  Then for $d\geq 1$, we set $2^{<d}:=\bigcup_{i=0}^{d-1}2^i$.  Given $0\leq i\leq j$, $\sigma\in 2^i$ and $\tau\in 2^j$, we write $\sigma \triangleleft \tau$ to denote that $\sigma$ is a proper initial segment of $\tau$, and $\sigma\trianglelefteq$ to denote that either $\sigma\triangleleft\tau$ or $\sigma=\tau$.  By convention $<>\trianglelefteq \sigma$, for any $\sigma\in 2^d$ and $d\geq 0$.  We think of $2^{<d}$ as a binary branching tree with root $<>$.  Given $\sigma\in 2^i$ and $\sigma'\in 2^j$, $\sigma \wedge \sigma'$ denotes the element of $2^{i+j}$ obtained by adjoining $\sigma'$ to the end of $\sigma$ (e.g. $(01)\wedge (001)=(01001)$). \label{not:wedge}

\begin{definition}\label{def:treerank}
Suppose $G=(V,E)$ is a graph.  Given $d\geq 1$, a \emph{$d$-tree} in $G$ is a tuple $(a_{\sigma})_{\sigma\in 2^{<d}}(b_{\rho})_{\rho\in 2^d}$ with the property that  for all $\sigma\triangleleft \eta$, $b_{\sigma}a_{\eta}\in E$ if and only if $\sigma \wedge 1\trianglelefteq \eta$.

The \emph{rank} of $G$ is 
$$
\rk(G):=\max\Big(\{0\}\cup \{d\geq 1:\text{ there is a $d$-tree in $G$}\}\Big).
$$
\end{definition}

Given a $d$-tree as in Definition \ref{def:treerank}, we refer to the elements $\{a_{\eta}:\eta\in 2^d\}$ as the \emph{leaves} and the elements $\{b_{\sigma}:\sigma\in 2^{<d}\}$ as the \emph{nodes}.    There is a fundamental relationship between tree-rank and the order property.

\begin{theorem}[Shelah \cite{Shelah.1990o5n}]\label{thm:treerank}
$\;$
\begin{enumerate}[label=\normalfont(\arabic*)]
\item For all $k\geq 1$ there is a $d=d(k)\geq 1$ such that the following holds.  Suppose $G=(V,E)$ is a graph with rank $d$, witnessed by $(a_{\eta})_{\eta\in 2^d}(b_{\sigma})_{\sigma\in 2^{<d}}$. Then there are $a_1,\ldots, a_k\subseteq \{a_{\eta}:\eta\in 2^d\}$ and $b_1,\ldots, b_k\in \{b_{\sigma}:\sigma\in 2^{<d}\}$ such that  $a_ib_j\in E$ if and only if $i\leq j$. 

\item For all $d\geq 1$ there is a $k=k(d)\geq 1$ such that the following holds.  Suppose $G=(V,E)$ is a graph with the $k$-order property, witnessed by $(a_i)_{i\in [k]}(b_i)_{i\in [k]}$.  Then there are $\{a_{\eta}:\eta\in 2^d\}\subseteq \{a_i:i\in [k]\}$ and $\{b_{\sigma}:\sigma\in 2^{<d}\}\subseteq \{b_i:i\in [k]\}$ such that   for all $\sigma\triangleleft \eta$, $a_{\eta}b_{\sigma}\in E$ if and only if $\sigma \wedge 1\trianglelefteq \eta$.
\end{enumerate}
\end{theorem}

It was shown by Hodges \cite{Hodges.1981} that in (1), one can take $d(k)=2^{k+2}-2$, and in (2), one can take $k(d)=2^{d+1}$.

Theorem \ref{thm:treerank} tells us that a graph with sufficiently large rank will have the $k$-order property, and on the other hand, a graph with a sufficiently long order property will have rank at least $d$.  We now define a notion of rank for subsets of vertices.  

\begin{definition}
Suppose $G=(V,E)$. The \emph{rank} of a set $A\subseteq V$ is 
$$
\rk(A):=\max(\{0\}\cup \{d\geq 1: \text{ there exists a $d$-tree in $G$ with leaves in $A$}\}).
$$ 
\end{definition}

Note that the rank of $A$ is $0$ if there do \emph{not} exist  $a_0,a_1\in A$ and $b_{<>}\in V$ such that $b_{<>}a_1\in E$ and $b_{<>}a_1\notin E$.  Thus, $\rk(A)=0$ if and only if $A$ is $0$-good in $G$.  

Observe that if $A$ has rank $d>0$, then for any $b\in V$, one of $N(b)\cap A$ or $(\neg N(b))\cap A$ must have rank less than $d$.   This observation will be crucial in our proofs.  We will also use the easy fact that if $A'\subseteq A$, then $\rk(A')\leq \rk(A)$. 

In our first lemma, we show that in any bipartite graph $G=(U\cup W, E)$ with tree rank $d$, we can partition $W$ into sets of rank less than $d$, and sets which are somewhat good.

\begin{lemma}\label{lem:goodsets1}
Suppose $d\geq 1$, $\e>0$ and  $f:\mathbb{N}\rightarrow (0,1]$ is a non-increasing function satisfying $\lim_{n\rightarrow \infty}f(n)=0$.  

Suppose $G=(U\cup W,E)$ is a bipartite graph and $\rk(W)\leq d$.  Then there is an integer $t\leq |W|$, a sequence $(s_1,\ldots, s_t)$, and a partition $W=\bigcup_{i\in [t]}\bigcup_{j\in [s_i]}W_{i,j}$ such that the following hold, where for each $1\leq i\leq t$,  $W_i:=\bigcup_{j\in [s_i]}W_{i,j}$.
\begin{enumerate}[label=\normalfont(\arabic*)]
\item For each $1\leq i\leq t$ and $1\leq j\leq s_i$, $\rk(W_{i,j})<d$, and
$$
|W_{i,j}|=f(i)\Big|W\setminus (\bigcup_{1\leq i'<i}W_{i'}\cup \bigcup_{j'=1}^{j-1}W_{i,j'})\Big|.
$$ 
\item For each $1\leq i\leq t$, $W\setminus (\bigcup_{i'=1}^i W_{i'})$ is $f(i)$-good.
\end{enumerate}
\end{lemma}
\begin{proof}
Suppose $G=(U\cup W, E)$ is a bipartite graph in which $\rk(W)\leq d$.  We define an integer $t$, and sequences $(s_0,\ldots, s_t)$, $Z_0,\ldots, Z_t$, and $W_1,\ldots, W_t$ as follows.  

\underline{Step $0$:} Set $s_0=0$, $Z_0=W$ and $W_0=\emptyset$.

\underline{Step $\alpha$:} Suppose $\alpha>0$ and we have chosen $Z_0,\ldots, Z_{\alpha-1},W_0,\ldots, W_{\alpha-1}$, and for all $0\leq \alpha'<\alpha$,  $s_{\alpha'}$ and $W_{\alpha',1},\ldots, W_{\alpha',s_{\alpha'}}$ so that the following hold.
\begin{itemize}
\item For all $0\leq \alpha'<\alpha$, $W_{\alpha'}=\bigcup_{u=1}^{s_{\alpha'}}W_{\alpha',u}$, and  $Z_{\alpha'}=W\setminus (\bigcup_{\beta<\alpha'}W_{\beta})$, 
\item For each $0\leq \alpha'<\alpha$ and $1\leq u<s_{\alpha'}$, $|W_{\alpha',u}|=f(\alpha')|W\setminus Z_{\alpha'-1}\cup (\bigcup_{v=1}^{u-1}W_{\alpha',v})|$, 
\item For each $0< \alpha'< \alpha$ and $Z_{\alpha'}$ if $f(\alpha')$-good.
\end{itemize}

\underline{Substep $(\alpha, 0)$:} Let $Z_{\alpha,0}=Z_{\alpha-1}=W\setminus (\bigcup_{\beta<\alpha}W_\beta)$ and $W_{\alpha,0}=\emptyset$.  If $Z_{\alpha,0}=\emptyset$, or $Z_{\alpha,0}$ is $0$-good, set $t=\alpha-1$ and end the construction.

\underline{Substep $(\alpha,v+1)$:} Suppose that $v\geq 0$ and we have defined $W_{\alpha,0},\ldots, W_{\alpha,v}$, $Z_{\alpha,0},\ldots, Z_{\alpha,v}$ so that for each $0\leq v'\leq v$, $Z_{\alpha,v'}=Z_{\alpha-1}\setminus (\bigcup_{v=0}^{v'}W_{\alpha,v})$, and for each $1\leq v'\leq v$, $\rk(W_{\alpha,v'})<d$, and $|W_{\alpha,v'}|= f(\alpha)|Z_{\alpha,v'-1}|$.  Define 
$$
U_{\alpha}^{v+1}=\{u\in U: \min\{|N(u)\cap Z_{\alpha,v}|,|\neg N(u)\cap Z^{\alpha,v}|\}\geq f(\alpha)|Z_{\alpha, v}|\}.
$$
If $U_{\alpha}^{v+1}=\emptyset$, set $s_{\alpha}=v$, $W_{\alpha}=\bigcup_{1\leq v'\leq v}W_{\alpha,v'}$,  $Z_{\alpha}=Z_{\alpha,v}$, and go to step $\alpha+1$.

If $U_{\alpha}^{v+1}\neq \emptyset$, choose any $u\in U^{\alpha}_{v+1}$, and let $X\in \{|N(u)\cap Z^{\alpha}_{i,v}|,|\neg N(u)\cap Z^{\alpha}_{i,v}|\}$ be a set rank less than $d$ (such an $X$ exists because $\rk(W)\leq d$).  Define $W_{\alpha,v+1}$ to be any subset of $X$ of size $f(\alpha)|Z^{\alpha}_{v}|$.  Since $\rk(X)<d$ and $W_{\alpha,v+1}\subseteq X$, we also have $\rk(W_{\alpha,v+1})<d$.  Set $Z_{\alpha, v+1}=Z_{\alpha,v}\setminus W_{\alpha,v+1}$ and go to step $(\alpha,v+2)$. 

Since $W$ is finite and $\lim_{n\rightarrow \infty}f(n)=0$, there will be some $\alpha=t$ where this ends. 
\end{proof}

We now prove a much stronger version of Lemma \ref{lem:goodsets1}.  This is the main technical lemma behind Theorem \ref{thm:functionstablereg}, and plays a crucial role in the proof of Theorem \ref{thm:FOPfinite}.  For the reader curious about equitability considerations, we note here that these will be addressed later on, in Theorem \ref{thm:goodstrongequitable}.

\begin{theorem}\label{thm:goodstrong}
Suppose $d\geq 0$.  For all $\e>0$ and non-increasing functions $f:\mathbb{N}\rightarrow (0,1]$ with $\lim_{n\rightarrow \infty}f(n)=0$, there is $M\geq 1$ such that the following hold.  

Suppose $m\geq 1$, $G=(U\cup W, E)$ is a bipartite graph, and $W_1\cup \ldots\cup W_m$ is a partition of $W$ so that $\rk(W_i)\leq d$ for all $i\in [m]$. 

Then there are $m'\leq M$, and a set $\Omega \subseteq [m]$, such that $|\bigcup_{u\in \Omega}W_u|\leq \e|W|$, and for each $u\in [m]\setminus \Omega$, there are $s_u\leq m'$ and a partition
$$
W_u=W_{u,0}\cup W_{u,1}\cup \ldots \cup W_{u,s_u},
$$
such that $|W_{u,0}|\leq \e |W_u|$, and such that $W_{u,i}$ is $f(m')$-good for each $1\leq i\leq s_u$.
\end{theorem}

Before we give the proof, we point out some important features of this theorem.  First, observe that the bound $M$ does not depend on $m$, the number of parts in the partition we start with.  Further, a graph need not be stable to satisfy the hypotheses of Theorem \ref{thm:goodstrong}.  Indeed, when $m$ is large, $G$ could have large rank even though each $W_i$ has small rank.  We will use Theorem \ref{thm:goodstrong} in such a scenario in Chapter \ref{sec:fop}.  

\vspace{2mm}
\begin{proofof}{Theorem \ref{thm:goodstrong}}
To ease notation, let $I$ be the set of non-increasing functions $f:\mathbb{N}\rightarrow (0,1]$ satisfying $\lim_{n\rightarrow \infty}f(n)=0$.  We proceed by induction on $d\geq 0$.  

\underline{Case $d=0$}:  Suppose $d=0$.  Given $\e>0$, and  $f\in I$, set $M=1$. Suppose $G=(U\cup W,E)$ is a bipartite graph and assume $W_1,\ldots, W_m$ partition $W$ such that $\rk(W_i)= 0$ for each $i\in [m]$. This means each $W_i$ is $0$-good in $G$.  Setting $m'=1$ and $W_{i,1}'=W_i$ for each $i\in [m]$, we are done. 

\underline{Case $d>0$}:  Suppose now that $d>0$ and assume by induction we have defined $M(d', f, \e)$ for all $0\leq d'<d$, $\e>0$, and $f\in I$.

Fix $\e>0$ and  $f\in I$.  Without loss of generality, assume $\e<1/4$.  We inductively define a new function $f':\mathbb{N}\rightarrow (0,1]$, and a sequence $(M_i,g_i)_{i\in \mathbb{N}}$, where each $M_i\in \mathbb{N}$, and each $g_i\in I$. 

\underline{Step $1$}: Define $g_1$ by setting $g_1(x)=\e^4f(2\e^{-4}x)$.  Then let $M_1=M(d-1,g_1,\e^8)$ and $f'(1)=g_1(M_1)/2$.

\underline{Step $i+1$}: Suppose $i\geq 1$, and assume that for all $1\leq j\leq i$, we have defined $g_j$, $M_j$,, and $f'(j)$.  Set $P_i=\sum_{j=1}^i \e^{-2} f'(j)^{-1}$ and define $g_{i+1}$ by setting $g_{i+1}(x)=\e^{4}f(\e^{-4}(x\cdot (P_i+1)))$.  Then set $M_{i+1}=M(d-1, g_{i+1}, \e^8)$ and $f'(i+1)=f'(i)(i+2)^{-1}g_{i+1}(M_{i+1})$.

This finishes the definition of $(M_i,g_i)_{i\in \mathbb{N}}$ and of $f'$.  It is straightforward to check that by construction, and because $f\in I$, we have that $f'\in I$, and $g_i\in I$ for all $i\geq 1$.   Let $S_1=\lceil \e^{-2}\rceil$ and choose $M(d,\e,f)\gg M(d-1,\e^8,f)\e^{-8}\prod_{i=1}^{S_1}M_iP_{i}$. 

Suppose $G=(U\cup W, E)$ is a bipartite graph and $W=W_1\cup \ldots\cup W_m$ is a partition so that for each $i\in [m]$, $\rk(W_i)\leq d$. Let $\ell$ be the number of $W_i$ with $\rk(W_i)=d$.  If $\ell=0$, then we are done by the inductive hypothesis.

Thus we may assume $\ell>0$.  By relabeling if necessary, we may assume $W_1,\ldots, W_{\ell}$ have  rank $d$ and $W_{\ell+1},\ldots, W_m$ have rank less than $d$.

For each $1\leq u\leq \ell$, apply Lemma \ref{lem:goodsets1} with parameters $d$, $\e^2$, and $f'$ to the set $W_u$ to obtain an integer $t(u)\leq |W_u|$, a sequence $(s_i(u))_{i\in t(u)}$, and a partition $W_u=\bigcup_{i\in [t(u)]}\bigcup_{j\in [s_i(u)]}V_{i,j}(u)$ such that the following hold, where for each $1\leq i\leq t(u)$, 
$$
V_i(u):=\bigcup_{j\in [s_i(u)]}V_{i,j}(u)\text{ and }Z_{i,j}(u)= W_u\setminus (\bigcup_{i'=1}^{i-1}V_{i'}(u)\cup \bigcup_{j'=1}^{j-1}V_{i,j'}(u)).
$$
\begin{enumerate}
\item[(i)] For each each $i\in [t(u)]$ and $j\in [s_i(u)]$, $|V_{i,j}(u)|=f'(i)|Z_{i,j}(u)|$ and $\rk(V_{i,j}(u))<d$, 
\item[(ii)] For each $i\in [t(u)]$, $W_u\setminus (\bigcup_{i'=1}^i V_{i'}(u))$ is $f'(i)$-good.
\end{enumerate}
Let $t=\max\{t(u): 1\leq u\leq \ell\}$.  Given $u\in [\ell]$, if $t(u)<i\leq t$, set $V_i(u)=\emptyset$.  Then for each $1\leq i\leq t$, define $Y_i=\bigcup_{u=1}^{\ell}V_i(u)$, and let $i_1$ be minimal such that $|Y_i|\leq \e^8 |W|$.  Clearly such an $i_1$ exists and is at most $\e^{-8}$.  Then define
$$
\Omega_0=\{u\in [\ell]:|V_{i_1}(u)|> \e^4|W_u|\}.
$$
Note by definition, $\e^4|\bigcup_{u\in \Omega_0}W_u|< |Y_{i_1}|\leq \e^8 |W|$, so $|\bigcup_{u\in \Omega_0}W_u|\leq \e^4 |W|$.  

Now fix $u\in [\ell]\setminus \Omega_0$.  Set $J(u)=\{j\in [i_1-1]:  |V_j(u)|\leq \e^4f'(i_1-1)|W_u|\}$.  Observe that by definition, 
$$
|\bigcup_{j\in J(u)}V_j(u)|\leq |J(u)|\e^4f'(i_1-1)|W_u|\leq (i_1-1)\e^4f'(i_1-1)|W_u|<\e^4|W_u|,
$$
 where the last inequality is because, by definition, $f'(i_1)<(i_1)^{-1}$.  Set $I(u)=[i_1-1]\setminus J(u)$, and for each $i\in I(u)$, define
$$
T_i(u)=\{j\in [s_i(u)]: |V_{i,j}(u)|\geq \e^2f'(i)|V_i(u)|\},
$$
 We then let $t_i(u)=|T_i(u)|$, and observe that $t_i(u)$ is at most $\e^{-2}f'(i)^{-1}$.  Now define 
 $$
V_{i,0}(u)=\bigcup_{j\in [s_i(u)]\setminus T_i(u)}V_{i,j}(u).
$$
We claim $|V_{i,0}(u)|\leq 2\e^2 |V_i(u)|$.  If $V_{i,0}(u)=\emptyset$, this is obvious.  If $V_{i,0}(u)\neq \emptyset$, consider $j_0=\min([s_i(u)]\setminus T_i(u))$. By (i), and since $j_0\notin [s_i(u)]\setminus T_i(u)$, 
$$
|V_{i,j_0}(u)|=f'(i)|Z_{i,j_0}(u)|\leq \e^2f'(i)|V_i(u)|.
$$
Consequently $|Z_{i,j_0}(u)|\leq \e^2|V_i(u)|$.  By construction, and minimality of $j_0$, we have $V_{i,0}(u)\subseteq V_{i,j_0}(u)\cup Z_{i,j_0}(u)$.  This shows $|V_{i,0}(u)|\leq 2\e^2 |V_i(u)|$.   We now set
$$
 V_0(u)=V_{i_1}(u)\cup (\bigcup_{i\in J(u)}V_j(u))\cup \bigcup_{i\in I(u)}V_{i,0}(u).
 $$
By above,
$$
|V_0(u)|\leq \e^4|W_u|+\e^4|W_u|+\sum_{i\in I(u)}|V_{i,0}(u)|\leq 2\e^4|W_u|+\sum_{i\in I(u)}2\e^2|V_i(u)|\leq  4\e^2|W_u|.
$$
Recall that by (ii), $Z_{i_1}(u)=W_u\setminus (\bigcup_{v=1}^{i_1}V_i(u))$ is $f'(i_1)$-good.  If $Z_{i_1}(u)<\e^2 |W_u|$, set $V_0'=V_0\cup Z_{i_1}$ and $V_{s+1}(u)=\emptyset$.  If $Z_{i_1}(u)<\e^2 |W_u|$, set $V_0'=V_0$ and $V_{s+1}(u)=Z_{i_1}$.  In either case, $|V_0'|\leq 5\e^2|W_u|$, and we now have a partition
$$
W_u=V'_0(u)\cup V_{s+1}(u)\cup \bigcup_{i\in I_i(u)}\bigcup_{j\in T_i(u)}V_{i,j}(u).
$$
Set $S(u):=\sum_{i\in I(u)} t_i(u)$.  Since each  $t_{i}(u)\leq (\e^2 f'(i))^{-1}$, we have 
$$
S(u)\leq \sum_{i=1}^{i_1-1}\e^{-2}f'(i)^{-1}=P_{i_1-1}.
$$
Let $R_1(u),\ldots, R_{S(u)}(u)$ be an enumeration of $\{V_{i,j}(u): i\in I(u), j\in T_i(u)\}$.

Consider the partition $\calP=\{R_i(u): u\in [\ell]\setminus \Omega_0, i\in [S(u)]\}\cup \{W_{\ell+1},\ldots, W_m\}$ of the set $W':=\bigcup_{X\in \calP}X$.  Note that each element in $\calP$ has rank at most $d-1$ in $G'=G[U\cup W']$.   By the induction hypothesis, applied to $G'=G[U\cup W']$, there are $\Omega_1\subseteq \calP$ and $m'\leq M_{i_1}=M(t-1,g_{i_1},\e^8)$ such that the following hold.
\begin{itemize}
\item $|\bigcup_{X\in \Omega_1}X|\leq \e^8 |W'|$,
\item For each $\ell+1\leq u\leq m$ with $W_u\notin \Omega_1$, there is an integer $c_u\leq m'$, and a partition 
$$
W_u=W_{u,0}\cup W_{u,1}\cup \ldots \cup W_{u,c_{u}},
$$
so that $|W_{u,0}|\leq \e^8 |W_u|$ and such that for each $v\in [c_u]$, $W_{u,v}$ is $g_{i_1}(m')$-good in $G'$,
\item For each $u\in [\ell]\setminus \Omega_0$, and each $i\in [S(u)]$ with $R_i(u)\notin \Omega_1$, there is an integer $b_i(u)\leq m'$ and a partition 
$$
R_i(u)=R_{i,0}(u)\cup R_{i,1}(u)\cup \ldots \cup R_{i,b_i(u)}(u),
$$
so that $|R_{i,0}(u)|\leq \e^8 |R_i(u)|$ and such that for each $v\in [b_i(u)]$, $R_{i,v}(u)$ is $g_{i_1}(m')$-good in $G'$.  
\end{itemize}

Note that since $G$ is bipartite, and $W'\subseteq W$, any $g_{i_1}(m')$-good subset of $W'$ in $G'$ is also $g_{i_1}(m')$-good in $G$.  Let $\Omega_2=\{u\in [\ell]: |W_u\cap (\bigcup_{X\in \Omega_1}X)|\geq \e^4 |W_u|\}$. Observe that 
$$
\e^4|\sum_{u\in \Omega_2}W_u|\leq |\bigcup_{X\in \Omega_1}X|\leq \e^8|W'|,
$$
and consequently, $|\sum_{u\in \Omega_2}W_u|\leq \e^4|W'|\leq \e^4|W|$.  Now define 
$$
\Omega=\{u\in [\ell]: u\in \Omega_0\cup \Omega_2\}\cup \{u\in [\ell+1,m]: W_u\in \Omega_1\}.
$$
Observe that $|\bigcup_{X\in \Omega}X|\leq \e^4|W|+\e^8|W'|<\e|W|$. For each $u\in [\ell]\setminus \Omega$, define
$$
W_{u,0}=V_0'(u)\cup \bigcup_{R_{i,v}(u)\in \Omega_1}R_{i,v}(u).
$$
Further, if $V_{s+1}(u)=\emptyset$, define $c_u:=|\{R_{i,v}(u): i\in [S(u)], v\in [b_i(u)]\}\setminus \Omega_1\}|$, and let $W_{u,1},\ldots, W_{u,c_{u}}$ be an enumeration of $\{R_{i,v}(u): i\in [S(u)], v\in [b_i(u)]\}\setminus \Omega_1$.  On the other hand, if $V_{s+1}(u)\neq \emptyset$, define 
$$
c_u:=1+|\{R_{i,v}(u): i\in [S(u)], v\in [b_i(u)]\}\setminus \Omega_1\}|,
$$
 and let $W_{u,1},\ldots, W_{u,c_{u}}$ enumerate $\{V_{s+1}(u)\}\cup (\{R_{i,v}(u): i\in [S(u)], v\in [b_i(u)]\}\setminus \Omega_1)$.  Since $u\notin \Omega_2$,
$$
|W_{u,0}|\leq |V'_0|+|\bigcup_{R_{i,v}(u)\in \Omega_1}R_{i,v}(u)|\leq 5\e^2|W_u|+\e^4|W_u|\leq 6\e^2|W_u|.
$$
Furthermore, by construction, for all $u\in [\ell]\setminus \Omega$ and all $1\leq v\leq c_u$, $W_{u,v}$ is $\max\{f'(i_1),g_{i_1}(m')\}$-good, and
$$
c_u\leq S(u)m'+1\leq P_{i_1-1}m'+1\leq m'(P_{i_1-1}+1)\leq M_{i_1}(P_{i_1-1}+1).
$$

We have now constructed, for all $u\in [m]\setminus \Omega$, a partition 
$$
W_u=W_{u,0}\cup W_{u,1}\cup \ldots \cup W_{u,c_u},
$$
where $|W_{u,0}|\leq 6\e^2 |W_u|<\e|W_u|$, and where for all $1\leq v\leq c_u$, $W_{u,v}$ is $\mu^*$-good, where $\mu^*:=\max\{g_{i_1}(m'), f'(i_1)\}$.  Further, for each $u\in [m]\setminus \Omega$, $c_u\leq M_{i_1}(P_{i_1-1}+1)\leq M$.  Since by construction, $\mu^*\leq f(M_{i_1}(P_{i_1-1}+1))$, this finishes the proof.
\end{proofof}
\vspace{2mm}

We can now give an easy proof of Theorem \ref{thm:functionstablereg}.

\vspace{2mm}

\begin{proofof}{Theorem \ref{thm:functionstablereg}}
Fix $k\geq 1$, $\e>0$ and $f:\mathbb{N}\rightarrow (0,1]$ non-increasing with $\lim_{n\rightarrow \infty}f(n)=0$.  Let $d=d(k)$ be as in part (2) of Theorem \ref{thm:treerank}. We will want to apply Theorem \ref{thm:goodstrong} using an auxiliary function, $\psi$, which will go to $0$ sufficiently fast compared to $f$ (this will be needed to later in the proof, specifically in displayed equation (\ref{lastline})).  In particular, we define $\psi: \mathbb{N}\rightarrow (0,1]$ by setting $\psi(x)=\e^4f(\e^{-4}x)^2/6x$.  Let $M=M(d,\psi, \e^2)$ be as in Theorem \ref{thm:goodstrong}.

Suppose $G=(V,E)$ is a $k$-stable graph on at least $N$ vertices. Consider $\bip(G)=(U\cup W, E')$. Note $\rk(W)\leq d$ in $\bip(G)$, by Theorem \ref{thm:treerank}.  By Theorem \ref{thm:goodstrong}, there some $m\leq M$ and a partition $W=W_0\cup \ldots\cup W_m$ such that $|W_0|\leq \e^2 |W|$ and such that for each $1\leq i\leq m$, $W_i$ is $\psi(m)$-good in $\bip(G)$.  By definition of $\bip(G)$, this yields a partition $V=V_0\cup \ldots \cup V_m$ so that $|V_0|\leq \e^2|V|$ and for each $i\in [m]$, $V_i$ is $\psi(m)$-good in $G$.   We now do a bit of work to obtain the desired equitability condition in Theorem \ref{thm:functionstablereg}.  Let $\Omega=\{i\in [m]: |V_i|<\e^2|V|/m\}$, and let $V_0'=V_0\cup \bigcup_{i\in \Omega}V_i$. Note that
$$
|V_0'|\leq \e^2 |V|+|\Omega|\e^2|V|/m\leq 2\e^2|V|.
$$
Set $K=\lfloor \e^4|V|/m\rfloor$.  For each $i\in [m]\setminus \Omega$, choose a partition $V_i=V_{i,0}\cup V_{i,1}\cup \ldots \cup V_{i,r_i}$ so that for each $1\leq v\leq r_i$, $|V_v|=K$ and $|V_{i,0}|<K$.  Now let $R_0=V_0'\cup \bigcup_{i\in [m]\setminus \Omega}V_{i,0}$, and let $R_1,\ldots, R_{m'}$ be an enumeration of $\{V_{i,j}: i\in [m]\setminus \Omega, j\in [r_i]\}$.  Then $m'\leq \e^{-4}m$, and
$$
|R_0|\leq 2\e^2|V|+mK\leq \e^2|V|+\e^4|V|<\e|V|.
$$
For each $1\leq i\leq m'$, $R_i$ is contained in some $V_{u,v}$ which is $\psi(m)$-good, so for all $a\in V$, there is some $\delta=\delta(a)$ so that $|N^{\delta}(a)\cap V_{u,v}|\leq \psi(m)|V_{u,v}|$.  Then
\begin{align}\label{lastline}
|N^{\delta}(a)\cap R_i|\leq \psi(m)|V_{u,v}|=\frac{\psi(m)|V_{u,v}|}{K}|R_i|\leq \e^{-4}m\psi(m)|R_i|<f(m')^2|R_i|/2,
\end{align}
where the second inequality is by definition of $K$, and since $|V_{u,v}|\leq |V|$, and the last inequality is by definition of $\psi$ and because $m'\leq \e^{-4}m$.  Thus for each $1\leq i\leq m'$, $R_i$ is $f(m')^2/2$-good.

By Corollary \ref{cor:goodpair}, for each $1\leq i,j\leq m'$, the pair $(R_i,R_j)$ has density in $[0,f(m'))\cup (1-f(m'),1]$.  By Proposition \ref{fact:homimpliesrandombinary}, $(R_i,R_j)$ satisfies $\disc_2(f(m'))$.
\end{proofof}

\section{Slicewise stable properties admit binary error}\label{subsec:wsbinary}

In this section we prove Proposition \ref{prop:wsbinary}, which shows that a slicewise stable hereditary $3$-graph property admits binary $\vdisc_3$-error and $\disc_{2,3}$-error.  

Recall that if $H=(V,E)$ is a $3$-graph and $a\in V$, then $H_a$ is the graph $(V, N(a))$, where $N(a)=\{bc\in {V\choose 2}: abc\in E\}$.  Recall further that $H$ is slicewise $k$-stable if and only if $H_a$ is $k$-stable for all $a\in V$.  Given a slicewise stable $H=(V,E)$, our strategy will be to find strong decompositions of $H_a$ for all $a\in V$, yielding a partition of $V\times V$ with special properties.  We will then apply the multicolored regularity lemma to the partition of $V\times V$.

The main technical lemma is Lemma \ref{lem:wsbinarylem1} below.  Its proof is closely based on the proof of Theorem \ref{thm:goodstrong} of the preceding section.   We will use the following notation in the proof.  Given a set $V$, $a\in V$, and $\Gamma\subseteq V\times V$, we will let $\Gamma(a):=\{b\in V: (a,b)\in \Gamma\}$.  Given a $3$-graph $H=(V,E)$, $a\in V$, and $C,B\subseteq V$, we let $H_a^{bip}[C,B]$ be the graph with vertex set $\{x_c: c\in C\}\cup B$ and edge set $\{x_cb: cb\in E(H_a)\}$.  Note that if $B'\subseteq B$, and $B'$ is $\e$-good in $H_a^{bip}[V,B]$, then $B'$ is also $\e$-good in $H_a$.  Further, the rank of $B'$ in $H_a^{bip}[V,B]$ is the same as the rank of $B'$ in $H_a$.

\begin{lemma}\label{lem:wsbinarylem1}
Suppose $d\geq 0$.  For all $\e>0$, and non-increasing functions $f:\mathbb{N}\rightarrow [0,1]$ with $\lim_{n\rightarrow \infty}f(n)=0$, there is $M\geq 1$ such that the following hold.  Suppose $H=(V,E)$ is a $3$-graph, $\Gamma\subseteq V\times V$, and $\Gamma=\Gamma_1\cup \ldots \cup \Gamma_m$ is a partition such that for all $a\in V$, and all $i\in [m]$, $\rk(\Gamma_i(a))\leq d$ in $H_a$. 

Then there exists $1\leq s\leq M$ and $\Omega\subseteq [m]$ so that $|\bigcup_{u\in \Omega}\Gamma_u|\leq \e |V|^2$, and such that for each $u\in [m]\setminus \Omega$, there is $s_u\leq s$, and a partition $\Gamma_u=\Gamma_{u,0}\cup \Gamma_{u,1}\cup \ldots \cup \Gamma_{u,s_u}$ with $|\Gamma_{u,0}|\leq \e |\Gamma_u|$, and so that for all $a\in V$ and $1\leq i\leq s_u$, $\Gamma_{u,i}(a)$ is $f(s)$-good in $H_a$.  
\end{lemma}
\begin{proof}
We proceed by induction on $d\geq 0$.  To ease notation, let $I$ be the set of non-increasing functions $f:\mathbb{N}\rightarrow [0,1]$ with $\lim_{n\rightarrow \infty}f(n)=0$.

\underline{Case $d=0$}:  Suppose $d=0$.  Given $\e>0$ and $f\in I$,  set $M=1$. Suppose $H=(V,E)$ is a $3$-graph, $\Gamma\subseteq V\times V$, and $\Gamma=\Gamma_1\cup \ldots \cup \Gamma_m$ is a partition such that for all $a\in V$ and $i\in [m]$, $\rk(\Gamma_i(a))\leq 0$ in $H_a$.   This means that for all $a\in V$, each of $\Gamma_1(a),\ldots, \Gamma_m(a)$ are $0$-good in $H_a$.  Setting $s=1$ and $\Gamma_{i,1}=\Gamma_i$ for each $1\leq i\leq m$, we are done.

\underline{Case $d>0$}:  Suppose now $d>0$ and assume by induction we have defined $M(d-1,  f, \e)$ for all $\e>0$ and $f\in I$.  

Fix $\e>0$ and $f\in I$.  We inductively define a new function $f'\in I$, and a sequence $(M_i,g_i)_{i\in \mathbb{N}}$, where each $M_i\in \mathbb{N}$ and each $g_i\in I$. 

\underline{Step $1$}: Define $g_1$ by setting $g_1(x)=\e^8f(2x+1)$. Choose $M_1$ by applying the inductive hypothesis to $(d-1,  g_1, \e^9)$, and set  $f'(1)=g(M_1)/2$.

\underline{Step $i+1$}: Suppose we have defined $g_j$, $M_j$, and $f'(j)$ for each $1\leq j\leq i$.  Let $P_i=\sum_{j=1}^i (\e f'(j))^{-1}$, and define $g_{i+1}$ by setting $g_{i+1}=g((P_i+1)\cdot x)$.  Choose $M_{i+1}$ by applying the inductive hypothesis to $(d-1,  g_{i+1}, \e^8)$, and set $f'(i+1)=g(M_{i+1})(i+2)^{-1}f'(i)$.

This finishes the definition of $f'$ and our sequence $(M_i,g_i)_{i\in \mathbb{N}}$.  Set 
$$
M=\max\{M_i: i\leq 2\e^{-1}\}.   
$$

Now suppose $H=(V,E)$ is a $3$-graph, $\Gamma\subseteq V\times V$, and $\Gamma=\Gamma_1\cup \ldots \cup \Gamma_m$ is a partition such that for all $a\in V$ and $i\in [m]$, $\rk(\Gamma_i(a))\leq d$ in $H_a$.  If $|\Gamma|\leq \e|V|^2$, then we are done by setting $s=1$ and $\Omega=[m]$.  So we may assume $|\Gamma|>\e|V|^2$.  

For each $1\leq \ell\leq m$, let
$$
V_\ell=\{a\in V: \text{ $\ell$-many elements in }\{\Gamma_1(a),\ldots, \Gamma_m(a)\}\text{ have rank $d$ in $H_a$}\}.
$$
For each $u\in [m]$, let 
$$
V_u^*=\{a\in V: \rk(\Gamma_u(a))=d\}\text{ and }\Gamma_u^*=\Gamma_u\cap (V_u^*\times V).
$$  

For each $u\in [m]$ and $a\in V_{u}^*$,   Apply Lemma \ref{lem:goodsets1} with parameters $d$, $\e^8$, and $f'$ to the set $\Gamma_u(a)$ in $H_a[V,\Gamma_u(a)]$  to obtain some $t^a(u)\leq |\Gamma_u(a)|$, a sequence $(s^a_i(u))_{i\in t^a(u)}$, and a partition $\Gamma_u(a)=\bigcup_{i\in [t^a(u)]}\bigcup_{j\in [s^a_i(u)]}V^a_{i,j}(u)$ such that the following hold, where for each $1\leq i\leq t^a(u)$, $V^a_i(u):=\bigcup_{j\in [s^a_i(u)]}V^a_{i,j}(u)$ and $Z^a_{i,j}(u)= \Gamma_u(a)\setminus (\bigcup_{i'=1}^{i-1}V^a_{i'}(u)\cup \bigcup_{j'=1}^{j-1}V^a_{i,j'}(u))$.
\begin{enumerate}
\item[(i)] For each each $i\in [t^a(u)]$ and $j\in [s^a_i(u)]$, $|V^a_{i,j}(u)|=f'(i)|Z^a_{i,j}(u)|$, and $\rk(V^a_{i,j}(u))<d$ in $H_a$,
\item[(ii)] For each $i\in [t^a(u)]$, $\Gamma_u(a)\setminus (\bigcup_{i'=1}^i V^a_{i'}(u))$ is $f'(i)$-good in $H_a$.
\end{enumerate}

For each $u\in [m]$, let $t(u)=\max\{t^a(u): a\in V_u^*\}$.  Then for each $i\leq t(u)$, let 
$$
s_i(u):=\max\{s^a_i(u): a\in V^*_{u}\text{ and }i\leq t^a(u)\}.
$$
For each $a\in V_u^*$, if $t^a(u)<i\leq t(u)$, set $V^a_i(u)=\emptyset$, and if $i\leq t^a(u)$ but $j>s_i^a(u)$, set $V^a_{i,j}(u)=\emptyset$. For each $i\in [t(u)]$ and $j\in [s_i(u)]$, we now define 
$$
\Sigma_{i,j}(u):=\{(a,b)\in V^2: a\in V_{u}^*\text{ and }b\in V_{i,j}^a(u)\}.
$$
For each $i\in [t(u)]$, set $\Sigma_i(u)=\bigcup_{j\in [s_i(u)]}\Sigma_{i,j}(u)$.  Finally, we set $t=\max\{t(u): u\in [m]\}$, and given $u\in [m]$, if $t(u)<i\leq t$, define $\Sigma_i(u)=\emptyset$.  

We can now define, for each $i\in [t]$,  $\Sigma_i=\bigcup_{u\in [m]}\Sigma_i(u)$.  Let $i_1$ be minimal such that $|\Sigma_{i_1}|\leq \e^8 |V|^2$.  Clearly such an $i_1$ exists and is at most $\e^{-8}$.  We now set
$$
\Omega_0=\{u\in [m]: |\Sigma_{i_1}(u)|\geq \e^4 |\Gamma_u|\}.
$$
Note that by definition, $\e^4|\bigcup_{u\in \Omega_0}\Gamma_u|\leq |\Sigma_{i_1}|\leq \e^8 |V|^2$, so $|\bigcup_{u\in \Omega_0}\Gamma_u|\leq \e^4 |V|^2$.  For each $u\in [m]\setminus \Omega_0$, set 
$$
J(u)=\{j\in [i_1-1]:  |\Sigma_j(u)|\leq \e^4f'(i_1-1)|\Gamma_u|\}.
$$
By definition, $|\bigcup_{j\in J(u)}\Sigma_j(u)|\leq |J(u)|\e^4f'(i_1-1)|\Gamma_u|\leq (i_1-1)\e^4f'(i_1-1)|\Gamma_u|<\e^4|\Gamma_u|$, where the last inequality is because, by construction, $f'(i_1-1)<(i_1)^{-1}$.  We then set $I(u)=[i_1-1]\setminus J(u)$, and for each $i\in I(u)$, let 
$$
T_i(u)=\{j\in [s_i(u)]: |\Sigma_{i,j}(u)|\geq \e^2f'(i)|\Sigma_i(u)|\}.
$$
Setting $t_i(u)=|T_i(u)|$, we have that $t_i(u)$ is at most $\e^{-2}f'(i)^{-1}$.  Now define 
 $$
\Sigma_{i,0}(u)=\bigcup_{j\in [s_i(u)]\setminus T_i(u)}\Sigma_{i,j}(u).
$$
We claim $|\Sigma_{i,0}(u)|\leq \e^2|\Sigma_i(u)|$.  If $[s_i(u)]\setminus T_i(u)=\emptyset$, then this is obvious.  Otherwise, consider $j_0=\min([s_i(u)]\setminus T_i(u))$.  Then for all $a\in V$, either $\Sigma_{j_0,0}(u)(a)=\emptyset$, or 
$$
|\Sigma_{i,j_0}(u)(a)|=|V^a_{i,j_0}(u)|=f'(i)|Z^a_{i,j_0}(u)|.
$$
Therefore if $\zeta_{i,j_0}(u):=\{(a,b): b\in Z^a_{i,j_0}(u)\}$, then $|\Sigma_{i,j_0}(u)|=f'(i)|\zeta_{i,j_0}(u)|\leq \e^2f'(i)|\Sigma_i(u)|$, where the last inequality follows from the fact that $j_0\notin T_i(u)$.  Consequently, $|\zeta_{i,j_0}(u)|\leq \e^2|\beta_i(u)|$.  By construction and the minimality of $j_0$, for all $a\in V$ and $j\in [s_i(u)]\setminus T_i(u)$, 
$$
V_{i,j}^a(u)\subseteq V^a_{i,j_0}\cup Z^a_{i,j_0}.
$$
This implies $\beta_{i,0}(u)\subseteq \Sigma_{i,j_0}(u)\cup \zeta_{i,j_0}$, and thus $|\Sigma_{i,0}(u)|\leq 2\e^2 |\Sigma_i(u)|$.  We now define $\Sigma_0(u)=\Sigma_{i_1}(u)\cup (\bigcup_{i\in J(u)}\Sigma_j(u))\cup \bigcup_{i\in I(u)}\Sigma_{i,0}$.   Then by above,
$$
|\Sigma_0(u)|\leq \e^4|\Gamma_u|+\e^4|\Gamma_u|+\sum_{i\in I(u)}|V_{i,0}|\leq 2\e^4|\Gamma_u|+\sum_{i\in I(u)}2\e^2|\Sigma_i(u)|\leq  3\e^2|\Gamma_u|.
$$
Recall that by (ii), $\zeta_{i_1}(u):=\Gamma_u\setminus (\bigcup_{v=1}^{i_1}\Sigma_i(u))$ has the property that for all $a\in V$, $\zeta_{i_1}(u)(a)$ is $f'(i_1)$-good in $H_a$.  If $|\zeta_{i_1}(u)|<\e^2 |\Gamma_u|$, set $\Sigma_0'(u)=\Sigma_0(u)\cup \zeta_{i_1}(u)$ and $\Sigma_{s+1}(u)=\emptyset$.  If $|\zeta_{i_1}(u)|\geq \e^2 |\Gamma_u|$, set $\Sigma_0'=\Sigma_0$ and $\Sigma_{s+1}(u)=\zeta_{i_1}(u)$.  In either case, $|\Sigma_0'(u)|\leq 5\e^2|\Gamma_u|$, and we now have a partition
$$
\Gamma^*_u=\Sigma'_0(u)\cup \Sigma_{s+1}(u)\cup \bigcup_{i\in I(u)}\bigcup_{j\in T_i(u)}\Sigma_{i,j}(u).
$$
Set $S(u):=\sum_{i\in I(u)} t_i(u)$.  Since each  $t_{i}(u)\leq (\e^2 f'(i))^{-1}$, we have that $S(u)\leq \sum_{i=1}^{i_1-1}\e^{-2}f'(i)^{-1}=P_{i_1-1}$.  Let $\gamma_1(u),\ldots, \gamma_{S(u)}(u)$ be an enumeration of $\{\Sigma_{i,j}(u): i\in I(u), j\in T_i(u)\}$.  

For each $u\in [m]$, define
$$
\Gamma'_u=\Gamma_u\cap ((V\setminus V_u^*)\times V),
$$
 and set $\Gamma'= \bigcup_{u\in [m]}\bigcup_{i\in [S(u)]}\gamma_i(u)$. Then we have the following partition of $\Gamma'$. 
 $$
 \calP=\{\gamma_i(u): u\in [m],  i\in [S(u)]\}\cup \{\Gamma_1',\ldots, \Gamma_m'\}
 $$
By construction, for every $a\in V$ and $\gamma\in \calP$, $\gamma(a)$ has rank at most $d-1$ in $H_a$.  By the induction hypothesis, there are $\Omega_1\subseteq \calP$ and $m'\leq M_{i_1}=M(d-1,g_{i_1},\e^8)$ such that the following hold.
\begin{itemize}
\item $|\bigcup_{\gamma\in \Omega_1}\gamma|\leq \e^8 |V|^2$,
\item For each $u\in [m]$ with $\Gamma_u'\notin \Omega_1$,  there is an integer $\ell_u\leq m'$, and a partition 
$$
\Gamma'_u=\Gamma'_{u,0}\cup \Gamma'_{u,1}\cup \ldots \cup \Gamma'_{u,\ell_{u}},
$$
so that $|\Gamma'_{u,0}|\leq \e^8 |\Gamma'_u|$ and for each $v\in [\ell_u]$, and $a\in V$, $\Gamma'_{u,v}(a)$ is $g_{i_1}(m')$-good in $H_a$,
\item For each $u\in [m]$ and $i\in [S(u)]$ with $\gamma_i(u)\notin \Omega_1$, there is an integer $b_i(u)\leq m'$ and a partition 
$$
\gamma_i(u)=\gamma_{i,0}(u)\cup \gamma_{i,1}(u)\cup \ldots \cup \gamma_{i,b_i(u)}(u),
$$
so that $|\gamma_{i,0}(u)|\leq \e^8 |\gamma_i(u)|$, and for each $v\in [b_i(u)]$ and $a\in V$, $\gamma_{i,v}(u)(a)$ is $g_{i_1}(m')$-good in $H_a$.  
\end{itemize}

Let $\Omega_2=\{u\in [m]: |\Gamma_u\cap (\bigcup_{\gamma\in \Omega_1}\gamma)|\geq \e^4 |\Gamma_u|\}$. Observe that 
$$
\e^4|\sum_{u\in \Omega_2}\Gamma_u|\leq |\bigcup_{\gamma\in \Omega_1}\gamma|\leq \e^8|V|^2.
$$
Consequently, $|\sum_{u\in \Omega_2}\Gamma_u|\leq \e^4|V|^2$.  For each $u\in [m]\setminus \Omega_2$, define
\[\Omega_1(u)=(\{\gamma_{i,v}(u): i\in [S(u)], v\in [b_i(u)]\}\cup \{\Gamma_{u,j}': j\in [\ell_u]\})\cap \Omega_1\]
and
\[\calP(u)=(\{\gamma_{i,v}(u): i\in [S(u)], v\in [b_i(u)]\}\cup \{\Gamma_{u,j}': j\in [\ell_u]\})\setminus \Omega_1.\]
Set $\tau_{u,0}=\Sigma_0'(u)\cup \bigcup_{\gamma\in \Omega_1(u)}\gamma$. If $\Sigma_{s+1}(u)=\emptyset$, set $d_u:=|\calP(u)|$ and let $\tau_{u,1},\ldots, \tau_{u,d_{u}}$ be an enumeration of $\calP(u)$. If $\Sigma_{s+1}(u)\neq \emptyset$, set $d_u:=1+|\calP(u)|$, and let $\tau_{u,1},\ldots, \tau_{u,d_{u}}$ enumerate $\{V_{s+1}(u)\}\cup \calP(u)$.     Since $u\notin \Omega_2$,
$$
|\tau_{u,0}|\leq  5\e^2|\Gamma_u|+\e^4|\Gamma_u|\leq 6\e^2|\Gamma_u|,
$$
and for all $1\leq v\leq d_u$, and $a\in V$, $\tau_{u,v}(a)$ is $\max\{f'(i_1),g_{i_1}(m')\}$-good in $H_a$.  Note 
$$
d_u\leq S(u)m'+1\leq P_{i_1-1}m'+1\leq m'(P_{i_1-1}+1)\leq M_{i_1}(P_{i_1-1}+1).
$$
We have now constructed, for all $u\in [m]\setminus \Omega_2$, a partition 
$$
\Gamma_u=\tau_{u,0}\cup \tau_{u,1}\cup \ldots \cup \tau_{u,d_u},
$$
where $|\tau_{u,0}|\leq 6\e^2 |\Gamma_u|<\e|\Gamma_u|$, and for all $1\leq v\leq d_u$, $\tau_{u,v}$ is $\mu^*$-good, where $\mu^*:=\max\{g_{i_1}(m'), f'(i_1)\}$.  Since by above, $d_u\leq M_{i_1}(P_{i_1-1}+1)\leq M$, and by construction, $\mu^*\leq f(M_{i_1}(P_{i_1-1}+1))$, this finishes the proof. 
\end{proof}

We now use Lemma \ref{lem:wsbinarylem1} to prove that any slicewise $k$-stable property admits binary $\vdisc_3$-error. 

\begin{corollary}
If $\calH$ is slicewise $k$-stable then it admits binary $\vdisc_3$-error.
\end{corollary}
\begin{proof}
Fix $k\geq 1$, $\e>0$, and $t_0\geq 1$.  Let $d=d(k)$ as in part (1) of Theorem \ref{thm:treerank}. Choose $\e_1>0$ and $\e_2:\mathbb{N}\rightarrow (0,1]$ as in Proposition \ref{prop:wnipdens} for $\e$ and $k$.  Choose $\e_1''\ll \e_1'\ll \e, \e_1, k^{-1}$, and define $\e'_2,\e_2'':\mathbb{N}\rightarrow (0,1]$ non-increasing so that for all $x$, $\e_2''(x)\ll \e_2'(x)\ll \e_2(x)$.  Then let $f$ and $g$ be as in Lemma \ref{lem:refinement} for $\e_1'',\e_2''$.  Now define $\psi:\mathbb{N}\rightarrow (0,1]$ to be any non-increasing function satisfying 
$$
\psi(x)\ll \e_1''f(1,1,t_02^{x^2}(\e_1'')^{-1},x+1)^{-8}g(1,1,t_02^{x^2}(\e_1'')^{-1},x+1)^{-8}x^{-1}.
$$
Note by definition $\lim_{x\rightarrow \infty}\psi(x)=0$.

Apply Lemma \ref{lem:wsbinarylem1} to the function $\psi$, the rank bound $d$, and $\e''_1$  to obtain $M$.  Then choose $T\gg (\e_1'')^{-1}t_0\e^{-1}M$ and $N\gg T$.

Suppose $H=([n],E)$ is a slicewise $k$-stable $3$-graph with $n\geq N$.  Apply Lemma \ref{lem:wsbinarylem1} to $\Gamma=V\times V$ to find $s\leq M$ and a partition $V\times V=\Gamma_0\cup \ldots\cup \Gamma_s$ so that $|\Gamma_0|\leq \e''_1 n^2$ and such that for all $1\leq i\leq s$ and $a\in V$, $\Gamma_i(a)$ is $\psi(s)$-good in $H_a$.  Observe that by Lemma \ref{lem:twosticks}, we know that for all $a\in V$ and $1\leq i,j\leq s$, there is some $\delta=\delta(i,j,a)$ so that 
\begin{align}\label{al:wspsi}
|E^{\delta}_a\cap K_2[\Gamma_i(a),\Gamma_j(a)]|\geq (1-2\sqrt{\psi(s)})|K_2[\Gamma_i(a),\Gamma_j(a)]|.
\end{align}
Let $\calP=\{V_1,\ldots,V_{t_1}\}$ be a minimal partition of $V$ so that if $a,a'$ are in the same part of $\calP$, then $\delta(i,j,a)=\delta(i,j,a')$ for all $1\leq i,j\leq s$.  Note  $t_1\leq 2^{s^2}$.

We now adjust $\calP$ to obtain an equipartition.  Set $t=(t_0t_1\e_1'')^{-1}$, and $K=\lfloor n/t\rfloor$.  For each $1\leq i\leq t_1$, choose a partition  $V_i=V_{i,0}\cup V_{i,1}\cup \ldots \cup V_{i,r_i}$ so that $|V_{i,0}|<K$ and for each $1\leq j\leq r_i$, $|V_{i,j}|=K$.   Define $V_0=\bigcup_{i=1}^{t_1}V_{i,0}$ and $V'=V\setminus V_0$. By construction, we have 
$$
|V_0|\leq \frac{t_1}{t}n\leq \e_1''n.
$$
Choose an equipartition $V_0=\bigcup_{i=1}^{t_1}\bigcup_{j=1}^{r_i}V^0_{i,j}$, and for each $1\leq i\leq t_1$ and $1\leq j\leq r_i$, set $V_{i,j}'=V_{i,j}\cup V_{0,i}$.  Clearly $|V_{i,j}'\setminus V_{i,j}|\leq \e_1''|V_{i,j}'|$.  We now define $\calQ_{vert}=\{V_{i,j}':i\in [t_1], j\in [r_i]\}$ and $t_2:=|\calQ_{vert}|$.  Note $t_2\leq t$.  Our next goal is to define a partition of $K_2[X,Y]$ for all $XY\in {\calQ_{vert}\choose 2}$.

Given $XY\in {\calQ_{vert}\choose 2}$, we know that $XY$ has the form $V_{i,j}'V_{i',j'}'$ for some $i,i'\in [t_1]$, $j\in [r_i]$, and $j'\in [r_{i'}]$.   If $1\leq i<i'\leq t_1$, then for each $1\leq j\leq r_i$, $1\leq j'\leq r_{i'}$, and $0\leq u\leq s$, define
$$
F_u[V_{i,j}', V_{i',j'}']=\{vv'\in K_2[V_{i,j}', V_{i',j'}']: v'\in \Gamma_u(v)\}.
$$
Similarly, if $1\leq i=i'\leq t_1$, then for each $1\leq j<j'\leq r_i$ and $0\leq u\leq s$, define
$$
F_u[V_{i,j}', V_{i,j'}']=\{vv'\in K_2[V_{i,j}', V_{i,j'}']: v'\in \Gamma_u(v)\}.
$$
Now set $\calQ_{edge}=\{F_u[X,Y]: XY\in {\calQ_{vert}\choose 2}, 0\leq u\leq s\}$.  Then $\calQ:=(\calQ_{vert},\calQ_{edge})$ is a $(t_2,s+1)$-decomposition of $V$.  Given $XYZ\in {\calQ_{vert}\choose 3}$ and $0\leq u,v,w\leq s$, let 
$$
Q_{XYZ}^{u,v,w}=(X\cup Y\cup Z, F_u[X,Y]\cup F_v[X,Z]\cup F_w[Y,Z]).
$$
We claim that for all $XYZ\in {\calQ_{vert}\choose 3}$ and $1\leq u,v,w\leq s$, there exists a value $\delta=\delta(X,Y,Z,u,v,w)\in \{0,1\}$ so that 
$$
|E^\delta\cap K_3^{(2)}(Q_{XYZ}^{u,v,w})|\geq (1-\psi(s)^{1/4})|K_3^{(2)}(Q_{XYZ}^{u,v,w})|.
$$
Fix $XYZ\in {\calQ_{vert}\choose 3}$.  By construction, there are $i_1,i_2,i_3\in [t_2]$, $j_1\in [r_{i_1}]$, $j_2\in [r_{i_2}]$, and $j_3\in [r_{i_3}]$ such that $X=V_{i_1,j_1}'$, $Y=V_{i_2,j_2}'$, and $Z=V_{i_3,j_3}'$.  

Suppose first $i_1<i_2,i_3$.  Then for all $xy\in F_u[X,Y]$, $y\in\Gamma_u(x)$ and for all $xz\in F_w[X,Z]$, $z\in \Gamma_w(x)$.  By (\ref{al:wspsi}), there is $\delta=\delta(u,w,x)\in \{0,1\}$ so that 
$$
|N_H^{\delta}(x)\cap K_2[\Gamma_u(x),\Gamma_w(x)]|\geq (1-2\sqrt{\psi(s)})|\Gamma_u(x)||\Gamma_w(x)|.
$$
Let $G'=(X\cup Y\cup Z, \{xyz\in K_3[X,Y,Z]: xy\in F_u[X,Y], xz\in F_w[X,Z]\})$.  Then $|E^{\delta}\cap K_3^{(2)}(Q_{XYZ}^{u,v,w})|$ is at least
\begin{align*}
|K_3^{(2)}(Q_{XYZ}^{u,v,w})|-\sum_{x\in X}2\sqrt{\psi(s)}|\Gamma_u(x)||\Gamma_w(x)|&\geq |K_3^{(2)}(Q_{XYZ}^{u,v,w})|(1-2\sqrt{\psi(s)}t^2)\\
&\geq |K_3^{(2)}(Q_{XYZ}^{u,v,w})|(1-\psi(s)^{1/4}),
\end{align*}
where the first inequality is since $|\Gamma_u(x)||\Gamma_w(x)|\leq t^2|Y||Z|$, and the second inequality is by the definition of $t$ and our choice of $\psi$ in the first paragraph of the proof.

A symmetric argument works if $i_2<i_1,i_3$ or $i_3<i_1,i_2$.  Suppose now $i_1=i_2<i_3$.  Since $X$ and $Y$ are distinct, we may assume after relabeling that $j_1<j_2$.  Then we have that for all $xy\in F_u[X,Y]$, $y\in \Gamma_u(x)$ and for all $xz\in F_w[X,Z]$, $z\in \Gamma_w(x)$.  We now finish this case as above, with $\delta=\delta(u,w,x)$.  A symmetric argument works if $i_1=i_3<i_2$ or $i_2=i_3<i_1$.  So we are now left with the case where   $i_1=i_2=i_3$.  Since $X,Y,Z$ are pairwise distinct, we must have the $j_1,j_2,j_3$ are pairwise distinct.  By relabeling if necessary, we may assume $j_1<j_2,j_3$.  Then we know that for all  $xy\in F_u[X,Y]$, $y\in \Gamma_u(x)$ and for all $xz\in F_w[X,Z]$, $z\in \Gamma_w(x)$.  We can thus finish this case as above with $\delta=\delta(u,w,x)$.

We now apply Lemma \ref{lem:refinement} to $\calQ$ and $\calP$, where $\calP$ is the trivial decomposition of $V$ consisting of $\calP_{vert}=\{V\}$ and $\calP_{edge}=\{K_2[V,V]\}$.  Then we obtain  $\ell'\leq f(1,1,t_2,s+1)$ and $t'\leq g(1,1,t_2,s+1)$, and a $(t',\ell',\e_1,',\e_2'(\ell'))$-decomposition $\calR$ of $V$ in which every element of $\calR_{edge}$ satisfies $\disc_2(\e_2'(\ell'))$ and which is an $(\e_1',\e_2'(\ell'))$-approximate refinement of both $\calQ$ and $\calP$ (see Definition \ref{def:refinement}). Let $\Sigma_1$ be as in the definition of $\calR$ being an $(\e_1',\e_2'(\ell'))$-approximate refinement of $\calQ$.  Then for all $XY\in {\calR_{vert}\choose 2}\setminus \Sigma_1$, there is $\phi(R_{XY}^{\alpha})\in \calQ_{edge}$ so that $|R_{XY}^{\alpha}\setminus \phi(R_{XY}^{\alpha})|\leq \e_1'|R_{XY}^{\alpha}|$, and so that $R_{XY}^{\alpha}\setminus \phi(R_{XY}^{\alpha})$ has $\disc_2(\e_2'(\ell))$.  Now define
$$
\Sigma_2=\Big\{XY\in {\calR_{vert}\choose 2}: F_0[X,Y]|\geq \sqrt{\e_1''}|X||Y|\Big\}.
$$

Since $|\Gamma_0|\leq \e_1'' n^2$ and $\calR_{vert}$ is an equipartition of $V$, we know that $|\Sigma_2|\leq \sqrt{\e_1''}(t')^2$.  We now set $\Sigma=\Sigma_1\cup \Sigma_2$.  Note 
$$
|\Sigma|\leq \e'_1 (t')^2+ \sqrt{\e_1''}(t')^2\leq \e_1 t^2.
$$
Now suppose we have $XYZ\in \calR_{vert}$ satisfies $XY, YZ, XZ\notin \Sigma$.  We claim that $XYZ$ is $\e'_1$-homogeneous with respect to $H$.     Define
$$
\Omega(XY)=\Big\{\alpha\in [\ell']: \phi(R_{XY}^{\alpha})\subseteq F_u[X',Y'] \text{ for some }u\neq 0 \text{ and } X'Y'\in {\calQ_{vert}\choose 2}\Big\}.
$$
Since $XY\notin \Sigma_2$, we can deduce that $|\Omega(XY)|\geq (1-(\e_1'')^{1/4})\ell'$.  We can then define $\Omega(XZ), \Omega(YZ)$ the same way and deduce that  $|\Omega(YZ)|,|\Omega(XZ)|\geq (1-(\e_1'')^{1/4})\ell'$.  Now fix some $\alpha\in \Omega(XY), \beta\in \Omega(XZ), \gamma\in \Omega(YZ)$.  Let $X',Y', Z'\in \calQ_{vert}$ and $1\leq u,v,w\leq s$ be such that $\phi(R_{XYZ}^{\alpha,\beta,\gamma})=Q_{X'Y'Z'}^{u,v,w}$.  By above, there is some $\delta=\delta(X',Y',Z',u,v,w)$ so that 
$$
|E^{\delta}\cap K_3^{(2)}(Q_{X'Y'Z'}^{u,v,w})|\geq (1-\psi(s)^{1/4})|K_3^{(2)}(Q_{X'Y'Z'}^{u,v,w})|.
$$
Let $A=R_{XY}^{\alpha}\setminus \phi(R_{XY}^{\alpha})$, $B=R_{XZ}^{\beta}\setminus \phi(R_{XY}^{\beta})$ and $C=R_{XZ}^{\gamma}\setminus \phi(R_{XZ}^{\gamma})$. Then each of $A,B,C$ satisfy $\disc_2(\e_2'(\ell'))$ and 
$$
|A|\leq (\e_1')^{1/4} |R_\alpha^{XY}|,\text{ }|B|\leq (\e_1')^{1/4} |R_{YZ}^{\beta}|,\text{ and }|C|\leq (\e_1')^{1/4} |R_{XZ}^{\gamma}|.
$$
Let $G_A=(X\cup Y\cup Z, A\cup  R_{YZ}^{\gamma}\cup R_{XZ}^{\gamma})$, $G_B=(X\cup Y\cup Z,   R_{XY}^{\alpha}\cup B\cup R_{XZ}^{\gamma})$, and $G_C=(X\cup Y\cup Z, R_{XY}^{\alpha}\cup R_{YZ}^{\beta}\cup C)$.  Then 

\begin{align*}
&|E^{\delta}\cap K_3^{(2)}(R_{XYZ}^{\alpha,\beta,\gamma})| \\
&\geq | K_3^{(2)}(R_{XYZ}^{\alpha,\beta,\gamma})\setminus K_3^{(2)}(Q_{X'Y'Z'}^{u,v,w})|-\psi(s)^{1/4}|K_3^{(2)}(Q_{X'Y'Z'}^{u,v,w})|.
\end{align*}
The latter expression equals
\[|K_3^{(2)}(R_{XYZ}^{\alpha,\beta,\gamma})|-|K_3^{(2)}(G_A)|-|K_3^{(2)}(G_B)|-|K_3^{(2)}(G_C)|- \psi(s)^{1/4}|K_3^{(2)}(Q_{X'Y'Z'}^{u,v,w})|,\]
which in turn, by Corollary \ref{cor:counting}, is at least
\[|K_3^{(2)}(R_{XYZ}^{\alpha,\beta,\gamma})|(1-3(\e_1')^{1/4}-\frac{(t')^2}{t^2}\psi(s)^{1/4})\geq  |K_3^{(2)}(R_{XYZ}^{\alpha,\beta,\gamma})|(1-(\e'_1)^{1/8}),\]
where the final inequality uses the definition of $\psi$.  Thus, we have that for every $\alpha\in \Omega(XY), \beta\in \Omega(XZ), \gamma\in \Omega(YZ)$, $R_{XYZ}^{\alpha,\beta,\gamma}$ is $\e_1'$-homogeneous with respect to $H$.  Combining this with Proposition \ref{prop:wnipdens}, the size bounds on $\Omega(XY)$, $\Omega(YZ)$, $\Omega(XZ)$ and since $\e_1'\ll \e_1$, we find that $XYZ$ is in fact $\e_1$-homogeneous with respect to $H$.   Consequently, $XYZ$ satisfies $\vdisc_3(\e_1)$ with respect to $H$ by Proposition \ref{prop:homimpliesrandomv}.  Since $|\Sigma|\leq \e_1 t^2$, this shows $\calP$ is $\vdisc_3(\e_1)$-homogeneous with respect to $H$ with binary error, as desired.
\end{proof}

Proposition \ref{prop:wsbinary} follows immediately from this and Proposition \ref{prop:simclasses2}.

%% file: chapter6.tex
\chapter{The functional order property and linear $\disc_{2,3}$-error}\label{sec:fop}

In this chapter we prove Theorem \ref{thm:FOP}, which says that a hereditary $3$-graph property admits linear $\disc_{2,3}$-error if and only if it is $\NFOP_2$.   We will begin in Section \ref{subsec:fopchar} by giving equivalent conditions for when a hereditary $3$-graph property has $\FOP_2$ or $\IP_2$.  In particular, we will show a hereditary $3$-graph property has $\IP_2$ (respectively $\FOP_2$) if and only if certain edge-colored auxiliary graphs contain copies of $\IP$ (respectively OP).  In Section \ref{subsec:stablerem}, we prove Theorem \ref{thm:stableremoval}, which is a strong removal lemma for graphs containing few $d$-trees.  In Section \ref{subsec:mainfop}, we then prove our main result of the chapter, which says that a hereditary $3$-graph property admits linear $\disc_{2,3}$-error if and only if it is $\NFOP_2$.  This proof will leverage the recharacterization from Section \ref{subsec:fopchar}, Theorem \ref{thm:goodstrong} of the previous chapter, and the removal lemma, Theorem \ref{thm:stableremoval}.  

\section{Encodings and recharacterizations of $\FOP_2$}\label{subsec:fopchar}
In this section we will give additional characterizations of when a hereditary property is $\NFOP_2$.  These recharacterizations will be used later in our proof of Theorem \ref{thm:FOP}.   Along the way, we will include some related results which demonstrate that $\FOP_2$ is analogous to the order property, in the same way that $\VC_2$-dimension is analogous to $\VC$-dimension.

The main recharacterization will be presented in terms of auxiliary graphs defined from regular triads in decompositions of $3$-graphs.

\begin{definition}
Suppose $\e_1>0$, $\e_2:\mathbb{N}\rightarrow \mathbb{N}$, $\ell, t\geq 1$, $V$ is a set, and $\calP$ is a $(t,\ell, \e_1,\e_2)$-decomposition for $V$ consisting of $\{V_i: i\in [t]\}$ and $\{P_{ij}^{\alpha}: ij\in {[t]\choose 2}, \alpha\in [\ell]\}$.    Define
\[\calP_{cnr}=\Big\{P_{jk}^\beta P_{ik}^{\gamma}: ijk\in {[t]\choose 3}, \beta,\gamma \in [\ell], \text{ and }P_{jk}^\beta, P_{ik}^{\gamma}\text{ satisfy }\disc_2(\e_2(\ell);1/\ell)\Big\}\]
and
\[\calP_{edge}^*=\Big\{P_{ij}^\alpha : ij\in {[t]\choose 2}, \alpha\in [\ell], \text{ and }P_{ij}^\alpha\text{ satisfies }\disc_2(\e_2(\ell);1/\ell)\Big\}.\]

\end{definition}

In the above, cnr stands for ``corner.''  The motivation for this  definition is that given an element $P_{ij}^{\alpha}\in \calP_{edge}$ and an element $P_{uv}^{\beta}P_{uw}^{\gamma}\in \calP_{cnr}$, if $\{v,w\}=\{i,j\}$, the pair $(P_{ij}^{\alpha}, P_{uv}^{\beta}P_{uw}^{\gamma})$ corresponds to a triad from $\calP$, in which case there is a natural way to define a reduced edge, or a reduced non-edge, between $P_{ij}^{\alpha}$ and $P_{uv}^{\beta}P_{uw}^{\gamma}$.  We will now define a bipartite edge colored graph associated to a regular decomposition of a $3$-graph, where a bipartite edge-colored graph is a tuple of the form $(A\cup B, E_1,\ldots, E_t)$ such that $K_2[A,B]=E_1\cup \ldots \cup E_t$.

\begin{definition}\label{def:reducedP}
Suppose $H=(V,E)$ is a $3$-graph, and $\calP$ is a $(t,\ell, \e_1,\e_2)$-decomposition for $V$. Given $G_{ijs}^{\alpha,\beta,\gamma}\in \triads(\calP)$, let $H_{ijs}^{\alpha,\beta,\gamma}=H|G_{ijs}^{\alpha,\beta,\gamma}$ and $d_{ijs}^{\alpha,\beta,\gamma}=|E\cap K_3^{(2)}(G_{ijs}^{\alpha\beta\gamma})|/|K_3^{(2)}(G_{ijs}^{\alpha\beta\gamma})|$.  Then define 
\begin{align*}
\mathbf{E}_1=\{P_{ij}^{\alpha}(P_{js}^{\beta} P_{is}^{\gamma})\in K_2[\calP^*_{edge}, \calP_{cnr}]:  (H_{ijs}^{\alpha,\beta,\gamma},G_{ijs}^{\alpha,\beta,\gamma})\text{ has }&\disc_{2,3}(\e_1,\e_2(\ell))\\ &\text{ and } d_{ijs}^{\alpha,\beta,\gamma}\geq 1/2\},
\end{align*}
\begin{align*}
\mathbf{E}_0=\{P_{ij}^{\alpha}(P_{js}^{\beta} P_{is}^{\gamma})\in K_2[\calP^*_{edge}, \calP_{cnr}]: (H_{ijs}^{\alpha,\beta,\gamma},G_{ijs}^{\alpha,\beta,\gamma})\text{ has }&\disc_{2,3}(\e_1,\e_2(\ell))\\ &\text{ and }d_{ijs}^{\alpha,\beta,\gamma}<1/2\},
\end{align*}
and
\[\mathbf{E}_2=K_2[\calP^*_{edge},\calP_{cnr}]\setminus (\mathbf{E}_1\cup \mathbf{E}_0).\]
\end{definition}
Note that $\mathbf{E}_2$ contains the pairs $P_{ij}^{\alpha}(P_{jk}^{\beta} P_{ik}^{\gamma})$ for which the corresponding triad $G_{ijs}^{\alpha,\beta,\gamma}$ fails  $\disc_{2,3}(\e_1,\e_2(\ell))$, as well as the ``degenerate" pairs $P_{ij}^{\alpha}(P_{uk}^{\beta} P_{vk}^{\gamma})$ which do not correspond to a triad because $uv\neq ij$.  We also observe that $\mathbf{E}_0,\mathbf{E}_1,\mathbf{E}_2$ each depend on $\calP$ and $H$, but we suppress these in the notation as they should be clear from context. Note Definition \ref{def:reducedP} associates to a decomposition $\calP$ of a $3$-graph $H$, an auxiliary edge-colored bipartite graph, $(\calP^*_{edge}\cup \calP_{cnr},\mathbf{E}_1,\mathbf{E}_0,\mathbf{E}_2)$. When the decomposition $\calP$ is $\disc_{2,3}$-homogeneous, elements in $\mathbf{E}_1$ will correspond to triads of density near $1$, and elements in $\mathbf{E}_0$ will correspond to triads of density near $0$. In that case, Definition \ref{def:reducedP} becomes an especially useful reflection of the structure in the original $3$-graph. This will be one of the main ideas behind our proofs in Section \ref{subsec:mainfop}.  We will also consider special copies of certain bipartite graphs made using $\mathbf{E}_1$ and $\mathbf{E}_0$. 

\begin{definition}\label{def:encodingR}
Suppose $G=(A\cup B,E_G)$ is a bipartite graph,  $H=(V,E)$ is a $3$-graph, and $\calP$ is a $(t,\ell,\e_1,\e_2)$-decomposition of $V$.  An \emph{$(A,B)$-encoding of $G$ in $(H,\calP)$} consists of a pair of functions $(g,f)$ where $g:A\rightarrow \calP_{cnr}$ and $f:B\rightarrow \calP^*_{edge}$ are such that the following hold for some $jk\in {[t]\choose 2}$.
 \begin{enumerate}
 \item $\mathrm{Im}(f)\subseteq \{P_{jk}^{\alpha}: \alpha\in [\ell]\}$, and $\mathrm{Im}(g)\subseteq \{P^{\beta}_{ij} P^\gamma_{ik}: i\in [t], \beta,\gamma\in [\ell]\}$, and
 \item for all $a\in A$ and $b\in B$, if $ab\in E_G$, then $g(a)f(b)\in \mathbf{E}_1$, and if $ab\notin E_G$, then $g(a)f(b)\in \mathbf{E}_0$.
 \end{enumerate}
 
If $\mathrm{Im}(g)=\{P^{\beta_1}_{i_1j} P^{\gamma_1}_{i_1k}, \ldots, P^{\beta_s}_{i_sj} P^{\gamma_s}_{i_sk}\}$, then we say $(V_{j}, V_{k}, \bigcup_{u=1}^sV_{i_u})$ is the \emph{partition corresponding to $(g,f)$}.
 \end{definition}
 
 Note that there is an inherent asymmetry in Definition \ref{def:encodingR}, since $A$ is mapped into $\calP_{cnr}$ while $B$ is mapped into $\calP^*_{edge}$.  We will be interested in the cases where $G=H(k)$ and $G=U(k)$, so it will be convenient to set notation for the labeling of their parts. Given $k\geq 1$, we set $U(k)=(A_k\cup C_{\calP([k])},E_{U(k)})$, where $A_k=\{a_i: i\in [k]\}$, $C_{\calP([k])}=\{c_S: S\subseteq [k]\}$, and $E_{U(k)}=\{a_ic_S: i\in S\}$.  Similarly, we set $H(k)=(A_k\cup B_k, E_{H(k)})$ where $A_k=\{a_i: i\in [k]\}$, $B_k=\{b_i: i\in [k]\}$, and $E_{H(k)}=\{a_ib_j: i\leq j\}$.   To ease notation, an \emph{encoding of $H(k)$} will always mean an $(A_k,B_k)$-encoding of $H(k)$, and an \emph{encoding of $U(k)$} will always mean an $(A_k,C_{\calP([k]})$-encoding of $U(k)$.

We now make a definition for when a hereditary $3$-graph property has encodings of a given bipartite graph $G$ using arbitrary regularity parameters.
 
\begin{definition}\label{def:binaryencode}
Suppose $G=(A\cup B, E_G)$ is a bipartite graph and $\calH$ is a hereditary $3$-graph property.  We say that $G$ is \emph{$(A,B)$-encoded in $\calH$} if for all $\e_1>0$ and $\e_2:\mathbb{N}\rightarrow (0,1]$,  there are $T$, $L$ such that for all $N$, the following holds.  

There exist a $3$-graph $H=(V,F)$ with $|V|\geq N$, some $t\leq T$, $\ell\leq L$, and a $(t,\ell,\e_1,\e_2(\ell))$-decomposition $\calP$ of $V$, such that there exists an $(A,B)$-encoding of $G$ in $(H,\calP)$.
\end{definition}
 
We will simply say that $H(k)$ is \emph{encoded in $\calH$} if it is $(A_k,B_k)$-encoded in $\calH$, and $U(k)$ is \emph{encoded in $\calH$} if it is $(A_k,C_{\calP([k])})$-encoded in $\calH$.  We now draw a connection between encodings of $H(k)$ and $U(k)$, and $k$-$\FOP_2$ and $k$-$\IP_2$, respectively.   Analogues of these results also appear in the arithmetic setting in \cite{Terry.2021a}.  

\begin{theorem}\label{thm:suffcond}
Suppose $\calH$ is a hereditary $3$-graph property. 
\begin{enumerate}[label=\normalfont(\arabic*)]
\item If $H(k)$ is encoded in $\calH$, then $\calH$ has $k$-$\FOP_2$.
\item If $U(k)$ is encoded in $\calH$, then $\calH$ has $k$-$\IP_2$.
\end{enumerate}
\end{theorem}

We will show later that the converses of Theorem \ref{thm:suffcond} also hold (see Theorems \ref{prop:fopequiv} and \ref{prop:vc2equiv}).  The proof of Theorem \ref{thm:suffcond} will rely on two facts, Propositions \ref{prop:suffvc2} and \ref{prop:sufffop}.  

\begin{proposition}\label{prop:suffvc2}
For all $k\geq 1$, there exist $\e_1>0$ and $\e_2:\mathbb{N}\rightarrow (0,1]$ such that for all $t,\ell\geq 1$, there is $N$ such that the following hold.  Suppose $H=(V,E)$ is a $3$-graph with $|V|\geq N$, and $\calP$ is a $(t,\ell, \e_1,\e_2(\ell))$-decomposition of $V$. If there exists an encoding of $U(k)$ in $(H,\calP)$, then $H$ has $k$-$\IP_2$. 

\end{proposition}
\begin{proof}
Fix $k\geq 1$. Choose $\delta_3$ as in Theorem \ref{thm:counting} for $m=2k+2^{k^2}$, $\xi=1/2$, and $d_3=1/2$.   Set $\e_1=\delta_3^2$, and let $\e_2(x)$ be the function sending $x$ to $\delta_2(x^{-1})$  from Theorem \ref{thm:counting}, applied with $\delta_3, m, \xi, d_3$.   Given $t,\ell\geq 1$, and choose $N\gg t,\ell, \e_1^{-1}, \e_2(\ell)^{-1}$.  

Suppose $H=(V,E)$ with $|V|=n\geq N$, and $\calP$ is a $(t,\ell,\e_1,\e_2(\ell))$-decomposition of $V$.  Say $\calP_{vert}=\{V_i: i\in [t]\}$ and $\calP_{edge}=\{P_{ij}^{\alpha}: ij\in {[t]\choose 2}, \alpha\leq \ell\}$. Assume there exists an encoding of $U(k)$ in $(H,\calP)$.   Then there are $1\leq j_0\neq k_0\leq t$, some $\{\alpha_S\in [\ell]: S\subseteq [k]\}$, $\{i_u\in [t]: u\in [k]\}$, and $ \{(\beta_u,\gamma_u)\in [\ell]^2: u\in [k]\}$ such that for each $u\in [k]$ and $S\subseteq [k]$, each of $P_{i_uj_0}^{\gamma_u}, P_{i_uk_0}^{\beta_u}, P_{j_0k_0}^{\alpha_S}$ satisfy $\disc_2(\e_2(\ell);1/\ell)$ and if $u\in S$ then $(H_{i_uj_0k_0}^{\alpha_S,\beta_u, \gamma_u},G_{i_uj_0k_0}^{\alpha_S,\beta_u, \gamma_u})$ has $\disc_{2,3}(\e_1,\e_2(\ell))$ and density at least $1/2$ and if $u\in S$, and density less than $1/2$ if $u\notin S$.  

Given $X\subseteq [k]^2$ and $i\in [k]$, let $N_X(i)=\{j\in [k]: (i,j)\in S\}$.  For each $u\in [k]$, set $A_u=V_{i_u}$ and $B_u=V_{j_0}$.  Then for each  $X\subseteq [k]$, let $C_X=V_{k_0}$. Given $X\subseteq [k]^2$ and $u\in [k]$, define $P_{C_XB_u}=P_{j_0k_0}^{\alpha_{N_X(u)}}$ and $P_{C_XA_u}=P_{i_uk_0}^{\beta_u}$.  For $(u,v)\in [k]^2$, define $P_{A_vB_u}=P_{i_uj_0}^{\gamma_u}$.

Given $(u,v)\in [k]^2$ and $X\subseteq [k]^2$, set $G_{uvX}:=(A_v\cup B_u\cup C_X,P_{A_vC_X}\cup P_{B_uC_X}\cup P_{A_vB_u})$ (in other words, $G_{uvX}=G_{i_vj_0k_0}^{\alpha_{N_X(u)},\beta_v,\gamma_v}$).  Then define $H_{uvX}=(A_v\cup B_u\cup C_X, E\cap K_3^{(2)}(G_{vuX}))$.   By assumption, $(H_{uvX},G_{uvX})$ has $\disc_{2,3}(\e_1,\e_2(\ell))$ and density at least $1/2$ if  $v\in N_X(u)$ (i.e. if $(u,v)\in X$) and less than $1/2$ if $v\notin N_X(u)$ (i.e. $(u,v)\notin X$).   By Theorem \ref{thm:counting}, there are at least $(1/\ell)^{m\choose 2}(1/2)^{m\choose 3}(n/t)^{m}>0$ many tuples $(a_u)_{u\in [k]}(b_v)_{v\in [k]}(c_X)_{X\subseteq [k]^2}$ such that $a_ub_vc_X\in E$ if and only if $(u,v)\in X$. This finishes the proof. 
\end{proof}

Note that in Proposition \ref{prop:suffvc2}, we did not just show that if there is an encoding of $U(k)$ in $(H,\calP)$, then $H$ has $k$-$\IP_2$, but rather the stronger statement that $\trip(H)$ contains many copies of $V(k)$ (recall Definition \ref{def:vk}).  

We now prove the analogous result for $\FOP_2$, where $U(k)$ is replaced by $H(k)$ and $V(k)$ by $F(k)$ (recall Definition \ref{def:fell}). Again, we will find that if there is an encoding of $H(k)$, then $\trip(H)$ contains many copies of $F(k)$. In fact, we will give an even more detailed conclusion, as we will need these details later on in the chapter.

\begin{proposition}\label{prop:sufffop}
For all $k\geq 1$, there exist $\e_1>0$ and $\e_2:\mathbb{N}\rightarrow (0,1]$ such that for all $t,\ell\geq 1$, there is $N$ such that the following hold.  Suppose $H=(V,E)$ is a $3$-graph with $|V|\geq N$, and $\calP$ is a $(t,\ell,\e_1,\e_2(\ell))$-decomposition of $V$.  Suppose $(g,f)$ is an encoding of $H(k)$ in $(H,\calP)$, with corresponding partition $(V_{j_0}, V_{k_0}, \bigcup_{w=1}^kV_{i_w})$.  Then
\begin{enumerate}[label=\normalfont(\arabic*)]
\item there exist $(c_w)_{w\in [k]}\in \prod_{w=1}^kV_{i_w}$,  $\{a^f_u: u\in [k]\}\subseteq V_{j_0}$, and $\{b_v: v\in [k]\}\subseteq V_{k_0}$, such that $a^f_ub_vc_w\in E$ if and only if $w\leq f(u,v)$, and
\item There exist $(c_w)_{w\in [k]}\in \prod_{w=1}^kV_{i_u}$, $\{a^f_u: u\in [k]\}\subseteq V_{k_0}$, and $\{b_v: v\in [k]\}\subseteq V_{j_0}$ such that $a^f_ub_vc_w\in E$ if and only if $w\leq f(u,v)$.
\end{enumerate}
In particular, $H$ has $k$-$\FOP_2$.
\end{proposition}
\begin{proof}
Fix $k\geq 1$. Choose $\delta_3$ as in Theorem \ref{thm:counting} for $m=2k+k^{k^2}$, $\xi=1/2$, and $d_3=1/2$.   Set $\e_1=\delta_3^2$.  Let $s(x)$ be the function sending $x$ to $\delta_2(x^{-1})$ from Theorem \ref{thm:counting} applied with $\delta_3, m, \xi, d_3$.  For all $x$, define $\e_2(x)\ll s(x)$. Given $t,\ell$, choose $N\gg t,\ell, \e_1^{-1}, \e_2(\ell)^{-1}$.  

Suppose $H=(V,E)$ with $|V|=n\geq N$, and $\calP$ is a $(t,\ell,\e_1,\e_2(\ell))$-decomposition of $V$.  Say $\calP$ consists of $\calP_{vert}=\{V_i: i\in [t]\}$ and $\calP_{edge}=\{P_{ij}^{\alpha}: ij\in {[t]\choose 2}, \alpha\leq \ell\}$.  Assume there exists an encoding of $H(k)$ in $(H,\calP)$.  Then there are $1\leq j_0\neq k_0\leq t$, some $\{(\alpha_u,\beta_u,\gamma_u)\in [\ell]^3: u\leq k\}$,  and $i_1,\ldots,i_k\leq t$ such that for each $u\in [k]$, each of $P_{i_uj_0}^{\alpha_u}, P_{i_uk_0}^{\beta_u}, P_{j_0k_0}^{\gamma_u}$ satisfy $\disc_2(\e_2(\ell);1/\ell)$ and if $u\leq v$ then $(H^{\gamma_v,\alpha_u,\beta_u}_{j_0k_0i_u},G_{j_0k_0i_u}^{\gamma_v,\alpha_u,\beta_u})$ has $\disc_{2,3}(\e_1,\e_2(\ell))$ and density $\geq 1/2$, and if $u> v$ then  $(H^{\gamma_v,\alpha_u,\beta_u}_{j_0k_0i_u},G_{j_0k_0i_u}^{\gamma_v,\alpha_u,\beta_u})$ has $\disc_{2,3}(\e_1,\e_2(\ell))$ and density $<1/2$.   

Let $I$ be the set of functions $f:[k]^2\rightarrow [k]$.  For each $u\in [k]$, set $C_u=V_{i_u}$, $B_u=V_{k_0}$, and for each $f\in I$, set $A_u^f=V_{j_0}$. Then for each $u,v,w\in [k]$ and $f\in I$, define $Q_{A_u^fC_w}=P_{j_0i_w}^{\alpha_w}$, $Q_{B_vC_w}=P_{k_0,i_w}^{\beta_w}$ and $Q_{A_u^fB_v}=P_{j_0k_0}^{\gamma_{f(u,v)}}$.  For each $u,v,w\in [k]$ and $f\in I$, let 
\[G_{uvw}^f=(A^f_u \cup B_v\cup C_w, Q_{B_vC_w}\cup Q_{A_u^fC_w}\cup Q_{A_u^fB_v})\]
and
\[H_{uvw}^f=(A^f_u\cup B_v\cup C_w, E\cap K_3^{(2)}(G_{uvw}^f)).\]

We know that for each $u,v,w\in [k]$ and each $f\in I$, $(H_{uvw}^f, G_{uvw}^f)$ satisfies $\disc_{2,3}(\e_1, \e_2(\ell))$, with density at least $1/2$ if $w\leq f(u,v)$ and less than $1/2$ if $w>f(u,v)$.  By Theorem \ref{thm:counting}, there are at least $(1/\ell)^{m\choose 2}(1/2)^{m\choose 3}(n/t)^m>0$ many tuples $(a^f_u)_{u\in [k], f\in I}(b_v)_{v\in [k]}(c_w)_{w\in [k]}\in \prod_{u\in [k], f\in I}A_u^f\times \prod_{v=1}^kB_v\times \prod_{w=1}^k C_w$ such that $a_u^fb_vc_w\in E$ if and only if $w\leq f(u,v)$. Note that for each such tuple, $c_w\in V_{i_w}$ for each $w\in [k]$, $\{a_u^f: u\in [k], f\in I\}\subseteq V_{j_0}$, and $\{b_v: v\in [k]\}\subseteq V_{k_0}$.  Thus (1) holds.

To prove (2), we proceed in the same way, but make the following small adjustments in the third paragraph.  For each $v\in [k]$, set $B_v=V_{j_0}$, and for each $u\in [k]$ and $f\in I$, set $A_u^f=V_{k_0}$.  As before, define $C_w=V_{i_w}$ for each $w\in [k]$.  Then for each $u,v,w\in [k]$ and $f\in I$, define $Q_{A_u^fC_w}=P_{k_0,i_w}^{\alpha_w}$, $Q_{B_vC_w}=P_{j_0,i_w}^{\beta_w}$ and $Q_{A_u^fB_v}=P_{j_0k_0}^{\gamma_{f(u,v)}}$.  One then finishes the proof in the same way to obtain (2).
\end{proof}

We are now able to prove Theorem \ref{thm:suffcond}, which shows that if $\calH$ encodes $H(k)$ (respectively $U(k)$) for all $k$ , then $\calH$ has $\FOP_2$ (respectively $\IP_2$).

\vspace{2mm}

\begin{proofof}{Theorem \ref{thm:suffcond}}
Suppose that $\calH$ encodes $U(k)$.  Let $\e_1>0$ and $\e_2:\mathbb{N}\rightarrow (0,1]$ be as in Proposition \ref{prop:suffvc2}.  Let $T=T(\e_1,\e_2)$ and $L=(\e_1,\e_2)$ be as in the definition of $U(k)$ being encoded in $\calH$.  Let $N=N(T,L)$ be as in Proposition \ref{prop:suffvc2}.  Since $U(k)$ is encoded in $\calH$ there is $H=(V,E)\in \calH$ with $|V|\geq N$, some $t\leq T$, $\ell\leq L$, a $(t,\ell,\e_1,\e_2(\ell))$-decomposition $\calP$ of $V$, and an encoding of $U(k)$ in $(H,\calP)$.  By Proposition \ref{prop:suffvc2}, $H$ has $k$-$\IP_2$, so $\calH$ has $k$-$\IP_2$.  

Suppose now that $\calH$ encodes $H(k)$.  Let $\e_1>0$ and $\e_2:\mathbb{N}\rightarrow (0,1]$ be as in Proposition \ref{prop:sufffop}.  Let $T=T(\e_1,\e_2)$ and $L=(\e_1,\e_2)$ be as in the definition of $H(k)$ being encoded in $\calH$.  Let $N=N(T,L)$ be as in Proposition \ref{prop:suffvc2}.  Since $H(k)$ is encoded in $\calH$ there are $H=(V,E)\in \calH$ with $|V|\geq N$, some $t\leq T$, $\ell\leq L$, a $(t,\ell,\e_1,\e_2(\ell))$-decomposition $\calP$ of $V$, and an encoding of $H(k)$ in $(H,\calP)$.  By Proposition \ref{prop:sufffop}, $H$ has $k$-$\FOP_2$, so $\calH$ has $k$-$\FOP_2$.  
\end{proofof}

\vspace{2mm}
 Our next goal is to show that the converse of Theorem \ref{thm:suffcond} is also true, i.e. if $\calH$ has $\FOP_2$  (respectively $\IP_2$), then it encodes $H(k)$ (respectively $U(k)$) for all $k$.  

We will first prove this for $\FOP_2$.  This will require Lemma \ref{lem:otherway} below, which shows that we can construct certain $3$-graphs from $F(n)$, when $n$ is large.

\begin{lemma}\label{lem:otherway}
For all $m,k\geq 1$, there is $N$ such that the following holds.  Suppose $n\geq N$, $U=\{u_1,\ldots, u_m\}$, $W=\{w_1,\ldots, w_m\}$ are sets $\bigcup_{\alpha\leq k} P^{\alpha}$ is a partition of $K_2[U,W]$, and $H=(A\cup B\cup C, E)$ is a clean copy of $F(n)$.

 Then there exist disjoint sets $Z_1,\ldots, Z_k\subseteq C$, each of size $m$, along with $\{a_1,\ldots, a_m\}\subseteq A$ and $\{b_1,\ldots, b_m\}\subseteq B$ such that the following hold, where $Q^{\alpha}=\{a_ib_j: u_iw_j\in P^{\alpha}\}$.  For all $c\in Z_\beta$ and $ab\in Q^{\alpha}$, $abc\in E$ if and only if $\beta\leq \alpha$.  
\end{lemma}
\begin{proof}
Fix $m,k\geq 1$.  Choose $N\gg m,k$.  Suppose $U=\{u_1,\ldots, u_m\}$ and $W=\{w_1,\ldots, w_m\}$ are sets of size $m$, and $\bigcup_{\alpha\leq k} P^{\alpha}$ is a partition of $K_2[U,W]$.  

Assume $n\geq N$ and $H=(A\cup B\cup C, E)$ is a clean copy of $F(n)$.  We may assume that $A=\{a_i^f: i\in [n], f:[n]^2\rightarrow [n]\}$, $B=\{b_i: i\in [n]\}$, and $C=\{c_i: i\in [n]\}$ are such that $a^f_ib_jc_s\in E$ if and only if $s\leq f(i,j)$.  

Define a function $g:[n]^2\rightarrow [n]$, as follows.  For each $(s,t)\in [m]^2$, set let $g(s,t)=m\alpha$, where $\alpha$ is such that $u_sw_t\in P_{UW}^{\alpha}$.  For all  $(s,t)\in [n]^2\setminus [m]^2$, define $g(s,t)=1$.  Let $X=\{a_1^g,\ldots, a_n^g\}$ and $Y=\{b_1,\ldots, b_m\}$.  For each $1\leq \alpha\leq k$, define $Z_\alpha=\{c_j: j\in ((\alpha-1) m, \alpha m] \}$.  For each $\alpha\leq k$, let $Q^{\alpha}=\{a_i^gb_j: u_iw_j\in P_{UW}^{\alpha}\}$.  Now consider some $1\leq i,j,s\leq m$, and $a_i^gb_jc_s\in K_3[X,Y,Z_1\cup \ldots \cup Z_k]$.  By construction there are some $\alpha,\beta\in [k]$ so that $c_s\in Z_{\beta}$ and $a^g_ib_j\in Q^{\alpha}$.  By assumption on $H$, $a^g_ib_jc_s\in E$ if and only if $s\leq g(i,j)$. By definition of $g$, $g(i,j)=m\alpha$, and by definition of $Z_{\beta}$, $(\beta-1)m<s\leq \beta m$.  Thus $s\leq g(i,j)$ if and only if $\beta\leq \alpha$, so $a_i^gb_jc_s\in E$ if and only if $\beta\leq \alpha$.
\end{proof}

Note that in Lemma \ref{lem:otherway}, the partition of $K_2[U,W]$ is completely arbitrary.  This will be used in the ``(2) $\Rightarrow$ (1)'' direction of the following theorem, which says that in the context of hereditary properties, $\FOP_2$ is equivalent to encoding $H(k)$ for all $k\geq 1$.

\begin{theorem}\label{prop:fopequiv}
Suppose $\calH$ is a hereditary $3$-graph property.  Then the following are equivalent.
\begin{enumerate}[label=\normalfont(\arabic*)]
\item $H(k)$ is encoded in $\calH$ for all $k\geq 1$, 
\item $\calH$ has $k$-$\FOP_2$ for all $k\geq 1$.
\end{enumerate}
\end{theorem}
\begin{proof}

That (1) implies (2) is immediate from Proposition \ref{prop:sufffop} and Definition \ref{def:binaryencode}.  Conversely, assume $\calH$ satisfies (2), and fix $k\geq 1$.   We show $H(k)$ is encoded in $\calH$.  Fix $\e_1,\rho>0$,  $\e_2:\mathbb{N}\rightarrow (0,1]$.  Let $N,T, L\geq 1$ be as in Theorem \ref{thm:reg2} for $\e_1$, $\e_2$, $t_0=2k+k^{2k}$, and $\ell_0=1$.  Let $n\gg_{\e_1,\e_2,T,L}N$.  Fix sets $U,W,V_1,\ldots, V_k$, each of size $n$, and any partitions $K_2[U,W]=\bigcup_{\alpha\leq k}P_{UW}^{\alpha}$  such that $P_{UW}^{\alpha}$ satisfies $\disc_2(\e_2(L);1/k)$.  Such a partition exists by Lemma \ref{lem:3.8}.

Let $M\gg n$.  Since $\calH$ has $\FOP_2$, by Lemma \ref{lem:clean}, it contains a clean copy of $F(M)$, say $H=(X\sqcup Y\sqcup Z,E)$.  By Lemma \ref{lem:otherway}, there are $X'=\{x_1,\ldots, x_n\}\subseteq X$, $Y'=\{y_1,\ldots, y_n\}\subseteq Y$, and $Z_1',\ldots, Z_k'\subseteq Z$ such that if $Q_{XY}^{\alpha}=\{x_iy_j: u_iw_j\in P_{WU}^{\alpha}\}$, then  the following holds:  for all $z\in Z_\beta$ and $x_iy_j\in Q_{XY}^{\alpha}$, $zx_iy_j\in E$ if and only if $\beta\leq \alpha$. 

For each $m\in [k]$, choose partitions (guaranteed by Lemma \ref{lem:3.8}) $K_2[X',Z'_m]=\bigcup_{\alpha\leq k} Q_{XZ_m}^{\alpha}$ and $K_2[Y',Z'_m]=\bigcup_{\alpha\leq k} Q_{YZ_m}^{\alpha}$ with the property that for each $\alpha\leq k$, $Q_{XZ_m}^{\alpha}$ and $Q_{YZ_m}^{\alpha}$ satisfy $\disc_2(\e_2(k);1/k)$.  Define a decomposition $\calQ$ of $V'=X'\cup Y'\cup \bigcup_{i=1}^kZ_i'$ by setting $\calQ_{vert}=\{X',Y',Z_1',\ldots, Z_k'\}$ and $\calQ_{edge}=\{Q_{XY}^{\alpha}: \alpha\leq k\}\cup \{Q^{\alpha}_{YZ_m},Q_{XZ_m}^{\alpha}:m\in [k], \alpha\leq k\}$.  By construction, $\calQ$ is a $(2+k,k,\e_1,\e_2(k))$-decomposition $V'$.   Let $H'=H[V']$. Clearly there is an encoding of $H(k)$ in $(H',\calQ)$, namely $(f,g)$, where $f:a_i\mapsto Q_{XY}^i$ and $g: b_j\mapsto Q_{XZ_j}^1Q_{XZ_j}^1$.  We have now shown (1) holds, as desired.
\end{proof}

We will now prove the analogous result for $U(k)$ and $\IP_2$.  

\begin{theorem}\label{prop:vc2equiv}
Suppose $\calH$ is a hereditary $3$-graph property.  Then the following are equivalent.
\begin{enumerate}[label=\normalfont(\arabic*)]
\item $U(k)$ is encoded in $\calH$ for all $k\geq 1$, 
\item $\calH$ has $k$-$\IP_2$ for all $k\geq 1$.
\end{enumerate}
\end{theorem}
\begin{proof}
That (1) implies (2) follows immediately from Proposition \ref{prop:suffvc2} and Definition \ref{def:binaryencode}.  Conversely, assume $\calH$ satisfies (2), and fix $k\geq 1$.  We show $U(k)$ is encoded in $\calH$.  Fix $\e_1>$, $\e_2:\mathbb{N}\rightarrow (0,1]$, and set $\ell=2^k$ and $t=2^k+2$.  Choose $m_0$ as in Lemma \ref{lem:3.8} for $\e=\delta=\e_2(\ell)$ and $\ell$.  Suppose $n\geq m_0$, and  $B$, $C$, and $A_1,\ldots, A_{\ell}$ are sets of size $n$. By Lemma \ref{lem:3.8}, there exists a partition $K_2[B, C]=\bigcup_{S\subseteq [k]}P_{BC}^{S}$ so that each $P_{BC}^{S}$ has $\disc_2(\e_2(\ell);1/\ell)$.  Let $H=(V,E)$ be the $3$-partite $3$-graph with vertex set $A_1\cup \ldots \cup A_{\ell} \cup B\cup C$, and edge set $\bigcup_{S\subseteq [k]}\bigcup_{i\in S}\bigcup_{a\in A_i}\{abc: bc\in P_{BC}^S\}$.

By Fact \ref{fact:vc2universal2}, there is a clean copy $H'=(V',E')$ of $H$ in $\calH$.  We may assume $V'=X_1\cup \ldots \cup X_{\ell}\cup Y\cup Z$ where $Y=\{y_b: b\in B\}$, $Z=\{z_c: c\in C\}$, for each $i\in [\ell]$, $X_i=\{x_a: a\in A_i\}$, such that $x_ay_bz_c\in E'$ if and only if $abc\in E$.  Use Lemma \ref{lem:3.8} to choose partitions $K_2[X_i, Y]=\bigcup_{\alpha\leq \ell}Q_{X_iY}^{\alpha}$ and $K_2[X_i, Z]=\bigcup_{\alpha\leq \ell}P_{X_iZ}^{\alpha}$ so that each $Q_{X_iY}^{\alpha}$, $Q_{X_iZ}^{\alpha}$ satisfies $\disc_2(\e_2(\ell);1/\ell)$.  For each $S\subseteq [k]$, let $Q_{YZ}^{\alpha}=\{y_bz_c: bc\in P_{BC}^S\}$.  Let $\calQ$ be the decomposition of $V'$ with $\calQ_{vert}=\{X_1,\ldots, X_{\ell}, Y,Z\}$ and 
$$
\calQ_{edge}=\{Q_{X_iY}^{\alpha}:i\in [k],\alpha\leq \ell\}\cup \{Q_{X_iZ}^{\alpha}:i\in [k],\alpha\leq \ell\}\cup \{Q_{YZ}^S:S\subseteq [k]\}.
$$  
By construction, $\calQ$ is a $(t,\ell,\e_1,\e_2(\ell))$-decomposition of $V'$.  Clearly $U(k)$ is encoded in $(H',\calQ)$ via $f:a_i\mapsto Q_{X_iY}^1Q_{X_iZ}^1$ and $g:c_S\mapsto Q_{YZ}^S$.  This shows $U(k)$ is encoded in $\calH$, so (1) holds.
\end{proof}

\section{Removal in almost stable bipartite graphs}\label{subsec:stablerem}

 In this section, we prove a removal type result, which roughly says the following. If a bipartite graph $G=(U\cup W,E)$ contains few $d$-trees with leaves in $W$, then it is close to a graph in which $\rk(W)\leq d$ (see Section \ref{subsec:strongstable}).  We begin with a definition which will help us keep tabs on the locations of the leaves and nodes of a $d$-tree.

\begin{definition}\label{def:abd}
Suppose $d\geq 1$, $G=(V, E)$ is a bipartite graph, and $A,B\subseteq V$.  An  $(A,B,d)$-tree in $G$ is a $d$-tree in $G$ with leaves in $B$ and nodes in $A$. 
\end{definition}

We now define a notion of closeness tailored to the bipartite setting.

\begin{definition}\label{def:bipclose}
Suppose $\delta>0$ and $G=(U\cup W, E)$, $G'=(U\cup W, E')$ are bipartite graphs.  Then $G$ and $G'$ are \emph{uniformly $(U,W,\delta)$-close} if for all $u\in U$, $|N_G(u)\Delta N_{G'}(u)|\leq \delta |W|$.
\end{definition}

Definition \ref{def:bipclose} is clearly stronger than the notion of $\delta$-closeness (recall $G$ and $G'$ on the same vertex set are $\delta$-close if $|E(G)\Delta E(G')|\leq \delta |V(G)|^2$).  We will need this extra strength in our proof of Theorem \ref{thm:FOPfinite}.  We can now state our removal lemma, Theorem \ref{thm:stableremoval} below.  Roughly speaking, Theorem \ref{thm:stableremoval} says that if one has a bipartite graph $G$ with few $d$-trees, then one can find a graph $G'$ which looks very similar to $G$ and which contains no $d$-trees. 
 
\begin{theorem}\label{thm:stableremoval}
For all $d\geq 1$ and $0<\e<1/9$, there is $\delta_0(\e,d)$ such that for all $0<\delta<\delta_0$ there is $\rho=\rho(\e,d,\delta)>0$ and $N=N(\e,  d,\rho)$ such that the following holds.  

Suppose $G=(U\cup W, E)$ is a bipartite graph with $|U|, |W|\geq N$ and the number of $(U,W,d)$-trees in $G$ is at most $\rho|U|^{2^d-1}|W|^{2^d}$.  Then there is $U_0\subseteq U$ and $W_0\subseteq W$ such that $|U_0|\leq \delta |U|$, $|W_0|\leq \e |W|$, and a bipartite graph $G'=((U\setminus U_0)\cup (W\setminus W_0), E')$ and so that $G[U\setminus U_0, W\setminus W_0]$ is uniformly $(U\setminus U_0, W\setminus W_0, \delta)$-close to $G'$, and so that $G'$ contains no $(U\setminus U_0,W\setminus W_0,d)$-trees.
\end{theorem}

In the proof of Theorem \ref{thm:FOPfinite}, we will apply Theorem \ref{thm:stableremoval} to a graph $G=(U\cup W,E)$ containing few $(U,W,d)$-trees and where $|W|\ll |U|$.  This imbalance in the size of the parts is another reason we require Theorem \ref{thm:stableremoval}, and cannot use an off-the-shelf removal lemma, such as Lemma \ref{lem:indremgraph}.

We begin with a lemma which shows that if a bipartite graph $G=(U\cup W, E)$ contains few $(U,W, d)$-trees, then $W$ can be almost partitioned into almost good sets.  This can be thought of as a version of Theorem \ref{thm:goodstrong} with a weaker conclusion, but which has the advantage of more general hypotheses.  In particular, Lemma \ref{lem:closetostableinWweak} applies when there are few $(U,W, d)$-trees (rather than none, as required in Theorem \ref{thm:goodstrong}).  

\begin{lemma}\label{lem:closetostableinWweak}
For all $d\geq 1$, there is $C=C(d)>0$ such that for all $0<\e<1/9$, there is $\mu_0(\e,d)>0$, such that for all $0<\mu<\mu_0$, there are $M=M(\e,d,\mu)$ and $N=N(\e,  d,\mu)$ such that the following holds.  

Suppose $G=(U\cup W, E)$ is a bipartite graph with $|U|, |W|\geq N$ and such that the number of $(U,W,d)$-trees in $G$ is at most $\mu^C |U|^{2^d-1}|W|^{2^d}$.  Then there is $m\leq M$, $U'\subseteq U$ of size at most $\mu |U|$, and a partition $W_0\cup W_1\cup \ldots \cup W_m$ such that $|W_0|\leq \e |W|$, and such that for each $1\leq i\leq m$, $|W_i|\geq  |W|/M$ and $W_i$ is $\mu$-good with respect to all $u\in U\setminus U'$.
\end{lemma}
\begin{proof}
We proceed by induction on $d$.  Suppose first $d=1$.  In this case, set $C=4$.  Given $0<\e<1/9$, set $N=M=1$ and let $\mu_0=\e$.  

Now assume $0<\mu<\mu_0$, and $G=(U\cup W,E)$ is a bipartite graph containing at most $\mu^C |U||W|^2$ many $(U,W,1)$-trees. Note that for all $u\in U$, if $uw_1\in E$ and $uw_0\notin E$, then $(u,w_0,w_1)$ is a $(U,W,1)$-tree.  Thus, for all $u\in U$, there are at least $|N(u)\cap W||\neg N(u)\cap W|$ many $(U,W,1)$-trees in $G$ with $a_{<>}=u$.  Consequently there must be a set $U'\subseteq U$ of size at most $\mu^{C/2}|U|$, such that for all $u\in U\setminus U'$, $|N(u)\cap W||\neg N(u)\cap W|\leq \mu^{C/2}|W|^2$ (since otherwise, $G$ would contain more than $\mu^C |U||W|^2$ many $(U,W,1)$-trees).  This means that for all $u\in U\setminus U'$, we must have that $\min\{|N(u)\cap W|, |\neg N(u)\cap W|\}\leq \mu^{C/4}|W|=\mu |W|$.  This shows that $W$ itself is $\mu$-good with respect to all $u\in U\setminus U'$. This finishes the case $d=1$.

Suppose now $d>1$, and assume that for all $1\leq d'<d$, we have defined $C(d')$, $\mu_0(d',\e')$, $N(d',\e',\mu)$, and $M(d',\e',\mu)$, for all all $0<\e'<1/4$ and all $0<\mu<\mu_0(d',\e')$.  Set $C'=C(d-1)$, and define $C(d)=8C'+2^d+2$.  Fix $0<\e<1/4$. Define $\mu_0(d,\e)=\e^{2^{d+1}}\mu_0(\e,d-1)/4$.  Given $0<\mu<\mu_0$, let $N'=N(d-1,\e^5,\mu^4)$ and $M'=M(d-1,\e^5,\mu^4)$.  Finally, we define $M=M(d,\e):=\e^{-100}M'$ and $N=N(d,\e):=\e^{-100}N'M'$. 

Suppose $G=(U\cup W, E)$ is a bipartite graph with $|U|,|W|\geq N$, and assume there are at most $\mu^C |U|^{2^d-1}|W|^{2^d}$ many $(U,W,d)$-trees in $G$. We define an integer $s$ and a sequence $W_0,\ldots, W_s$ and $Z_0,\ldots, Z_s$ as follows.

\underline{Step 0:} Let $W_0=\emptyset$ and $Z_0=W$.

\underline{Step $i+1$:} Suppose we have defined $W_0,\ldots, W_i$, $Z_0,\ldots, Z_i$ so that for each $0\leq i'\leq i$, $Z_{i'}=W\setminus (\bigcup_{u=0}^{i'}W_u)$, and for all $1\leq i'\leq i$,  there are at most $\mu^{4C'}|U|^{2^{d-1}-1}|W_{i'}|^{2^{d-1}}$ many $(U, W_{i'},d-1)$-trees in $G$, and $|W_{i'}|\geq \e^2\mu|W|$.  Let 
$$
U_{i+1}=\{u\in U: \min\{|N(u)\cap Z_i|,|\neg N(u)\cap Z_i|\}\geq \mu|Z_i|\}.
$$
If either $|Z_i|\leq \e^2 |W|$ or $|U_{i+1}|\leq \mu |U|$, define $s=i$ and end the construction.  

Otherwise, for each $u\in U_{i+1}$ set  $W_{i,1}(u)=N(u)\cap Z_i$ and $W_{i,0}(u)=\neg N(u)\cap Z_i$.  By definition of $U_{i+1}$ and by assumption on the size of $Z_i$, we know that for all $u\in U_{i+1}$, both $W_{i,1}(u)$ and $W_{i,0}(u)$ have size at least $\mu |Z_i|\geq \mu \e^2|W|$. 

Define $U_{i+1}'$ to be the set of $u\in U_{i+1}$ for which the number of $(U,W_{i,1}(u),d-1)$-trees in $G$ exceeds $\mu^{4C'} |U|^{2^{d-1}-1}|W_{i,1}(u)|^{2^{d-1}}$ and the number of $(U,W_{i,0}(u),d-1)$-trees in $G$ exceeds $\mu^{4C'}|U|^{2^{d-1}-1}|W_{i,0}(u)|^{2^{d-1}}$.  We claim that $U_{i+1}\setminus U_{i+1}'\neq \emptyset$.  Indeed, suppose not.  Then we construct many copies of $T(d)$ as follows.
\begin{itemize}
\item Choose $u\in U_{i+1}'$ and set $u_{<>}=u$.  There are $|U'_{i+1}|=|U_{i+1}|> \mu |U|$ choices.
\item Choose $(u^1_{\sigma})_{\sigma\in 2^{<d-1}}(w^1_{\rho})_{\rho\in 2^{d-1}}$ a $(U,W_{i,1}(u),d-1)$-tree in $G$, and for each  $\sigma\in 2^{<d-1}$ and $\rho\in 2^{d-1}$, set $u_{1\wedge \sigma}=u_{\sigma}^1$ and $w_{1\wedge \rho}=w_{\rho}^1$.  By assumption, there more than 
$$
\mu^{4C'} |U|^{2^{d-1}-1}|W_{i,1}(u)|^{2^{d-1}}\geq \mu^{4C'}|U|^{2^{d-1}-1}(\mu\e^2|W|)^{2^{d-1}}
$$
 ways to do this, where the inequality is because $u\in U_{i+1}$.
\item Choose $(u^0_{\sigma})_{\sigma\in 2^{<d-1}}(w^0_{\rho})_{\rho\in 2^{d-1}}$ a $(U,W_{i,0}(u),d-1)$-tree in $G$, and for each $\sigma\in 2^{<d-1}$ and $\rho\in 2^{d-1}$, set $u_{0\wedge \sigma}=u_{\sigma}^1$ and $w_{0\wedge \rho}=w_{\rho}^1$.  By assumption, there are more than 
$$
\mu^{4C'} |U|^{2^{d-1}-1}|W_{i,1}(u)|^{2^{d-1}}\geq \mu^{4C'}|U|^{2^{d-1}-1}(\mu\e^2|W|)^{2^{d-1}}
$$
 ways to do this, where the inequality is because $u\in U_{i+1}$.
\end{itemize}
Every such construction produces a distinct tuple $(u_{\sigma})_{\sigma\in 2^{<d}}(w_{\rho})_{\rho\in 2^d}$ which is a $(U,W,d)$-tree in $G$.  But this means the number of $(U,W,d)$-trees in $G$ is at least 
$$
(\mu |U|)(\mu^{4C'} |U|^{2^{d-1}-1}( \mu\e^2|W|)^{2^{d-1}})^2,$$ 
which equals
$$
\e^{2^{d+1}}\mu^{8C'+2^d}|U|^{2^d-1}|W|^{2^d}>\mu^{C}|U|^{2^d-1}|W|^{2^d},
$$
where the inequality uses the fact that $\mu<\mu_0$ and the definition of $C$.  This is a contradiction, so we must have that $U_{i+1}\setminus U_{i+1}'\neq \emptyset$.  

Thus we may choose some $u\in U_{i+1}\setminus U_{i+1}'$.  By definition of $U'_{i+1}$, there is some $W_{i+1}\in \{W_{i,1}(u),W_{i,0}(u)\}$ with the property that the number of $(U,W_{i+1},d-1)$-trees in $G$ is at most $\mu^{4C'} |U|^{2^{d-1}-1}|W_{i+1}|^{2^{d-1}}$. Set $Z_{i+1}=W\setminus (\bigcup_{j=1}^{i+1}W_j)$. This ends step $i+1$.

Clearly this construction will end in some $s\leq \e^{-2}\mu^{-1}$ steps (since otherwise $|\bigcup_{i=1}^s W_i|>|W|$).  At the end, we will have a partition $W=W_1\cup \ldots\cup W_s\cup Z_s$, where $|Z_s|\leq \e^2|W|$, and where for each $1\leq i\leq  s$,  there are at most $\mu^{4C'}|U|^{2^{<d-1}}|W_i|^{2^{d-1}}$ many $(U,W_i,d-1)$-trees in $G$ and $|W_i|\geq \mu\e^2|W|\geq N'=N(\e^2,d-1,\mu^4)$.  For each $i\in [s]$, the  inductive hypothesis applied to $G[U\cup W_i]$ implies there is some $s_i\leq M'=M(\e^2, d-1,\mu^4)$, a partition $W_i=W_{i,0}\cup W_{i,1}\cup \ldots \cup W_{i,s_i}$, and $U_i\subseteq U$ such that $|U_i|\leq \mu^4 |U|$, $|W_{i,0}|\leq \e^2|W_i|$, and such that for each $j\in [s_i]$, $W_{i,j}$ is $\mu^4$-good in $G[U\cup W_i]$ with respect to every $u\in U\setminus U_i$.  Clearly this implies each $W_{i,j}$ is $\mu^4$-good with respect to every $u\in U\setminus U_i$ in $G$ as well.  Observe that each $|W_{i,j}|\geq \mu \e^2 |W|/M'\geq |W|/M$.

Set $U_0=\bigcup_{i=1}^s U_i$ and $W_0=Z_s\cup \bigcup_{i=1}^sW_{i,0}$.  Since $s\leq \e^{-2}\mu^{-1}$, $|U_0|\leq \mu^4\e^{-2}\mu^{-1}|U|< \mu^2|U|$.  Further, 
$$
|W_0|\leq \e^2|W|+ \e^2\sum_{i=1}^s |W_i|\leq 2\e^2 |W|\leq \e |W|.
$$
Let $W_1',\ldots, W_{m'}$ enumerate $ \{W_{i,j}: 1\leq i\leq s, 1\leq j\leq s_i\}$.  Note 
$$
m'=\sum_{i=1}^s s_i\leq sM'\leq \e^{-2}\mu^{-1}M'\leq M.
$$
 Thus we have a partition $W=W_0\cup W_1'\cup \ldots \cup W_{m'}$, where $m'\leq M$, $|W_0|\leq \e |W|$, and such that for each $1\leq u\leq m'$, $W_u'$ is  $\mu$-good to every $u\in U\setminus U_0$ and satisfies $|W_u'|\geq \mu |W|/M$.  This finishes the proof.
\end{proof}

The last remaining ingredient is a bipartite counting lemma. This is likely implicit in work of Alon, Fischer and Newman \cite{Alon.2007}. The authors also used Lemma \ref{lem:graphcounting1} in \cite{Terry.2021a}, where a proof appears (see Appendix C of \cite{Terry.2021a}). We note that Lemma \ref{lem:graphcounting1} assumes the pairs involved have density near $0$ or $1$, which is stronger than assuming they are merely regular (see Fact \ref{fact:homimpliesrandombinary}). 

\begin{lemma}\label{lem:graphcounting1}
For all $k\geq 1$ and $\e>0$ there are $\eta>0$ and $N$ such that the following hold.  Suppose $F=(U\cup W,R)$ is a bipartite graph with $|U|+|W|= k$.  Assume $G=(V,E)$ is a graph, and $V=\bigcup_{u\in U}V_u\cup \bigcup_{w\in W}V_w$ is a partition of $V$ into parts of size at least $N$.  Suppose that for each $uw\in R$, $|E\cap K_2[V_u,V_w]|\geq (1-\eta)|V_u||V_w|$, and for each $uw\notin R$, $|E\cap K_2[V_u,V_w]|\leq \eta |V_u||V_w|$.  Then there are at least $(1-\e)\prod_{x\in U\cup W}|V_x|$ many tuples $(v_u)_{u\in U}(v_w)_{w\in W}\in \prod_{u\in U}V_u\times \prod_{w\in W}V_w$ such that $v_uv_w\in E$ if and only if $uw\in R$.
\end{lemma}

We are now ready to prove Theorem \ref{thm:stableremoval}, the strong removal lemma for the ``trees" of Definition \ref{def:abd}.

\vspace{2mm}

\begin{proofof}{Theorem \ref{thm:stableremoval}}
Fix $d\geq 1$ and $0<\e<1/9$.  Let $\eta=\eta(2^d+2^d-1,1/2)$ and $N_2=N_2(2^d+2^d-1, 1/2)$ be from Lemma \ref{lem:graphcounting1}.  Let $C=C(d)$ and $\mu_0(\e,d)$ be from Lemma \ref{lem:closetostableinWweak}.  Set $\delta_0=\min\{\mu_0(\e,d), \eta\}$.  Given $0<\delta<\delta_0$, let $M=M(\e,\delta,d)$, $N_1=N_1(\e,\delta,d)$ be from Lemma \ref{lem:closetostableinWweak}.  Then define $\rho=\frac{1}{4}(\delta/M2^M)^C$ and let $N=2^M\mu^{-1}\max\{N_1,N_2\}$.  

Suppose $G=(U\cup W, E)$ is a bipartite graph with $|U|, |W|\geq N$ such that the number of $(U,W,d)$-trees in $G$ is at most $\rho |U|^{2^d-1}|W|^{2^d}$.  Since $\rho<\delta^C$, Lemma \ref{lem:closetostableinWweak} implies there exist $m\leq M$, $U_0\subseteq U$ of size at most $\delta |U|$, a partition $W_0\cup W_1\cup \ldots \cup W_m$ such that $|W_0|\leq \e |W|$, and such that for each $1\leq i\leq m$, $|W_i|\geq  |W|/M\geq N_2$ and $W_i$ is $\delta$-good with respect to all $u\in U\setminus U_0$.

For each $1\leq i\leq m$ and $u\in U\setminus U_0$, let $\xi(u,W_i)$ be such that $|N^{\xi(u,W_i)}(u)\cap W_i|\geq (1-\delta)|W_i|$.  Partition $U\setminus U_0$ into a minimal number of parts $U_{1}\cup \ldots \cup U_t$ so that if $u$ and $u'$ are in the same part, then $\xi(u,W_j)=\xi(u',W_j)$ for all $j\in [m]$.  Clearly $t\leq 2^m$.   Then define $I=\{ i\in [t]:|U_i|<\mu|U|/2^t\}$, and set $U_0'=U_0\cup \bigcup_{i\in I}U_i$.  Clearly $|U_0'|\leq 2\delta |U|$.  Let $U_1',\ldots, U'_{t'}$ be a reindexing of $\{U_i: i\in [t]\setminus I\}$.   For each $i\in [t']$, let $\xi(U'_i,W_j)$ be such that $\xi(u,W_j)=\xi(U'_i,W_j)$ for all $u\in U'_i$. 

Set $\mathbf{U}=\{U'_1,\ldots, U'_{t'}\}$ and $\mathbf{W}=\{W_1,\ldots, W_m\}$, and define the bipartite graph 
$$
\mathbf{R}=(\mathbf{U}\cup \mathbf{W}, \{U'_iW_j: \xi(U_i',W_j)=1\}).
$$
 Note that each $|W_i|\geq  |W|/M\geq N/M\geq N_2$, and each $|U_i'|\geq \delta |U|/2^t\geq \delta N/2^t\geq N_2$.  We claim there are no $(\mathbf{U},\mathbf{W},d)$-trees in $\mathbf{R}$. 

Suppose towards a contradiction there exists $(U_{\sigma})_{\sigma\in 2^{<d}}(W_{\lambda})_{\lambda\in 2^d}$ a $(\mathbf{U},\mathbf{W},d)$-tree in $\mathbf{R}$.  Then for each $\sigma \in 2^{<d}$ and $\lambda\in 2^d$, $|E\cap K_2[U_{\sigma},W_{\lambda}]|/K_2[U_{\sigma},W_{\lambda}]|\in [0,\delta)\cup (1-\delta,1]$.  By Lemma \ref{lem:graphcounting1}, the number of $(U,W,d)$-trees $(u_{\sigma})\in 2^{<d}(w_{\lambda})_{\lambda\in 2^d}\in \prod_{\sigma\in 2^{<d}}U_{\sigma}\times \prod_{\lambda\in 2^d}W_{\lambda}$ is at least 
$$
\frac{1}{2}\prod_{\sigma\in 2^{<d}}|U_{\sigma}|\prod_{\lambda\in 2^d}|W_{\lambda}|\geq \frac{1}{2}(\delta |U|/2^t)^{2^d-1}(|W|/M)^{2^d}> \rho |U|^{2^d-1}|W|^{2^d},
$$
a contradiction.  Thus there exist no $(\mathbf{U},\mathbf{W},d)$-trees in $\mathbf{R}$. Define
$$
G'=((U\setminus U_0')\cup (W\setminus W_0), E')\text{ where }E'=\bigcup_{U_iW_j\in E(\mathbf{R})} K_2[U_i,W_j].
$$
We claim that $G'$ contains no $(U\setminus U_0', W\setminus W_0', d)$-trees.  Suppose towards a contradiction there exists some $(u_{\sigma})_{\sigma\in 2^{<d}}(w_{\lambda})_{\lambda\in 2^d}$ which is a $(U\setminus U_0', W\setminus W_0', d)$-tree in $G'$.  By the definition of $G'$ and the fact that there are no there are no $(\mathbf{U},\mathbf{W},d)$-trees in $\mathbf{R}$, we know there must be some $\lambda_1\neq \lambda_2$ in $2^d$ and a $W_i\in \mathbf{W}$ with $w_{\lambda_1},w_{\lambda_2}\in W_i$.  Since $\lambda_1\neq \lambda_1$, there is some $\tau\in 2^{<d}$ such that $\tau\wedge 0\trianglelefteq \lambda_1$ and $\tau\wedge 1\trianglelefteq \lambda_2$, or vice versa.  Without loss of generality, assume $\tau\wedge 0\trianglelefteq \lambda_1$ and $\tau\wedge 1\trianglelefteq \lambda_2$.  Thus there is $u_{\tau}\in U\setminus U_0'$ so that $u_{\tau}w_{\lambda_1}\in E$ and $u_{\tau}w_{\lambda_2}\notin E$.  But there is some $U_j'\in \mathbf{U}$ with $u_{\tau}\in U_j'$.  But by definition, since $w_{\lambda_1},w_{\lambda_2}\in W_i$, we have that $u_{\tau}w_{\lambda_1}$ and $u_{\tau}w_{\lambda_2}$ are both in $E^{\delta(U_j',W_i)}$, a contradiction.  This finishes our verification that $G'$ contains no $(U\setminus U_0', W\setminus W_0', d)$-trees.

It is easy to see that for all $u\in U\setminus U_0'$, 
$$
|(N_{G'}(u)\Delta  N_G(u))\cap (W\setminus W_0'))|\leq \sum_{W_j\in \mathbf{W}}|(N_{G'}(u)\Delta N_G(u))\cap W_j|,
$$
which is at most
$$
\sum_{W_j\in \mathbf{W}}\delta |W_j|\leq \delta |W\setminus W_0'|.
$$
Thus $G[U\setminus U_0', W\setminus W_0]$ is uniformly $(U\setminus U_0', W\setminus W_0, \delta)$-close to $G'$.
\end{proofof}

\section{A property admits linear $\disc_{2,3}$-error if and only if it is $\NFOP_2$}\label{subsec:mainfop}
In this section we prove our main result about $\NFOP_2$ properties, which is that they are exactly the properties admitting linear $\disc_{2,3}$-error.  The harder direction is Theorem \ref{thm:FOPfinite}, which implies that an $\NFOP_2$ property admits linear $\disc_{2,3}$-error.   The proof will apply Theorem \ref{thm:stableremoval}  along with a refinement of Theorem \ref{thm:goodstrong}  to auxiliary graph structures like the ones defined in Section \ref{subsec:fopchar} (see Definition \ref{def:reducedP}).

Before we start the proof of Theorem \ref{thm:FOPfinite}, we prove two additional lemmas.  The first is a corollary of Theorem \ref{thm:goodstrong}, which contains additional equitability requirements in both its assumptions and conclusions.

 \begin{theorem}\label{thm:goodstrongequitable}
Suppose $d\geq 0$.  For all $\e>0$ and non-increasing functions $f:\mathbb{N}\rightarrow (0,1]$ with $\lim_{n\rightarrow \infty}f(n)=0$, there is $M\geq 1$ such that the following hold.  

Suppose $G=(U\cup W, E)$ is a bipartite graph, and $W=W_1\cup \ldots\cup W_m$ is a partition satisfying $|W_1|=\ldots=|W_m|$, and such that for each $i\in [m]$, $G$ contains no $(U,W_i,d)$-trees.

Then there are some $m'\leq M$, $n'\leq |W|$, and a set $\Omega \subseteq [m]$, such that $|\bigcup_{u\in \Omega}W_u|\leq \e|W|$, and such that for all $u\in [m]\setminus \Omega$, there is a partition
$$
W_u=W_{u,0}\cup W_{u,1}\cup \ldots \cup W_{u,m'},
$$
so that $|W_{u,0}|\leq \e |W_u|$, and so that for each $1\leq i\leq m'$, $W_{u,i}$ is $f(m')$-good in $G$ and $|W_{i,u}|=n'$.
\end{theorem}

The difference between the above and Theorem \ref{thm:goodstrong} is that Theorem \ref{thm:goodstrongequitable} begins with an equipartition of $W$ (rather than just a partition), and ends up with an equipartition of most of $W$ into good sets (rather than just a partition).   The proof of  Theorem \ref{thm:goodstrongequitable} consists of applying  Theorem \ref{thm:goodstrong} with appropriately chosen functions, then showing that the desired level of ``goodness" is preserved after refining things into an equipartition.  We have chosen to present Theorem \ref{thm:goodstrongequitable} as a separate result in order to isolate the important combinatorial content (the proof of Theorem \ref{thm:goodstrong}) from the rather routine computations used below to deduce Theorem \ref{thm:goodstrongequitable} from Theorem \ref{thm:goodstrong}. 

\vspace{2mm}

\begin{proofof}{Theorem \ref{thm:goodstrongequitable}}
Fix $d\geq 0$, $\e>0$, and a non-increasing function $f:\mathbb{N}\rightarrow (0,1]$ with $\lim_{n\rightarrow \infty}f(n)=0$.  Define $\psi(x)=\e^4f(\e^{-4}x)/x^2$.  Let $M_1$ be as in Theorem \ref{thm:goodstrong} for $d$, $\e^8$, and $\psi$.  Then define $M=\e^{-2}(M_1)^2$.

Suppose $G=(U\cup W, E)$ is a bipartite graph, and $W=W_1\cup\ldots\cup W_m$ is a partition so that $|W_1|=\ldots=|W_m|$, and so that $G$ contains no $(U,W_i,d)$-trees for each $i\in [m]$.  

By Theorem \ref{thm:goodstrong}, there is some $m_1\leq M_1$ and a set $\Omega \subseteq [m]$, such that $|\bigcup_{u\in \Omega}W_u|\leq \e^8 |W|$ and so that for each $u\in [m]\setminus \Omega$, there are $s_u\leq m_1$ and a partition
$$
W_u=W_{u,0}\cup W_{u,1}\cup \ldots \cup W_{u,s_u},
$$
so that $|W_{u,0}|\leq \e^8 |W_u|$, and so that for each $1\leq i\leq s_u$, $W_{u,i}$ is $\psi(m_1)$-good.  By reindexing if necessary, we may assume that for each $u\in [m]\setminus \Omega$, $|W_{u,1}|\geq \ldots \geq |W_{u,s_u}|$.  Set  $K=\lceil  \e^2 |W|/m_1 \rceil $, and let $\Omega'=\Omega\cup \{u\in [m]\setminus \Omega: |W_{u,1}|<K\}$.

For each $u\in [m]\setminus \Omega'$, let $p_u$ be maximal so that $|W_{u,p_u}|\geq K$, and for each $1\leq v\leq p_u$, choose a partition 
$$
W_{u,v}=W_{u,v,0}\cup W_{u,v,1}\cup \ldots \cup W_{u,v,s_{u,v}},
$$
with the property that $|W_{u,v,0}|<K$ and for each $1\leq w\leq s_{u,v}$, $|W_{u,v,w}|=K$.  Then set $q_u=\sum_{v\in [p_u]}s_{u,v}$, and fix an enumeration 
$$
\{W_{u,v,w}: v\in [p_u], w\in [s_{u,v}]\}=\{X_{u,1},\ldots, X_{u,q_u}\}\text{ so that }|X_{u,1}|\geq \ldots \geq |X_{u,q_u}|.
$$
Then observe that 
\begin{align*}
|W_u\setminus \Big(\bigcup_{v=1}^{q_u}X_{u,v}\Big)|\leq |W_{u,0}|+|\bigcup_{v=p_u+1}^{r_u}W_{u,v}|+K &\leq \e^2|W_u|+\frac{\e^2 |W_u|r_u}{m_1}+\frac{\e^2|W_u|}{m_1},
\end{align*}
which is at most  $4\e^2|W_u|$.
Note that for all $u\in [m]\setminus \Omega'$, $|(\bigcup_{v=1}^{q_u}X_{u,v})|=q_uK\geq (1-4\e^2)|W_u|=(1-4\e^2)|W_1|$.  Therefore, setting $q=\min\{q_u: u\in [m]\}$, we know $|W_1|\geq qK\geq (1-4\e^2)|W_1|$. For each $u\in [m]\setminus \Omega'$, define 
$$
X_{u,0}=W_{u,0}\cup \Big(\bigcup_{v=p_u+1}^{r_u}W_{u,v}\Big)\cup \Big(\bigcup_{v=1}^{p_u}W_{u,v,0}\Big)\cup \Big(\bigcup_{q<q'\leq q_u}X_{u,v}\Big).
$$
Then we note that 
\[|X_{u,0}|\leq \e^2 |W_u|+\frac{2(r_u-p_u)\e^2 |W_u|}{m_1}+\frac{2p_u\e^4 |W_u|}{m_1}+(|W_u|-qK_2)\]
which in turn is at most 
\[5\e^2|W_u|+4\e^2|W_1| = 9\e^2|W_u|< \e |W_u|,\]
where the last inequality is because $\e< 1/4$.  For all $u\in [m]\setminus \Omega'$, we now have a new partition of $W_u$, namely  
$$
W_u=X_{u,0}\cup X_{u,1}\cup \ldots \cup X_{u,q}.
$$
where $|X_{u,0}|\leq \e |W_u|$ and for each $1\leq v\leq q$, $|X_{u,v}|=K$.  Further, note that since $qK\leq |W_1|$, $q\leq \e^{-4}m_1\leq M$.

We now show that for all $u\in [m]\setminus \Omega'$ and $q'\in [q]$, $X_{u,q'}$ is $f(q)$-good.  Fix $u\in [m]\setminus \Omega'$ and $1\leq q'\leq q$.  By construction, there is some $v\in [p_u]$ and $w\in [s_{u,v}]$ so that $X_{u,q'}=W_{u,v,w}$, where $W_{u,v,w}$ has size $K$, and is contained in $W_{u,v}$, a set which is $\psi(m_1)$-good.  Thus for all $u\in U$ there is some $\delta\in \{0,1\}$ so that $|N^{\delta}(u)\cap W_{u,v}|\leq \psi(m_1)|W_{u,v}|$.  By definition of $\psi$, we have that
\[|N^{\delta}(u)\cap W_{u,v,w}| \leq \psi(m_1)|W_{u,v}|\leq \frac{\psi(m_1)|W_{u}|}{K}|W_{u,v,w}|,\]
which, by definition of $K$, is at most
\[\psi(m_1)\e^{-4}m_1|W_{u,v,w}|\leq  f(q)|W_{u,v,w}|.\]
This shows $X_{u,q'}$ is $f(q)$-good, which finishes the proof.
\end{proofof}

\vspace{2mm}

Our second lemma, Lemma \ref{lem:divis} below, ensures that certain divisibility requirements hold in regular decompositions.  It is proved by refining sufficiently regular decompositions, which are guaranteed to exist by Theorem \ref{thm:disc3}. 

\begin{lemma}\label{lem:divis}
For all $\e_1>0$, $\e_2:\mathbb{N}\rightarrow (0,1]$, and $k\geq 1$, there are $\mu_1>0$ and $\mu_2:\mathbb{N}\rightarrow (0,1]$, such that for all $L,T\geq 1$, there is $N\geq 1$ such that the following holds.

Suppose $H=(V,E)$ is a $3$-graph and $|V|=n\geq N$.  Suppose $1\leq t\leq T$, $1\leq \ell\leq L$, and $\calP$ is a $(\mu_1,\mu_2(\ell), t,\ell)$-decomposition of $V$ which is $(\mu_1,\mu_2(\ell),t,\ell)$-regular and $\mu_1$-homogeneous with respect to $H$.  Then there exists an $(\e_1,\e_2(\ell'), t,\ell')$-decomposition $\calP'$ of $V$ such that $\ell'=k\ell$, and such that $\calP'$ is $(\e_1,\e_2(\ell'),t,\ell')$-regular and $\e_1$-homogeneous with respect to $H$. 
\end{lemma}
\begin{proof}
Choose $\mu'_1,\mu'_2$ as in Lemma \ref{lem:intersecting} for $\e_1$ and $\e_2$.  Define $\mu_2:\mathbb{N}\rightarrow (0,1]$ by choosing $\mu_2(x)\ll_k \mu_2'(x)$ for all $x$.  Given $T, L\geq 1$, choose $N\gg N_1, T, L, \e_2(L)^{-1}$. 

Now suppose $H=(V,E)$ is a $3$-graph with $|V|\geq N$.  Suppose $1\leq t\leq T$, $1\leq \ell\leq L$, and $\calP$ is an $(\mu'_1,\mu_2(\ell), t,\ell)$ decomposition of $V$ which is $(\mu'_1,\mu_2(\ell),t,\ell)$-regular and $\mu'_1$-homogeneous with respect to $H$.  For each $P_{ij}^{\alpha}\in \calP_{edge}$ satisfying $\disc_2(\mu_2(\ell); 1/\ell)$, apply Lemma \ref{lem:3.8} with $\e=\delta=\mu_2(\ell)$ to obtain a partition $P_{ij}^{\alpha}=P_{ij}^{\alpha,1}\cup \ldots \cup P_{ij}^{\alpha, \ell k}$, so that each $P_{ij}^{\alpha,u}$ satisfies $\disc_2(\mu_2(\ell))$ and has density $(1\pm 2\mu_2(\ell))\frac{1}{k\ell}$.  By our choice of $\mu_2$, each such $P_{ij}^{\alpha,\beta}$ satisfies $\disc_2(\mu_2'(k\ell);1/k\ell)$.  Let $\calP'$ have $\calP'_{vert}=\calP_{vert}$ and 
\begin{align*}
\calP'_{edge}=\{P_{ij}^{\alpha,\beta}&: P_{ij}^{\alpha}\text{ satisfies }\disc_2(\mu_2(\ell);1/\ell), 1\leq \beta\leq k\}\\
&\cup \{P_{ij}^{\alpha}: P_{ij}^{\alpha}\text{ fails }\disc_2(\mu_2(\ell);1/\ell)\}.
\end{align*}
Then $\calP'$ is a $(\mu'_1,\mu'_2(\ell'), t,\ell')$, decomposition of $V$, where $t'=t$ and $\ell'=\ell k$. It is straightforward to see that $\calP'$ is a $(0,0)$-approximate refinement of $\calP$ (see Definition \ref{def:refinement}).  By Lemma \ref{lem:intersecting}, $\calP'$ is a $(t,\ell',\e_1,\e_2(\ell'))$-decomposition of $V$ which is $(\e_1,\e_2(\ell'))$-regular and $\e_1$-homogeneous with respect to $H$.
\end{proof}

We are now ready to prove Theorem \ref{thm:FOPfinite}, which says that if a hereditary $3$-graph property $\calH$ is $\NFOP_2$, then $\calH$ admits linear $\disc_{2,3}$-error.  At a high level, the idea of the proof is to define certain auxiliary graphs which are ``close" to being stable, then apply a form of stable \emph{graph} regularity to the auxiliary graphs, and finally show that  these partitions can be used to build an efficient decomposition of the $3$-graph.  Before beginning the proof, we give a detailed outline of the steps below.  The reader may wish to refer back to this outline as they read the full proof.  

Suppose $\calH$ is hereditary $3$-graph property which is $\NFOP_2$, meaning there is some $k\geq 1$ so that every $H\in \calH$ does not have $k$-$\FOP_2$.  Given a sufficiently large $3$-graph $H$ from $\calH$, our goal is to find a $\disc_{2,3}$-regular decomposition for $H$ with linear error, which we will do using the following steps. 
\begin{enumerate}
\item Because $H$ does not have $k$-$\FOP_2$, it also has $\VC_2$-dimension less than $k$ (Fact \ref{fact:vc2universal}).  Since $3$-graphs of small $\VC_2$-dimension admit homogeneous $\disc_{2,3}$-decompositions (Theorem \ref{thm:vc2finite}), we can find a $\disc_{2,3}$-homogeneous decomposition $\calP_1$ for $H$.
\item Refine $\calP_1$ to obtain a decomposition $\calP$, which is still homogeneous with respect to $H$, and which satisfies additional equitability conditions.  
\item Because $\calP$ is homogeneous, there exists a partition $\Gamma_0\cup \Gamma_1\cup \Gamma_{err}=\triads(\calP)$, where $\Gamma_0/\Gamma_1$ correspond to the triads of density close to $0/1$ respectively, and $\Gamma_{err}$ corresponds to the irregular triads.
\item Define $\Sigma_0$ to be the set of triples $V_iV_jV_s\in {\calP_{vert}\choose 3}$ with non-trivial intersection with $\Gamma_{err}$.  The set $\Sigma_0$ will be a small proportion of ${\calP_{vert}\choose 3}$ and will play a crucial role in defining the desired ``linear" error set.  It is worth noting here that there will be nothing we can say about $\Sigma_0$ other than the fact it is small. As such, it can be seen as the limiting factor in this proof, and the reason we produce a ``linear" error set (as opposed to something more structured).
\item Define $\Psi$ to be the set of pairs $V_iV_j$ for which a non-trivial number of vertex parts $V_s$ satisfy $V_iV_jV_s\in \Sigma_0$. We think of $\Psi$ as a collection of ``bad pairs."
\item For each pair $V_jV_k$ from ${\calP_{vert}\choose 2}\setminus \Psi$, define a bipartite, edge-colored auxiliary \emph{graph} $H_{jk}=(W_{jk}\cup U_{jk}, E_{jk}^1,E_{jk}^0,E_{jk}^2)$ (basically what appears in Definition \ref{def:reducedP}).  In particular, $W_{jk}$ will be the set $\{P_{jk}^{\alpha}: \alpha\in [\ell]\}$.
\item By Proposition \ref{prop:sufffop}, for each $V_jV_k\in {\calP_{vert}\choose 2}\setminus \Psi$, $H_{jk}$ contains no copy of the $k$-order property using the edge colors $E_{jk}^1,E_{jk}^0$.  
\item Consequently, by Theorem \ref{thm:stableremoval}, for each $V_jV_k\in {\calP_{vert}\choose 2}\setminus \Psi$, $H_{jk}$  is ``close" to a bipartite graph  $H_{jk}^*$ which contain no trees (in the sense of Definition \ref{def:abd}). 
\item We then union the auxiliary graphs of the form $H_{jk}^*$ for $V_jV_k\in {\calP_{vert}\choose 2}\setminus \Psi$ to obtain an auxiliary graph $G^*$ satisfying the hypotheses of Theorem \ref{thm:goodstrong}. 
\item Apply Theorem \ref{thm:goodstrong} to $G^*$.  After some unpacking, one see this yields a slightly larger (but still small) collection of ``bad pairs"  $\Psi'\supseteq \Psi$ so that the following holds. For all $V_jV_k\in {\calP_{vert}\choose 2}\setminus \Psi'$, there is a partition $R_{jk}^0\cup \ldots \cup R_{jk}^{m_2}$ of $W_{jk}$, so that each of $R_{jk}^1,\ldots, R_{jk}^{m_2}$ are good sets in $H_{jk}$.
\item For each $V_jV_k\in {\calP_{vert}\choose 2}\setminus \Psi'$, define a partition $K_2[V_j,V_k]=\mathbf{R}_{jk}^0\cup \ldots \cup \mathbf{R}_{jk}^{m_2}$ by letting $\mathbf{R}_{jk}^{\alpha}=\bigcup_{P\in R_{ij}^{\alpha}}P$ (recall each $R_{jk}^{\alpha}\subseteq \{P_{jk}^{\beta}: \beta\in [\ell]\}$).  We note $m_2$ will be significantly smaller than $\ell$.
\item We now define two additional small error sets, $\Sigma_1,\Sigma_3\subseteq {\calP_{vert}\choose 3}$. First, $\Sigma_1$ contains $\Sigma_0$ in addition to all triples $V_iV_jV_k$ where $\{V_iV_j,V_iV_k,V_jV_k\}\cap \Psi'\neq \emptyset$.  The definition of $\Sigma_3$ is more complicated, but roughly speaking, it contains all triples $V_iV_jV_k$ which have large interaction with certain error sets from Step (8).  We then define our final error set to be
$$
\Sigma=\Sigma_1\cup \Sigma_3.
$$
Since $\Sigma_1$ and $\Sigma_3$ are small, so is $\Sigma$.  
\item We then show that for all $V_iV_jV_k\in {\calP_{vert}\choose 3}\setminus \Sigma$ and all choices of $1\leq \alpha,\beta,\gamma\leq m_2$, the triad
$$
(V_i\cup V_j\cup V_k, \mathbf{R}_{ij}^{\alpha}\cup \mathbf{R}_{jk}^{\beta}\cup \mathbf{R}_{ik}^{\gamma}),
$$
is homogeneous.
\item Finally, for each $V_jV_k\in {\calP_{vert}\choose 2}\setminus \Psi'$, we use standard arguments to redistribute the ``error" $\mathbf{R}_{jk}^0$. This yields a slightly different partition $K_2[V_j,V_k]={\bf P}_{jk}^1\cup \ldots \cup {\bf P}_{jk}^m$, for some $m$ which is still significantly smaller than the original $\ell$.  For $V_jV_k\in \Psi'$, we define a somewhat arbitrary partition $K_2[V_j,V_k]=\bigcup_{\alpha\leq m}\mathbf{P}_{jk}^{\alpha}$.  We then define  $\calQ$ to be the decomposition with $\calQ_{vert}=\calP_{vert}$ and $\calQ_{edge}=\{\mathbf{P}^{\alpha}_{jk}\in {\calP_{vert}\choose 2}: \alpha\in [m]\}$.  
\item We then verify that certain triads of $\calQ$ inherit homogeneity  from step (13). In particular, we show that for all $V_iV_jV_k\in {\calQ_{vert}\choose 3}\setminus \Sigma$ and all $1\leq \alpha,\beta,\gamma\leq m$, the triad 
$$
(V_i\cup V_j\cup V_k, \mathbf{P}_{ij}^{\alpha}\cup \mathbf{P}_{jk}^{\beta}\cup \mathbf{P}_{ik}^{\gamma}),
$$
 is homogeneous. This will finish the proof.
\end{enumerate}

\vspace{2mm}
\begin{proofof}{Theorem \ref{thm:FOPfinite}}
\label{proof:FOPfinite} Fix $k\geq 1$, $\e_1>0$, $\e_2:\mathbb{N}\rightarrow (0,1]$, and $\ell_0,t_0\geq 1$.  We now fix a series of parameters.  Choose $0<\e\ll \e_1$, and let $d=d(k)$ as in part (1) in Theorem \ref{thm:treerank}. 

Let $\e^*_1>0$ and $\e^*_2:\mathbb{N}\rightarrow (0,1]$ be as in Proposition \ref{prop:sufffop} for $k$.   Define a function $\psi:\mathbb{N}\rightarrow (0,1]$ by setting $\psi(x)=\e^{100}/x^3$ for each $x\geq 1$.  Then let $M_1=M_1(d,\e^{100},f)$ be as in Theorem \ref{thm:goodstrongequitable}, where $f:\mathbb{N}\rightarrow (0,1]$ is defined by $f(x)=\e^{100} \psi(2\e^{-100}x(x+1))/(3x(x+1))$.  Let $M\gg \e^{-100}M_1$, and choose $\e_0\ll \e/(M!)^3$.  Let $\delta_0=\delta_0(d,\e_0)$ and $N_1=N_1(d,\e_0)$ be as in Theorem \ref{thm:stableremoval}.  Then choose $\e_1'$ and $\delta$ such that 
$$
\e_1'\ll \delta\ll f((M!)^3)/(M!)^3, \e_1^*, \e_1, \delta_0, \e^5, d^{-1}. 
$$
Define $\e'_2:\mathbb{N}\rightarrow (0,1]$ by choosing $\e'_2(x)\ll \e_1'\e_2(x)/xM!$.  Set $t_0'=t_0$ and $\ell_0'=\ell_0N_1\delta^{-1}$.   

Now choose $\e_1''=\e_1''(\e_1',\e_2')$ and $\e_2''=\e_2''(\e_1',\e_2')$ as in Lemma \ref{lem:disc2}.  Let $\mu_1\ll \e_1'''\ll \e_1''$, and define $\mu_2,\e'''_2:\mathbb{N}\rightarrow (0,1]$ by choosing $\mu_2(x)\ll \e_2'''(x)\ll \e_2''(x)$ for all $x$.  Then let $\mu'_1$ and $\mu'_2$ be as in Lemma \ref{lem:divis} for $\mu_1, \mu_2$ and $(M!)^3$.  Finally, let $T$, $L$, and $N_2$ be as in Theorem \ref{thm:vc2finite}  for $t_0'$, $\ell_0'$, $\mu'_1$ and $\mu'_2$, and let $N\gg N_1N_2,T, M!, (\mu'_1)^{-1}, (\mu'_2(M!L))^{-1}, \delta^{-1}$.  

Suppose $H=(V,E)$ is a $3$-graph with $|V|\geq N$ and which does not have $k$-$\FOP_2$.  Since $H$ has $\VC_2$-dimension less than $k$ (see Fact \ref{fact:vc2universal}), Theorem \ref{thm:vc2finite} implies there exist $\ell'_0\leq \ell_1\leq L$, $t_0'\leq t\leq T$, and $\calP_1$ a  $(t,\ell_1,\mu'_1,\mu'_2(\ell_1))$-decomposition of $V$ which is $(\mu'_1,\mu'_2(\ell_1))$-regular and $\mu'_1$-homogeneous with respect to $H$.  By Lemma \ref{lem:divis}, there is a $(t,\ell,\mu_1,\mu_2(\ell))$-decomposition $\calP_2$ of $V$, where $\ell=(M!)^3\ell_1$, and so that $\calP_2$ is  $(\mu_1,\mu_2(\ell))$-regular and $\mu_1$-homogeneous with respect to $H$.  By applying Lemma \ref{lem:disc2} to $\calP_2$, we obtain a $(t,\ell,\e_1''',\e_2'''(\ell))$-decomposition $\calP$ of $V$,  so that $\calP$ is  $(\e'''_1,\e'''_2(\ell))$-regular and $\e_1'''$-homogeneous with respect to $H$, and so that $\calP$ has no $\disc_2$-irregular triads.  Say $\calP_{vert}=\{V_1,\ldots, V_t\}$ and $\calP_{edge}=\{P_{ij}^\alpha: ij\in {[t]\choose 2}, \alpha \in [\ell]\}$.  

Given $ijs\in {[t]\choose 3}$ and $1\leq \alpha,\beta,\gamma\leq \ell$, set $G_{ijs}^{\alpha,\beta,\gamma}=(V_i\cup V_j\cup V_s, P_{ij}^\alpha\cup P_{js}^\beta\cup P_{is}^\gamma)$, $H_{ijs}^{\alpha,\beta,\gamma}=(V_i\cup V_j\cup V_s, E\cap K_3^{(2)}(G_{ijs}^{\alpha,\beta,\gamma}))$ and let $d_{ijs}^{\alpha,\beta,\gamma}$ be such that $| E\cap K_3^{(2)}(G_{ijs}^{\alpha,\beta,\gamma})|=d_{ijs}^{\alpha,\beta,\gamma}|K_3^{(2)}(G_{ijs}^{\alpha,\beta,\gamma})|$.  We will use throughout that since $\e_2'''(x)\ll \e_2''(x)$, Corollary \ref{cor:counting} implies that for all $ijs\in {[t]\choose 3}$ and $\alpha,\beta,\gamma\in [\ell]$, 
\begin{align}\label{align:foptri}
|K_3^{(2)}(G_{ijs}^{\alpha,\beta,\gamma})|=(1\pm \e''_2(\ell))(n/\ell t)^3.
\end{align}
We will use $\calP$ to construct a different decomposition, of $V$, which will be called $\calQ$, and which will have the same vertex partition as $\calP$, but different edge partitions. In other words, we will have $\calQ_{vert}=\calP_{vert}$ but $\calQ_{edge}\neq \calP_{edge}$.

We begin defining several sets which will help us control the location of the error triads of $\calP$.  First, set 
\begin{align*}
\Gamma_1=\{G_{ijs}^{\alpha,\beta,\gamma} \in \triads(\calP):(H_{ijs}^{\alpha,\beta,\gamma},G_{ijs}^{\alpha,\beta,\gamma}) \text{ has }&\disc_{2,3}(\e_1'',\e''_2(\ell)) \\&\text{ and } d_{ijs}^{\alpha,\beta,\gamma}\geq 1-\e_1''\},\\
\end{align*}
\begin{align*}
\Gamma_0=\{G_{ijs}^{\alpha,\beta,\gamma} \in \triads(\calP): (H_{ijs}^{\alpha,\beta,\gamma},G_{ijs}^{\alpha,\beta,\gamma}) \text{ has }&\disc_{2,3}(\e_1'',\e''_2(\ell))\\& \text{ and } d_{ijs}^{\alpha,\beta,\gamma}\leq \e_1''\},
\end{align*}
and 
\[\Gamma_{err}=\triads(\calP)\setminus (\Gamma_0\cup \Gamma_1).\]

By assumption, $\triads(\calP)=\Gamma_{err}\sqcup \Gamma_1\sqcup \Gamma_0$, and at most $\e_1'''n^3$ triples $xyz\in {V\choose 3}$ are in an element of $\Gamma_{err}$, i.e.
\begin{align}\label{k32}
\Big|\bigcup_{G\in \Gamma_{err}}K_3^{(2)}(G)\Big|\leq \e_1'''n^3.
\end{align}
 By (\ref{align:foptri}), this implies $|\Gamma_{err}|\leq \e'''_1 n^3/((n^3/t^3\ell^3)(1-\e_2''(\ell))\leq \e_1' t^3\ell^3$.  We now define the set of triples from $\calP_{vert}$ which have non-trivial intersection with $\Gamma_{err}$.
\begin{align*}
\Sigma_0&=\{V_iV_jV_s:ijs\in {[t]\choose 3} \text{ and }|\{(\alpha,\beta,\gamma)\in [\ell]^3: G_{ijk}^{\alpha,\beta,\gamma}\in \Gamma_{err}\}|\geq 2(\e_1')^{1/2} \ell^3 \}.
\end{align*}
From (\ref{align:foptri}) and  (\ref{k32}), we have that $|\Sigma_0|\leq (\e_1')^{1/2} t^3$.  We next define the set of pairs $V_iV_j$ contained in a non-trivial number of triples from $\Sigma_0$.
\begin{align*}
\Psi&=\{V_iV_j: |\{s\in [t]: V_iV_jV_s\in \Sigma_0\}| \geq (\e_1')^{1/4} t \}.
\end{align*}
Since $|\Sigma_0|\leq 2(\e_1')^{1/2}t^3$, we have that $|\Psi|\leq (\e_1')^{1/4} t^2$.  We now define, for each $V_iV_j\notin \Psi$, an auxiliary edge-colored graph $(U_{ij}\cup W_{ij},E^1_{ij},E_{ij}^0,E_{ij}^2)$, where 
\begin{align*}
W_{ij}&=\{P_{ij}^\alpha:\alpha\leq \ell \}\text{ and }U_{ij}=\{P_{is}^\beta P_{js}^\gamma: s\in [t]\setminus\{i,j\}, \beta,\gamma \in [\ell]\},\\
E^1_{ij}&=\{P_{ij}^\alpha (P_{is}^\beta P_{js}^\gamma)\in K_2[W_{ij},U_{ij}]: G_{ijk}^{\alpha,\beta,\gamma}\in \Gamma_1\text{ and }V_iV_jV_s\notin \Sigma_0\},\\
E^0_{ij}&=\{P_{ij}^\alpha (P_{is}^\beta P_{js}^\gamma)\in K_2[W_{ij},U_{ij}]: G_{ijk}^{\alpha,\beta,\gamma}\in \Gamma_0\text{ and }V_iV_jV_s\notin \Sigma_0\},\\
E^2_{ij}&=K_2[W_{ij},U_{ij}]\setminus (E_{ij}^1\cup E_{ij}^0). 
\end{align*}
Define $H_{ij}=(U_{ij}\cup W_{ij},E_{ij}^1)$.  Note $H_{ij}$ is a bipartite graph with parts $W_{ij},U_{ij}$, where $|W_{ij}|=\ell$ and $|U_{ij}|=(t-2)\ell^2$.   Our next step is to prove Claim \ref{cl:stable} below, which shows that for all $V_iV_j\notin \Psi$, there is a large subset $U^*_{ij}\cup W_{ij}^*$ of $V(H_{ij})$, on which $H_{ij}$ is uniformly $\delta$-close to some $H_{ij}^*$ which contains no $d$-trees with leaves in $W^*_{ij}$.    The idea behind the proof is as follows. First, by Proposition \ref{prop:sufffop}, we know there are no encodings of $H(k)$ in $(H,\calP)$.  Using Theorem \ref{thm:treerank}, we will show that this implies that any $(U_{ij},W_{ij},d)$-tree in $H_{ij}$ must involve some triad $G_{ijs}^{\alpha,\beta,\gamma}\in\Gamma_{err}$ where $V_iV_jV_s\notin \Sigma_0$.  If $V_iV_j\notin \Psi$, then there are few such triads, so there will be few such $(U_{ij},W_{ij},d)$-trees in $H_{ij}$.  We then apply Theorem \ref{thm:stableremoval} to obtain the desired $H_{ij}^*$.

\begin{claim}\label{cl:stable}
Suppose $V_iV_j\notin \Psi$.  Then there are $U_{ij}^{bad}\subseteq U_{ij}$ and $W_{ij}^{bad}\subseteq W_{ij}$ such that the following hold, where $U_{ij}^*=U_{ij}\setminus U_{ij}^{bad}$ and $W_{ij}^*=W_{ij}\setminus W_{ij}^{bad}$.
\begin{itemize}
\item $|W_{ij}^{bad}|\leq \e_0|W_{ij}|$ and $|U^{bad}_{ij}|\leq \delta |U_{ij}|$,
\item $H_{ij}[U_{ij}^*,W^*_{ij}]$ is uniformly $(U_{ij}^*,W_{ij}^*,\delta)$-close to a bipartite graph $H_{ij}^*=(U_{ij}^*\cup W_{ij}^*, E_{ij}^*)$ containing no $(U_{ij}^*,W^*_{ij},d)$-trees. 
\end{itemize} 
\end{claim}
\begin{proofof}{Claim \ref{cl:stable}}
Fix $V_iV_j\notin \Psi$.  We first show that the number of $(W_{ij},U_{ij},d)$-trees in $H_{ij}$ is at most $(\e_1')^{1/4} |U_{ij}|^{2^{<d}}|W_{ij}|^{2^{d}}$. Define $\Omega_{ij}=\{G_{ijs}^{\alpha,\beta,\gamma}\in \Gamma_{err}:V_iV_jV_s\notin \Sigma_0\}$.  Note that if $V_s$ satisfies $V_iV_jV_s\notin \Sigma_0$, then by definition of $\Sigma_0$,
$|\{(\alpha,\beta,\gamma)\in [\ell]^3: G_{ijs}^{\alpha,\beta,\gamma}\in \Gamma_{err}\}|\leq 2(\e_1')^{1/2} \ell^3$.  Therefore, since $V_iV_j\notin \Psi$,
\begin{align*}
|\Omega_{ij}|&\leq (\e_1')^{1/4}t\ell^3+2(\e_1')^{1/2}t\ell^3\leq 3(\e_1')^{1/4}t\ell^3.
\end{align*}  
We now make some observations about any $(U_{ij},W_{ij},d)$-tree in $H_{ij}$.  Suppose 
$$
(P_{is_{\sigma}}^{\alpha_{\sigma}}P_{js_{\sigma}}^{\beta_{\sigma}})_{\sigma\in 2^{<d}} (P_{ij}^{\gamma_{\rho}})_{\rho\in 2^{d}} 
$$
is a $(U_{ij},W_{ij},d)$-tree in $H_{ij}$.  Note that for all $\sigma\in 2^{<d}$, there is $\rho\in 2^{d}$ so that $\sigma\wedge 1\trianglelefteq \rho$, and for this $\rho$, $(P_{is_{\sigma}}^{\alpha_{\sigma}}P_{js_{\sigma}}^{\beta_{\sigma}})P_{ij}^{\gamma_{\rho}}\in E^1_{ij}$.  By definition of $E^1_{ij}$, this implies $V_iV_jV_{s_{\sigma}}\notin \Sigma_0$.  Thus for all $\sigma\in 2^{<d}$, $V_iV_jV_{s_{\sigma}}\notin \Sigma_0$.

By part (1) of Theorem \ref{thm:treerank}, there are  $\{P_{ij}^{\gamma'_u}:u\leq k\}\subseteq \{P_{ij}^{\gamma_{\rho}}: \rho\in 2^{d+1}\}$ and $\{P_{is'_u}^{\alpha'_u}P_{js_u'}^{\beta_u'}\}\subseteq \{P_{is_{\sigma}}^{\alpha_{\sigma}}P_{js_{\sigma}}^{\beta_{\sigma}}: \sigma\in 2^{<d}\}$ so that $P_{is'_u}^{\alpha'_u}P_{js_u'}^{\beta_u'}P_{ij}^{\gamma'_v}\in E^1_{ij}$ if and only if $u\leq v$.  Translating, this means that for all $u\leq v$, $G_{ijs_u'}^{\alpha_u'\beta_u'\gamma_v'}\in \Gamma_1$ and for all $v<u$, $G_{ijs_u'}^{\alpha_u'\beta_u'\gamma_v'}\notin \Gamma_1$ (i.e. $G_{ijs_u'}^{\alpha_u'\beta_u'\gamma_v'}\in \Gamma_{err}\cup \Gamma_0$).  By above, for all $u\in [k]$, $V_iV_jV_{s_u'}\notin \Sigma_0$, so for all $v<u$, we have that $G_{ijs_u'}^{\alpha_u'\beta_u'\gamma_v'}\in \Omega_{ij}\cup \Gamma_0$.  If it were true that $G_{ijs_u'}^{\alpha_u'\beta_u'\gamma_v'}\in  \Gamma_0$ for all $v<u$, then we would have an encoding of $H(k)$ in $(H,\calP)$ via $f:a_u\mapsto P_{ij}^{\gamma'_u}$ and $g:b_v\mapsto P_{is'_v}^{\alpha'_v}P_{js'_v}^{\beta'_v}$.  By Proposition \ref{prop:sufffop}, this would imply $H$ has $k$-$\FOP_2$, a contradiction.  Thus there is some $v<u$ so that $G_{ijs_u'}^{\alpha_u'\beta_u'\gamma_v'}\in  \Omega_{ij}$.

The above argument shows that for every $(U_{ij},W_{ij},d)$-tree in $H_{ij}$ of the form $(P_{is_{\sigma}}^{\alpha_{\sigma}}P_{js_{\sigma}}^{\beta_{\sigma}})_{\sigma\in 2^{<d}} (P_{ij}^{\gamma_{\rho}})_{\rho\in 2^{d}}$, there is some  $\sigma\in 2^{<d}$ and $\rho\in 2^{d}$ such that $P_{ijs_{\sigma}}^{\alpha_{\sigma},\beta_{\sigma},\gamma_{\rho}}\in \Omega_{ij}$.  Thus, any $(U_{ij},W_{ij},d)$-tree in $H_{ij}$ will be constructed in the following process.
\begin{enumerate}
\item Choose some $\sigma\in 2^{<d}$ and $\rho\in 2^{d}$. There are $(2^d-1)2^d$ choices.
\item Choose an element $G_{ijs}^{\alpha,\beta,\gamma}\in \Omega_{ij}$, and set $P_{is_{\sigma}}^{\alpha_{\sigma}}P_{js_{\sigma}}^{\beta_{\sigma}}=P_{is}^{\alpha}P_{js}^{\beta}$ and $P_{ij}^{\gamma_{\rho}}=P_{ij}^{\gamma}$.  There are at most $|\Omega_{ij}|\leq 2(\e_1')^{1/2}t \ell^3$ ways to do this.

\item For all $\sigma'\in 2^{<d}\setminus \{\sigma\}$ and $\rho'\in 2^{d}\setminus \{\rho\}$, choose any $P_{is_{\sigma'}}^{\alpha_{\sigma'}}P_{js_{\sigma'}}^{\beta_{\sigma'}}\in U_{ij}$ and any $P_{ij}^{\gamma_{\rho'}}\in W_{ij}$.  Then number of choices is at most 
$$
|U_{ij}|^{2^{d}-2}|W_{ij}|^{2^{d}-1}=((t-2)\ell^2)^{2^{d}-2}\ell^{2^{d}-1}.
$$
\end{enumerate}
This yields that the number of  $(U_{ij},W_{ij},d)$-trees in $H_{ij}$ is at most 
$$
(2^d-1)2^d2(\e_1')^{1/2} t\ell^3 ((t-2)\ell^2)^{2^{d}-2}\ell^{2^{d}-1}<(\e_1')^{1/4}|U_{ij}|^{2^{d}-1}|W_{ij}|^{2^{d}},
$$
where the inequality is due to our choice of $\e_1'$.  Now, since $\e_1'\ll \delta\ll \delta_0$, Theorem \ref{thm:stableremoval} implies that there is some $W_{ij}^{bad}\subseteq W_{ij}$ of size at most $\e_0|W_{ij}|$, and $U_{ij}^{bad}\subseteq U_{ij}$ of size at most $\delta |U_{ij}|$ such that, setting $U^*_{ij}=U_{ij}\setminus U_{ij}^{bad}$ and $W_{ij}^*=W_{ij}\setminus W_{ij}^{bad}$, $H_{ij}[U^*_{ij}, W^*_{ij}]$ is uniformly $(U_{ij}^*,W_{ij}^*,\delta)$-close to a bipartite graph $H_{ij}^*=(U_{ij}^*\cup W^*, E_{ij}^*)$ containing no $(U_{ij}^*,W_{ij}^*,d)$-trees.  This finishes the proof of Claim \ref{cl:stable}.
\end{proofof}

\vspace{2mm}
Our next goal is to apply Theorem \ref{thm:goodstrongequitable}.  Set $X=\min\{|W_{ij}^*|: V_iV_j\notin \Psi\}$.  By Claim \ref{cl:stable}, $X\geq (1-\e_0)\ell$.  For each $V_iV_j\notin \Psi$, let $W_{ij}^{**}$ be any $X$-element subset of $W_{ij}^*$. Note that $H_{ij}^*[U^*_{ij},W_{ij}^{**}]$ contains no $(U^*_{ij},W_{ij}^{**}, d)$-trees.  Set $W=\bigcup_{V_iV_j\notin \Psi}W_{ij}^{**}$ and $U=\bigcup_{V_iV_j\in \Psi}U^{*}_{ij}$, and define 
$$
G^*=(U\cup W,\bigcup_{V_iV_j\notin \Psi}E(H_{ij}^*[U\cup W]).
$$
Note $G^*$ comes equipped with an equipartition $\calS=\{W_{ij}^{**}:V_iV_j\notin \Psi\}$ of $W$ such that for all $W_{ij}^{**}\in \calS$, $G^*$ contains no $(U, W_{ij}^{**},d)$-trees.  Consequently, by Theorem \ref{thm:goodstrongequitable}, there are $1\leq m_1\leq M$, $L_1\leq X$, and a set $\Omega\subseteq \calS$  such that $|\bigcup_{W_{ij}^{**}\in \Omega}W^{**}_{ij}|\leq \e^{100}|W|$, and for all  $W^{**}_{ij}\in \calS\setminus \Omega$, there is a partition 
$$
W_{ij}^{**}=B_{ij}^0\cup B_{ij}^1\cup \ldots \cup B_{ij}^{m_1},
$$
so that $|B_{ij}^0|\leq \e^{100}|W_{ij}^{**}|$, and for each $u\in [m_1]$, $|B_{ij}^u|=L_1$ and $B_{ij}^u$ is $f(m_1)$-good in $G^*$.  Set 
$$
\Psi':=\Psi\cup \Big\{V_iV_j\in {\calP_{vert}\choose 2}\setminus \Psi: W_{ij}^{**}\in \Omega\Big\}.\label{psiprime}
$$
Note $|\Psi'|\leq (\e_1')^{1/4}t^2+\e^{100}t^2\leq 2\e^{100}t^2$.  Set $m=\lceil \e^{-100}m_1(m_1+1)\rceil $. By construction, $(1-\e_0)\frac{\ell}{m_1}\leq L_1\leq \frac{\ell}{m_1}$, so since $m_1\leq M$, our choice of $\e_0$ implies $L_1>\frac{\ell}{m}$.   Let $r$ be such that $L_1=rm+c$ for some $c<m$.  Then for all $V_iV_j\notin \Psi'$, and each $u\in [m_1]$, we can choose a partition $B_{ij}^u=B_{ij}^{u,0}\cup B_{ij}^{u,1}\cup \ldots \cup B_{ij}^{u,r}$ so that $|B_{ij}^{u,0}|=c<\frac{\ell}{m}$ and for all $1\leq v\leq r$, $|B_{ij}^{u,v}|=\frac{\ell}{m}$ (recall that $\ell$ is divisible by $(M!)^3$, so we know that $\ell/m\in \mathbb{N}$).  Define 
$$
R_{ij}^0=B_{ij}^0\cup \bigcup_{v=1}^{m_1}B_{ij}^{u,0}.
$$
Note $|R_{ij}^0|\leq \e^{100}\ell+\frac{m_1\ell}{m}<2\e^{100}\ell\leq \e^{75}\ell$.  We claim that for all $W_{ij}^{**}\in \calS\setminus \Omega$ and all $u\in [m_1]$, and $1\leq v\leq r$, $B_{ij}^{u,v}$ is $\psi(m)$-good in $H_{ij}$ with respect to every $u\in U^*_{ij}$.  Fix $W_{ij}^{**}\in \calS\setminus \Omega$, $u\in [m_1]$, $1\leq v\leq r$, and $u\in U_{ij}^*$.  We know that $B_{ij}^u$ is $f(m_1)$-good in $G^*$, so there is some $\tau\in \{0,1\}$ so that 
$$
|N^{\tau}_{G^*}(u)\cap B_{ij}^u|=|N^{\tau}_{H_{ij}^*}(u)\cap B_{ij}^u|\leq f(m_1)|B_{ij}^u|.
$$
Then since $H_{ij}^*$ is uniformly $\delta$-close to $H_{ij}[U_{ij}^*, W_{ij}^*]$, 
\[|N^{\tau}_{H_{ij}}(u)\cap B_{ij}^{u,v}|\leq f(m_1)|B_{ij}^u|+\delta |W_{ij}^*|=|B_{ij}^{u,v}|\Big(\frac{f(m_1)|B_{ij}^u|}{|B_{ij}^{u,v}|}+\frac{\delta|W_{ij}^*|}{|B_{ij}^{u,v}|}\Big),\]
which by definition of $m$, since $|B_{ij}^u|=L_1$, and since $|B_{ij}^{u,v}|=\ell/m$, is at most
\[|B_{ij}^{u,v}|(f(m)m+\delta m)\leq |B_{ij}^{u,v}|(f(m_1)+\delta)(2\e^{-100}m_1(m_1+1))\leq \psi(m)|B_{ij}^{u,v}|,\]
where the last inequality is by definition of $f$, since $m\leq 2\e^{-100}m_1(m_1+1)$, and our choice of $\delta$.  We point out to the reader that this step (namely deducing goodness of the $B_{ij}^{u,v}$ in $H_{ij}$ from the goodness of $W_{ij}^*$ in $H_{ij}^*$) is where it is important to have the degree of goodness in Theorem \ref{thm:goodstrongequitable} be small as a function of the number of parts.  Now, choose a reindexing 
$$
\{R_{ij}^1,\ldots, R_{ij}^{m_2}\}=\{B_{ij}^{u,v}: u\in [m_1], v\in [r]\},
$$
where $m_2=m_1r$.  Note $m_1r\leq \e^{-100}m_1(m_1+1)<(M!)^3$, so $\ell$ is divisible by $m_2$.

We have now defined, for all $V_iV_j\notin \Psi'$, a partition
$$
W_{ij}=W_{ij}^{bad}\cup R_{ij}^0\cup R_{ij}^1\cup \ldots \cup R_{ij}^{m_2}.
$$
This naturally translates into a new partition of $K_2[V_i,V_j]$.  In particular, for each  $1\leq u\leq m_2$, we let 
$$
\mathbf{R}_{ij}^u=\{xy\in K_2[V_i,V_j]: xy\in P_{ij}^{\alpha}\text{ for some }P_{ij}^{\alpha}\in R_{ij}^u\}.
$$
Setting $\mathbf{R}_{ij}^0=K_2[V_i, V_j]\setminus \bigcup_{u=1}^{m_2}\mathbf{R}_{ij}^u$, we have a partition $K_2[V_i,V_j]=\mathbf{R}_{ij}^0\cup \mathbf{R}_{ij}^1\cup \ldots \cup \mathbf{R}_{ij}^{m_2}$.  By construction, for each $u\in [m_2]$, $\mathbf{R}_{ij}^u$ is a union of $\ell/m$ many $P_{ij}^{\alpha}\in \calP_{edge}$, so by Fact \ref{fact:adding}, $\mathbf{R}_{ij}^u$ has $\disc_2((\ell/m)\e_2'''(\ell);(1/\ell)(\ell/m))$.  By our choice of $\e_2'''(\ell)$, we can conclude $\mathbf{R}_{ij}^u$ has $\disc_2(\e''_2(1/m);1/m)$.   By (\ref{align:foptri}), Claim \ref{cl:stable}, and since $\e_0\ll \e$, we have that
$$
|\mathbf{R}_{ij}^0|\leq (1+\e_2''(\ell))\frac{n^2}{t^2\ell}(|R_{ij}^0|+|W_{ij}^{bad}|)\leq (1+\e_2''(\ell))\frac{n^2}{t^2}(\e^{75}+\e_0)\leq \e^{50}\frac{n^2}{t^2}.
$$
Given $ijs\in {[t]\choose 3}$ with $V_iV_j, V_jV_s, V_iV_s\notin \Psi'$, and $(u,v,w)\in [m_2]^3$, we define
$$
\mathbf{R}_{ijs}^{uvw}:=(V_i\cup V_j\cup V_s, \mathbf{R}^u_{ij}\cup \mathbf{R}^v_{js}\cup \mathbf{R}^w_{is}).
$$
  By Corollary \ref{cor:counting}, for each such $\mathbf{R}_{ijs}^{uvw}$,
\begin{align}\label{align:size}
|K_3^{(2)}(\mathbf{R}_{ijs}^{uvw})|=(1\pm \e'_2(m))\frac{n^3}{t^3m^3}.
\end{align}
Our next goal is to show most of these $\mathbf{R}_{ijs}^{uvw}$ have density near $0$ or $1$ with respect to $H$.  In order to do this, we will need a few more sets to help us avoid the error triads.  Set
\begin{align*}
\Sigma_1&=\Sigma_0\cup \Big\{V_iV_jV_s:ijs\in {[t]\choose 3}\text{ and }{\{V_i,V_j,V_s\}\choose 2}\cap \Psi'\neq \emptyset\Big\}.
\end{align*}
By the bounds given above, $|\Sigma_1|\leq (\e_1')^{1/2}t^3+2\e^{100}t^2\leq 3\e^{100}t^3$. Given $ijs\in {[t]\choose 3}$, we define 
$$
U_{ij}(s):=\{P_{is}^\beta P_{js}^\gamma: \beta,\gamma\in [\ell]\}.
$$
Note that for all $ij\in {[t]\choose 2}$, $U_{ij}=\bigcup_{s\in [t]}U_{ij}(s)$.  Then define
\begin{align*}
\Sigma_3=\Big\{V_{i}V_{j}V_{s}&\in {\calP_{vert}\choose 3}  \setminus \Sigma_1:\\
&\max\{|U_{ij}(s)\cap U^{bad}_{ij}|, |U_{js}(i)\cap U^{bad}_{js}|, |U_{is}(j)\cap  U^{bad}_{is}|\geq \delta^{1/2} \ell^2\Big\}.
\end{align*}

Observe that for all $V_iV_j\notin \Psi'$, the number of $V_s$ with $V_iV_jV_s\in \Sigma_3$ is at most $\delta^{1/2} t$.  Consequently, $|\Sigma_3| \leq \delta t^3$. We now set 
$$
\Sigma:=\Sigma_1\cup \Sigma_3.
$$
 This will end up being the desired set $\Sigma$ given in the outline of the proof.  For this reason, the following bound on its size is crucial.
 
 \begin{claim}\label{cl:sigmasize}
$|\Sigma|\leq 4\e^{100} t^3$.
\end{claim}
\begin{proof}
By definition of $\Sigma$, $|\Sigma|\leq |\Sigma_1|+|\Sigma_3|\leq 3\e^{100}t^3+\delta t^3\leq 4\e^{100} t^3$.
\end{proof}

  We will show that if $V_iV_jV_s\notin \Sigma$, then the density of edges on $K_3^{(2)}(\mathbf{R}_{ijs}^{uvw})$ is close to $0$ or $1$ for all $u,v,w\in [m_2]$ (see Claim \ref{cl:fophom} below).

First, we need prove several facts about triples outside $\Sigma$.  Given  $V_iV_jV_s\notin \Sigma$, and  $u\in [m_2]$, define
\[U^1_{ij}(s)=\{P_{is}^\beta P_{js}^\gamma\in U_{ij}(s): |\{P_{ij}^\alpha\in R_{ij}^u:  G_{ijs}^{\alpha\beta\gamma}\in  \Gamma_1\} |\geq (1-\psi(m))|R_{ij}^u|\}\]
and
\[U^0_{ij}(s)=\{P_{is}^\beta P_{js}^\gamma\in U_{ij}(s): |\{P_{ij}^\alpha\in R_{ij}^u:  G_{ijs}^{\alpha\beta\gamma}\notin\Gamma_1\} |\geq (1-\psi(m))|R_{ij}^u|\}.\]
By Claim \ref{cl:stable}, we know that $U_{ij}(s)=U_{ij}^{bad}\cup U^1_{ij}(s)\cup U^0_{ij}(s)$.  However, we are more interested in the following set, rather than $U_{ij}^0(s)$.
\[U^{00}_{ij}(s)=\{P_{is}^\beta P_{js}^\gamma\in U_{ij}(s): |\{P_{ij}^\alpha\in R_{ij}^u:  (P_{is}^\beta P_{js}^\gamma)P_{ij}^{\alpha}\in E_{ij}^2 \} |\geq (1-\psi(m)^{1/2})|R_{ij}^u|\}.\]

Setting $U_{ij}^{err}(s)=U_{ij}(s)\setminus (U^1_{ij}(s)\cup U^{00}_{ij}(s))$, we wish to show that $|U_{ij}^{err}(s)|\leq \e^9|U_{ij}(s)|$.  Observe that $|K_2[U_{ij}(s),W_{ij}]\cap E_{ij}^2|$ is at least
$$
(\psi(m)^{1/2}-\psi(m))|U^{0}_{ij}(s)\setminus U^{00}_{ij}(s)|R_{ij}^u|\geq \psi(m)|U^{0}_{ij}(s)\setminus U^{00}_{ij}(s)|\ell/m,
$$
where the last inequality uses the fact that $|R_{ij}^u|=\ell/m$ and $\psi(m)<\e<1/4$. Because $ijs\notin \Sigma_1$, $|K_2[U_{ij}(s),W_{ij}]\cap E_{ij}^2|\leq 2(\e_1')^{1/2}\ell^3$.  Consequently, this implies
\begin{align*}
|U^{0}_{ij}(s)\setminus U^{00}_{ij}(s)|\leq 2(\e_1')^{1/2}\ell^3\psi(m)^{-1}\frac{m}{\ell}=2(\e_1')^{1/2}m\psi(m)^{-1}\ell^2\leq (\e_1')^{1/4}\ell^2,
\end{align*}
where the last inequality is because $\e_1'\ll \psi((M!)^3)/(M!)^3$. Further, since $ijs\notin \Sigma_3$, $|U_{ij}(s)\cap U_{ij}^{bad}|\leq \delta^{1/2} \ell^2 $, so $|U^1_{ij}(s)\cup U^{00}_{ij}(s)|$ is at least
\begin{align*}
|U_{ij}(s)|- |U^{0}_{ij}(s)\setminus U^{00}_{ij}(s)|-|U_{ij}^{bad}(s)|&\geq (1-\delta^{1/2})\ell^2 -(\e_1')^{1/4}\ell^2\\
&\geq \ell^2(1-\e^{50}),
\end{align*}
where the last inequality is because $\e'_1\ll \delta \ll\e$.  This shows $|U_{ij}^{err}(s)|\leq \e^{50}|U_{ij}(s)|$.  We will use this below.

\begin{claim}\label{cl:fophom}
For all  $ijs\in {[t]\choose 3}$ with $V_iV_jV_s\notin\Sigma$, and all $(u,v,w)\in [m_2]^3$, 
$$
\frac{|E\cap K_3^{(2)}(\mathbf{R}_{ijs}^{uvw})|}{|K_3^{(2)}(\mathbf{R}_{ijs}^{uvw})|}\in [0,\e)\cup (1-\e,1].
$$
\end{claim}
\begin{proof}
Fix $ijs\in {[t]\choose 3}$ with $V_iV_jV_s\notin \Sigma$, and $(u,v,w)\in [m_2]^3$.  Suppose towards a contradiction that 
\begin{align}\label{align:cont}
\frac{|E\cap K_3^{(2)}(\mathbf{R}_{ijs}^{uvw})|}{|K_3^{(2)}(\mathbf{R}_{ijs}^{uvw})|}\in [\e,1-\e].
\end{align}

To ease notation, define $A=R_{is}^w$, $B= R_{js}^v$, $C=R_{ij}^u$, and let 
\begin{align*}
\mathbf{F}_0&=K_2[A,B]\cap U_{ij}^{00}(s)\text{ and }\mathbf{F}_1=K_2[A,B]\cap U_{ij}^1(s).
\end{align*}
Then let $\mathbf{F}_{err}=K_2[A,B]\cap U_{ij}^{err}(s)$.  Note $|A|=|B|=|C|=\ell/m$, and our observations above show that $|\mathbf{F}_{err}|\leq \e^{50}|A||B|$.  We claim that 
\begin{align}\label{align:fopmain}
\min \{ |\mathbf{F}_1|, |K_2[A,B]\setminus \mathbf{F}_1| \}\geq \e^{3} |A||B|.
\end{align}
Suppose towards a contradiction that (\ref{align:fopmain}) is false.  Then we have that either $ |\mathbf{F}_1|<\e^3|A||B|$ or $|K_2[A,B]\setminus \mathbf{F}_1|< \e^{3} |A||B|$. Assume first that $ |\mathbf{F}_1|<\e^{3} |A||B|$.  Combining this with the above, we have that 
$$
|K_2[A,B]\setminus \mathbf{F}_0|\leq (\e^{9}+\e^3)|A||B|\leq 2\e^3|A||B|.
$$ 
Therefore, using (\ref{align:foptri}),
\begin{align*}
\sum_{\{P_{is}^{\alpha}P_{js}^{\beta}\in K_2[A,B]\setminus \mathbf{F}_0\}}\;\;\sum_{P_{ij}^{\gamma}\in C}|K_3(G_{ijk}^{\alpha,\beta,\gamma})|&\leq (1+\e''_2(\ell))\frac{2\e^3 n^3}{t^3\ell^3}|A||B||C|\leq 3\e^3\frac{n^3}{t^3m^3}.
\end{align*}
Note that if $P_{is}^{\alpha}P_{js}^{\beta}\in \mathbf{F}_0$, then at least $(1-\psi(m)^{1/2})|C|$ elements in $P_{ij}^{\gamma}\in C$ are such that $G_{ijs}^{\alpha,\beta,\gamma}\in \Gamma_0$.  Consequently,
\begin{align*}
\sum_{P_{ij}^{\gamma}\in C}|E\cap K_3(G_{ijk}^{\alpha,\beta,\gamma})|&\leq \sum_{\{P_{ij}^{\gamma}:  G_{ijs}^{\alpha,\beta,\gamma}\in \Gamma_0\}}\e_1''|K_3(G_{ijk}^{\alpha,\beta,\gamma})|+ \sum_{\{P_{ij}^{\gamma}:  G_{ijs}^{\alpha,\beta,\gamma}\notin \Gamma_0\}}|K_3(G_{ijk}^{\alpha,\beta,\gamma})|\\
&\leq (1+\e_2(\ell))\frac{n^3}{t^3\ell^3}(\e_1''|C|(1-\psi(m)^{1/2})+\psi(m)^{1/2}|C|)\\
&\leq \e^5\frac{n^3}{t^3m\ell^2},
\end{align*}
where the last inequality uses that $|C|=\ell/m$ and that $\psi(m), \e_1'(\ell),\e_1''\ll \e$.  We therefore have that
\begin{align*}
|E\cap K_3^{(2)}(\mathbf{R}_{ijs}^{uvw})|&\leq 2\e^3\frac{n^3}{t^3m^3} +\e^5|A||B|\frac{n^3}{t^3m\ell^2}\leq 3\e^3 \Big(\frac{n}{t}\Big)^3\frac{1}{m^3}.
\end{align*}
  By (\ref{align:size}), this shows that $|E\cap K_3^{(2)}(\mathbf{R}_{ijs}^{uvw})|< \e |K_3^{(2)}(\mathbf{R}_{ijs}^{uvw})|$, contradicting (\ref{align:cont}).  Thus we must have that $|K_2[A,B]\setminus \mathbf{F}_1|<\e^{3}K_2[A,B]$. In this case, 
  $$
|\mathbf{F}_0\cup \mathbf{F}_{err}|\leq \e^3|A||B|= \e^3 \frac{\ell^2}{m^2}.
$$ 
Consequently, using (\ref{align:foptri}),
  $$
  \sum_{\{P_{is}^{\alpha}P_{js}^{\beta}\in \mathbf{F}_0\cup \mathbf{F}_{err}\}}\;\;\sum_{P_{ij}^{\gamma}\in C}|K_3(G_{ijk}^{\alpha,\beta,\gamma})|\leq (1+\e''_2(\ell))\frac{\e^3n^3}{t^2m^3}\leq \frac{2\e^3n^3}{t^3m^3}.
  $$
  
  On the other hand, for each $P_{is}^{\alpha}P_{js}^{\beta}\in \mathbf{F}_1$, we know that at least $(1-\psi(m))|C|$ many elements $P_{ij}^{\gamma}\in C$ are such that $G_{isj}^{\alpha,\beta,\gamma}\in \Gamma_1$.  Consequently,
\[\sum_{\{P_{is}^{\alpha}P_{js}^{\beta}\in \mathbf{F}_0\}}\;\;\sum_{P_{ij}^{\gamma}\in C}|K_3(G_{ijk}^{\alpha,\beta,\gamma})\setminus E|\]
is at most
\[(1+\e_2''(\ell))\Big(\frac{n}{t\ell}\Big)^3|A||B|( \psi(m)|C| +(1-\psi(m))\e_1''|C|)\leq \e^5\Big(\frac{n}{t}\Big)^3\frac{1}{m^3},\]
where the last inequality uses that $|C|=\ell/m$ and that $\psi(m), \e_1'(\ell),\e_1''\ll \e$.  Combining these, we obtain that
\begin{align*}
|K_3^{(2)}(\mathbf{R}_{ijs}^{uvw})\setminus E|&\leq \frac{2\e^3n^3}{m^3}+\e^5\frac{n^3}{t^3m^3}\leq 3\e^3\Big(\frac{n}{t}\Big)^3\frac{1}{m^3}.
\end{align*}
  By (\ref{align:size}), this shows that $|K_3^{(2)}(\mathbf{R}_{ijs}^{uvw})\setminus E| <\e |K_3^{(2)}(\mathbf{R}_{ijs}^{uvw})|$, contradicting (\ref{align:cont}).  Thus we must have that $|K_2[A,B]\setminus \mathbf{F}_1|\geq \e^{3}|A||B|$, and therefore (\ref{align:fopmain}) is true. 
  
Consider the bipartite graph $G_{AB}=(A\cup B, \mathbf{F}_1)$.  Since (\ref{align:fopmain}) is true, Lemma \ref{lem:twosticks} implies one of the following holds.  
 \begin{enumerate}
 \item There is a set $S\subseteq B$ of size at least $\e^{6} |B|/2$ such that for all $b\in S$ 
$$
\min\{|N_{G_{AB}}(b)|, |A\setminus N_{G_{AB}}(b)| \}\geq \e^{6} |A|/2.
$$
\item There is a set $S\subseteq A$ of size at least $\e^{6} |A|/2$ such that for all $a\in S$ 
$$
\min\{|N_{G_{AB}}(a)|, |A\setminus N_{G_{AB}}(b)| \}\geq \e^{6} |B|/2.
$$
\end{enumerate}
Let us assume (1) holds (case (2) is similar).  So there is some $S\subseteq B$ of size at least $\e^{6} |B|/2$ such that for all $b\in S$ 
$$
\min\{|N_{G_{AB}}(b)|, |A\setminus N_{G_{AB}}(b)| \}\geq \e^{6} |A|/2.
$$
Recall we showed $|\mathbf{F}_{err}|\leq \e^{50}|A||B|$, so if we let 
$$
S'=\{b\in S: |\{a\in A: ab\in \mathbf{F}_{err}\}|\leq \e^{25}|A|\},
$$
then $|S'|\geq |S|-\e^{25}|B| \geq \e^6|B|/2 - \e^{25}|B|\geq \e^6|B|/4$, where the inequality is because $\e< 1/4$.  For each $b\in S'$, set
$$
A_0(b)=\{a\in A: ab \in \mathbf{F}_0\}\text{ and }A_1(b)=\{a\in A: ab \in \mathbf{F}_1\}.
$$
By definition of $S'$, $|A_1(b)|\geq \e^6 |A|/2$, and $|A_0(b)|\geq \e^6 |A|/2- \e^{25}|A|\geq \e^6 |A|/4$.   Now consider the bipartite graphs  $G_0(b)=(A_0(b)\cup C, E_0(b))$ and $G_1(b)=(A_1(b)\cup C, E_1(b))$, where 
\begin{align*}
E_0(b)&=\{ca\in K_2[C,A_0(b)] : G_{ijs}^{\alpha,\beta,\gamma}\in \Gamma_0\text{ where }c=P_{ij}^\alpha, a= P_{is}^\beta, \text{ and } b=P_{js}^{\gamma} \}\\
E_1(b)&=\{ca\in K_2[C,A_1(b)] : G_{ijs}^{\alpha,\beta,\gamma}\in \Gamma_1\text{ where }c=P_{ij}^\alpha, a= P_{is}^\beta, \text{ and } b=P_{js}^{\gamma} \}.
\end{align*}

We claim that for all $b\in S'$ and $a\in A_0(b)$, $|N_{G_0(b)}(a)|\geq (1-\psi(m)^{1/2})|C|$.  Say $b=P_{js}^{\gamma}\in S'$ and $a=P_{is}^{\beta}\in A_0(b)$.  Since $ab\in \mathbf{F}_0$, we know that $P_{js}^{\gamma}P_{is}^{\beta}\in W_{ij}^{00}(s)$.  By definition, this means that 
$$
|N_{G_0(b)}(a)|=|\{c=P_{ij}^{\alpha}\in C: G_{ijs}^{\alpha,\beta,\gamma}\in \Gamma_0\}|\geq (1-\psi(m)^{1/2})|C|,
$$
as desired.  Consequently, $|E_0(b)|\geq (1-\psi(m)^{1/2})|C||A_0(b)|$.  Standard arguments show that this implies there is a set $C_0(b)\subseteq C$ with $|C_0(b)|\geq (1-\psi(m)^{1/4})|C|$ so that for all $c\in C_0(b)$,
$$
|N_{G_0(b)}(c)|\geq (1-\psi(m)^{1/4})|A_0(b)|\geq (1-\psi(m)^{1/4})\e^6|A|/2\geq \e^6|A|/4.
$$
A similar argument shows that for all $b\in S'$ and $a\in A_1(b)$, $|N_{G_1(b)}(a)|\geq (1-\psi(m)^{1/2})|C|$, so $|E_1(b)|\geq (1-\psi(m)^{1/4})|C||A_1(b)|$.  Thus there exists a set $C_1(b)\subseteq C$ with $|C_1(b)|\geq (1-\psi(m)^{1/4})|C|$ so that for all $c\in C_1(b)$,
$$
|N_{G_1(b)}(c)|\geq (1-\psi(m)^{1/4})|A_1(b)|\geq (1-\psi(m)^{1/4})\e^6|A|/2\geq \e^6|A|/4. 
$$
Combining these, we have that for all $b\in S'$, there is a set $C(b)=C_0(b)\cap C_1(b)\subseteq C$, with size $|C(b)|\geq (1-2\psi(m)^{1/4})|C|$, such that for all $c\in C(b)$, 
$$
\min\{|N_{G_0(b)}(c)|, |N_{G_1(b)}(c)|\}\geq \e^6|A|/4.
$$
We claim that for all $b\in S'$ and $c\in C(b)$, $bc\in U_{is}^{bad}$.  To this end, fix $b\in S'$ and $c\in C(b)$.  Then $b=P_{js}^{\beta}$ and $c=P_{ij}^{\alpha}$ for some $\beta,\alpha\in [\ell]$.  By above and the definition of $H_{is}$,
\begin{align*}
\e^6|A|/4&\leq \min\{|N_{G_0(b)}(c)|, |N_{G_1(b)}(c)|\}\\
&=\min\{|P_{is}^{\gamma}\in A: G_{ijs}^{\alpha,\beta,\gamma}\in \Gamma_0\}|, \{|P_{is}^{\gamma}\in A: G_{ijs}^{\alpha,\beta,\gamma}\in \Gamma_1\}|\}\\
&\leq \min\{|N_{H_{is}}(P_{js}^{\beta}P_{ij}^{\alpha})\cap A|, |A\setminus N_{H_{is}}(P_{js}^{\beta}P_{ij}^{\alpha})|\}
\end{align*}

This shows $A=R_{is}^w$ is not $\e^6/4$-good in $H_{is}$ with respect to $bc$. Since $V_iV_jV_s\notin \Sigma$, we know that $R_{is}^w$ is $\psi(m)$-good with respect to every $u\in U_{is}\setminus U_{is}^{bad}$, so since $\psi(m)\ll \e$, this implies that $bc\in U_{is}^{bad}$.  Thus we have shown that for all $b\in S'$ and $c\in C(b)$, $bc\in U_{is}^{bad}$. Combining this with the bounds above, we have that 
$$
|U_{is}^{bad}|\geq |S'|(1-2\psi(m)^{1/4})|C|\geq \e^6|B||C|/4=\e^6\ell^2/4m^2.
$$
However, by Claim \ref{cl:stable}, $|U_{is}^{bad}|\leq \delta \ell^2<\e^6\ell^2/4m^2$, where the last inequality is since $\delta\ll \e/M$ and $m\leq M^3$.  Thus, we have arrived at a contradiction, which finishes the proof of Claim \ref{cl:fophom}.
\end{proof}

We are nearly done.  The last main task is to redistribute $\mathbf{R}_{ij}^0$ for each $V_iV_j\notin \Psi'$ (recall the definition of $\Psi'$ on page \pageref{psiprime}).  Fix $V_iV_j\notin \Psi'$.  Recall that $|W_{ij}|=\ell$, and $|R_{ij}^0|=|W_{ij}|-|B_{ij}^1\cup \ldots \cup B_{ij}^{m_2}|=\ell -\ell m_2/m$.  Since $m_2\leq M$, and $(M!)^3$ divides $\ell$, we know $m_2$ divides $\ell$, say $\ell=m_3m_2$. Note this implies $|R_{ij}^0|=m_2m_3-m_2\ell/m$.    Similarly, since $m\leq (M!)^3$, $m$ also divides $\ell$, so we can choose a partition $R_{ij}^0=S_{ij}^1\cup \ldots \cup S_{ij}^{m_2}$  into sets of size $(m_2m_3-m_2\ell/m)/m_2=m_3-\ell/m$.  For each $u\in [m_2]$, define
$$
\mathbf{S}^u_{ij}=\{xy\in K_2[V_i,V_j]: xy\in P_{ij}^{\alpha}\text{ for some }P_{ij}^{\alpha}\in S_{ij}^u\}.
$$
By Fact \ref{fact:adding}, $\mathbf{S}_{ij}^u$ has $\disc_2(m\e_2'''(\ell))$.  By our choice of $\e_2''$, this implies $\mathbf{S}_{ij}^u$ has $\disc_2(\e_2'(\ell))$.  Now set $\mathbf{P}_{ij}^u=\mathbf{R}_{ij}^u\cup \mathbf{S}_{ij}^u$ for each $u\in [m_2]$.  Note each $\mathbf{P}_{ij}^u$ is a union of $m_3-\ell/m+\ell/m=m_3$ elements of $\calP_{edge}$, so by Fact \ref{fact:adding}, $\mathbf{P}_{ij}^u$ satisfies $\disc_2(m\e''_2(\ell))$ and has size $m_3\frac{n^2}{\ell}(1\pm m_3\e_2'(\ell))$.  Combining this with our choice of $\e_2''$ and the fact that $m_3/\ell=1/m_2$, we have that each $\mathbf{P}_{ij}^u$ satisfies $\disc_2(\e_2(m_2); 1/m_2)$.  By construction, we know that each $|S_{ij}^u|\leq |R_{ij}^0|/m_2\leq \e^{75}\ell/m_2$.  Therefore, $|\mathbf{S}_{ij}^u|\leq \e^{25}\frac{n^2}{t^2m_2^2}\leq \e^{10}|\mathbf{P}_{ij}^u|$.  

Now for each $V_iV_j\in \Psi'$, choose an arbitrary partition $K_2[V_i,V_j]=\bigcup_{\alpha\leq m_2}\mathbf{P}_{ij}^{\alpha}$.  Since $|\Psi'|\leq 3\e^{100}t^2$, and every $\mathbf{P}_{ij}^u$ satisfies $\disc_2(\e_2(m_2);1/m_2)$, we obtain a $(t,m_2,\e^{50},\e_2(m_2))$-decomposition $\calQ$ with $\calQ_{vert}=\{V_1,\ldots, V_t\}$ and 
$$
\calQ_{edge}=\Big\{\mathbf{P}_{ij}^u: ij\in {[t]\choose 2}, u\in \{1,\ldots, m_2\}\Big\}.
$$
We claim that for all $V_iV_jV_s\notin \Sigma$ and $\alpha,\beta,\gamma\in [m_2]$, $\mathbf{P}_{ijs}^{\alpha,\beta,\gamma}=(V_i\cup V_j\cup V_s, \mathbf{P}_{ij}^{\alpha}\cup \mathbf{P}_{is}^{\beta}\cup \mathbf{P}_{js}^{\gamma})$ is $2\e$-homogeneous with respect to $H$.  

Fix $V_iV_jV_s\notin \Sigma$, and $\alpha,\beta,\gamma\in [m_2]$.  By construction, there are $u,v,w\in [m_2]$ so that 
\begin{align*}
\mathbf{P}_{ij}^{\alpha}=\mathbf{R}_{ij}^u\cup \mathbf{S}_{ij}^u, \text{ }\mathbf{P}_{is}^{\beta}=\mathbf{R}_{is}^v\cup \mathbf{S}_{is}^v,\text{ and }
\mathbf{P}_{is}^{\gamma}=\mathbf{R}_{js}^w\cup \mathbf{S}_{js}^w.
\end{align*}
To ease notation in what follows, let $E^1=E$ and $E^0={V(H)\choose 3}\setminus E$. Observe that Claim \ref{cl:fophom} implies there exists some $\tau\in \{0,1\}$ such that $|E^{\tau}\cap K_3^{(2)}(\mathbf{R}_{ijs}^{u,v,w})|\leq  \e n^3/(t^3m_2^3)$.  Define
\begin{align*}
T_{is}&=(V_i\cup V_j\cup V_s, \mathbf{P}_{ij}^{\alpha}\cup \mathbf{S}_{is}^v\cup \mathbf{P}_{is}^{\gamma}),\\
T_{ij}&=(V_i\cup V_j\cup V_s, \mathbf{S}_{ij}^u\cup \mathbf{P}_{is}^{\beta}\cup \mathbf{P}_{is}),\\
T_{js}&=(V_i\cup V_j\cup V_s,  \mathbf{P}_{ij}^{\alpha}\cup \mathbf{P}_{is}^{\beta}\cup \mathbf{S}_{js}^w).
\end{align*}
By Corollary \ref{cor:counting}, and the bound above on the size of each $\mathbf{S}_{ij}^u, \mathbf{S}_{is}^v,  \mathbf{S}_{js}^w$, we have that 
$$
\max\{|K_3^{(2)}(T_{ij})|, |K_3^{(2)}(T_{ij})|,|K_3^{(2)}(T_{ij})|\}\leq \e^{10}\frac{n^3}{t^3m_2^3}\leq \e^{5}|K_3^{(2)}(\mathbf{P}_{ijs}^{\alpha,\beta,\gamma})|,
$$
where the last inequality uses Corollary \ref{cor:counting}.  Thus we have that
\begin{align*}
|E^{\tau}\cap K_3^{(2)}(\mathbf{P}_{ijs}^{\alpha,\beta,\gamma})|&\leq \e \frac{n^3}{t^3m_2^3}+\e^{10}|K_3^{(2)}(\mathbf{P}_{ijs}^{\alpha,\beta,\gamma})|\leq 2\e|K_3^{(2)}(\mathbf{P}_{ijs}^{\alpha,\beta,\gamma})|.
\end{align*}
Since $\e\ll \e_1$, this implies that $\mathbf{P}_{ijs}^{\alpha,\beta,\gamma}$ is $\e_1$-homogeneous with respect to $H$ and moreover, by Proposition \ref{prop:homimpliesrandome}, satisfies $\disc_{2,3}(\e_1,\e_2(m_2))$ with respect to $H$.  This finishes the proof of Theorem \ref{thm:FOPfinite}.
\end{proofof}
\vspace{2mm}

We now prove Theorem \ref{thm:FOP}, which says that a hereditary $3$-graph property admits linear $\disc_{2,3}$-error if and only if it is $\NFOP_2$.

\vspace{2mm}

\label{proof:FOP}
\begin{proofof}{Theorem \ref{thm:FOP}}
 Suppose $\calH$ is an $\NFOP_2$, hereditary $3$-graph property.  That $\calH$ admits linear $\disc_{2,3}$-error is the content of Theorem \ref{thm:FOPfinite}.  

Conversely, suppose that $\calH$ has $\FOP_2$. We show that it requires non-linear $\disc_{2,3}$-error.  Choose $0<\e_1'\ll \e_1\ll 1$, and choose $\e_2,\e_2':\mathbb{N}\rightarrow (0,1]$ so that for all $x\in \mathbb{N}$, $0<\e_2'(x)\ll \e_2(x)\leq 1/4x$.  Fix $T,L,N\geq 1$.  We show that there is some $H\in \calH$ with $|V(H)|\geq N$ so that there exists no $(t,\ell,\e_1', \e_2'(\ell))$-decomposition $\calP$ of $V(H)$ which is $\disc_{2,3}$-regular with respect to $H$ with linear error.

Choose $m\gg L, T, N, (\e_1')^{-1}, (\e_2'(L))^{-1}$ and $n\gg m$.  By Lemma \ref{lem:3.8}, there exist sets $B=\{b_i: i\in [n]\}$ and $C=\{c_i: i\in [n]\}$ and $\gamma_1,\ldots, \gamma_m\subseteq K_2[B,C]$ so that for each $i\in [m]$, $\gamma_i$ has $\disc_2(\e_2(m)/T^2; 1/m)$.  Let $A=\{a_i: i\in [n]\}$ be a set of size $n$.

By Lemma \ref{lem:otherway}, there exists $H=(V,E)\in \calH$ with $V=A\cup B\cup C$, such that for all $i,j,s\in [n]$, $a_ib_jc_s\in E$ if and only if $b_jc_s\in \gamma_u$  for some $u\geq \lfloor i n/\ell \rfloor$. We show that there do not exist $1\leq \ell\leq L$, $1\leq t\leq T$ and a  $(t,\ell,\e'_1,\e'_2)$-decomposition of $V$ which is $(\e'_1,\e'_2(\ell))$-regular with respect to $H$ with linear $\disc_{2,3}$-error.  

Suppose towards a contradiction that there exist some $1\leq \ell\leq L$, $1\leq t\leq T$ and a  $(t,\ell,\e'_1,\e'_2)$-decomposition $\calP$ of $V$ which is $(\e'_1,\e'_2(\ell))$-regular with respect to $H$, and a set $\Sigma\subseteq {[t]\choose 3}$ with $|\Sigma|\leq \e_1't^3$ so that if $G\in \triads(\calP)$ fails $\disc_{2,3}(\e_1',\e_2'(\ell))$ with respect to $(H,\calP)$, then $V(G)\subseteq V_i\cup V_j\cup V_k$ for some $ijk\in \Sigma$.  Say $\calP_{vert}=\{V_i: i\in [t]\}$ and $\calP_{edge}=\{P_{ij}^{\alpha}: ij\in {[t]\choose 2}, \alpha \leq \ell\}$.  

For each $i\in [T]$, and $X\in \{A,B,C\}$, set $X_i=X\cap V_i$.  For each $X\neq Y$, $i\neq j$, and $\alpha\in [\ell]$ we let $P_{X_iY_j}^{\alpha}=P_{ij}^{\alpha}\cap K_2[X_i,Y_j]$.

Given $i\in [t]$ and $X\in \{A,B,C\}$, we say that $X_i$ is \emph{trivial} if $|X_i|< \e_1 |V_i|$, and that a vertex $x\in V$ is trivial if $x\in X_i$ for some trivial $X_i$.  For each $i\in [t]$, at most two elements of $\{A_i,B_i,C_i\}$ are trivial, and therefore, $V_i$ contains at most $2\e_1|V_i|$ many trivial vertices.  Thus, at most $2\e_1(3n)$ vertices are trivial.  For each $X\in \{A,B,C\}$, let $\calI_X$ be the set of $i\in [T]$ such that $X_i$ is non-trivial.  Then for each $X\in \{A,B,C\}$, we have $|X\cap (\bigcup_{i\in \calI_X}X_i)|\geq (1-12\e_1)|X|$, and consequently
$$
|\calI_X|\geq (1-12\e_1)|X|/(3n/t)=(1-12\e_1)t/3.
$$
Now define
$$
\calS_1=\{X_iY_j:X\neq Y \in \{A,B,C\}\text{ and either } i=j, \text{ or }i\in \calI_X\text{ or }j\in \calI_Y \}.
$$
By the above, $|\calS_1|\leq 3t+3(12\e_1)t^2/3\leq 40\e_1t^2$.  Now let 
$$
\calS_2=\{X_iY_j: X\neq Y \in \{A,B,C\}, ij\notin \calS_1\text{ and }|\{Z_s: X_iY_jZ_s\in \Sigma\}|\geq \sqrt{\e'_1}t\}.
$$
Since $|\Sigma|\leq \e_1't^3$, $|\calS_2|\leq \sqrt{\e_1'}t^2$.  Now let $\calS=\calS_1\cup \calS_2$. Then $|\calS|\leq 2\sqrt{\e_1} t^2$ (recall $\e_1'\ll \e_1\ll 1$).  By above, $|\calI_B||\calI_C|\geq (1-12\e_1)^2t^2/3$, so there must exist some $B_jC_s\notin \calS$ with $B_j\in \calI_B$ and $C_s\in \calI_C$.  Define $\calI_A'=\{i\in \calI_A: A_iB_jC_s\notin \Sigma\}$, and set $A'=\bigcup_{i\in \calI_A'}A_i$.  By construction, and since $\e_1'\ll \e_1$,
$$
|A'|\geq (1-12\e_1)|A|-\sqrt{\e_1'}t(3n/t)\geq (1-13 \e_1)|A|.
$$
 Given $\alpha\in [\ell]$, define
\begin{align*}
I^\alpha_{B_jC_s}=\{u\in [\ell]: |P^\alpha_{B_jC_s}\cap \gamma_u|\geq \e_1 |P^\alpha_{B_jC_s}|/\ell\},
\end{align*}
and set $R^\alpha_{B_jC_s}=\bigcup_{u\in I^\alpha_{B_jC_s}}(P^\alpha_{B_jC_s}\cap \gamma_u)$.  Note
$$
|P^\alpha_{B_jC_s}\setminus R^\alpha_{B_jC_s}|\leq \ell (\e_1 |P^\alpha_{B_jC_s}|/\ell)<\e_1|P^\alpha_{B_jC_s}|.
$$
Given $v\in [m]$, let $\gamma_v^+=\bigcup_{u> v}\gamma_u$ and $\gamma_v^-=\bigcup_{v<u}\gamma_u$.  Then define 
\begin{align*}
I^\alpha_{B_jC_s}(\e_1)&=\{ u\in I^\alpha_{B_jC_s}: \min\{|P^\alpha_{B_jC_s}\cap \gamma_u^+|,  |P^\alpha_{B_jC_s}\cap \gamma_u^-|\}\geq \e_1 |P^\alpha_{B_jC_s}|\}.
\end{align*}
Clearly $I^\alpha_{B_jC_s}(\e_1)$ is an interval, say $I^\alpha_{B_jC_s}(\e_1)=[u_1,u_2]$. Define 
$$
R^\alpha_{B_jC_s}(\e_1)=\bigcup_{u\in I^\alpha_{B_jC_s}(\e_1)}(P^{\alpha}_{B_jC_s}\cap \gamma_u).
$$
This implies that $|P^\alpha_{B_jC_s}\setminus R^\alpha_{B_jC_s}(\e)|$ is at most
\[|P^\alpha_{B_jC_s}\setminus R^\alpha_{B_jC_s}|+|R^\alpha_{B_jC_s}\cap [u_2+1,n]|+|P^\alpha_{B_jC_s}\cap [1,u_1-1]|\leq 3\e_1|P^\alpha_{B_jC_s}|.\]
Consequently,
\begin{align*}
|\bigcup_{\alpha\in [\ell]} R^\alpha_{B_jC_s}(\e_1)|\geq (1-3\e_1)|B_j||C_s|,
\end{align*}
which implies that $|\bigcup_{\alpha\in [\ell]} I^\alpha_{B_jC_s}(\e_1)|\geq (1-3\e_1)|B_j||C_s|/(|B_j||C_s|/\ell)\geq (1-3\e_1)\ell$.

Given $u\in [\ell]$, let $u^+=[\lceil (u-1)n/\ell\rceil ,n]$, $u^-=[1,\lfloor (u+1)n/\ell \rfloor ]$, and set
$$
u^*=[ \lfloor (u+1)n/\ell \rfloor,\lceil (u-1)n/\ell\rceil].
$$
Given $i\in [t]$, we then define
$$
I_i=\{u\in [\ell]: |A_i\cap u^*|\geq \e_1|A_i|\},
$$
and we set $I_{i}(\e)=\{u\in [\ell]:\min\{|A_i\cap u^+|, |A_i\cap u^-|\} \geq \e_1 |A_{i}|\}$.  Clearly $I_i(\e_1)$ is an interval, say $I_i(\e_1)=[v_1,v_2]$. Set $A_{i}(\e_1)=\bigcup_{u\in I_i(\e_1)}(A_i\cap u^*)$, and observe that
$$
|A_i\setminus A_i(\e_1)|\leq |A_i\cap v_1^-|+|A_i\cap v_2^+ |\leq 2\e_1 |A_i|.
$$
Consequently, $|\bigcup_{i\in \calI'_A}A_i(\e_1)|\geq  (1-2\e_1)|A'|$, which implies that 
$$
|\bigcup_{i\in \calI_A'}I_i(\e_1)|\geq (1-2\e_1)|A'|/(n/\ell)\geq (1-2\e_1)^2\ell(1-13\e_1)\ell\geq (1-20\e_1)\ell.
$$
This implies, with the above, that
\begin{align*}
|\Big(\bigcup_{i\in \calI_A'}I_i(\e_1)\Big)\cap \Big( \bigcup_{\alpha\in [\ell]} I^\alpha_{B_jC_s}(\e_1) \Big)|&\geq (1-23\e_1)\ell>0.
\end{align*}
Choose any $u\in (\bigcup_{i\in \calI_A'}I_i(\e_1)\Big)\cap \Big( \bigcup_{\alpha\in [\ell]} I^\alpha_{B_jC_s}(\e_1) )$ and define
\begin{align*}
G_1&=((A_i\cap u^+)\cup B_j\cup C_S, (P_{B_jC_s}^{\alpha}\cap \gamma_u^-)\cup P_{A_iC_s}^{\beta}\cup P_{B_jA_i}^{\gamma}),\\
G_0&=((A_i\cap u^-)\cup B_j\cup C_S, (P_{B_jC_s}^{\alpha}\cap \gamma_u^+)\cup P_{A_iC_s}^{\beta}\cup P_{B_jA_i}^{\gamma}).
\end{align*}
By construction, $K_3^{(2)}(G_1)\subseteq E$ and $K_3^{(2)}(G_0)\cap E\neq \emptyset$.  On the other hand, since $V_iV_jV_s\notin \Sigma$, we must have that $(H_{ijs}^{\alpha,\beta,\gamma},G_{ijs}^{\alpha,\beta,\gamma})$ satisfies $\disc_{2,3}(\e_1',\e_2'(\ell))$.  Let  $d\in [0,1]$ be such that $|E\cap K_3^{(2)}(G_{ijs}^{\alpha,\beta,\gamma})|=d|K_3^{(2)}(G_{ijs}^{\alpha,\beta,\gamma})|$.  Without loss of generality, say $d\leq 1/2$ (the other case is similar). By above, and by definition of $\disc_{2,3}(\e_1',\e_2'(\ell))$,
$$
\frac{1}{2}|K_3^{(2)}(G_1)|\leq |K_3^{(2)}(G_1)|(1-d)=\Big||E\cap K_3^{(2)}(G_1)|-d|K_3^{(2)}(G_1)|\Big|\leq \e'_1\frac{n^3}{\ell^3t^3}.
$$
Thus $ |K_3^{(2)}(G_1)|\leq 2\e_1' n^3/(\ell^3t^3)$.

Given $b\in B_j$, let $X_b^-=\{a\in A_i\cap u^-: ab\in P_{ij}^{\gamma}\}$.  Because $P_{ij}^{\gamma}$ satisfies $\disc_2(\e_2'(\ell);1/\ell)$, Lemma \ref{lem:sl} implies $P_{A_iB_j}^{\gamma}[A_i\cap u^-, B_j]$ satisfies $\disc_2(\e_2(\ell);1/\ell)$.  Therefore, a standard argument shows that there is a set $B_j'\subseteq V_j$ of size at least $(1-\sqrt{\e_2(\ell)})|B_j|$ such that for all $b\in B_j'$, 
$$
|X_b^-|=\frac{|A_i\cap u^-|}{\ell}(1\pm \sqrt{\e_2'(\ell)})\geq \frac{|A_i\cap u^-|}{2\ell}\geq \e_1|V_i|/2\ell.
$$

Similarly, let $Y_b^+=\{c\in C_s: vw\in P_{B_jC_s}^{\alpha}\cap \gamma_u^+\}$.  Because $P_{js}^{\alpha}$ satisfies $\disc_2(\e_2'(\ell);1/\ell)$, Lemma \ref{lem:sl} implies that $P_{B_jC_s}^{\gamma}$ satisfies $\disc_2(\e_2(\ell);1/\ell)$.  Thus, there is a set $B''_j\subseteq B_j$ of size at least $(1-\sqrt{\e_2(\ell)})|B_j|$ so that for all $v\in B_j''$, $|Y_b^+|=\frac{|C_s|}{\ell}(1\pm \sqrt{\e_2(\ell)})\geq \e_1|V_s|/2\ell$. Combining the above with the fact that $B_j$ is non-trivial, we have 
$$
|B_j''\cap B_j'|\geq (1-2\sqrt{\e_2(\ell)})|B_j|\geq \e_1^2|V_j|/2.
$$

Now given $b\in B_j'\cap B_j''$, since $|X^-_b|\geq \frac{1}{2\ell}|A_i|\geq \frac{\e_1}{2\ell}|V_i|$ and $|Y^+_b|\geq \frac{\e_1^2}{2\ell}|C_s|\geq \frac{\e_1^3}{2\ell}|V_s|$, and since $P_{is}^{\beta}$ satisfies $\disc_2(\e_2'(\ell);1/\ell)$, we must have that $P_{is}^{\beta}[X^-_b,Y^+_b]\geq \frac{|X_b||Y_b|}{2\ell}$.  Note that for all $ac\in P_{is}^{\beta}[X_b,Y_b]$, $bac\in K_3^{(2)}(G_1)$.  Thus we have shown
$$
|K_3^{(2)}(G_1)|\geq |B_j'\cap B_j''|\frac{\e_1|V_s|}{2\ell}\frac{\e_1|V_i|}{2\ell}\geq \e_1^4|V_i||V_j||V_s|/4\ell^2.
$$
Since $\e_1'\ll \e_1$, this contradicts that $ |K_3^{(2)}(G_1)|\leq 2\e_1'\frac{n^3}{\ell^3t^3}$, which finishes the proof. \end{proofof}

%% file: chapter7.tex
\chapter{Binary error}\label{sec:binary}

In this chapter we switch gears and consider binary error.   First, in Section \ref{subsec:rbinary}, we will give several examples of properties which admit binary $\disc_{2,3}$-error but contain $\overline{U}(k)$ for all $k$.  These examples will be used to prove Proposition \ref{prop:disc3vsVC}, which says that there exist properties which are far from SNIP but which admit $\disc_{2,3}$-error.   In Section \ref{subsec:rbinaryvdisc}, we will give an example which admits binary $\vdisc_3$-error but which is still far from SNIP, showing that being close to SNIP cannot characterize when a property admits binary $\vdisc_3$-error.  We will also show in Section \ref{subsec:quad}, that both $\calH_{\HP}$ and $\calH_{\GS_p}$ have finite $\VC$-dimension (we recall these are the minimal hereditary properties containing the examples appearing in Definitions \ref{def:hop3} and \ref{def:gsexamples}, respectively). This implies, by Proposition \ref{prop:wnipreduction}, that to show $\calH_{\HP}$ and $\calH_{\GS_p}$ require non-binary $\disc_{2,3}$-error, it suffices to show that they both require non-binary $\vdisc_3$-error.   

The rest of the chapter will then be concerned with showing that $\calH_{\HP}$ and $\calH_{\GS_p}$ require non-binary $\vdisc_3$-error.   There are two ways one could do this.  The most straightforward option is to present separate proofs for each.  However, the proofs are so similar that this becomes quite repetitive.  The second option is to isolate the important commonalities between the two examples, and then give a unified proof for both at once.  Although it is a bit more work, we have chosen  the second option, as we think it is of independent interest to know this can be done.  Thus, in Section \ref{subsec:special} we make a general definition, namely that of a $(p,\mu,\tau,\alpha,\rho)$-special $3$-graph (see Definition \ref{def:special}).  Roughly speaking, this is a $3$-partite $3$-graph equipped with a metric on each of its parts, and which satisfies certain axioms relating the edge relation and the metric structure.  We will show that both $\HP(N)$ and $\GS_p(n)$ are  $(p,\mu,\tau,\alpha,\rho)$-special, for appropriately chosen parameters. 

In Section \ref{subsec:binaryproof}, we show that any $\vdisc_3$-regular partition of a sufficiently large $(p,\mu,\tau,\alpha,\rho)$-special $3$-graph requires non-binary $\vdisc_3$-error (see Theorem \ref{thm:irreg}).  From this we will deduce Theorems \ref{thm:ternarytriads} and \ref{thm:ternarytriads1}.  We will end the section with a proof of Theorem \ref{thm:disc3}.

\section[Binary $\disc_{2,3}$-error is distinct from zero $\vdisc_3$-error]{Admitting binary $\disc_{2,3}$-error is distinct from admitting zero $\vdisc_3$-error}\label{subsec:rbinary} 
In this section we give a class of examples admitting binary $\disc_{2,3}$-error.  We will use these examples to show that admitting binary $\disc_{2,3}$-error is distinct from admitting zero $\vdisc_3$-error, and that there exist properties far from SNIP which still admit binary $\disc_{2,3}$-error.  The idea behind these examples is that if a property is ``essentially binary'' (in the sense of Definition \ref{def:binary} below), then it admits binary $\disc_{2,3}$-error. 

\begin{definition}\label{def:binary}
Suppose $G=(V\cup U\cup W,E)$ is a finite $3$-partite, $3$-graph.  We say that $G$ is \emph{$\ell$-binary} if the following holds.  There are partitions $\calP_{UV}$ of $K_2(U, V)$, $\calP_{UW}$ of $K_2(U, W)$, and $\calP_{VW}$ of $K_2(V,W)$, such that 
\begin{itemize}
\item $|\calP_{UV}|,|\calP_{UW}|,|\calP_{VW}|\leq \ell$,
\item for all $(\alpha,\beta,\gamma)\in \calP_{UV}\times\calP_{UW}\times\calP_{VW}$, if $G_{\alpha,\beta,\gamma}=(V\cup U\cup W, \alpha\cup \beta\cup \gamma)$, then either $K_3^{(2)}(G_{\alpha,\beta,\gamma})\subseteq E$ or $K_3^{(2)}(G_{\alpha,\beta,\gamma})\cap E=\emptyset$.
\end{itemize}
\end{definition}

Here are some examples of $\ell$-binary $3$-graphs.
\begin{itemize}
\item Given $n$, let $A_n=\{a_1,\ldots, a_n\}$, $B_n=\{b_1,\ldots, b_n\}$, and $C_n=\{c_1,\ldots, c_n\}$.  Then $G=(A_n\cup B_n\cup C_n,\{a_ib_jc_k: i<j<k\})$ is $2$-binary.  Indeed, set 
$$
\calP_{AB}=\{\{a_ib_j: i< j\}, \{a_ib_j:i\geq j\}\},
$$
 and define $\calP_{BC}$ and $\calP_{AC}$ analogously.
\item For any $3$-partite graph $G=(V\cup U\cup W,E)$, the hypergraph with vertex set $V\cup U\cup W$ and edge set $\{(x,y,z)\in VUW: xyz\text{ is a triangle in $G$}\}$ is $2$-binary.  Indeed, define
$$
\calP_{UV}=\{E\cap K_2(U,V), K_2(U,V)\setminus E\},
$$
 and define $\calP_{UW}$ and $\calP_{VW}$ analogously.
\end{itemize}

We now show that $\ell$-binary hypergraphs admit $\disc_{2,3}$-regular decompositions with binary $\disc_{2,3}$-error.   

\begin{proposition}\label{prop:binary-disctame}
For all $\e_1>0$, $\e_2:\mathbb{N}\rightarrow (0,1]$ and $\ell\geq 1$, there is $n_0$ and $T,L, N$ such that the following holds.  If $G=(V\cup U\cup W,E)$ is an $\ell$-binary, $3$-partite, $3$-graph with $|U\cup V\cup W|=n\geq N$, then there is $t_0\leq t\leq T$ and $\ell'\leq L$ such that $G$ has an $(\e_1,\e_2(\ell'))$-regular $(t,\ell, \e_2(\ell'),\e_1)$-decomposition with binary $\disc_{2,3}$-error.  
\end{proposition}
\begin{proof}
Fix $\e_1>0$,  $\e_2:\mathbb{N}\rightarrow (0,1]$ and $\ell\geq 1$.  Choose $\e_1'\ll \mu_1\ll \e_1$, and  $\e_2',\mu:\mathbb{N}\rightarrow (0,1]$ with the property that $\e_2'(x)\ll \mu_2(x)\ll \e_2(x)$ for all $x\geq 1$.  Choose $f$ and $g$ as in Lemma \ref{lem:refinement} for $\e_1'$ and $\e_2'$. Then set $L=f(3,\ell, 3,1)$, $T=f(3,\ell, 3,1)$, and let $N_1=N_1(3,\ell,3,1)$ from Lemma \ref{lem:refinement}.  Finally, choose some $N\gg N_1,(\e_1')^{-1}, \e_2'(L)^{-1}$.

 Now suppose $H=(V\cup U\cup W,E)$ is a $3$-partite $3$-graph which is $\ell$-binary, and where $|U\cup W\cup V|=n\geq N$.  Choose partitions $K_2[U,V]=\bigcup_{\alpha\leq \ell_{UV}}P^{UV}_\alpha$, $K_2[U,W]=\bigcup_{\alpha\leq \ell_{UW}}P^{UW}_\alpha$, and $K_2[V,W]=\bigcup_{\alpha\leq \ell_{VW}}P^{VW}_\alpha$, with $\ell_{UV},\ell_{UW},\ell_{VW}\leq \ell$ witnessing that $G$ is $\ell$-binary.  Note that each triple $(P^{UV}_{\alpha}, P^{UW}_{\beta}, P^{VW}_{\gamma})$ comes  equipped with a $\delta=\delta(\alpha,\beta,\gamma)\in \{0,1\}$ so that $K_3^{(2)}(U\cup V\cup W, P^{UV}_{\alpha}\cup P^{UW}_{\beta}\cup P^{VW}_{\gamma})\subseteq E^{\delta}$.
 
Let $\calP$ be the $(3,\ell)$-decomposition of $V$ with $\calP_{vert}=\{U,V,W\}$ and $P_{edge}=\{P^{XY}_{\alpha}: XY\in \{UV,UW,VW\}, \alpha\leq \ell_{XY}\}$.  Let $\calQ$ be the $(3,1)$ decomposition of $V$ consisting of $\calQ_{vert}=\{U,V,W\}$ and $\calQ_{edge}=\{K_2[U,V], K_2[V,W], K_2[U,W]\}$. Apply Lemma \ref{lem:refinement} to $\calP$ and $\calQ$ to obtain $\ell'\leq L$, $t\leq T$, and a $(t,\ell',\e_1',\e'_2(\ell'))$-decomposition $\calR$ of $V$ which has no $\disc_2$-irregular triads and which is an $(\e_1',\e_2'(\ell'))$-approximate refinement of $\calP$.  Let $\Sigma\subseteq {\calR_{vert}\choose 2}$ be as in the definition of an $(\e_1',\e_2'(\ell'))$-approximate refinement of $\calP$.  Let $\calR_{vert}=\{V_1,\ldots, V_t\}$ and $\calR_{edge}=\{R_{ij}^{\alpha}: ij\in {[t]\choose2}, \alpha\leq \ell'\}$.  
 
Suppose $R^{\alpha}_{ij}, R^{\beta}_{ik},R^{\gamma}_{jk}\in \calR_{edge}$ are such that $V_iV_j, V_jV_k, V_iV_k\notin \Sigma$.  Let $R_{ijk}^{\alpha,\beta,\gamma}$ denote the triad $(V_i\cup V_j\cup V_k, R^{\alpha}_{ij} \cup R^{\beta}_{ik}\cup R^{\gamma}_{jk})$.  By definition of an $(\e_1',\e_2'(\ell'))$ refinement, there are $X,Y,Z\in \{U,V,W\}$ with $|V_i\setminus X|<\e_1'|V_i|$, $|V_j\setminus Y|<\e_1'|V_j|$, and $|V_k\setminus Z|<\e_1'|V_k|$.  If $X=Y$, $X=Z$, or $Y=Z$, then $|E\cap K_3^{(2)}(R_{ijk}^{\alpha,\beta,\gamma})|=0$.  By Proposition \ref{prop:homimpliesrandome}, $R_{ijk}^{\alpha,\beta,\gamma}$ satisfies $\disc_{2,3}(\e_1,\e_2(\ell'))$ with respect to $H$.  
 
Assume now that $X,Y,Z$ are pairwise distinct.  Then by definition of an $(\e_1',\e_2'(\ell'))$-approximate refinement, there are $\alpha',\beta',\gamma'$ so that the following hold (possibly after relabeling),
 $$
 |R^{\alpha}_{ij}\setminus P_{\alpha'}^{UV}|<\e'_1|R^{\alpha}_{ij}|, \text{ }  |R^{\beta}_{ik}\setminus P_{\beta'}^{UW}|<\e'_1|R^{\beta}_{ik}|, \text{ and } | R^{\gamma}_{jk}\setminus P_{\gamma'}^{VW}|<\e_1'|R^{\gamma}_{jk}|,
 $$
and further, each of $R^{\alpha}_{ij}\setminus P_{\alpha'}^{UV}$, $R^{\beta}_{ik}\setminus P_{\beta'}^{UW}$, and $R^{\gamma}_{jk}\setminus P_{\gamma'}^{VW}$ have $\disc_2(\e'_2(\ell'))$.  By Corollary \ref{cor:counting}, $K_3^{(2)}(R_{ijk}^{\alpha,\beta,\gamma})=(1\pm \mu_2(\ell'))\frac{|V_i||V_j||V_k|}{(\ell')^3}$, and the number of triples in $K_3^{(2)}(R_{ijk}^{\alpha,\beta,\gamma})$ involving a pair from one of $R^{\alpha}_{ij}\setminus P_{\alpha'}^{UV}$, $R^{\beta}_{ik}\setminus P_{\beta'}^{UW}$, or $R^{\gamma}_{jk}\setminus P_{\gamma'}^{VW}$ is at most
$$
3\e_1'(1+\mu_2(\ell'))\frac{1}{(\ell')^3}|V_i||V_j||V_k|\leq 4\mu_1 |K_3^{(2)}(R_{ijk}^{\alpha,\beta,\gamma})|.
$$
Thus, if $\delta=\delta(\alpha',\beta',\gamma')$, then
$$
|E^{\delta}\cap K_3^{(2)}(R_{ijk}^{\alpha,\beta,\gamma})|\geq (1-4\mu_1)|K_3^{(2)}(R_{ijk}^{\alpha,\beta,\gamma})|.
$$
By Proposition \ref{prop:homimpliesrandome}, this shows $R_{ijk}^{\alpha,\beta,\gamma}$ has $\disc_{2,3}(\e_1,\e_2(\ell'))$ with respect to $H$.  By definition, this shows that $\calR$ has binary $\disc_{2,3}$-error with respect to $H$.
 \end{proof}

As a consequence we can now show that many examples admit binary $\disc_{2,3}$-error. In particular, we can now prove Propopistion \ref{prop:disc3vsVC}, which shows that there is a property containing $\Ubar(k)$ for every $k$ and which admits binary $\disc_{2,3}$-error.  Recall that given an infinite $3$-graph $H$, $\age(H)$ denotes the class of all finite $3$-graphs which are isomorphic to an induced sub-$3$-graph of $H$.  It is easy to see that $\age(H)$ is always a hereditary $3$-graph property.

\begin{definition}\label{def:W}
Define $\calW$ to be the $3$-graph 
$$
\calW=(\{a_i,b_i: i\in \mathbb{N}\}\cup \{c_S: S\subseteq\mathbb{N}\}, \{a_ib_jc_S:j\in S\}).
$$  
Let $\calH_{\calW}:=\age(\calW)$.
\end{definition}

\vspace{2mm}

\begin{proofof}{Proposition \ref{prop:disc3vsVC}} Clearly every element in $\calH_{\calW}$ is $2$-binary.  Consequently, by Proposition \ref{prop:binary-disctame}, $\calH_{\calW}$ admits binary $\disc_{2,3}$-error.  On the other hand, $\calH_{\calW}$ contains $\Ubar(k)\in \calW$ for all $k$ by definition. 
\end{proofof}

\vspace{2mm}

We can use $\calH_{\calW}$ as an example to demonstrate several other distinctions among the definitions from the introduction.  Recall that in Corollary \ref{cor:wnipcor}, we showed that a property $\calH$ admitting binary $\vdisc_3$-error must be $\vdisc_3$-homogeneous. For a property $\calH$ to be $\vdisc_3$-homogeneous, by Theorem \ref{thm:vdiscubark}, there must be a $k\geq 1$ such that $\Ubar(k)\notin \trip( \calH)$.   Consequently, $\calH_{\calW}$ does not admit binary $\vdisc_3$-error, but does admit binary $\disc_{2,3}$-error.  Moreover, $\calH_{\calW}$ has unbounded slicewise $\VC$-dimension (which implies unbounded $\VC$-dimension) and is not slicewise stable (and thus not stable).  Thus $\calH_{\calW}$ shows that admitting binary $\disc_{2,3}$-error is distinct from admitting binary $\vdisc_3$-error, from being NIP, from  being close to slicewise NIP, from being stable and from being close to slicewise stable.  

\section[Binary $\vdisc_3$-error is not characterized by being close to slicewise stable]{Admitting binary $\vdisc_3$-error is not characterized by being close to slicewise stable}\label{subsec:rbinaryvdisc}

We will now also give an example to show that being close to slicewise stable cannot characterize when a hereditary $3$-graph property admits binary $\vdisc_3$-error.  

\begin{definition}\label{def:V}
Define $\calV$ to be the $3$-graph 
$$
\calV=(\{a_i,b_i, c_i: i\in \mathbb{N}\}, \{a_ib_jc_k:i\leq j\}).
$$  
Let $\calH_{\calV}:=\age(\calV)$.
\end{definition}

\begin{proposition}\label{prop:slvdiscbin}
$\calH_{\calV}$ is far from slicewise stable but admits binary $\vdisc_3$-error.
\end{proposition}

\begin{proof}
Since $\overline{H}(k)\in \trip(\calH_{\calV})$ for all $k\in \mathbb{N}$, Theorem \ref{thm:simclassws} implies $\calH_{\calV}$ is far from slicewise stable.  On the other hand, it is easy to check that $\calH_{\calV}$ is $2$-binary and thus admits binary $\disc_{2,3}$-error by Proposition \ref{prop:binary-disctame}.  It is an exercise to show $\calH_{\calV}$ has finite slicewise VC-dimension. As an immediate consequence of this and Proposition \ref{prop:wnipreduction}, we obtain that $\calH_{\calV}$ admits binary $\vdisc_3$-error.
\end{proof}

\section{Proof of Proposition \lowercase{\ref{prop:finitevcnotdisc3tame}}}\label{subsec:quad}

In this section we  show that both $\calH_{\HP}$ and $\calH_{\GS_p}$ have bounded $\VC$-dimension, which will allow us to apply Proposition \ref{prop:wnipreduction} to both.   We recall some terminology from the study of $\VC$-dimension.

Suppose $X$ is a set. A \emph{set system} is a pair $(X,\calF)$ where  $\calF\subseteq \calP(X)$.  Given $Y\subseteq X$, define $\calF\cap Y:=\{Y\cap F: F\in \calF\}$.  We then say that $\calF$ \emph{shatters $Y$} if $|\calF\cap Y|=2^{|Y|}$.  

\begin{definition}
The \emph{$\VC$-dimension of $(X,\calF)$}, $\VC(X,\calF)$, is the largest size of a shattered subset of $X$.
\end{definition}

Note that given a $3$-graph $H=(V,E)$, by definition,
\begin{align*}
\VC(H)=\max\Big\{\VC\Big(V, \Big\{N_H(ab): ab\in {V\choose 2}\Big\}\Big), \VC\Big({V\choose 2}, \{N_H(a): a\in V\}\Big)\Big\}.
\end{align*}

We will now show that $\calH_{\GS_p}$ has bounded $\VC$-dimension.  The key ingredient is a fact proved in \cite[Corollary A.10]{Terry.2021a}.  

\begin{theorem}\label{thm:gsvc}
For all $p\geq 3$, the set system $\calF_{\GS_p}=(\F_p^n, \{A(p,n)-g:g\in \F_p^n\})$ has $\VC$-dimension at most $3$.  
\end{theorem}

\vspace{2mm}

\begin{corollary}\label{cor:gsvc}
For all $p\geq 3$, $\calH_{\GS_p}$ has $\VC$-dimension at most $3$.
\end{corollary}
\begin{proof}
Fix $p\geq 3$. We show $\calH_{\GS_p}$ has $\VC$-dimension at most $3$.  Suppose towards a contradiction there is $G=(V,E)\in \calH_{\GS_p}$ with $\VC$-dimension at least $4$.  By definition of $\calH_{\GS_p}$ there is $n\geq 1$ such that $G$ is isomorphic to an induced sub-$3$-graph of $\GS_{p}(n)$.  Thus, $\GS_{p}(n)$ also has $\VC$-dimension at least $4$.  Let $A=\{a_g: g\in \mathbb{F}_p^n\}$, $B=\{b_g: g\in \mathbb{F}_p^n\}$, $C=\{c_g: g\in \mathbb{F}_p^n\}$.  Then $\GS_{p}(n)$ has vertex set $V'=A\cup B\cup C$ and edge set $E'=\{a_gb_{g'}c_{g''}: g+g'+g''\in A(p,n)\}$.  Since $\VC(\GS_{p}(n))\geq 4$, by definition, one of the following hold.
\begin{enumerate}
\item There is $X\subseteq V'$ a set of size $4$, and for each $Y\subseteq X$, $e_Y\in {V'\choose 2}$ such that for all $x\in X$, $\{x\}\cup e_Y\in E'$ holds if and only if $x\in Y$.
\item There is $X\subseteq {V'\choose 2}$ of size $4$ and for each $Y\subseteq X$, $v_Y\in A\cup B\cup C$ such that for all $e\in X$, $\{v_Y\}\cup e\in E'$ if and only if $e\in Y$.  
\end{enumerate}

Suppose first we are in case (1).  Because there exists $e\in {V'\choose 2}$ with $X=N(e)\cap X$, we must have that either $X\subseteq A$, $X\subseteq B$, or $X\subseteq C$.  Without loss of generality, let us assume $X\subseteq A$.  Say $X=\{a_x,a_y,a_z,a_w\}$.  By definition of $\GS_{p}(n)$, for any $\emptyset\neq Y\subseteq X$, since $N(e_Y)\cap X\neq \emptyset$, must have that $e_Y=b_{g_Y}c_{h_Y}$ for some $g_Y,h_Y\in \mathbb{F}_p^n$.  But then $\{x,y,z,w\}$ is shattered by $\calF_{\GS_p}$, contradicting Theorem \ref{thm:gsvc}.  

Suppose now we are in case (2).  Because there is some $v\in V$ such that $X\cap N(v)=X$, we must have that either $X\subseteq AB$, $X\subseteq AC$ or $X\subseteq BC$.  Without loss of generality, assume $X\subseteq AB$, say $X=\{a_{g_1}b_{h_1},\ldots, a_{g_4}b_{h_4}\}$. Note that for any $\emptyset \neq Y\subseteq X$, since $N(e_Y)\cap X'\neq \emptyset$, we must have that $v_Y=c_{\nu_Y}$ for some $\nu_Y\in \mathbb{F}_p^n$.  Again, we see that $\{g_1+h_1,\ldots, g_4+h_4\}$ is shattered by $\calF_{\GS_p}$,  contradicting Theorem \ref{thm:gsvc}.
\end{proof}

We now turn to showing that $\calH_{\HP}$ has bounded $\VC$-dimension.  This will be clear to the model theorist, as the edge relation of $\calH_{\HP}$ is uniformly definable in Presburger arithmetic, for instance.  On the other hand, this fact will also be known to many combinatorialists, as the $3$-graphs of $\calH_{\HP}$ are examples of semi-algebraic hypergraphs (see e.g. \cite{Alon.20058t}), which are well known to have bounded $\VC$-dimension. However, we include a proof here for the sake of completeness.  

\begin{fact}
$\calH_{\HP}$ has $\VC$-dimension $1$. 
\end{fact}
\begin{proof}
By definition of $\calH_{\HP}$, it suffices to show that $\VC(\HP(N))\leq 1$ for all $N$.  Fix $N\geq 1$, and consider the $3$-graph $\HP(N)=(V,  E)$ where $V=A\cup B\cup C$, $A=\{a_i:i\in [N]\}$, $B=\{b_i:i\in [N]\}$, $C=\{c_i:i\in [N]\}$, and $E=\{a_ib_jc_k: i+j+k\geq N+2\}$.   Suppose towards a contradiction that $\HP(N)$ has $\VC$-dimension as least $2$.  This means there is either a set $X\subseteq V$ of size $2$ shattered in $(V, \calF_1)$ where $\calF_1=\{N_{\HP(N)}(uv): uv\in {V\choose 2}\}$, or a there is a set $X\subseteq {V\choose 2}$ of size $2$ shattered $({V\choose 2},\calF_2)$, where $\calF_2=\{N_{\HP(N)}(u): u\in V\}$.  

Suppose first there is  $X\subseteq {V\choose 2}$ of size $2$ shattered in $({V\choose 2}, \calF_1)$.  Arguing as in in the proof of Corollary \ref{cor:gsvc}, it is not difficult to show that because $\HP(N)$ is $3$-partite, $X\subseteq K_2[A,B]$, or $X\subseteq K_2[A,C]$, or $X\subseteq K_2[B,C]$.  Since $A$, $B$, and $C$ behave symmetrically in $\HP(N)$, let us assume that $X\subseteq K_2[A,B]$.  Say $X=\{a_{i_1}b_{j_1},a_{i_2}b_{j_2}\}$.  By relabeling if necessary, assume $i_1+j_1\leq i_2+j_2$.  It is not difficult to see that because $X$ is shattered in $\HP(N)$, there is some $c_k\in C$ so that $c_ka_{i_1}b_{j_1}\in E$ but $c_ka_{i_2}b_{i_2}\notin E$.  But this means that $i_1+j_1+k\geq N+2$ while $i_2+j_2+k<N+2$, which is not possible, as $i_1+j_1\leq i_2+j_2$.

Suppose now that there is $X\subseteq V$ of size $2$ shattered $(V,\calF_2)$.  Again, it is not hard to see that since $\HP(N)$ is $3$-partite, we must have $X\subseteq A$, $X\subseteq B$, or $X\subseteq C$.  As the roles of $A,B,C$ are symmetric, we may assume, say, $X=\{a_{i_1},a_{i_2}\}\subseteq A$, where $i_1\leq i_2$.  Since $X$ is shattered, there is some $b_jc_k$ so that $b_jc_ka_{i_1}\in E$ but $b_jc_ka_{i_2}\notin A$.  But this means $j+k+i_1\geq N+2$ while $j+k+i_2<N+2$, which is not possible, as $i_1\leq i_2$.
\end{proof}

\vspace{2mm}

\section{Special $3$-graphs and the proof of Theorems \lowercase{\ref{thm:ternarytriads}} and \lowercase{\ref{thm:ternarytriads1}}}\label{subsec:special}

In this section we define $(p,\mu,\tau, \alpha,\rho)$-special $3$-graphs (see Definition \ref{def:special}).  The purpose of Definition \ref{def:special}  is to isolate common properties our main examples, $\GS_p(n)$ and $\HP(n)$, specifically those which are required to prove both examples require non-binary $\vdisc_3$-error (which we will do in Section \ref{subsub:binaryproof}).   

The main idea is that both examples come equipped with natural metrics on their set of vertices, which interact nicely with their edge sets.  As Definition \ref{def:special} is quite technical, we will begin by examining each of the two examples in turn by way of motivation.  We will begin with $\GS_p(n)$, where the picture is somewhat simpler.

\subsection{Special properties of $\GS_p(n)$}\label{subsub:gs} In this subsection we will go through the main properties of $\GS_p(n)$ used in the proof that it requires non-binary $\vdisc_3$-error.  We first recall the definition of $\GS_p(n)$, and prove some basic facts about it.  Given $x\in \F_p^n$ and $i\in [n]$, $x_i$ denotes the $i$th coordinate of $x$.  Given $1\leq i\leq n$, we let $H_i$ denote the subspace of $\mathbb{F}_p^n$ with the property that for all $x\in H_i$, $x_j=0$ for all $1\leq j\leq i$.  By convention, we set $H_0=\mathbb{F}_p^n$.  

Recall that $A(p,n)=\bigcup_{i=1}^nH_i+e^i$, where $e^i$ is the $i$th standard basis vector in $\F_p^n$.  In other words, a vector $x\in \F_p^n$ is in $A(p,n)$ if the first non-zero coordinate in $x$ is $1$.  Then $\GS_p(n)$ is defined as follows.
$$
\GS_p(n)=(\{a_g,b_g,c_g: g\in \F_p^n\}, \{a_gb_{g'}c_{g''}: g+g'+g''\in A(p,n)\}).
$$
To ease notation, while in Section \ref{subsub:gs}, we will let $A=\{a_g: g\in \F_p^n\}$, $B=\{b_g: g\in \F_p^n\}$, $C=\{c_g: g\in \F_p^n\}$, $E=E(\GS_p(n))$, and $N=p^n$.  Give $x, y\in \mathbb{F}_p^n$ we set 
$$
\lambda(x,y)=\max\{i\in \{0,\ldots, n\}: \text{ $x$ and $y$ are in the same coset of $H_i$}\}.
$$
Note that by definition, $x_{\lambda(x,y)+1}\neq y_{\lambda(x,y)+1}$ and for all $1\leq j\leq \lambda(x,y)$, $y_j=x_j$.  

\begin{lemma}\label{lem:gsmetric}
Given $x,y\in \mathbb{F}_p^n$, set  $d(x,y)=0$ if $x=y$ and set $d(x,y)=p^{-\lambda(x,y)}$ if $x\neq y$.  Then $d$ is an ultra metric on $\mathbb{F}_p^n$, with distances in $[0,1]$, and for all $x\in \F_p^n$, $B_{p^{-i}}(x)$ is the coset of $H_{i+2}$ containing $x$.   
\end{lemma}
\begin{proof}
If $x\neq y$, then $x$ and $y$ are in different cosets of $H_n$, so $d(x,y)\geq p^{-n}>0$.  Given $x\in \F_p^n$, it is clear that $d(x,y)=p^{-i}$ if and only $y$ is in the same coset of $H_i$ as $x$, but not the same coset of $H_{i+1}$.  Thus $B_{p^{-i}}(x)=\{y\in \F_p^n: d(x,y)<p^{-i}\}$ is the coset of $H_{i+2}$ containing $x$. This shows that $d(x,y)=d(y,x)$ for all $x,y\in \F_p^n$ and $d(x,y)=0$ if and only if $x=y$.

Consider $x,y,z\in \F_p^n$.  Suppose $d(x,y)=p^{-i}$, so $x$, $y$ are in the same coset of $H_i$, say $H_i+g$, but not the same coset of $H_{i+1}$, say $x\in H_{i+1}+g+g'$ and $y\in H_{i+1}+g+g''$.  If $z\notin H_{i}+g$, then $d(x,z)=d(y,z)=p^{-i-1}=pd(x,y)$.  If $z\in H_{i}+g$, then it cannot be in both $H_{i+1}+g+g'$ and $H_{i+1}+g+g'$, so $\max\{d(z,x),d(z,y)\}=p^{-i}=d(\xbar,\ybar)$.  This shows $d(x,y)\leq \max\{d(x,z),d(y,z)\}$, and thus $d$ is an ultrametric.
\end{proof}

The fact that $d$ is an ultrametric on $\F_p^n$ makes things particularly easy to analyze in $\GS_p(n)$.  For example, for any $r\in (0,1]$ and $x\in \F_p^n$, $|B_r(x)|=|H_{m+1}|=p^{-m-1}N$, where $m=\max\{i\leq n: r\leq p^{-i}\}=\lceil \log_p(r^{-1})\rceil$.  Thus every ball of radius $r$ contains the same number of points, and $p^{-1}rN\leq |B_r(x)|\leq prN$.  It is also very simple to partition balls into smaller ones, as our next claim shows.  

\begin{claim}[Balls can be partitioned]\label{cl:gsball}
For any $x\in \F_p^n$ and $0<r_1<r_2\leq 1$, there is a partition of $B_{r_2}(x)$ into at most $pr_2r_1^{-1}$ disjoint open balls of radius $r_1$, one of which is centered at $x$.  
\end{claim}
\begin{proof}
Let $m_1=\lceil \log_p(r_1^{-1})\rceil$, $m_2=\lceil \log_p(r_2^{-1})\rceil$.  Then $B_{r_2}(x)$ is the coset $H_{m_2+2}+g$ which contains $x$.  We then partition $H_{m_2+2}+g$ into cosets of $H_{m_1+2}$.  All of these cosets are balls of radius $r_1$, and one of them must contain $x$.  We can easily compute that the number of balls in this partition is $p^{m_1-m_2}\leq pr_2r_1^{-1}$.
\end{proof}

Another nice feature of $\GS_p(n)$ is that it is vertex regular (i.e. all of its vertices have the same degree).   

\begin{claim}[$\GS_p(n)$ is vertex regular]\label{cl:gsdegree}
For all $u,v\in A\cup B\cup C$,  $|N(v)|=|N(u)|=N^2(\frac{1}{p-1}\pm p^{-n})$.  
\end{claim}
\begin{proof}
Given $a_g, b_{g'}\in \GS_p(n)$, the neighborhood of $a_gb_{g'}$ is 
$$
N(a_gb_{g'})=\{c_{g''}: g''\in A(p,n)-g-g'\}.
$$
Clearly $|N(a_gb_{g'})|=|A(p,n)|$. Observe that 
$$
|A(p,n)|=\sum_{i=1}^n |H_i|=\sum_{i=1}^n p^{n-i}=N\sum_{i=1}^np^{-i}=N(\frac{1}{p-1}-\frac{p^{-n}}{p-1})\geq \frac{1}{p-1}p^n-\frac{1}{p-1}.
$$
Thus $\frac{1}{p-1}p^n-1\leq |A(p,n)|=|N(a_gb_{g'})|\leq \frac{1}{p-1}p^n$.   Then for any $a_g$, we compute the size of its neighborhood as $|N(a_g)|=\sum_{g'\in \F_p^n}N(a_gb_{g'})=p^n|A(p,n)|$.   Clearly any $b_g$ or $c_g$ has the same degree.  
\end{proof}

It is also important to note that we showed in Claim \ref{cl:gsdegree} that the degree of every vertex has the form $\alpha N^2$ for some $\alpha$ bounded away from $0$ and $1$ (e.g. $\alpha=1/2p$ would work).  We now consider some important interactions between the metric and hypergraph structure.  Equip each of the sets $A, B, C$ with the natural metrics coming from Lemma \ref{lem:gsmetric} (i.e. for $g,g'\in \F_p^n$, set $d(b_g,b_{g'})=d(c_g,c_{g'})=d(a_g,a_{g'})=d(g,g')$).  In our next claim, we use the vertex-regularity in $\GS_p(n)$ to show that given open balls $B_1\subseteq A$ and $B_2\subseteq B$, there are many edges and nonedges in $K_3[B_1,B_2,C]$.  This can be thought of as a beefed-up version of the fact that for every $a_gb_{g'}$, $N(a_gb_{g'})$ and $\neg N(a_gb_{g'})$ are large  (see Claim \ref{cl:gsdegree}).

\begin{claim}\label{cl:gsdegreeball}
For all $a_g\in A$, $b_{g'}\in B$, and $r\in (0,1]$, there are elements $f_0(r,g,g')$ and $f_1(r,g,g')$ of $\F_p^n$  such that 
\[K_3[B_r(a_g),B_r(b_{g'}), B_r(c_{f_1(r,g,g')})]\subseteq E\]
and
\[K_3[B_r(a_g),B_r(b_{g'}), B_r(c_{f_0(r,g,g')})]\cap E=\emptyset.\]
\end{claim}
\begin{proof}
Set $m=\lceil \log_p(r^{-1})\rceil$, so $B_r(a_{g})=B_{p^{-m}}(a_g)$ and $B_r(b_{g'})=B_{p^{-m}}(b_g')$.  Let $h,h'\in \F_p^n$ be such that in $\F_p^n$, $B_{p^{-m}}(g)=H_{m+1}+g$ and $B_{p^{-m}}(g')=H_{m+1}+g'$.  Define $s=\min\{m, \lambda(g,g')\}$.  

Let $f_0(r,g,g')=-\sum_{i=1}^{s+1}(g_i+g_i')e^i+2e^{s+1}$ and let $f_1(r,g,g')=-\sum_{i=1}^{s+1}(g_i+g_i')e^i+e^{s+1}$.  Note $H_{s+1} +g+g'+f_0(g,g')=H_{s+1}+2e^{s+1}$ is disjoint from $A(p,n)$ and $H_{s+1}+g+g'+f_1(g,g')=H_{s+1}+e^{s+1}$ is contained in $A(p,n)$.  Thus $K_3[B_r(a_g),B_r(b_{g'}), B_s(c_{f_1(r,g,g')})]\subseteq E$ and $K_3[B_r(a_g),B_r(b_{g'}), B_s(c_{f_0(r,g,g')})]\cap E=\emptyset$, as desired.  Note that by construction, 
\[d(f_0(r,g,g'),f_1(r,g,g'))\geq \max\{r, d(g,g')\},\]
which completes the proof of the claim.
\end{proof}

It will be important for us to know that the $f_0(r,g,g')$, $f_1(r,g,g')$, $f_0(r',g'',g''')$, $f_1(r',g'',g''')$ in Claim \ref{cl:gsdegreeball} are far apart for certain choices of $g,g',g'',g''',r,r'$. 

\begin{claim}\label{cl:gsdegreeballmove}
If $g'\neq g''$, then for each $u\in \{0,1\}$, $d(c_{f_u(r,g,g')}), c_{f_u(r,g,g'')})=d(b_{g'},b_{g''})$.  Further, if  $h\neq h'$, $r'\geq r$, and $d(h,h')\geq \max\{d(g,g'),r'\}$.  Then
$$
\min\{d(c_{f_u(r,g,h)}, c_{f_v(r',g',h')}): u,v\in \{0,1\}\geq d(h,h').
$$
\end{claim}
\begin{proof}
Note that if $g''\neq g'$, then 
\begin{align*}
\lambda(c_{f_0(r,g,g')}, c_{f_0(r,g,g'')})=\lambda(-\sum_{i=1}^{s+1}(g_i+g_i')e^i+2e^{s+1},-\sum_{i=1}^{s'+1}(g_i+g_i'')e^i+2e^{s'+1})),
\end{align*}
where $s'=\min\{m, \lambda(g,g'')\}$.  This shows $d(c_{f_0(g,g')}, c_{f_0(g,g'')})=d(b_{g'},b_{g''})$.  A similar argument shows that $d(c_{f_1(g,g')}, c_{f_1(g,g'')})=d(b_{g'},b_{g''})$.

Suppose that $g\neq g'$, $h\neq h'$, $r'\geq r$, and that $d(h,h')\geq \max\{d(g,g'),r'\}$.  Set $m'=\lceil \log ((r')^{-1})\rceil$, $s=\min\{m,\lambda(g,h)\}$ and $s'=\min\{m', \lambda(g',h')\}$.  By assumption,  $\lambda(h,h')\leq \min\{\lambda(g,g'),m'\}$.  Given $u\in \{0,1\}$ let $u\dotplus 1$ be $2$ if $u=1$ and $1$ if $u=0$.  Now fix $(u,v)\in \{0,1\}^2$ and consider $\lambda(f_u(r,h,g), f_v(r', h',g'))$.  This is 
\[\lambda(-\sum_{i=1}^{s+1}(h_i+g_i)e^i+(u\dotplus 1)e^{s+1}, -\sum_{i=1}^{s'+1}(g'_i+h_i')e^i+(v\dotplus 1)e^{s'+1}),\]
which equals
\[\lambda(-\sum_{i=1}^{s}2h_i+(u\dotplus1-h_{s+1}-g_{s+1})e^{s+1}, -\sum_{i=1}^{s'}2h'_i+(v\dotplus1-h_{s'+1}'-g_{s'+1}')e^{s'+1}).\]
By definition of $s$ and $s'$, and because $\lambda(h,h')\leq \min\{\lambda(g,g'),m'\}$, this is nothing other than $\lambda(h,h')$.
Thus $d(c_{f_u(r,h,g)}, c_{f_v(r', h',g')})\geq d(h,h')$, as desired.
\end{proof}

We now show that the intersections of certain neighbhorhoods and non neighborhoods are contained in open balls.

\begin{claim}\label{cl:gsint}
Given $a_g,b_{g'},b_{g''}\in V(\GS_p(n))$, $N(a_g,b_{g'})\cap \neg N(a,b_{g''})$ is contained in an open ball in $C$ ball of radius $d(b_{g'},b_{g''})$.
\end{claim}
\begin{proof}
 Given $a\in \mathbb{F}_p^n$ and $\alpha\in \mathbb{F}_p$, define $\tau_i^\alpha(a):=\alpha e^i-a$, and note that 
$$
A(p,n)-a=\bigcup_{i=1}^n H_i+\tau_i^1(a).
$$
Fix $a_{g},b_{g'},b_{b''}$ and consider $C'=N(a_g,b_{g'})\cap \neg N(a,b_{g''})$.  Let $h=g+g'$ and $f=g+g''$.  Then let $C'=\{z_{w}: w\in (A(p,n)-h)\cap (\neg A(p,n)-f)\}$, and set $m=\lambda(h,f)$. Note
$$
(A(p,n)-h)\cap (\neg A(p,n)-f)=\Big(\bigcup_{i=1}^nH_i+\tau_i^1(h)\Big)\cap \Big(\bigcup_{i=1}^nH_i+\tau_i^0(f)\Big).
$$
We prove in \cite[Lemma A.7]{Terry.2021a} that the set above is contained in a single coset $K$ of $H_{m-1}$. Consequently, $C'\subseteq \{c_{v}: v\in K\}$.  We are now done since any coset of $H_{m-1}$  is an open ball of radius $p^{-m}=p^{-\lambda(h,f)}=d(b_{g'},b_{g''})$.
\end{proof}

Our last lemma shows that sufficiently small balls around distinct points can be ``split apart'' using the edge relation in $\GS_p(n)$.

\begin{claim}\label{cl:gssplit}
Suppose $a_g,a_{g'}\in A$, and $d(a_g,a_{g'})=r>0$. Then for all $b_{g''}\in B$, there is $c_h\in C$ such that 
\[K_3[B_{r}(c_h),B_{r}(b_{g''}),B_{r}(a_{g})]\subseteq E\] 
and 
\[K_3[B_{r}(c_h),B_{r}(b_{g''}),B_{r'}(a_{g'})]\cap E=\emptyset,\] 
or vice versa.  
\end{claim}
\begin{proof}
We begin by defining an element $f(a,b)\in \F_p^n$ for all distinct $a,b\in \F_p$ with the property that for all $a\neq b$ in $\F_p$, 
$$
\{a+f(a,b), b+f(a,b)\}=\{1,u\}\text{ for some $u\notin\{ 0,1\}$}.
$$
Given $0<a<b$, let $f(a,b)=a-1$. Then set $f(0,p-1)=f(p-1,0)=2$, for each $b\notin \{0,p-1\}$, and set $f(0,b)=f(b,0)=1$.  Finally, set $f(1,b)=f(b,1)=0$ for all $b\notin \{0,1\}$.  
 
Now suppose  $a_g,a_{g'}\in A$, and $d(a_g,a_{g'})=r>0$ and $b_{g''}\in B$.  Let $m=\lambda(g,g')$, so $r=p^{-m}$.  Note $g_i=g_i'$ for all $i\leq m$, by definition.  Set 
\begin{align*}
h=(\sum_{i=1}^m-g_i-g''_i)+f(g_{m+1}+g''_{m+1}, g'_{m+1}+g''_{m+1})e^{m+1}.
\end{align*}
Then 
\[H_{m+1}+g+g''+h=H_{m+1}+(f(g_{m+1}+g''_{m+1}, g'_{m+1}+g''_{m+1})+g_{m+1}+g''_{m+1} )e^{d+1}\]
and
\[H_{m+1}+g+g''+h=H_{m+1}+(f(g_{m+1}+g''_{m+1}, g'_{m+1}+g''_{m+1})+g'_{m+1}+g''_{m+1} )e^{d+1}.\]
By our definition of $f$, one of these has the form $H_{m+1}+e^{m+1}$ and one of them has the form $H_{m+1}+ue^{m+1}$ for some $u\notin \{0,1\}$.  Thus one is completely contained in $A(p,n)$, and the other is completely disjoint from $A(p,n)$.  Consequently, $K_3[B_{r}(c_h),B_{r}(b_{g''}),B_{r}(a_{g})]\subseteq E$ and $K_3[B_{r}(c_h),B_{r}(b_{g''}),B_{r'}(a_{g'})]\cap E=\emptyset$, or vice versa.  
\end{proof}

The existence of a metric which satisfies Claims \ref{cl:gsball}, \ref{cl:gsdegree}, \ref{cl:gsdegreeball}, \ref{cl:gsdegreeballmove}, \ref{cl:gsint}, \ref{cl:gssplit} is sufficient to show that $\GS_p(n)$ requires non-binary $\vdisc_3$-error.  However, we note that this does require a proof.  We do not provide this proof here, as our goal is instead to prove a more general result, which will apply to both $\GS_p(n)$ and $\HP(N)$.  In the next subsection, we will show that  $\HP(N)$ satisfies claims similar to Claims \ref{cl:gsball}, \ref{cl:gsdegree}, \ref{cl:gsdegreeball}, \ref{cl:gsdegreeballmove}, \ref{cl:gsint}, \ref{cl:gssplit}, although some weakening will be required.   Based on these,  we will later define a single set of axioms, which hold in both $\GS_p(n)$ and $\HP(N)$, and which imply that a property requires non-binary $\vdisc_3$-error. 

\subsection{Special properties of $\HP(N)$}\label{subsub:hp}  In this subsection we show that $\HP(N)$ satisfies weakened versions of the properties of the preceding subsection.  We begin by equipping $\HP(N)$ with a metric in the obvious way.  Recall that
$$
\HP(N)=(\{a_i, b_i, c_i: i\in [N]\},  \{a_ib_jc_k: i+j+k\geq N+2\}).
$$
To ease notation, in Section \ref{subsub:hp}, we let $A=\{a_i: i\in [N]\}$, $B=\{b_i: i\in [N]\}$, $C=\{c_i: i\in [N]\}$, and $E=E(\HP(N))$.

\begin{lemma}\label{lem:hpmetric}
For each $X\in \{A,B,C\}$ and $x_i,x_j\in X$, set $d(x_i,x_j)=\frac{|i-j|}{N}$.  This makes each of $A,B,C$ into a metric space with distances in $[0,1]$.
\end{lemma}

In an ideal world, we would like to prove all the same facts here for $\HP(N)$ as we showed held for $\GS_p(n)$ in the previous subsection.   Indeed, some of the facts will go through largely unchanged.  For example, the following is a straightforward analogue of Claim \ref{cl:gsint}.

\begin{claim}\label{cl:hpint}
Given $a_i\in A$ and $b_j,b_{j'}\in B$, $N(a_i,b_{j'})\cap \neg N(a_i,b_{j'})$ is contained in a closed ball in $C$ ball of radius at most $d(b_{j'},b_{j'})$.
\end{claim}
\begin{proof}
Without loss of generality, assume $j'<j$.  Note that
$$
N(a_i,b_j)\cap \neg N(a_{i},b_{j'})=\{c_s: s\in [N+2-i-j, N]\cap [0, N+2-i-j')\}.
$$
If this is empty, let $C$ be the ball of radius $0$ centered at $c_1$.  Otherwise, the above is equal to 
$$
\{c_s: s\in [N+2-i-j, N+2-i-j')\},
$$
 which has size $|j-j'|=d(b_j,b_{j'})N$.  This set is clearly contained in a closed ball of radius $ |j-j'| /2N$ in $C$.
\end{proof}

Other properties from Section \ref{subsub:gs} must be tweaked in a more serious way before they will hold in $\HP(N)$.  For instance, in $\HP(N)$, a ball $B_r(a_i)$ in $A$ is of the form 
$$
\{a_j : j\in (\min\{i-rN, 1\}, \max\{i+rN, N\})\}.
$$
This could have size about $rN/2$ (e.g.  if  $i=N$), or it could have size $2rN$ (e.g. if  $i=\lfloor N/2\rfloor$, and $r<1/2$).  Consequently, not all balls with the same radius contain the same number of vertices.  However, we can clearly give a range of $\lfloor \frac{r}{2}N\rfloor \leq |B_r(a_i)|\leq 2rN$.  

Similarly, not every vertex in $\HP(N)$ has the same degree.  For instance, if $i,j\geq N/2$, then $|N(a_ib_j)|=N$, while if $i=j=1$, then $|N(a_ib_j)|=1$.  However, we can show that most vertices have degree within a certain non-trivial range, and that many (but not all) pairs of vertices are in many edges and non-edges.

\begin{claim}\label{cl:hpdegree}
Suppose $0<\tau<1/4$ and $0<3\mu<\tau$.  The following hold in $\HP(N)$.
\begin{enumerate}[label=\normalfont(\arabic*)]
\item For all $i\in (\mu N, 2\mu N)$ and $j\in (\tau N,(1-\tau)N)$, $|N(a_ib_j)| \in [\tau N/2, (1-\tau /2)N]$.
\item For all $i\in (\mu N, 2\mu N)$, $\min\{|N(a_i)|, |\neg N(a_i)|\}\geq N^2/4$.
\item For all $j\in (\tau N,(1-\tau)N)$, $|N(b_j)|\in [\mu\tau N^2/4, (1-\mu\tau/4)N^2]$.
\end{enumerate}
\end{claim}
\begin{proof}
Given $i\in (\mu N, 2\mu N)$ and $j\in (\tau N,(1-\tau)N)$, note 
$$
N(a_ib_j)=\{c_k: k\in [N+2-i-j,N]\}.
$$ 
 Therefore, 
$$
\tau N/2 <\mu N+\tau N-2\leq |N(a_ib_j)|=i+j-2<N(1-\tau+2\mu)-2<N(1-\tau/2).
$$
This shows (1) holds.  Further, for all $j\in (\tau N, (1-\tau)N)$, we have 
$$
|N(b_j)|\geq |(\mu N, 2\mu N)|\tau N/2\geq \mu \tau N^2/4,
$$
and $|N(b_j)|\leq N^2-|(\mu N, 2\mu N)|\tau N/2\leq N^2 -\mu\tau N^2/4$, so (3) holds.  For (3), we have that for all $i\in (\mu N, 2\mu N)$, 
$$
N^2/4\leq \sum_{j=1}^N\sum_{k=N+2-j-2\mu N}k\leq \sum_{j=1}^N \sum_{k=N+2-j-i}^Nk\leq  |N(a_i)|\leq \sum_{j=1}^N\sum_{k=1}^{N-j}k\leq N^2/2
$$  
as claimed.
\end{proof}

Claim \ref{cl:hpdegree} shows that almost all vertices (namely those with index in $(\tau N,(1-\tau)N)$) have degree $\alpha N^2$ for some $\alpha$ bounded away from $0$ and $1$.  Further, there are many pairs of vertices with degree $\beta N$ for some $\beta$ bounded away from $0$ and $1$ (namely those where one vertex has index in $(\mu N, 2\mu N)$ and the other has index in $(\tau N,(1-\tau)N)$).  

This lemma, and the observations preceding it, illustrate a general strategy for how to deal with $\HP(N)$ (as compared to $\GS_p(n)$).  In particular, we must make the statements more approximate, and also relativize some of the quantifications to certain  ``special small index sets" (e.g. the vertices with indices in $(\mu N, 2\mu N)$) and ``special non-small index sets" (e.g. the vertices with indices in $(\tau N, (1-\tau)N)$).   To ease notation and emphasize this idea, we define, for $\tau, \mu>0$, 
\[A_{lg}(\tau)=\{a_i: i\in (\tau N,(1-\tau)N)\}\;\text{ and }\;A_{sm}(\mu)=\{a_i: i\in (\mu N,2\mu N)\}.\]
We define $B_{lg}(\tau)$, $B_{sm}(\mu)$, $C_{lg}(\tau)$, and $B_{sm}(\mu)$ in the same way.  We will always use these notions with $\mu$ being significantly smaller than $\tau$, in which case $A_{lg}(\tau)\cap A_{sm}(\mu)=\emptyset$.  Note that Claim \ref{cl:hpdegree}(1) says that every $ab\in A_{sm}(\mu)\times B_{lg}(\tau)$ has degree $\beta N$ for some $\beta$ bounded away from $0$ and $1$, and Claim \ref{cl:hpdegree}(2) says that every $b\in B_{lg}(\tau)\cup B_{sm}(\mu)$ has degree of $\alpha N^2$, for some $\alpha$ bounded away from $0$ and $1$.  

When it comes to partitioning open balls into smaller ones (i.e. proving an analogue of Claim \ref{cl:gsball}), things are again more complicated.  For instance, unlike in the ultrametric space $\GS_p(n)$, we cannot partition an open ball $B_r(a_i)$ into disjoint open sub-balls.  This is because any open ball has the form $\{a_j: j\in I\}$, for some open interval $I$. Any attempt to partition $I$ into disjoint open intervals will always miss a few points of $I$, namely the endpoints of the sub-intervals.  Problems can also arise due to divisibility issues near the end points of intervals. The following lemma will be used to deal with these issues, and says that we can almost partition an interval into almost evenly sized sub-intervals.  Below, $d_N$ denotes the usual metric $d_N(i,j)=|i-j|$ on $[N]$.  A ball in this metric space has the form $(x-d, x+d)\cap [N]$ for some $x\in [N]$ and $d\in [N]$.  

\begin{lemma}\label{lem:divide}
For all $0<r_1<1$, there is $n_0$ so that if $N\geq n_0$, the following holds.  Suppose $0<r_1<r_2$, and $I=(\alpha,\beta)$ is an open ball of radius $r_2N$ in $[N]$. Then there is $m\leq 4r_2/r_1$ and a collection $\calP$ of $m$ open balls contained in $I$, each of radius at most $r_1N$ and at least $r_1N/3$, such that $|I\setminus (\bigcup \calP)|\leq 2m$.
\end{lemma}
\begin{proof}

Assume $N\gg r_1^{-1},r_2^{-1}$.  Set $d_2=\lceil r_2 N\rceil$ and $d_1=\lceil r_1 N/3\rceil$.  Let $i\in [N]$ be such that 
$$
I=(\alpha,\beta)=(i-d_2,i+d_2)\cap [1,N].
$$
  Note $d_2+1\leq |I|=\beta-\alpha+1\leq 2d_2+1$.  Set $s=\lfloor |\beta-\alpha|/(2d_1)\rfloor$, and for each $j\in [s]$, set $x_j=\alpha+(2j-1)d_1$.  Note that for each $j\in [s]$, $x_j+d_1<\beta$, so the ball of radius $d_1$ around $x_j$ is contained in $I$.  However, $I\setminus (\bigcup_{j=1}^sB_{d_1}(x_j))$ could be large if $(x_s+d_1, \beta)$ is large.  To make sure this is not an issue, we will combine the last two balls into one ball of a slightly different radius.

By definition of $x_{s-1}$, $2d_1<\beta-(x_{s-1}-d_1)\leq 4d_1$. Set $d_1'=\lfloor (\beta-(x_{s-1}-d_1))/2\rfloor$, and let $x_{s-1}'=x_{s-1}+d_1'$.  Set
$$
\calP=\{B_{d_1}(x_i): 1\leq i\leq s-1\}\cup \{B_{d_1'}(x'_{s-1})\}.
$$
By construction, every element in $\calP$ is an open ball of radius at most $r_1N$ and at least $r_1N/3$, and 
$$
|I\setminus (\bigcup \calP)|\leq |\{\alpha,x_1,\ldots, x_{s-2},x_{s-1}'\}\cup \{\beta,\beta-1\}|\leq s+2.
$$
By definition $s\leq (2r_2N+2)/2d_1N\leq (2r_2N+2)/(2r_1N/3)\leq 4r_2/r_1$, where the last inequality uses that $N$ is large.
\end{proof} 

We can now give a weak analogue of Claim \ref{cl:gsball}.  

\begin{claim}\label{cl:hpcover}
Suppose $0<r_1<\mu^2$, $r_1<r_2\leq 1$, $N$ is sufficiently large compared to $r_1$, and $a\in A$.  Then there is some $m\leq 4r_1/r_2$ and a set $\calP$ of $m$ disjoint open balls, each of radius at least $r_1/3$ and at most $r_1$, all contained in $B_{r_2}(a)$, such that $|B_{r_2}(a)\setminus (\bigcup \calP)|\leq 2m$. 
\end{claim}
\begin{proof}
Note $B_{r_2}(a)=\{a_i: i\in I\}$ for some open ball $I$ in $[N]$.  The conclusion is thus an immediate corollary of Lemma \ref{lem:divide}.
\end{proof}

We recall that in Claim \ref{cl:gsball}, there was an additional conclusion beyond what we have shown above, namely that we were also able to choose one of the balls in the cover ahead of time.  Due to divisibility issues, we cannot always do this in $\HP(N)$, but it turns out that Claim \ref{cl:hpcover} along with Claim \ref{cl:hpball1} below will suffice. In particular, Claim \ref{cl:hpball1} says that given an open ball $B_r(a)$ in $A$, whose center is not too close to the endpoints of $[N]$ (where problems arise), we can almost partition $A$ into balls, one of which is $B_r(a)$.  To help deal with this requirement, we define
\[A_{sm}(\mu)^+=\{a_i: i\in ((\mu-\mu^2)N, (\mu+\mu^2)N)\}\]
and
\[A_{lg}(\tau)^+=\{a_i: i\in ((\tau-\mu^2)N, (1-\tau+\mu^2)N)\},\]
which are slightly enlarged versions of our special sets. We similarly define the sets $B_{sm}(\mu)^+,B_{lg}(\mu)^+$ and $C_{sm}(\mu)^+,C_{lg}(\mu)^+$. 

\begin{claim}\label{cl:hpball1}
Assume $0<\tau<1/4$, $0<\mu<\tau/3$, and $0<r<\mu^2$.  For $N$ sufficiently large, the following holds.  For all $a_i\in A_{lg}(\tau)^+\cup A_{sm}(\mu)^+$, there are $m\leq 8r^{-1}$ and a set $\calP$ of $m$ disjoint open balls in $A$, each of radius at most $r$ and at least $r/3$, such that  $B_r(a_i)\in \calP$ and  $|A\setminus (\bigcup\calP)|\leq 2m$.  
\end{claim}
\begin{proof}
Assume $N\gg r^{-1},\mu^{-1}$.  Suppose $a_i\in A_{lg}(\tau)^+\cup A_{sm}(\mu)^+$ and $0<r<\mu^2$. By assumption, $i\in ((\tau-\mu^2) N, (1-\tau+\mu^2)N)\cup ((\mu-\mu^2)N, (2\mu+\mu^2)N)$.   Let $d=\lceil rN\rceil$.  Note that $B_r(a_i)$ has the form $\{a_i: i\in I\}$ for some open ball $I=(i-d,i+d)\cap [N]$ in $[N]$.  Since $a_i\in A_{lg}^+(\tau)\cup A_{sm}^+(\mu)$ and $r<\mu^2$, we know that $i-d>\mu-2\mu^2$ and $i+d<\tau+2\mu^2$.  Consequently, both $i-d$ and $i+d$ are in $[N]$.  Let $I_1=[1,i-d)$ and $I_2=(i+d,N]$.  These can also be written as open balls in $[N]$, say of radii $d_1$ and $d_2$, respectively.  Note that $d_1,d_2\geq (\mu-2\mu^2)N/2$. 

By Lemma \ref{lem:divide} there are $m_1\leq 4\frac{r_1}{N}r^{-1}$ and a collection $\calP_1$ of open balls in $([N],d_N)$, each of radius at least $rN/3$ and at most $rN$ such that each element of $\calP_1$ is contained in $(i+d,N)$ and $|(i+d,N)\setminus (\bigcup\calP_1)|\leq 2m_1$.   Similarly, there are $m_2\leq 4r_2r^{-1}$  and a collection $\calP_1$ of open balls in $([N],d_N)$, each of radius at least $rN/3$ and at most $rN$ such that each element of $\calP_2$ is contained in $(1,i-d)$ and $|(1,i-d)\setminus (\bigcup\calP_2)|\leq 2m_2$. We then let 
$$
\calP=\{B_r(a_i)\}\cup \calP_1\cup \calP_2.
$$
Note $|\calP|=m_1+m_2+1\leq 4r^{-1}$.  By construction, each element of $\calP$ is a ball of radius at least $r/3$ and at most $r$, and $|A\setminus (\bigcup \calP)|\leq 2m_1+2m_2+2\leq 8r^{-1}$, as desired.
\end{proof}

We will need one more partition fact, namely that given a ball $B_r(a_i)$ whose center $a_i$ is in $A_{sm}(\mu)$, the approximate subcover obtained in Claim  \ref{cl:hpcover} will have the property that $B_r(a_i)\cap A_{sm}(\mu)$ is mostly covered by balls whose centers are in $A_{sm}(\mu)^+$.  
 
\begin{claim}\label{cl:hpcover2}
In Claim \ref{cl:hpcover}, if $\calP'$ are the elements of $\calP$ centered at elements in $A^+_{sm}(\mu)$, and $\calP''$ are those centered at $A_{lg}(\tau)^+$, then 
\[\max\Big\{|(B_{r_2}(a)\cap A_{sm}(\mu))\setminus (\bigcup\calP')|, |(B_{r_2}(a)\cap A_{lg}(\tau))\setminus (\bigcup\calP'')|\Big\}\leq 2m.\]
\end{claim}
\begin{proof}
This is immediate from the definitions of $A_{sm}(\mu)^+$ and $A_{lg}(\tau)^+$.
\end{proof}

We now turn to the analogues of Claims \ref{cl:gsdegreeball} and \ref{cl:gsdegreeballmove}.  For these we must restrict ourselves to our special sets $A_{sm}(\mu)^+$ and $A_{lg}(\mu)^+$.

\begin{claim}\label{cl:hpdegreeball}
Suppose $0<\tau<1/4$, $0<\mu<\tau/3$, and $0<r<\mu^2$.  For all sufficiently large $N$ the following holds.  For every $a\in A_{lg}(\tau)^+$ and $b\in B_{sm}(\mu)^+$, there are $f_0(r,a,b),f_1(r,a,b)\in [N]$ such that 
\[K_3[B_{r/2}(c_{f_0(r,a,b)}), B_r(a), B_{r}(b)]\subseteq E\]
and 
\[K_3[B_{r/2}(c_{f_1(r,a,b)}), B_r(a), B_{r}(b)]\cap E=\emptyset,\]
and such that $d(c_{f_0(r,a,b)},c_{f_1(r,a,b)})\leq 7r$.  
\end{claim}
\begin{proof}
Suppose $(\tau-\mu^2)N< i<(1-\tau+\mu^2)N$, $(\mu-\mu^2)N<j< (2\mu+\mu^2)N$, and $0<r<\mu^2$.  Set $d_1=\lceil rN\rceil$, and define 
\begin{align*}
f_0(r,a_i,b_j)=N+2-i-j-3d_1\text{ and }f_1(r,a_i,b_j)=N+2-i-j+3d_1.
\end{align*}
  Note each of these are in $[N]$ due to the assumptions on $i$, $j$, and $r$.  Suppose $a_wb_uc_v\in B_r(a_i)\times B_r(b_j)\times B_{r/2}(c_{f_0(r,a_i,b_j)})$.  Then
\begin{align*}
w+u+v&<(i+d_1)+(j+d_1)+(N+2-i-j-3d_1)+\lceil d_1/2\rceil\\
&=N+2-d_1+\lceil rN/2\rceil \\
&<N+2.
\end{align*}
On the other hand, if $a_wb_uc_v\in B_r(a_i)\times B_r(b_j)\times B_{r/2}(c_{f_1(r,a_i,b_j)})$.  Then 
\begin{align*}
w+u+v&> (i-\lfloor rN\rfloor)+(j-\lfloor rN\rfloor)+(N+2-i-j+3d_1)-\lfloor rN/2\rfloor\\
&\geq N+d_1-\lfloor rN/2\rfloor\\
&>N+2,
\end{align*}
where the last inequality uses the $N$ is large. Note that by construction, 
\[d(c_{f_1(r,a_i,b_j)},c_{f_0(r,a_i,b_j)})=6d_1/N\leq 7r,\] 
which completes the proof of the claim.
\end{proof}

As in Section \ref{subsub:gs}, it is important to know the $f_u(r,a,b)$, $f_v(r',a',b')$ obtained in Claim \ref{cl:hpdegreeball} are far for certain choices of $a,b,r, a',b',r'$.

\begin{claim}\label{cl:hpdegreeballmove}
Assume $0<\tau<1/4$, $0<\mu<\tau/3$, and $0<r\leq r'<\mu^2$.  For any  $a,a'\in A_{lg}(\tau)^+$ and $b,b'\in B_{sm}(\mu)^+$ the following hold. 
\begin{enumerate}[label=\normalfont(\arabic*)]
\item For each $u\in \{0,1\}$, we have that $d(f_u(r,a,b), f_u(r,a,b'))\geq d(b,b')$ and $d(f_u(r,a',b), f_u(r,a,b))\geq d(a,a')$. 
\item If  $d(a,a')\geq d(b,b')+6r'$, then 
$$
\min\{d(c_{f_u(r,a,b)}, c_{f_v(r',a',b')}): u,v\in \{0,1\}\}\geq d(a,a')-d(b,b')-7r'.
$$
\item If $d(b,b')\geq d(a,a')+6r'$, then 
$$
\min\{d(c_{f_u(r,a,b)}, c_{f_v(r',a',b')}): u,v\in \{0,1\}\}\geq d(b,b')-d(a,a')-7r'.
$$
\end{enumerate}
\end{claim}
\begin{proof}
Suppose $a_i,a_{i'}\in A_{lg}(\tau)^+$ and $b_j,b_{j'}\in B_{sm}(\mu)^+$.  By definition, for each $u\in \{0,1\}$, 
$$
d(f_u(r,a_i,b_j),f_u(r,a_{i'},b_j))=|i-i'|. 
$$
Thus $d(c_{f_u(r,a_i,b_j)},d_{f_u(r,a_{i'},b_j)})=|i-i'|/N=d(a,a')$.  A symmetric argument shows that $d(f_u(r,a,b), f_u(r,a,b'))= d(b,b')$, so (1) holds.  For (2), suppose that $d(b_j,b_{j'})\geq d(a_i,a_{i'})+6r'$. Without loss of generality, say $j'>j$.  Then for any $(u,v)\in \{0,1\}^2$, $d(f_u(r,a,b), f_v(r',a',b'))$ equals
\[ |N+2-i-j+(-1)^u3\lceil rN\rceil-(N+2-i-j+(-1)^v3\lceil r'N\rceil)|,\]
which, since $d(b_j,b_{j'})\geq d(a_i,a_{i'})+6r'$ and $r'\geq r$, is bounded below by
\[|j'-j|-|i'-i|-6\lceil r'N\rceil.\]
It follows that
\begin{align*}
d(c_{f_u(r,a_i,b_j)}, c_{f_v(r',a_{i'}, b_{j'})})&\geq |j'-j|-|i'-i|-6\lceil r'N\rceil/N\\&\geq d(b_j,b_{j'})-d(a_i,a_{i'})-7r'.
\end{align*}
The argument for (3) is similar. 
\end{proof}

We now give the analogue of Claim \ref{cl:gssplit}.  The version here is a bit weaker because the balls are smaller, and we will also require things to be sufficiently spaced out.
\begin{claim}\label{cl:hpsplit}
Suppose $0<\tau<1/4$, $0<3\mu<\tau$, $0<r,r'<\mu^2$ and $r''\geq \max\{r,r'\}$.  Then if $N$ is sufficiently large, the following hold.  

For all $a_i\in A^+_{lg}(\tau)$ and $a_{i'}\in A$ with $d(a,a')=r''+r'+r$, there is $b\in B_{sm}(\mu)$ and $c\in C$ such that if $r_0=\min\{r,r',r''\}/3$, then $K_3[B_{r_0/2}(c),B_{r_0/2}(b),B_{r'}(a_{i'})]\subseteq E$ and $K_3[B_{r_0/2}(c),B_{r_0/2}(b),B_{r}(a_i)]\cap E=\emptyset$, or vice versa.  
\end{claim}
\begin{proof}
Fix $a_i\in A_{lg}^+(\tau)$ and $a_{i'}\in A$ satisfying $d(a,a')=r''+r'+r$, where $r''>2\max\{r',r\}$ and $0<r,r'<\mu^2$.  Set $r_0=\min\{r,r',r''\}/3$, and choose any $b_j\in B_{sm}(\mu)$ (note this implies $\mu N< j< 2\mu N$).

Suppose first that $i'\in ((\tau-\mu^2)N, N(1-\tau+\mu^2))$ and  $i \leq i'$. Let $k=N+2-i'-j+\lceil r'N\rceil+\lceil r_0N\rceil$.  Our assumptions imply $i'+j\leq (1-\tau+\mu)N$ so $k\geq N+2-i'-j\geq (\tau-\mu^2-\mu)N\geq 1$.  On the other hand, 
$$
k\leq N+2-(\tau -\mu^2)N-\mu N+2+r'N+r_0N\leq N,
$$
where the last inequality is because $r',r_0<\mu^2$ and $N$ is large.  Thus $k\in [N]$, so it makes sense to define $c=c_k$.  Suppose $b_uc_v\in B_{r_0/2}(b_j)\times B_{r_0/2}(c)$, $a_s\in B_{r'}(a_{i'})$ and $a_t\in B_r(a_i)$.  Then using that $i'-i=r+r'+r''$ by assumption, we have
\begin{align*}
t+u+v&< (i+rN)+(j+r_0N/2)+(N+2-i'-j+\lceil r'N\rceil+\lceil r_0N\rceil +r_0N/2)\\
&=N+2 +r_0N-(i'-i)+\lceil r'N\rceil+\lceil r_0N\rceil\\
&\leq N+7-(r+r'+r'')N-r'N+r_0N\\
&<N+2,
\end{align*}
where the last inequality is because $r''>2\max\{r,r'\}$ and $N$ is large.  Thus $a_tb_uc_v\notin E$.  On the other hand, 
\begin{align*}
s+u+v&> (i'-r'N)+(j-r_0N/2)+(N+2-i'-j+\lceil r'N\rceil+\lceil r_0N\rceil-r_0N/2)\\
&=N+2-r_0N+\lceil r'N\rceil+\lceil r_0N\rceil\\
&\geq N+2.
\end{align*}
Thus $a_sb_uc_v\in E$.  

If $i'\in ((\tau-\mu^2)N, N(1-\tau+\mu^2))$ and  $i' \leq i$, then proceed as above with the roles of $i$ and $i'$ switched.

Suppose now that $i'<(\tau-\mu^2)N$.   In this case, set $k=N+2-i-j+\lceil rN\rceil+\lceil r_0N\rceil$.  It is straightforward to check that $k\in [N]$ since $i\in ((\tau-\mu^2)N, (1-\tau+\mu^2)N)$ and $j\in (\mu N, 2\mu N)$. Thus it makes sense to define $c=c_k$.  Suppose $b_uc_v\in B_{r_0/2}(b_j)\times B_{r_0/2}(c)$, $a_s\in B_{r'}(a_{i'})$ and $a_t\in B_r(a_i)$. Then
\begin{align*}
t+u+v&>(i-rN)+(j-r_0N/2)+(N+2-i-j+\lceil rN\rceil +\lceil r_0N\rceil-r_0N/2)\\
&=N+2-rN+\lceil rN\rceil -r_0N+\lceil r_0N\rceil\\
&\geq N+2.
\end{align*} 
Thus $a_tb_uc_v\in E$.  On the other hand
\begin{align*}
s+u+v&<(i'+r'N)+(j+r_0N/2)+(N+2-i-j+\lceil rN\rceil +\lceil r_0N\rceil+r_0N/2)\\
&=N+2-(i-i')+r_0N+\lceil rN\rceil +\lceil r_0N\rceil\\
&\leq N+2-\lfloor (r+r'+r'')N\rfloor +\lceil rN\rceil +\lceil r_0N\rceil\\
&\leq N+2,
\end{align*}
where the last inequality is since $N$ is large and $r''>2\max\{r,r'\}\geq 6r_0$.  Thus $a_sb_uc_v\notin E$. 

Finally, consider the case where $i'>N(1-\tau+\mu^2)$.  In this case, set 
$$
k=N+2-i-j-\lceil rN\rceil-\lceil r_0 N\rceil.
$$
As above, it is straightforward to check that $k\in [N]$.  Thus it makes sense to define $c=c_k$.  Suppose $b_uc_v\in B_{r_0/2}(b_j)\times B_{r_0/2}(c)$, $a_s\in B_{r'}(a_{i'})$ and $a_t\in B_r(a_i)$. Then
\begin{align*}
s+u+v&>(i'-r'N)+(j-r_0N/2)+(N+2-i-j-\lceil rN\rceil -\lceil r_0N\rceil-r_0N/2)\\
&=N+2+(i'-i)-r'N-\lceil rN\rceil -r_0N-\lceil r_0N\rceil\\
&\geq N+2+\lfloor (r+r'+r'')N\rfloor -r'N-\lceil rN\rceil -r_0N-\lceil r_0N\rceil\\
&\geq N+2,
\end{align*} 
where the last inequality is because $N$ is large and $r''>2\max\{r,r'\}\geq 6r_0$. Thus $a_sb_uc_v\notin E$.  On the other hand
\begin{align*}
t+u+v&<(i+rN)+(j+r_0N/2)+(N+2-i-j-\lceil rN\rceil -\lceil r_0N\rceil+r_0N/2)\\
&=N+2+rN-\lceil rN\rceil +r_0N-\lceil r_0N\rceil\\
&\leq N+2,
\end{align*}
thus $a_tb_uc_v\notin E$.  This finishes the proof.
\end{proof}

We will see that the properties described in this subsection are sufficient for showing that $\HP(N)$ requires non-binary $\vdisc_3$-error.  Further,  they are strictly weaker than the properties laid out in Section \ref{subsub:gs}.  For this reason, the axioms in Definition \ref{def:special} will be based on what holds in $\HP(N)$.  

\subsection{Special $3$-graphs} In this subsection we define special $3$-graphs, and show that both $\HP(N)$ and $\GS_p(n)$ are special, for appropriately chosen parameters.

\begin{definition}\label{def:special}
Suppose $p\geq 2$ and $G=(A\cup B\cup C, E)$ is a $3$-partite $3$-graph where $|A|=|B|=|C|=N$.  We say $G$ is \emph{$(p,\mu,\tau,\alpha,\rho)$-special} if there is a metric $d$ on $A$, $B$, and $C$ with distances in $[0,1]$, and distinguished subsets $A_{sm}\subseteq A_{sm}^+$ and $A_{lg}\subseteq A_{lg}^+$ of $A$, $B_{sm}\subseteq B_{sm}^+$ and $B_{lg}\subseteq B_{lg}^+$ of $B$, and $C_{sm}\subseteq C_{sm}^+$, and $C_{lg}\subseteq C_{lg}^+$ of $C$, satisfying the following axioms  for all permutations $XYZ$ of $ABC$, where given $X\in \{A,B,C\}$, $B_r(x)=\{x'\in X: d(x,x')<r\}$.  
\begin{enumerate}
\item (Special Subsets Axiom) There are $r_{sm}\geq \mu/2$, $r_{lg}\geq \tau/2$, and $x_{sm},x_{lg}\in X$ such that $X_{sm}=B_{r_{sm}}(x_{sm})$, $X_{lg}=B_{r_{lg}}(x_{lg})$, $X_{sm}^+=B_{r_{sm}+\mu^2}(x_{sm})$, and $X_{lg}^+=B_{r_{lg}+\mu^2}(x_{lg})$.
\item (Size Axiom) For each $r\in (0,1]$ and $x\in X$, $r|X|/2p\leq |B_r(x)|\leq pr|X|$. 
\item (Degree Axiom) For all $x\in X_{lg}^+\cup X_{sm}^+$, 
$$
\min\{|N(x)\cap K_2[Y,Z]|,|\neg N(x)\cap K_2[Y,Z]|\}\geq \alpha n^2.
$$
\item (Intersection Axiom) For all $y\in Y$ and $z,z'\in Z$, there is a closed ball of radius at most $d(z,z')$ containing $ N(y,z)\cap \neg N(y,z')$.
\item (Partition Axiom) For any $\rho<r<\mu^2$ and $x\in X_{lg}^+\cup X_{sm}^+$, there are $m\leq 4p(r^{-1})$ and a set $\calP$ of $m$  disjoint open balls in $X$ of radius at least $r/3$ and at most $r$, one of which is $B_r(x)$, such that  $|X\setminus (\bigcup \calP)|\leq 2m$.  
\item (Covering Axiom) Suppose $\rho<r_1<\mu^2$ and $r_1<r_2\leq 1$, and $B$ is a ball of radius $r_2$ in $X$.  Then there is some $m\leq 2pr_1/r_2$, and a set $\calP$ of $m$ disjoint open balls, each of radius at least $r_1/3$ and at most $r_1$, all contained in $B$, such that $|B\setminus (\bigcup \calP)|\leq 2m$. Moreover, the following hold.
\begin{enumerate}
\item[(a)] If $\calP'=\{B_r(x)\in \calP: x\in X_{lg}^+\}$, then $|(B\cap X_{lg})\setminus (\bigcup \calP')|\leq 2m$ 
\item[(b)]  If $\calP''=\{B_r(x)\in \calP: x\in X_{sm}^+\}$, then $|(B\cap X_{sm})\setminus (\bigcup \calP'')|\leq 2m$.
\end{enumerate}
\item (Splitting Axiom 1) For every $x\in X_{lg}^+$ and $y\in Y_{sm}^+$ and $\rho<r<\mu^2$, there are $f_0(r,x,y)$ and $f_1(r,x,y)$ such that $B_{r/2}(f_u(r,x,y))\times B_r(x)\times B_{r}(y)\subseteq E^u$ for each $u\in \{0,1\}$, and such that $d(f_0(r,x,y),f_1(r,x,y))\leq 3r$.  Further, the following hold. 
\begin{enumerate}
\item[(a)] If $y\neq y'$, then for each $u\in\{0,1\}$, $d(f_u(r,x,y),f_u(r,x,y'))\geq d(y,y')$,
\item[(b)] If $d(y,y')\geq d(x,x')+6r'$, then
$$
\min\{d(f_u(r,x,y), f_v(r',x',y')) : u,v\in \{0,1\}\}\geq d(y,y')-d(x,x')-7r'.
$$
\item[(c)] If $d(x,x')\geq d(y',y)+6r'$, then
$$
\min\{d(f_u(r,x,y), f_v(r',x',y')) : u,v\in \{0,1\}\}\geq d(x,x')-d(y,y')-7r'.
$$
\end{enumerate}
\item (Splitting Axiom 2) For all $x\in X^+_{lg}$ and $x'\in X$ satisfying $d(x,x')=r''+r'+r$, where $\rho\leq r,r'<\mu^2$ and $r''>2\max\{r',r\}$, and all $y\in Y_{sm}$ there is $z=f(x,x',y)\in Z$ such that if $r_0=\min\{r,r',r''\}/3$, then $K_3[B_{r_0/2}(y),B_{r_0/2}(z),B_{r'}(x')]\subseteq E$ and $K_3[B_{r_0/2}(y),B_{r_0/2}(z),B_{r}(x)]\cap E=\emptyset$, or vice versa.  
\item (Disjointness Axiom) The sets $X_{sm}^+$ and $X_{lg}^+$ are disjoint.
\end{enumerate}
\end{definition}

It seems likely that with a bit more work, one could get rid of Splitting Axiom 2 by deducing it from a version of Splitting Axiom. However, optimizing the list of axioms is not our main goal. We can now easily verify that our main examples are special.

\begin{proposition}\label{prop:gsspecial}
For all $p\geq 2$, there is $n_0$ such that for all $\rho>0$ and all $n\geq n_0$, $\GS_p(n)$ is $(p,1/p^2,1/p^2,\frac{1}{2(p-1)},\rho)$-special.
\end{proposition}
\begin{proof}
Equip $A$, $B$, and $C$ with the metrics  defined in Section \ref{subsub:gs}.  For each $X\in \{A,B,C\}$, set 
\begin{align*}
X_{lg}&=X_{lg}^+=\{x_g: g\in H_1+e_1\}\text{ and }X_{sm}=X_{sm}^+=\{x_g: g\in H_1+e_2\}.
\end{align*}
Clearly these are balls of radius at least $p^{-2}$ (see Claim \ref{cl:gsball}), and $X_{sm}^+\cap X_{lg}^+=\emptyset$.   Thus Axiom (1) holds. Axiom  (2) holds by the comments following \ref{lem:gsmetric}.  Axioms (3)-(7) hold immediately from Claims \ref{cl:gsball}, \ref{cl:gsdegree}, \ref{cl:gsdegreeball}, \ref{cl:gsdegreeballmove}, \ref{cl:gsint}, \ref{cl:gssplit}, and the fact that the roles of $A$, $B$ and $C$ are symmetric in $\GS_p(n)$. 
\end{proof}

\begin{proposition}\label{prop:hpspecial}
For all $0<\tau<1/4$, $0<\mu<\tau/3$, $0<\rho<\mu^3$, and $p\geq 2$, there is $n_0$ such that for all $N\geq n_0$, $\HP(N)$ is $(p,\mu,1-\tau,\mu^2,\rho)$-special.
\end{proposition}
\begin{proof}
Choose $n_0$ sufficiently large for $\tau$, $\mu$, and $r,r',r''\geq \rho$ in Claims \ref{cl:hpcover} and \ref{cl:hpsplit}.  Assume $N\geq n_0$ and equip each of $A,B,C$ in $\HP(N)$ the metrics described in Section \ref{subsub:hp}.  Set $x_{lg}=x_{\lfloor n/2\rfloor}$, $x_{sm}=x_{\lfloor 3\mu N/2\rfloor}$, $r_{sm}=\mu/2$, and $r_{lg}=(1-\tau)/2$.  Then define $X_{sm},X_{sm}^+,X_{lg},X_{lg}^+$ to satisfy Axiom of 1 of Definition \ref{def:special}.  These will each differ from the sets defined as $X_{sm}(\mu),X_{sm}^+(\mu),X_{lg}(\tau),X_{lg}(\tau)^+$  in Section \ref{subsub:hp} by at most two points each (note not every interval in $[N]$ is an open ball due to divisibility issues).  This difference will not affect the results proved in Section \ref{subsub:hp} as $N$ is large.   Note Axiom (9) holds by the definition of $X_{lg}^+, X_{sm}^+$ and the assumptions on $\tau$ and $\mu$.  Axiom (2) clearly holds by the discussion following Claim \ref{cl:hpdegree}.  Then (3)-(7) hold by Claims  \ref{cl:hpdegree}, \ref{cl:hpdegreeball}, \ref{cl:hpdegreeballmove}, \ref{cl:hpcover}, \ref{cl:hpcover2}, \ref{cl:hpint}, \ref{cl:hpsplit}, along with the fact that the roles of $A$, $B$ and $C$ are symmetric in $\HP(N)$. 
\end{proof}

\subsection{Proof of Theorems \lowercase{\ref{thm:ternarytriads}} and \lowercase{\ref{thm:ternarytriads1}}}\label{subsub:binaryproof}

In this subsection we give the proofs of  Theorems \ref{thm:ternarytriads} and \ref{thm:ternarytriads1}.  We begin with Theorem \ref{thm:irreg} below, which shows that special $3$-graphs will not admit $\vdisc_3$-homogeneous partitions with binary error.

\begin{theorem}\label{thm:irreg}
 Suppose $p\geq 2$, $0<\mu<1/16$, $0<\tau<1-4\mu$, and $\alpha>0$.  For all $0<\e<\alpha \mu/p^{16}$, there is $t_0\geq 1$, so that for all $T\geq t_0$, the following holds. There is $\rho>0$ such that for all sufficiently large $n$, if $G=(A\cup B\cup C,E)$ is a $(p,\mu,\tau,\alpha,\rho)$-special $3$-partite $3$-graph with $|A|=|B|=|C|=n$, then $G$ does not admit an $\e$-homogeneous partition with binary error and $t$ parts, for any $t_0\leq t\leq T$. 
\end{theorem}
\begin{proof}
Fix $p\geq 2$, $0<\mu<1/16$, $0<\tau<1-4\mu$, and $\alpha>0$.  Choose $0<\e\ll \alpha, \mu, \tau$ and $t_0\gg \e^{-1}$.  Fix $T\geq t_0$.  Then choose  $0<\rho\ll T^{-1}p^{-100}\mu\tau$,  $n\gg T\rho^{-1}\e^{-1}p^{-1}$ and set $N=3n$.  

Suppose $G=(A\cup B\cup C,E)$ is a $3$-partite $3$-graph with $|A|=|B|=|C|=n$, and which is $(p,\alpha,\mu,\tau,\rho)$-special.  For ease of notation, set $V=A\cup B\cup C$.  Suppose towards a contradiction there are $t_0\leq t\leq T$ and an equipartition $\calP=\{V_1\cup \ldots \cup V_t\}$ of $V$ and $\Sigma\subseteq {[t]\choose 2}$ so that $|\Sigma|\leq \e t^2$, and for all $ijk\in {[t]\choose 2}$, if $ij, jk, ik\notin \Sigma$, then $|E\cap K_3[V_i,V_j,V_k]|/|K_3[V_i,V_j,V_k]\in [0,\e)\cup (1-\e,1]$.

Enumerate the sets $A=\{a_i: i\in [n]\}$, $B=\{b_i: i\in [n]\}$, and $C=\{c_i: i\in [n]\}$. For each $X\in \{A,B,C\}$, let $X_{sm},X_{sm}^+,X_{lg},X_{lg}^+$ be as in Definition \ref{def:special}.  Set $\e_3=\e^{1/8}$.  We will use throughout that $t\geq t_0\gg \e^{-1}$, and that each $V_i$ is very large (since $|V_i|=(3n/t)\pm 1$ and $t\leq T$).

We say that a pair $xy\in {V\choose 2}$ is a \emph{non-$\Sigma$ pair} if $xy\in K_2[V_i,V_j]$ for some $V_iV_j\in {\calP\choose 2}\setminus  \Sigma$.  On the other hand, we say $xy$ is a \emph{$\Sigma$-pair} if $xy\in K_2[V_i,V_j]$ for some $V_iV_j\in \Sigma$.  Since $|\Sigma|\leq \e t^2$, the number of non-$\Sigma$ pairs is at least 
$$
{|V|\choose 2}-\sum_{i=1}^t{|V_i|\choose 2}-\sum_{ij\in \Sigma}|V_i||V_j|\geq {|V|\choose 2}-|V|^2/t-\e |V|^2\geq (1-2\e ){|V|\choose 2},
$$
where the last inequality uses that $t\geq t_0$, and that $|V|$ is large.  For each $1\leq i\leq t$, set 
 $$
 A_i=A^*\cap V_i,\text{ }B_i=B^*\cap V_i,\text{ and }C_i=C^*\cap V_i.
 $$
 Given $i\in [t]$ and $X\in \{A,B,C\}$, we say $X_i$ is \emph{trivial} if $|X_i|< \e_3 |V_i|$, and a vertex $x\in V$ is \emph{trivial}, if $x\in X_i$ for some trivial $X_i$.  For each $i\in [t]$, at most two elements of $\{A_i,B_i,C_i\}$ are trivial, and therefore, the number of trivial vertices in $V_i$ is at most $2\e_3|V_i|$.  Consequently, at most $4\e_3 N$ vertices are trivial.  For each $X\in \{A,B,C\}$, let 
  $$
  \calI_X:=\{i\in [t]: X_i\text{ is non-trivial}\}.
  $$
Since all but at most $4\e_3 N$ vertices are non-trivial, 
\begin{align}\label{al:triv}
\Big|\bigcup_{i\in \calI_X}X_i\Big|\geq (1-12\e_3)|X|, 
\end{align}
and consequently, $|\calI_X|\geq (1-12\e_3)|X|/(N/t)\geq (1-12\e_3)t/3$, for each $X\in \{A,B,C\}$.  We now turn to the sets of pairs which are in some sense ``error pairs'' with respect to the partition $\{V_1,\ldots, V_t\}$. In particular, we define
\begin{align*}
\calS_{1}=\{X_iY_j: i,j\in [t], XY\in {\{A,B,C\}\choose 2}& \text{ and either }i=j, \\&\text{ or one of }X_i,Y_j\text{ is trivial}\}
\end{align*}
and
\[\calS_{2}=\{X_iY_j\notin \calS_1: V_iV_j\in \Sigma\}.\]
Set $\calS=\calS_{1}\cup\calS_{2}$.  Note that $|\calS_{2}|\leq 3|\Sigma|\leq 3\e t^2$.  Using (\ref{al:triv}), we have that 
\begin{align}\label{al:sigmacont}
|\calS|\leq 3|\Sigma|+3t+t\Big|\bigcup_{i\notin \calI_X}X_i\Big|/(N/t)\leq 13\e_3 t^2,
\end{align}
where the last inequality uses that $t\gg t_0\gg \e^{-1}$.  We think of the set $\calS$ as controlling all potential error triples.  Indeed, we will show that any triple $A_iB_jC_k$ avoiding $\calS$  will be somewhat homogenous with respect to $H$.  Specifically, we claim that if $A_iB_jC_k$ has $A_iB_j,A_iC_k, B_jC_k\notin \calS$, then $A_iB_jC_k$ is $\sqrt{\e}\mu^8$-homogeneous with respect to $H$.  Suppose towards a contradiction that $A_iB_j,A_iC_k, B_jC_k\notin \calS$, but $|E\cap K_3[A_i,B_j,C_k]|/|A_i||B_j||C_k|\in (\sqrt{\e}\mu^8,1-\sqrt{\e}\mu^8)$.  Since $A_iB_j,A_iC_k, B_jC_k\notin \Sigma$, by assumption, there is $u\in \{0,1\}$ such that $|E^u\cap K_3[V_i,V_j,V_k]|\geq (1-\e)|V_i||V_j||V_k|$.  But we also have
$$
|K_3[V_i,V_j,V_k]\setminus E^u|\geq \sqrt{\e}|K_3[A_i,B_j,C_k]|\geq \e^{1/2}\mu^8\e_3^3|V_i||V_j||V_k|>\e|V_i||V_j||V_k|,
$$
where the last inequality is by our choice of $\e\ll \mu$.  This contradicts that $|E^u\cap K_3[V_i,V_j,V_k]|\geq (1-\e)|V_i||V_j||V_k|$.  Thus any triple $A_iB_jC_k$ avoiding $\calS$ is $\sqrt{\e}\mu^8$-homogenous with respect to $H$.  

Eventually, we will contradict (\ref{al:sigmacont}).  We will use the following definition.  Given $X\in \{A,B,C\}$ and $i\in [t]$,  set 
\begin{align*}
\Irr(X_i)=\{Y_jZ_k: &X_iY_j, X_iZ_k\notin \calS, \\ &|E\cap K_3[X_i,Y_j,Z_k]|/|K_3[X_i,Y_j,Z_k]|\in (\e^{1/2}\mu^8,1-\e^{1/2}\mu^8)\}.
\end{align*}
If $X_i$ is non-trivial, then our arguments above imply that if $Y_jZ_k\in \Irr(X_i)$, then $Y_jZ_k\in \calS$. Consequently, for any non-trivial $X_i$, $|\Irr(X_i)|\leq 13\e_3 t^2$, and thus $|\bigcup_{Y_jZ_k\in \Irr(X_i)}Y_jZ_k|\leq 39\e_3 n^2$.

It will be important to know that, based on our assumptions, the pairs in $\calS$ are somewhat constrained in their distibution.  This is the purpose of our next series of definitions. First, set
\begin{align*}
S_1=\bigcup_{X_iY_j\in \calS_1}K_2[X_i,Y_j],\text{ and }S_2=\bigcup_{X_iY_j\in \calS_2}K_2[X_i,Y_j].
\end{align*}
By (\ref{al:triv}), we have that
\[|S_1|\leq \sum_{i=1}^t t{V_i\choose 2}+\sum_{X\in \{A,B,C\}}\sum_{i\in \calI_X}|X_i||V\setminus X|\leq t(3n/t)^2+2n\sum_{X\in \{A,B,C\}}\sum_{i\in \calI_X}|X_i|,\]
which is bounded above by
\[9n^2/t+6(12\e_3)n^2\leq 61\e_3 n^2,\]
where the last inequality uses that $t\gg \e^{-1}$.  Given $X_i$, set
$$
\calS_1(X_i):=\{Y_j: X_iY_j\in \calS_1\}\; \text{ and }\;S_1(X_i)=\bigcup_{Y_j\in \calS_1(X_i)}Y_j.
$$
If $X_i$ is non-trivial, then using the above and our bound on the number of non-trivial vertices, we have
$$
|S_1(X_i)|\leq |V_i\setminus X_i|+|\bigcup_{j\notin \calI_Y}Y_j\cup \bigcup_{j\notin \calI_Z}Z_j| \leq |V_i|+4\e_3N\leq 7\e_3(2n)=14\e_3 n.
$$
Because $|\calS_2|\leq 3\e t^2$, $|S_2|\leq 9\e n^2$.  Similar to the above, given $X_i$, we set 
$$
\calS_2(X_i)=\{Y_j: X_iY_j\in \calS_2\}\;\text{ and }\;S_2(X_i)=\bigcup_{Y_j\in \calS_2(X_i)}Y_j.
$$
To deduce that $S_2(X_i)$ is small, it is not enough to just know that $X_i$ is non-trivial.  This is the motivation of the next definition.  Given $X\in \{A,B,C\}$, define
$$
Gd(X):=\{x\in X: x\text{ is non-trivial and }|\{y\in V\setminus X : xy\in S_2\}|\leq 3\sqrt{\e}(2n)\}.
$$
We say $x\in X$ is \emph{good} if $x\in Gd(X)$.  Since $|S_2|\leq 9\e n^2$, we must have that for any $X\in \{A,B,C\}$, $|Gd(X)|\geq (1-3\e^{1/2})n$.    Given $X\in \{A,B,C\}$ and $i\in [t]$, we say that $X_i$ is \emph{good} if $X_i$ is non-trivial and $|X_i\cap Gd(X)|\geq (1-(3\e)^{1/3})|X_i|$.  Note that for each $X\in \{A,B,C\}$, $|\bigcup_{X_i\text{ not good}}X_i|$ is at most
$$
|\bigcup_{i\notin \calI_X}X_i|+\frac{|X\setminus Gd(X)|}{(3\e)^{1/3}\e_3(3n/t)}\leq 12\e_3|X|+3\e^{1/2-1/3-1/8}|X|\leq 4\e_3^{1/3}|X|.
$$
Thus, most of the elements in $X$ are in a good $X_i$.  We claim that if $X_i$ is good, then $|\calS_2(X_i)|\leq 4\e^{5/24}t$.  Suppose $X_i$ is good. By definition, $|\{xy: x\in X, xy\in S_2\}|$ is at most
\begin{align*}
(3\sqrt{\e})2n|X_i\cap Gd(X)|+2n|X_i\setminus Gd(X)|&\leq 6\sqrt{\e}n|X_i|+2n(3\e)^{1/3}|X_i|\\
&\leq 12\e^{1/3}|X_i|n.
\end{align*}
Note that if $Y_j\in \calS_2(X_i)$, then $Y_j$ is non-trivial so $K_2[X_i,Y_j]$ contains at least $|X_i|\e_3(3n/t)$ many elements of $S_2(X_i)$.  Consequently,
$$
|\calS_2(X_i)|\leq (12\e^{1/3}|X_i|n)/(|X_i|\e_3(3n/t))=4\e^{5/24} t.
$$
This also shows that if $X_i$ is good, then $|S_2(X_i)|\leq (3n/t)4\e^{5/24} t=12 \e^{5/24}n$.  We have now establishes the required properties about the distribution of $\calS$.  We now turn to considering the implications of the axioms from Definition \ref{def:special}.

We say that $X_i$ is \emph{mixed-degree} if it is good, and moreover,
$$
|\{x\in X_i: \min\{|N(x)|, |\neg N(x)|\}\geq \alpha n^2\}|\geq \e_3^{1/10}|X_i|.
$$
Note that by the Degree Axiom, if a set $X_i$ is good and satisfies $|X_i\cap (X_{sm}^+\cup X_{lg}^+)|\geq \e_3^{1/10}|X_i|$, then $X_i$ is mixed degree.  Our next goal is to show that every mixed-degree $X_i$ is mostly contained in an open ball of small radius.   Given $\eta\in (0,1)$, define
$$
\Mix(\eta,X_i)=\{yz\in YZ: |N_H(yz)\cap X_i|/|X_i|\in (\eta,1-\eta)\}.
$$
We can already show that if $X_i$ is good, then $\Mix(\eta,X_i)$ is not too large for certain values of $\eta$.
\begin{claim}\label{cl:mixbound}
If $X_i$ is good then $|\Mix(\e^{1/4}\mu^4,X_i)|\leq 2\e^{1/4}\mu^4 n^2$.
\end{claim}
\begin{proof}
Recall that since $X_i$ is non-trivial, for any $X_iY_j, X_iY_j\notin \calS$ and $Y_jZ_k\notin \Irr(X_i)$, the triple $X_iY_jZ_k$ is $\e^{1/2}\mu^8$-homogeneous.  Thus, for any such $Y_j,Z_k$, there is some $\eta\in \{0,1\}$ so that at least $(1-\e^{1/4}\mu^4)|Y_j||Z_k|$ many elements $yz\in Y_jZ_k$ have $|N^{\eta}_H(yz)\cap X_i|\geq (1-\e^{1/4}\mu^4)|X_i|$.  Using the size estimates from earlier in the proof, and the fact that $X_i$ is good, we have
\begin{align*}
|\Mix(\e^{1/4}\mu^4,X_i)| &\leq \e^{1/4}\mu^4\Big|\bigcup_{Y_jZ_k\notin ( \calS_1(X_i)\cup\calS_2(X_i))\setminus \Irr(X_i)}K_2[Y_j,Z_k]\Big|\\
 &\hspace{20pt}+\Big|\bigcup_{Y_jZ_k\in \Irr(X_i)}K_2[Y_j,Z_k]\Big|+n|S_1(X_i)|+n|S_2(X_i)|\\
&\leq \e^{1/4}\mu^4 n^2+39\e_3 n^2+14 \e_3 n^2+12\e^{1/3} n^2\\
&\leq 2\e^{1/4}\mu^4 n^2,
\end{align*}
where the last inequality is because $\e\ll \mu$.   
\end{proof}

Claim \ref{cl:mixbound} tells us that, for any mixed degree $X_i$, there is only a small number of pairs $yz$ with $N(yz)\cap X_i$ and $\neg N(yz)\cap X_i$ both being large.  On the other hand, because $X_i$ has mixed degree, its vertices are involved in many edges as well as many non-edges.  This suggests that there must be many pairs $yz$ for which $N(yz)\cap X_i$ is most of $X_i$, as well as many pairs $yz$ for which $\neg N(yz)\cap X_i$ is almost all of $X_i$.   This idea is important for the rest of the proof.  

Consider the following sets.
\begin{align*}
\Gamma^0(X_i)&:=\{yz\in K_2[Y, Z]: |(\neg N(yz))\cap X_i|\geq (1-2\e^{1/4})|X_i|\}\text{ and }\\
\Gamma^1(X_i)&:=\{yz\in K_2[Y, Z]: |N(yz)\cap X_i|\geq (1-2\e^{1/4})|X_i|\}.
\end{align*}
Note that $\Mix(2\e^{1/4},X_i)=K_2[Y, Z]\setminus (\Gamma^0(X_i)\cup \Gamma^1(X_i))$.  Define
\begin{align*}
Q(X_j)=\{yzz'\in K_3[Y,Z,Z]: &yz\in \Gamma^1(X_i), yz'\in \Gamma^0(X_i)\}\\
&\cup \{yy'z\in K_3[Y,Y,Z]: yz\in \Gamma^1(X_i), y'z\in \Gamma^0(X_i)\}.
\end{align*}

We show in our next claim that if $X_i$ is mixed-degree, then $Q(X_i)$ is nonempty.  Later, we will use this, along with the Intersection Axiom, to deduce that every mixed-degree $X_i$ is mostly contained in an open ball of small radius.

\begin{claim}\label{cl:intQ} If $X_i$ is mixed-degree, then $Q(X_i)$ is nonempty. 
\end{claim}
\begin{proof}
 Fix $X\in \{A,B,C\}$ and suppose $X_i$ is mixed-degree.  Let $Y,Z$ be such that $\{X,Y,Z\}=\{A,B,C\}$.  Define 
 $$
 J_Y=\{Y_j: X_iY_j\notin \calS\}\text{ and }J_Z=\{Z_k: X_iZ_k\notin \calS\}.
 $$
Note that since $X_i$ is non-trivial, $Y_j\in J_Y$ if and only if $j\neq i$, $Y_j$ is non-trivial and $Y_j\notin \calS_2(X_i)$.  Consequently, since $X_i$ is good,
$$
|J_Y|\geq |\calI_Y|-1-|\calS_2(X_i)|\geq |\calI_Y|-1- 4\e^{5/24}t\geq (1-13 \e_3)|\calI_Y|\geq(1-13 \e_3)^2t/3,
$$
where the last inequality uses the fact that $|\calI_Y|\geq (1-12\e_3)t/3$.  Similarly $Z_k\in J_Z$ if and only if $k\neq i$, $Z_k$ is non-trivial and $Z_k\notin \calS_2(X_i)$.  A similar argument also shows $|J_Z|\geq (1-13 \e_3)|\calI_Z|\geq (1-13 \e_3)^2t/3$.  Define 
\[F_0=\{Y_jZ_k\in K_2[J_Y, J_Z]: |E\cap K_3[X_i,Y_j,Z_k]|/|X_i||Y_j||Z_k|\leq \e^{1/2}\}\]
and
\[F_1=\{Y_jZ_k\in K_2[J_Y, J_Z]: |E\cap K_3[X_i,Y_j,Z_k]|/|X_i||Y_j||Z_k|\geq 1-\e^{1/2}\}.\]
Note that by the definition of $J_Y,J_Z$, we have that if $Y_jZ_k\in K_2[J_Y,J_Z]\setminus (F_1\cup F_0)$, then $Y_jZ_k\in \Irr(X_i)$. 

We claim that it suffices to show that either there is some $Y_jZ_kZ_{k'}$ with $Y_jZ_k\in F_0$ and $Y_jZ_{k'}\in F_1$, or there is some $Y_jY_{j'}Z_k$ with $Y_jZ_k\in F_0$ and $Y_{j'}Z_{k}\in F_1$.  Indeed, suppose there exists some $Y_jZ_kZ_{k'}$ with $Y_jZ_k\in F_0$ and $Y_jZ_{k'}\in F_1$.  Since $Y_jZ_k\in F_0$, standard counting arguments imply that there is a set $\Omega_0\subseteq Y_j$ of size at least $(1-\e^{1/4})|Y_j|$ such that for all $y\in \Omega_0$,
$$
|\{zx\in K_2[Z_k,X_i]: xyz\notin E\}|\geq (1-\e^{1/4})|X_i||Z_k|.
$$
Similarly, since $Y_jZ'_k\in F_1$,  there is a set $\Omega_1\subseteq Y_j$ of size at least $(1-\e^{1/4})|Y_j|$ such that for all $y\in \Omega_1$,
$$
|\{zx\in K_2[Z_k,X_i]: xyz\in E\}|\geq (1-\e^{1/4})|X_i||Z_k|.
$$
Choose any $y_{\star}\in \Omega_1\cap \Omega_0$ (clearly $\Omega_0\cap \Omega_1$ is nonempty, as its size is at least $(1-2\e^{1/4})|Y_j|$).  Consider the bipartite graphs 
$$
G_0=(Z_k\cup X_i, N_H(y_{\star})\cap K_2[Z_k,X_i])\text{ and }G_1=(Z_{k'}\cup X_i, N_H(y_{\star})\cap K_2[Z_{k'},X_i]).
$$
  By construction, the density of edges in $G_0$ is at most $\e^{1/4}$, and the density of edges in $G_1$ is at least $1-\e^{1/4}$.  Thus, there exists a $z\in Z_k$ of degree at most $\e^{1/4}|X_i|$ in $G_0$, and there exists a $z'\in Z_{k'}$ with degree at least $(1-\e^{1/4})|X_i|$ in $G_1$.  Now  
$$
\min\{|\neg N_H(y_{\star}z)\cap X_i|,|N_H(y_{\star}z')\cap X_i|\}= \min\{|N_{G_1}(z')|,|N_{G_0}(z)|\}\geq (1-\e^{1/4})|X_i|,
$$
 so $y_{\star}zz'\in Q(X_i)$.   An symmetric argument shows that if there is some $Y_jY_{j'}Z_k$ with $Y_jZ_k\in F_0$ and $Y_{j'}Z_{k}\in F_1$, then there exists some $yy'z_{\star}\in Q(X_i)$. Thus we just need to show that either there is some $Y_jY_{j'}Z_k$ with $Y_jZ_k\in F_0$ and $Y_jZ_{k'}\in F_1$, or there is some $Y_jZ_kZ_{k'}$  with $Y_jZ_k\in F_0$ and $Y_{j'}Z_{k}\in F_1$. 

Suppose towards a contradiction that no such $Y_jY_{j'}Z_k$ or $Y_jZ_kZ_{k'}$ exist.  Let $\calG$ be the bipartite graph $\calG$ with vertex set $J_Y\cup J_Z$ and edge set $F_1$.  Recall that for all $Y_jZ_k\in (J_Y\times J_Z)\setminus (F_0\cup F_1)$, we must have $Y_jZ_k\in \Irr(X_i)$.  Consequently, 
$$
|K_2[J_Y,J_Z]\setminus (F_0\cup F_1)|\leq 13\e_3 t^2\leq 100\e_3 |J_Y||J_Z|,
$$
where the last inequality uses our lower bounds on $|J_Y|,|J_Z|$.  Therefore, if we set 
\[J_Y'=\{Y_j\in J_Y: |\{Z_k\in J_Z: Y_jZ_k\notin F_0\cup F_1\}|\leq \sqrt{100\e_3}|J_Z|\}\]
and
\[J_Z'=\{Z_k\in J_Z: |\{Y_j\in J_Y: Y_jZ_k\notin F_0\cup F_1\}|\leq \sqrt{100\e_3}|J_Y|\},\]
then  $|J_Y'|\geq \sqrt{100\e_3}|J_Y|$ and $|J_Z'|\geq \sqrt{100\e_3}|J_Z|$.  Because there does not exist $Y_j,Z_k,Z_{k'}$ with $Y_jZ_k\in F_0$ and $Y_jZ_{k'}\in F_1$, we have that for all $Y_j\in J'_Y$, either $|N_{\calG}(Y_j)|\geq (1-\sqrt{100\e_3})|J_Z|$ or $|N_{\calG}(Y_j)|=0$.  Similarly, because there does not exist $Y_j,Y_{j'},Z_{k}$ with $Y_jZ_k\in F_0$ and $Y_{j'}Z_{k}\in F_1$, we have that for all $Z_k\in J_Z'$, either $|N_{\calG}(Z_k)|\geq (1-\sqrt{100\e_3})|J_Z|$ or $|N_{\calG}(Z_k)|=0$.  By Lemma \ref{lem:twosticks}, either $|F_1|\leq 2(100\e_3)^{1/4}|J_Y||J_Z|$ or $|F_1|\geq (1-2(100\e_3)^{1/4})|J_Y||J_Z|$.  

Suppose first that  $|F_1|\leq 2(100\e_3)^{1/4}|J_Y||J_Z|$.  Then $|\bigcup_{x\in X_i}N_H(x)|$ is at most
\[|X_i|\Big|\bigcup_{Y_jZ_k\notin F_0\cup F_1}K_2[Y_j,Z_k] \Big| +\sum_{Y_jZ_k\in F_1}|X_i||Y_j||Z_k|+\sum_{Y_jZ_k\in F_0}|E\cap K_3[X_i,Y_j,Z_k]|,\]
which in turn is bounded above by
\begin{align*}
 |X_i|&39\e_3 n^2+|F_1|(3n/t)^2+\e^{1/2}|J_Y||J_Z|(3n/t)^2\\
&\leq  |X_i|\Big(39\e_3 n^2+3(100\e_3)^{1/4}t^2(3n/t)^2+\e^{1/2}t^2(3n/t)^2\Big)\\
&=|X_i|n^2\Big(39\e_3 n^2+27(100\e_3)^{1/4}+9\e^{1/2} \Big)\\
&\leq \e^{1/5}_3 |X_i|n^2,
\end{align*}  
where the last inequality is because $\e$ is small.  This implies that there is a set $X_i'\subseteq X_i$ such that $|X_i'|\geq (1-\e_3^{1/10}) |X_i|$, and such that for all $x\in X_i'$, $|N_H(x)|\leq \e_3^{1/10}n^2<\alpha n^2$, but this contradicts that $X_i$ is mixed degree.  

Suppose on the other hand that $|F_1|\geq (1-2(100\e_3)^{1/4})|J_Y||J_Z|$.  Note that 
\begin{align*}
\Big|\bigcup_{Y_jZ_k\notin F_1}K_2[Y_j,Z_k]\Big|&\leq \Big|\bigcup_{Y_jZ_k\notin F_0\cup F_1}K_2[Y_j,Z_k] \Big|+\Big|\bigcup_{Y_jZ_k\in F_0}K_2[Y_j,Z_k]\Big|\\
&\leq 39\e_3 n^2 +2(100\e_3)^{1/4}t^2\Big(\frac{3n}{t}\Big)^2\\
&\leq \e_3^{1/5} n^2/2,
\end{align*} 
where the last inequality is because $\e$ is small.  Consequently, we have that 
$$
|\bigcup_{x\in X_i}N_H(x)|\geq (1-2\e^{1/4})\sum_{Y_jZ_k\in F_1}|Y_j||Z_k|\geq (1-\e^{1/2})n^2(1-\e_3^{1/5}/2)\geq n^2(1-\e^{1/5}_3).
$$
  From this we deduce there is a set $X_i'\subseteq X_i$ of size at least $(1-\e_3^{1/10})|X_i|$, so that for all $x\in X_i'$, $|N_H(x)|\geq (1-\e_3^{1/10})n^2>(1-\alpha)n^2$, contradicting that $X_i$ is mixed-degree. This finishes the proof of Claim \ref{cl:intQ}.
\end{proof}

We now deduce that every mixed-degree $X_i$ is mostly contained in an open ball.  Given a mixed-degree $X_i$, define
\begin{align*}
D(X_i)=\{d>0: \text{there is }& yzz'\in Q(X_i) \text{ with } d(z,z')=d
\\& \text{ or } yy'z\in Q(X_i) \text{ with } d(y,y')=d\}.
\end{align*}
Since $D(X_i)$ is nonempty by Claim \ref{cl:intQ} (and $V$ is finite), we may choose $d_{X_i}:=\min D(X_i)$ and define a ball $B_{X_i}$ as follows.  If $d_{X_i}$ is witnessed by some $yy'z'\in Q(X_i)$, let $\rho_{X_i}\leq d_{X_i}$ be minimal such that for some $x$, $N(yz)\cap \neg N(y',z)\subseteq \overline{B_{\rho_{X_i}}(x)}$ (such an $x$ and $\rho_{X_i}$ exist by the Intersection Axiom).  If on the other hand, $d_{X_i}$ is witnessed by some $yzz'\in Q(X_i)$, let $\rho_{X_i}\leq d_{X_i}$ be minimal such that for some $x$, $N(yz)\cap \neg N(y,z')\subseteq \overline{B_{\rho_{X_i}}(x)}$ (again such an $x$ and $\rho_{X_i}$ exists by the Intersection Axiom).  In either case, set $B_{X_i}=\overline{B_{\rho_{X_i}}(x)}$, and note $|X_i\cap B_{X_i}|\geq (1-2\e^{1/4})|X_i|$.  Combining this with the Size Axiom and the fact that $X_i$ is non-trivial (since it is good) we have that 
$$
(1-2\e^{1/4})\e_3(3n/t)\leq |B_{X_i}|\leq p\rho_{X_i}n,
$$
and thus, $\rho_{X_i}\geq (1-2\e^{1/4})t^{-1}p^{-1}\e_3>50\rho p^{50}$, where the last inequality is because $t\leq T$ and by our choice of $\rho$.

We say $X_i$ is \emph{sm-interesting} if it is good and $|X_i\cap X_{sm}|\geq \e_3^{1/10}|X_i|$.  Similarly, we say $X_i$ is \emph{lg-interesting} if it is  good and $|X_i\cap X_{lg}|\geq \e_3^{1/10}|X_i|$.  We will say $X_i$ is simply \emph{interesting} if it is either sm-interesting or lg-interesting.  Note that by the Degree Axiom, any interesting $X_i$ is automatically mixed-degree.  Therefore, when $X_i$ is sm-interesting, we can apply Claim \ref{cl:intQ} to find $B_{X_i}$.  In this case, we have that $|X_i\cap B_{X_i}|\geq (1-2\e^{1/4})|X_i|$ and $|X_i\cap X_{sm}|\geq \e_3^{1/10}|X_i|$.  This implies that 
$$
|X_i\cap B_{X_i}'\cap X_{sm}|\geq (\e_3^{1/10}-2\e^{1/4})|X_i|\geq \e_3^{1/10}|X_i|/2>0.
$$
Thus, we may choose some $x_i\in X_i\cap B_{X_i}\cap X_{sm}$, and define $B_{X_i}'=B_{2\rho_{X_i}}(x_i)$.  

Similarly, when $X_i$ lg-interesting, $|X_i\cap B_{X_i}|\geq (1-2\e^{1/4})|X_i|$ and $|X_i\cap X_{lg}|\geq \e_3|X_i|$ imply 
$$
|X_i\cap B_{X_i}'\cap X_{lg}|\geq (\e_3^{1/10}-2\e^{1/4})|X_i|\geq \e_3^{1/0}|X_i|/2>0.
$$
Thus, we may choose some $x_i\in X_i\cap B_{X_i}\cap X_{lg}$, and define $B_{X_i}'=B_{2\rho_{X_i}}(x_i)$.  

We have now shown that when $X_i$ is sm-interesting (respectively lg-interesting), there is an open ball $B_{X_i}'$ of radius $\rho_{X_i}':=2\rho_{X_i}\leq 2d_{X_i}$, with center $x_i\in X_{sm}$ (respectively $x_i\in X_{lg}$), and so that $|X_i\cap B'_{X_i}|\geq (1-2\e^{1/4})|X_i|$.  Our next goal is to show that this $\rho'_{X_i}$ is not too big.  

\begin{claim}\label{cl:rho}
If $X_i$ is interesting, then $\rho'_{X_i}\leq \mu^2/2$.
\end{claim}
\begin{proof}
Suppose towards a contradiction there exists an interesting $X_i$ with $\rho'_{X_i}>\mu^2/2$.  If $X_i$ is sm-interesting, set $X_*=X_{sm}$ and $Y_*=Y_{lg}$ if $X_i$ is lg-interesting, set $X_*=X_{lg}$ and $Y_*=Y_{sm}$.

By the Covering Axiom, there exists $\calQ_1$ a collection of $m_1$ open balls of radius, each of radius at most $\mu^2/32$, and at least $\mu^2/96$, each contained in $B_{X_i}'$, such that $m_1\leq2p (\mu^{2}/32)(\rho'_{X_i})^{-1}$, and such that if $\calQ_1'$ is set of elements of $\calQ_1$ centered at a point in $X_*^+$, then $|(B_{X_i}'\cap X_*)\setminus (\bigcup \calQ_1')|\leq 2m_1$.  Note that by assumption, $m_1\leq 2p (\mu^{2}/32)(\rho'_{X_i})^{-1}\leq 2p$.  By the pigeonhole principle, there is $B_{r_1}(x')\in \calQ_1'$ such that 
$$
|B_{r_1}(x')\cap X_*\cap X_i|\geq \Big(\e_3^{1/10}|X_i|/2-2m_1\Big)/m_1\geq \e_3^{1/10}|X_i|/4p, 
$$
where the last inequality is because $m_1\leq 2p\ll |X_i|$.  We claim that
$$
|B_{\mu^2/16}(x')\cap X_i|<(1-2\e^{1/4})|X_i|.
$$
Suppose towards a contradiction $|B_{\mu^2/16}(x')\cap X_i|\geq (1-2\e^{1/4})|X_i|$. Choose any $y\in Y_{*}$.  By Splitting Axiom 1, there exist $z_0,z_1\in Z$ such that $d(z_0,z_1)\leq 3\mu^2/16$ such that $B_{\mu^2/16}(x')\subseteq N_H(z_1y)\cap \neg N_H(z_0y)$.  But now $yz_1z_0\in Q(X_i)$, so $d(z_1,z_0)\in D(X_i)$.  But $d(z_0,z_1)<\mu^2/4<\rho'_{X_i}/2=\rho_{X_i}$, a contradiction to the definition of $\rho_{X_i}$.

Thus  $|B_{\mu^2/16}(x')\cap X_i|<(1-2\e^{1/4})|X_i|$.  By the Partition Axiom, there is a set $\calQ_2$ of $m_2$ open balls, each of radius at most $\mu^2/16$, and at least $\mu^2/48$, one of which is $B_{\mu^2/16}(x')$, such that $m_2\leq 16\cdot 4p\mu^{-2}$, and such that $|X\setminus (\bigcup \calP)|\leq 2m_2$.  Since $|B_{\mu^2/16}(x')\cap X_i|<(1-2\e^{1/4})|X_i|$, by pigeonhole, there is some $B_{r_2}(x'')\in \calQ_2$, distinct from $B_{\mu^2/16}(x')$, such that 
$$
|X_i\cap B_{r_2}(x'')|\geq (2\e^{1/4}|X_i|-2m_2)/m_2\geq \frac{\e^{1/4}\mu^2|X_i|}{16\cdot 2p}-2\geq \frac{\e^{1/4}\mu^2|X_i|}{16\cdot 4p},
$$
where the last inequalities use the upper bound on $m_2$, and the fact that $X_i$ is large.  By the Covering Axiom, there exists a collection $\calQ_3$ of $m_3$ disjoint open balls, each of radius at most $r_2/4$ and at least $r_2/12$, each contained in $B_{r_2}(x'')$, such that $m_3\leq 8p$, and such that $|B_{r_2}(x'')\setminus (\bigcup \calQ_3)|\leq 2m_3$.  By the pigeonhole principle, there is $B_{r_3}(x''')\in \calQ_3$ such that 
$$
|X_i\cap B_{r_3}(x''')|\geq (\e^{1/4}\mu^2|X_i|/(16\cdot 4p)-2m_3)/m_3\geq \frac{\e^{1/4}\mu^2|X_i|}{16\cdot 32 p^2}-2\geq \e^{1/4}\mu^4|X_i|,
$$
where the last inequalities use the upper bound on $m_3$, the fact that $\e\ll \mu\ll 1$, and the fact that $X_i$ is large.  Note that $r_2>r_3\geq \rho$. 

Now, since $B_{\mu^2/16}(x')$ is disjoint from $B_{r_2}(x'')$ and $B_{r_3}(x''')\subseteq B_{r_2}(x'')$, we have that $d(x''',x')\geq \mu^2/16+r_3$.  Set $r''=d(x''',x')-r_1-r_3$, and $r_0=\min\{r_1,r_3\}/3$.  Recall $r_1\leq \mu^2/32$  and $\mu^2/48\cdot 12\leq r_2/12 \leq r_3\leq r_2/4\leq \mu^2/16\cdot 4$.  Consequently, 
$$
d(x''',x')\geq \mu^2/16+r_3>2\max\{r_1,r_3\}. 
$$
Therefore, by Splitting Axiom 2, there are is an open ball $B_1$ in $Y$ and an open ball $B_2$ in $Z$, each of radius $r_0/3$, such that 
\[K_3[B_1,B_2,B_{r_1}(x')]\subseteq E \text{ and } K_3[B_1,B_2,B_{r_3}(x'')]\subseteq \neg E,\] 
or vice versa.  Since
$$
|B_{r_1}(x')\cap X_i|\geq \e_3^{1/10} |X_i|/4p\text{ and }|  B_{r_3}(x'')\cap X_i|\geq \e^{1/4}\mu^4|X_i|,
$$
we must have that $K_2[B_1,B_2]\subseteq \Mix(\e^{1/4}\mu^4,X_i)$.  Since $r_1\geq \mu^2/96$ and $r_3\geq \mu^2/(48\cdot 12)$, we have $r_0/3\geq \mu^2/(48\cdot 36)$, so the Size Axiom implies we have 
$$
\min\{|B_1|,|B_2|\}\geq \mu^2n/(48\cdot 36)(2p) =\mu^2 n/(96\cdot 36 p).
$$
Consequently 
$$
|\Mix(\e^{1/4}\mu^4,X_i)|\geq |K_2[B_1,B_2]|\geq \frac{\mu^4 n^2}{(96\cdot 36 p)^2}.
$$
However, since $X_i$ is good, Claim \ref{cl:mixbound} implies $|\Mix(\e^{1/4}\mu^4,X_i)|\leq 2\e^{1/4}\mu^4 n^2$.  Because we chose $\e\ll \mu$, $2\e^{1/4}\mu^4 n^2< \mu^4 n^2/(96\cdot 36 p)^2$, so this is a contradiction.  This finishes the proof of Claim \ref{cl:rho}.
\end{proof}

Note that Claim \ref{cl:rho} tells us that if $X_i$ is sm-interesting, then since the center of $B_{X_i}'$ is in $X_{sm}$ and the radius of $B_{X_i}'$ is less than $\mu^2/2$, $B_{X_i}'\subseteq X_{sm}^+$.  Consequently $|B_{X_i}'\cap X_i\cap X_{sm}^+|\geq (1-2\e^{1/4})|X_i|$.  Similarly if $X_i$ is lg-interesting, then $B_{X_i}'\subseteq X_{lg}^+$, and consequently $|B_{X_i}'\cap X_i\cap X_{lg}^+|\geq (1-2\e^{1/4})|X_i|$.  

Let  $K^{lg}_X$ be the set of lg-interesting $X_i$, and let $K^{sm}_X$ be the set of sm-interesting $X_i$.    Using the above, we have that 
\begin{align*}
|X_{sm}\setminus (\bigcup_{X_i\notin K^{sm}_X}X_i)|&\leq |\bigcup_{X_i\text{ not good}}X_i|+\sum_{X_i\text{ good},X_i\notin K^{sm}_X}|X_i\cap X_{sm}|\\
&\leq 2\e^{1/3}_3|X|+2\e^{1/4}|X|\\
&\leq 3\e^{1/3}_3 |X_{sm}|,
\end{align*}
where the last inequality is because  $\e\ll \mu$, and because, by the Size and Special Sets Axioms, $|X_{sm}|\geq \mu |X|/3p$.  Consequently,
\begin{align}\label{al:specialsm}
|\bigcup_{X_i\in K^{sm}_X}X_i|\geq |X_{sm}|(1-4\e^{1/3}_3).
\end{align}
A similar argument shows 
\begin{align}\label{al:speciallg}
|\bigcup_{X_i\in K^{lg}_X}X_i|\geq |X_{lg}|(1-4\e^{1/3}_3).
\end{align}

Thus, for each $X\in \{A,B,C\}$, most of $X_{sm}$ is in an sm-interesting $X_i$, and most of $X_{lg}$ is in an lg-interesting $X_i$.   

\begin{claim}\label{cl:balls}
For each $X\in \{A,B,C\}$ and all interesting $X_i$,  and for each $u\in [3]$, there is a set $\phi_u(X_i)$ of open balls in $X$ of radius at most $\rho'_{X_i}/p^{2+u}$, and at least $\rho'_{X_i}/3p^{2+u}$ such that  the following holds for each $u\in [3]$.
\begin{enumerate}[label=\normalfont(\arabic*)]
\item For all $B\in \phi_u(X_i)$, $B\subseteq B_{X_i}'$,
\item For each $B\in \phi_u(X_i)$, $|X_i\cap B|\geq (1-2\e^{1/4})|X_i|/p^{3+u}$,
\item For each $u\in [3]$, $|X_i\cap (\bigcup_{B\in \phi_u(X_i)}B)|\geq (1-3\e^{1/4})|X_i|$.
\end{enumerate}
\end{claim}
\begin{proof}
Suppose $X_i$ is interesting.  If $X_i$ is lg-interesting, let $X_*=X_{lg}$ and if sm-interesting, let $X_*=X_{sm}$.  

By the Covering Axiom, for each $u\in [3]$, there is a set $\calP_u$ of $m_u\leq 2p^{3+u}$ disjoint open balls, each of some radius $r_u\leq \rho'_{X_i}/p^{2+u}$, such that $|B_{X_i}'\setminus (\bigcup\calP_u)|\leq 2m_u$, and such that each ball in $\calP_u$ is contained in $B_{X_i}'$.  For each $u\in [3]$, let 
\begin{align*}
\phi_u(X_i)=\{B_{\rho'_{X_i}/p^{2+u}}(x'): &B_{\rho'_{X_i}/p^{2+u}}(x')\in \calP_u \text{ and } 
\\ &|B_{\rho'_{X_i}/p^{2+u}}(x')\cap X_i|\geq (1-2\e^{1/4})|X_i|/m_u\}.
\end{align*}
Note (1) and (2) hold by definition of $\phi_u(X_i)$.  We show (3). By, assumption, $|B'_{{X_i}}\setminus (\bigcup \calP_u)|\leq 2 m_u$.  Thus,
$$
|X_i\cap (\bigcup_{B'\in \phi_u(X_i)}B')|\geq |X_i\cap B_{X_i}'|-2m_u\geq (1-2\e^{1/4})|X_i|-4p^{3+u}\geq (1-3\e^{1/4})|X_i|,
$$
where the last inequality is because $X_i$ is large.
\end{proof}

For ease of notation, given an interesting $X_i$ and $u\in \{1,2,3\}$, let $P_u(X_i)$ be the set of centers of the balls in $\phi_u(X_i)$, and let $\Psi_u(X_i)=\bigcup_{B'\in \phi_u(X_i)}B'$. Note that $P_u(X_i)\subseteq X_{sm}^+$ when $X_i$ is sm-interesting, and $P_u(X_i)\subseteq X_{lg}^+$ when $X_i$ is lg-interesting.

Our next claim is crucial to the overall argument.  The basic idea is that for certain sm-ineresting $Y_j$, and any lg-interesting $X_i$, then, given any $B'\in \phi_2(X_i)$ and $B''\in \phi_2(Y_j)$, there is a ball $B'''=B'''(B',B'')$  in $Z$, which is not too small, and such that $K_2[B'', B''']$ contains many $xz$ which ``split'' $Y_j\cap B''$, in the sense that both $N(xz)\cap (Y_j\cap B')$ and $\neg N(xz)\cap (Y_j\cap B'')$ are large.   It will also be important to know this produces different $B'''$ as we let $B'$ and $B''$ vary (see (i) and (ii) below).

\begin{claim}\label{cl:splitting}
Suppose $X\neq Y\in \{A,B,C\}$, $Y_j$ is sm-interesting and $X_i$ is lg-interesting.
\begin{enumerate}[label=\normalfont(\arabic*)]
\item Suppose $\rho'_{X_j}\geq \rho'_{Y_i}$.  Then for each $x_0\in P_2(X_i)$, and each $y_0\in P_2(Y_j)$,  there is $f_{X_iY_jZ}(x_0y_0)\in Z$ so that for at least half the $x\in B_{\rho'_{X_i}/2p^4}(x_0)\cap X_i$, and all the $z\in B_{\rho'_{X_i}/2p^4}(f_{Y_jX_iZ}(x_0y_0))$, $xz\in \Mix(2\e^{1/4}/p^5,Y_j)$.  
\item Suppose $\rho'_{Y_j}\geq \rho'_{X_i}$.  Then for each $x_0\in P_2(X_i)$, and each $y_0\in P_2(Y_j)$,  there is $f_{Y_jX_iZ}(y_0x_0)\in Z$ so that for at least half the $y\in B_{\rho'_{Y_j}/2p^4}(y_0)\cap Y_j$, and all the $z\in B_{\rho'_{Y_j}/2p^4}(f_{X_iY_jZ}(y_0x_0))$, $yz\in \Mix(2\e^{1/4}/p^5,X_i)$. 
\end{enumerate}
Moreover, the following hold.
\begin{enumerate}[label=\normalfont(\roman*)]
\item  If $X_i, X_i$ are lg-interesting, and $\rho'_{X_i}\geq \rho'_{X_{i'}}\geq \rho'_{Y_j}$,  then for all $x\in P_2(X_i)$, $x'\in P_2(X_{i'})$ with $d(x,x')\geq 10\rho_{X_{i}}$, and all $y,y'\in P_2(Y_j)$,  
$$
d(f_{X_iY_jZ}(xy),f_{X_{i'}Y_{j}Z}(x'y'))\geq 5\rho'_{X_i}.
$$ 
\item  If $Y_j, Y_{j'}$ are sm-interesting, and $\rho'_{Y_j}\geq \rho'_{Y_{j'}}\geq \rho'_{X_i}$,  then for all $y\in P_2(Y_j)$, $y'\in P_2(Y_{j'})$ with $d(y,y')\geq 10\rho'_{Y_{j}}$, and all $x,x'\in P_2(X_i)$,  
$$
d(f_{Y_jX_iZ}(yx),f_{Y_{j'}X_{i}Z}(y'x'))\geq 5\rho'_{Y_j}.
$$ 
\end{enumerate}
\end{claim}
\begin{proof}
Assume first that $\rho'_{X_i}\geq \rho'_{Y_j}$.  Suppose $y_0\in P_2(Y_j)$ and $x_0\in P_2(X_i)$.  By Claim \ref{cl:rho}, $\rho'_{X_i}\leq \mu^2/2$.  Set $r=\rho_{X_i}'/4p^2$ and let $B'=B_r(y_0)$ and $B''=B_r(x_0)$.  Apply Splitting Axiom 1  to obtain $z_0=f_0(r,x_0,y_0)$ and $z_1=f_1(r,x_0,y_0)$, so that so that for each $\eta\in \{0,1\}$, $B_{r/2}(z_{\eta})\times B'\times B''\subseteq E^{\eta}$ and $d(z_0,z_1)\leq 3r$.  Suppose there were $y\in Y$, $z\in B_{r/2}(z_0)$, and $z'\in B_{r/2}(z_1)$ such that $(y,z)\in \Gamma_1(X_i)$ and $(y,z')\in \Gamma_0(X_i)$, or vice versa.  Then we would have $yzz'\in Q(X_i)$, contradicting  our choice of $\rho_{X_i}$ (since $d(z,z')<4r\leq \rho_{X_i}$).  Thus, for all $y\in Y$, there is $\xi=\xi(u)\in \{0,1\}$ such that 
$$
(\{y\}\times B_{r/2}(z_1))\cup (\{y\}\times B_{r/2}(z_0))\subseteq \Mix(2\e^{1/4},X_i)\cup \Gamma_\xi(X_i).
$$
Let $B(1)=\{y\in B'\cap Y_j: \xi(y)=1\}$ and $B(0)=\{y\in B'\cap Y_j: \xi(y)=0\}$.  Then 
$$
\max\{|B(0)|, |B(1)|\}\geq |B'\cap Y_j|/2.
$$
Choose $\xi^*\in \{0,1\}$ so that $|B(\xi^*)\cap Y_j|\geq |B'\cap Y_j|/2$.  Then for all $y\in B(\xi^*)\cap Y_j$, and all $z\in B_{|1-\xi^*|}(z_{|1-\xi^*|})$, one of the following holds.
\begin{itemize}
\item $yz\in \Mix(2\e^{1/4}, X_i)$. In this case, by definition, 
\[\min\{|N(yz)\cap X_i|, |\neg N(yz)\cap X_i|\}\geq 2\e^{1/4}|X_i|.\] 
\item $yz\in \Gamma_{\xi^*}(X_i)$. In this case $|N^{\xi^*}(yz)|\geq (1-2\e^{1/4})|X_i|$ (by definition of $\Gamma_{\xi^*}(X_i)$).  Since $z\in B_{|1-\xi^*|}(z_{|1-\xi^*|})$, we also know $|\neg N(yz)\cap X_i|\geq |B''\cap X_i|\geq 2\e^{1/4}|X_i|/p^5$, where the last inequality is because $B''\in \phi_2(X_i)$.   
\end{itemize}
This shows $(B(\xi^*)\cap Y_j)\times B_{|1-\xi^*|}(z_{|1-\xi^*|})\subseteq \Mix(2\e^{1/4}/p^5,X_i)$.  Define 
\[f_{X_iY_jZ}(y_0x_0)=z_{|1-\xi^*|}.\]  
This finishes the proof of (1).  The proof of (2) is symmetric, with the roles of $Y$ and $X$ switched everywhere.

For (i), suppose $X_i$ and $X_{i'}$ are lg-interesting, and $\rho'_{X_i}\geq \rho'_{X_{i'}}\geq \rho'_{Y_j}$.  Suppose $x\in P_2(X_i)$, $x'\in P_2(X_{i'})$ and $y,y'\in P_2(Y_j)$ are such that $d(x,x')\geq 10\rho'_{X_i}$.  Since we chose $f_{X_iY_jZ}(xy)$ and $f_{X_{i'}Y_jZ}(x'y')$ using Splitting Axiom 1 with $r=\rho_{X_i}'/4p^2$ and $r=\rho_{X_{i'}}'/4p^2$, respectively, we have by part (c) of that axiom that, since $d(y,y')\leq \rho'_{Y_j}\leq \rho'_{X_i}$, 
$$
d(f_{X_iY_jZ}(xy),f_{X_{i'}Y_{j}Z}(x'y'))\geq d(x,x')-d(y,y')-7(\rho'_{X_i}/4p^4)\geq 5\rho'_{X_i}.
$$

The proof of (ii) is symmetric, with the roles of $X$ and $Y$ switched everywhere, and with Splitting Axiom part (b) in place of part (c).
\end{proof}

We are now ready to embark on the main argument of the proof.  Given $X\in \{A,B,C\}$, let $I_X(0)=K_X^{sm}\cup K_X^{lg}$.  Then for $X\neq Y$, define
\begin{align*}
M^Y(X_i)&=\{Y_j\in I_Y(0): X_iY_j\notin \calS\text{ and }\rho'_{Y_j}\leq \rho'_{X_i}\}.
\end{align*}
We will construct an integer $t^*$, and for each $1\leq u\leq t^*$, a collection of objects $w_u$, $W_u$, $\calS_{W_u}$ and $L_X(u), I_X(u), J_X(u)$ for each  $X\in \{A,B,C\}$ so that the following hold. 
\begin{enumerate}
\item For each $1\leq u\leq t^*$, $W_u\in I_A(0)\cup I_B(0)\cup I_C(0)$, and $w_u\in P_2(W_u)$,
\item For each $1\leq u\leq t^*$ and $X\in \{A,B,C\}$, $L_X(u)=\{W_1,\ldots, W_u\}\cap I_X(0)$,
\begin{align*}
\hspace{20pt}I_X(u)=\{X_s\in I_X(0)\setminus L_X(u): \; & X_s\in M^X(W_u)\\
&\text{ for all }W_u\in L_Y(u)\text{ with } Y\neq X\},
\end{align*}
and
\[\hspace{10pt}J_X(u)=\{X_s\in I_X(u): \text{ for all }X_i\in L_X(u), d_X(x_s,x_i)\geq 10\rho'_{X_i}\}.\]
\item For each $1\leq u\leq t^*$ and $X\in \{A,B,C\}$, if $W_u\in I_X(0)$ and $YZ$ lists $\{A,B,C\}\setminus \{X\}$ in alphabetical order, the following hold.  If $W_u$ is lg-interesting, then
\begin{align*}
\hspace{20pt}\calS_{W_{u}}=\{Y_jZ_k &\in K_2[J_Z(u), J_Y(u)]:  Y_j \in K_Y^{sm},\text{ and } \\
& \exists y\in P_2(Y_j) \text{ such that } \Psi_3(Z_k)\cap B_{\rho'_{W_u}/4}(f_{W_{u}Y_jZ}(w_uy))\neq \emptyset\},
 \end{align*}
 and if $W_u$ is sm-interesting, then 
 \begin{align*}
\hspace{20pt}\calS_{W_{u}}=\{Y_jZ_k &\in K_2[J_Z(u), J_Y(u)]:  Y_j \in K_Y^{lg},\text{ and } \\
& \exists y\in P_2(Y_j) \text{ such that } \Psi_3(Z_k)\cap B_{\rho'_{W_u}/4}(f_{W_{u}Y_jZ}(w_uy))\neq \emptyset\},
 \end{align*}
\item For each $1\leq u\leq t^*$, $\calS_{W_u}\subseteq \calS$ and $|\calS_{W_u}|\geq (1-\sqrt{\e})\e_3^{1/4}\rho'_{W_u}t/2$.
\item  For all $1\leq u\neq u'\leq t^*$, $\calS_{W_u}\cap \calS_{W_{u'}}=\emptyset$.
\item For some $X\in \{A,B,C\}$, $|\bigcup_{X_v\in L_X(t_0)}B'_{X_v}|\geq \e^{1/16}n$. 
\end{enumerate}

We will now construct the objects described above via an inductive argument. For each $i\geq 0$, the $i$th step will ensure that we have defined the desired objects $w_{u}$, $W_{u}$, $\calS_{W_{u}}$ and $L_X(u), I_X(u), J_X(u)$ for each  $X\in \{A,B,C\}$ and all $1\leq u\leq i$. 
\vspace{2mm}

\underline{Step 0:} There is nothing to do, as this case is vacuous.

\underline{Step $i+1$:} Suppose $i\geq 0$, and for each $1\leq u\leq i$, we have chosen $W_u$, $\calS_{W_u}$, and sets $L_U(u), I_U(u), J_U(u)$ for each $U\in \{A,B,C\}$, such that (1)-(4) hold for each $1\leq u\leq i$, such that (5) holds for each $1\leq u\neq u'\leq i$, and such that (6) holds for all $1\leq u<i$ (note the hypotheses say nothing if $i=0$).

If $|\bigcup_{U_s\in L_U(i)}B'_{U_s}|\geq \e^{1/16} n$ for some $U\in \{A,B,C\}$, set $t^*=i$ and end the construction (clearly this will not happen if $i=0$).  Otherwise, we have that for each $U\in \{A,B,C\}$, $|\bigcup_{U_s\in L_U(i)}B_{U_s}|<\e^{1/16}n$.  By the Size Axiom, this implies that for each $U_{v}\in L_U(i)$, $\rho'_{U_{v}}\leq 2p\e^{1/16}n$. Consequently, using (\ref{al:speciallg}),
\begin{align}\label{al:binarylb1}
\nonumber\Big|\bigcup_{U_s\in I_U(i)}U_s\Big|&\geq \Big|\bigcup_{U_s\in I_U(0)}U_s\Big|-|\bigcup_{U_s\in L_U(i)}B'_{U_s}|\\&\nonumber\geq |U_{lg}^+\cup U_{sm}^+|(1-4\e_3^{1/3}) -\e^{1/16}n\\
&\geq |U_{lg}^+\cup U_{sm}^+|(1-5\e_3^{1/3}).
\end{align}
By our assumptions, for each $U\in \{A,B,C\}$ and each $U_{i_s},U_{i_v}\in L_U(i)$, 
$$
B_{10\rho'_{U_{i_s}}}(u_{i_s})\cap B_{10\rho'_{U_{i_v}}}(u_{i_v})=\emptyset.
$$
  Consequently, by the Size Axiom and the definition of $J_U(i)$, 
\begin{align}\label{al:binarylb2}
\nonumber\Big|\bigcup_{U_s\in J_U(i)}U_s\Big|&\geq \Big|\bigcup_{U_s\in I_U(0)}U_s\Big|-|\bigcup_{U_{i_s}\in L_U(i)}B_{10\rho'_{U_{i_s}}}(u_{i_s})|\\
&\nonumber \geq |U_{lg}^+\cup U_{sm}^+|(1- 5\e_3^{1/3}) -10p\e^{1/16}n\\
&\geq |U_{lg}^+\cup U_{sm}^+|(1-6\e_3^{1/3}).
\end{align}
For each $U\in \{A,B,C\}$, set set $r_{i+1}^U=\max \{\rho'_{U_j}: U_j\in J_U(i)\}$.  Then define
$$
r_{i+1}=\max\{r^A_{i+1},r_{i+1}^B, r^C_{i+1}\}.
$$
Choose $X\in \{A,B,C\}$ such that $r_{i+1}^X=r_{i+1}$, and any $j_{i+1}\in J_X(i)$  such that $\rho'_{X_{j_{i+1}}}=r_{i+1}$.  Then define $W_{i+1}:=X_{j_{i+1}}$, $L_X(i+1)=L_X(i)\cup \{W_{i+1}\}$, and for each $U\neq X$, set $L_U(i+1)=L_U(i)$.  For each $U\in \{A,B,C\}\setminus \{X\}$, define  $I_U(i+1)=I_U(i)\cap M^U(W_{i+1})$ and $J_U(i+1)=J_U(i)\cap M^U(W_{i+1})$.  Then set $I_X(i+1)=I_X(i)\setminus \{W_{i+1}\}$ and 
\[J_X(i+1)=\{X_s\in I_X(i+1):  d(u_s,x_{j_{i+1}})\geq 10\rho'_{W_{i+1}}\}.\]
Let $YZ$ list the elements of $\{A,B,C\}\setminus X$ in alphabetical order, and define 
\[K_Y^{sm}(i+1)=\{Y_j\in J_Y(i+1):\text{ $Y_j$ is sm-interesting}\}\]
and
\[K_Y^{lg}(i+1)=\{Y_j\in J_Y(i+1):\text{ $Y_j$ is lg-interesting}\}.\]
Using (\ref{al:specialsm}) and (\ref{al:binarylb2}), and since $J_Y(i+1)=J_Y(i)$,
\begin{align*}
\Big|\bigcup_{Y_j\in K_Y^{sm}(i+1)}Y_j\Big|\geq \Big|\bigcup_{Y_j\in K^{sm}_Y(i+1)\setminus J_Y(i+1)}Y_j\Big|\geq (1-6\e_3^{1/3})|Y_{sm}|.
\end{align*}
A similar argument shows that 
\begin{align*}
\Big|\bigcup_{Y_j\in K_Y^{lg}(i+1)}Y_j\Big|\geq (1-6\e_3^{1/3})|Y_{lg}|.
\end{align*}
Consequently, combining this with the Size and Special Sets Axioms, and the fact that $\e\ll \mu, \tau$, we have
$$
|K_Y^{sm}(i+1)|\geq (\mu n/2p) (1-6\e_3^{1/3})/(3n/t)\geq \mu (1-6\e_3^{1/3})t/6p \geq \mu t/7p,
$$
and similarly $|K_Y^{lg}(i+1)|\geq \tau t/7p$.  By construction, $W_{i+1}$ is either lg-interesting or sm-interesting.  If $W_{i+1}$ is lg-interesting, let $Y_*=Y_{sm}$ and $K^*=K_Y^{sm}(i+1)$.  If $W_{i+1}$ is sm-interesting, let $Y_*=Y_{lg}$ and $K^*=K_Y^{lg}(i+1)$.

Fix $w_{i+1}\in P_2(W_{i+1})$.  Given $Y_j\in K^*$ and $y\in P_2(Y_j)$, consider the ball $B_{r_{i+1}/4}(f_{W_{i+1}Y_jZ}(w_{i+1}y))$ as in Claim \ref{cl:splitting}.  Note that since by construction, $r_{i+1}=\rho'_{W_{i+1}}\geq \rho'_{Z_k}$, if $Z_k\in J_Z(i+1)$, $B''\in \phi_3(Z_k)$, and $B''\cap B_{r_{i+1}/4}(f_{W_{i+1}Y_jZ}(w_{i+1}y))\neq \emptyset$, then $B''\subseteq B_{r_{i+1}/2}(f_{W_{i+1}Y_jZ}(w_{i+1}y))$. Thus, for all such $Z_k$, we have by Claim \ref{cl:splitting} that
$$
|E\cap K_3[Z_k,W_{i+1},Y_j]|/|K_3[Z_k,W_{i+1},Y_j]|\in (\e^{1/4}\mu^4/2, 1-\e^{1/4}\mu^4/2).
$$
Since we assumed $Y_j\in J_Y(i+1)\subseteq M^Y(W_{i+1})$, $W_{i+1}Y_j\notin \calS$.  Therefore, if $W_{i+1}Z_k\notin\calS$, and $Z_k$ is as above, then we would have to have $Y_jZ_k\in \Irr(W_{i+1})$.   Thus, if $Z_k\in J_Z(i+1)$ and $\Psi_3(Z_k) \cap B_{r_{i+1}/4}(f_{W_{i+1}Y_jZ}(w_{i+1}y))\neq \emptyset$, then $Y_jZ_k\in \Irr(W_{i+1})$.  Let $\calS_{W_{i+1}Y_j}$ be the set of pairs arising in this way for $Y_j$, i.e. 
 \begin{align*}
\calS_{W_{i+1}Y_j}=\{Y_jZ_k: Z_k\in J_Z(i+1),&\text{ and for some }y\in P_2(Y_j),\\
& \Psi_3(Z_k)\cap B_{r_{i+1}/4}(f_{W_{i+1}Y_jZ}(w_{i+1}y))\neq \emptyset\}
 \end{align*}
We have shown that for all $Y_j\in K^*$, $\calS_{W_{i+1}Y_j}\subseteq \Irr(W_{i+1})$. Clearly $|\calS_{W_{i+1}Y_j}|\geq s_j$, where we define $s_j$ as follows. 
 \begin{align*}
s_{j}:=\frac{|(\bigcup_{y\in P_2(Y_j)}B_{r_{i+1}/4}(f_{W_{i+1}Y_jZ}(w_{i+1}y))\cap (\bigcup_{Z_k\in J_Z(i+1)}\Psi_3(Z_k))|}{|N/t|}.
 \end{align*}
Note that for any $Y_j\in K^*$, 
\begin{align*}
\Big|\Big(\bigcup_{y\in P_2(Y_j)}B_{r_{i+1}/4}&(f_{W_{i+1}Y_jZ}(w_{i+1}y))\setminus \Big( \bigcup_{Z_k\in J_Z(i+1)}\Psi_3(Z_k)\Big)\Big|\\
&\geq \Big|\bigcup_{y\in P_2(Y_j)}B_{r_{i+1}/4}(f_{W_{i+1}Y_jZ}(w_{i+1}y))\Big|-s_j(N/t).
\end{align*}
By the Size Lemma, for each $y\in P_2(Y_j)$, 
\begin{align*}
&|B_{r_{i+1}/4}(f_{W_{i+1}Y_jZ}(w_{i+1}y))|\geq r_{i+1}n/8p=\rho'_{W_{i+1}}n/8p.
\end{align*}
Thus
\begin{align*}
\Big|\Big(\bigcup_{y\in P_2(Y_j)}B_{r_{i+1}/4}(f_{W_{i+1}Y_jZ}(w_{i+1}y)\Big)\setminus \Big( \bigcup_{Z_k\in J_Z(i+1)}&\Psi_3(Z_k)\Big)\Big|
\\ &\geq \rho'_{W_{i+1}}n/8p-s_j(N/t).
\end{align*}
Now define
$$
\calG^*=\{Y_j\in K^*: s_{j}\geq \rho'_{W_{i+1}}(1-\sqrt{\e})t/2p^3\}.
$$
We claim $|\bigcup_{Y_j\in \calG^*}Y_j|\geq \e_3^{1/4}n$.  Suppose towards a contradiction this is not the case, i.e. that $|\bigcup_{Y_j\in \calG^*}Y_j|< \e_3^{1/4}n$.  Note that if $Y_j\in  K^*\setminus \calG^*$, then by definition of $\calG^*$ and the above,
\[\Big|\Big(\bigcup_{y\in P_2(Y_j)}B_{r_{i+1}/4}(f_{W_{i+1}Y_jZ}(w_{i+1}y)\Big)\setminus \Big( \bigcup_{Z_k\in J_Z(i+1)}\Psi_3(Z_k)\Big)\Big|\]
is at least
\[\rho'_{W_{i+1}}n/8p-\rho'_{W_{i+1}}(1-\sqrt{\e})N/4p^3 \geq \rho'_{W_{i+1}}n(\frac{1}{8p}-\frac{3(1-\sqrt{\e})}{4p^3})\geq \rho'_{W_{i+1}}n/16p^3.\]
By assumption $|\bigcup_{Y_j\in \calG^*}Y_j |\leq \e_3^{1/4}n$, so 
\begin{align}\label{align:lb1}
|\bigcup_{Y_j\in K^*\setminus\calG^* }Y_j|\geq |\bigcup_{Y_j\in K^* }Y_j|-\e_3^{1/4}n\geq |Y_{*}|(1-7\e_3^{1/4}),
\end{align}
where the last inequality also uses the Size Axiom, the Special Sets Axiom, and the fact that $\e\ll \mu$ (if $Y_*=Y_{sm}$) and $\e\ll \tau$ (if $Y_*=Y_{lg}$).  Let $K$ be a maximal $10\rho'_{W_{i+1}}$-separated subset of $\bigcup_{K^*\setminus \calG^*}P_2(Y_j)$.  Since 
$$
\bigcup_{Y_j\in K^*\setminus \calG^*}Y_j\subseteq \bigcup_{y\in K}B_{10\rho'_{X_i}}(y),
$$
we know that by (\ref{align:lb1}), the Size Axiom, and the Special Sets Axiom,  
$$
\min\{\mu,\tau\}(1-7\e_3^{1/4}) n/2p\leq |Y_{*}|(1-7\e_3^{1/4})\leq |K|10p\rho'_{W_{i+1}}n.
$$
Thus $|K|\geq (\rho'_{W_{i+1}})^{-1}(1-7\e_3^{1/4})\min\{\mu,\tau\} /20p^2$.  By definition, for all $y\neq y'\in K$, $d(y,y')>10\rho'_{W_{i+1}}$, and therefore, by Claim \ref{cl:splitting}(i) or (ii),  
\[d(f_{W_{i+1}Y_jZ}(w_{i+1}y), f_{W_{i+1}Y_{j'}Z}(w_{i+1}y'))\geq 5\rho'_{W_{i+1}}.\]  
Consequently,  
$$
B_{r_{i+1}/4}(f_{W_{i+1}Y_jZ}(w_{i+1}y))\cap B_{r_{i+1}/4}(f_{W_{i+1}Y_{j'}Z}(w_{i+1}y'))=\emptyset.
$$
 Thus 
 \begin{align*}
&\hspace{-20pt} |\bigcup_{Z_k\in J_Z(i+1)}\Psi_3(Z_k)|\\
& \leq n-|(\bigcup_{Y_j\in K^*\setminus \calG^*}\bigcup_{y\in K}B_{r_{i+1}/4}(f_{W_{i+1}Y_jZ}(w_{i+1}y)))\setminus \bigcup_{Z_k\in J_Z(i+1)}\Psi_3(Z_k)|\\
&\leq  n-|K|\rho'_{W_{i+1}}n/16 p^3 \\
&\leq  n-(\rho'_{W_{i+1}})^{-1}(1-7\e_3^{1/4})\min\{\mu,\tau\} (20p^2)^{-1}\rho'_{W_{i+1}}n(16 p^3)^{-1}\\
&\leq  n(1-\min\{\mu,\tau\}/240 p^6).
 \end{align*}
However, by the definition of $\Psi_3(Z_k)$ in Claim \ref{cl:balls}, we know that for all $Z_k\in J_Z(i+1)$, $|\Psi(Z_k)|\geq (1-3\e^{1/4})|Z_k|$.  Consequently by (\ref{al:binarylb2}), 
 \begin{align*}
|\bigcup_{Z_k\in J_Z(i+1)}\Psi_3(Z_k)| &\geq (1-6\e^{1/3})|\bigcup_{Z_k\in J_Z(i+1)}Z_k|\\
&\geq (1-6\e^{1/3})(1-\e_3^{1/2})|Z_{lg}^+\cup Z_{sm}^+|\\
&>(1-\min\{\mu,\tau\}/240 p^6)n,
 \end{align*}
where the last inequality uses the fact that $\e\ll \mu,\tau$.  This is a contradiction, so we must have that  $|\bigcup_{Y_j\in \calG^*}Y_j|\geq \e_3^{1/4}n$, so $|\calG^*|\geq \e_3^{1/4}t/3$.  Now define 
$$
\calS_{W_{i+1}}=\bigcup_{Y_j\in \calG^*}\calS_{W_{i+1}Y_j}
$$
Then we have shown that $\calS_{W_{i+1}}\subseteq \calS$ and 
$$
|\calS_{W_{i+1}}|\geq |\calG^*|\rho'_{W_{i+1}}(1-\sqrt{\e})t/2p^3\geq \rho'_{W_{i+1}}\e_3^{1/4}(1-\sqrt{\e})^2t^2/12p^3.
$$
This finishes the definition of $\calS_{W_{i+1}}$.
 
We now show that  for all $u\leq i$,  $\calS_{W_u}\cap \calS_{W_{i+1}}=\emptyset$.  Suppose towards a contradiction, there is some $u\leq i$ with $\calS_{W_u}\cap \calS_{W_{i+1}}\neq \emptyset$, say $Y_jZ_k\in \calS_{W_u}\cap \calS_{W_{i+1}}$.  By construction, this implies $W_u\in L_X(u)$, say $W_u=X_{j_u}$ for some $j_u\neq j_{i+1}$.  We claim that either both $W_u,W_{i+1}\in K_X^{lg}$ or both $W_u,W_{i+1}\in K_X^{sm}$.  Suppose $W_u\in K_X^{sm}$ and $W_{i+1}\in K_X^{lg}$.  Then $W_u\in K_X^{sm}$ implies $Y_j\in K_Y^{lg}$, while $W_{i+1}\in K_X^{lg}$ implies $Y_{j}\in K_Y^{sm}$.  But now $Y_j\in K_Y^{sm}\cap K_Y^{lg}$.  This  implies that $|Y_j\cap Y^+_{sm}|>|Y_j|/2$ and $|Y_j\cap Y^+_{lg}|>|Y_j|/2$, which implies $Y^+_{sm}\cap Y^+_{lg}\neq \emptyset$, contradicting the Disjointness Axiom.   A similar argument shows we cannot have $W_{i+1}\in K_X^{sm}$ and $W_u\in K_X^{lg}$.  

Therefore we have that either $W_u,W_{i+1}\in K_X^{lg}$ or $W_u,W_{i+1}\in K_X^{sm}$. By construction, there is $w_u\in P_2(W_u)$ such that one of the following holds.
\begin{enumerate}
\item Both $W_u,W_{i+1}\in K_X^{lg}$ and $Y_jZ_k\in K_2[K_Y^{sm}(u),J_Z(u)]\cap K_2[K_Y^{sm}(i+1),J_Z(i+1)]$, and 
\[\Psi_3(Z_k)\cap B_{r_{u}/4}(f_{W_{u}Y_jZ}(w_uy))\neq \emptyset\text{ for some }y\in P_2(Y_j)\]
and
\[\Psi_3(Z_k)\cap B_{r_{i+1}/4}(f_{W_{i+1}Y_jZ}(w_{i+1}y))\neq \emptyset\text{ some }y\in P_2(Y_j).\]
 \item Both $W_u,W_{i+1}\in K_X^{sm}$ and $Y_jZ_k\in K_2[K_Y^{lg},J_Z(u)]\cap K_2[K_Y^{lg}(i+1),J_Z(i+1)]$, and 
\[ \Psi_3(Z_k)\cap B_{r_{u}/4}(f_{W_{u}Y_jZ}(w_uy))\neq \emptyset\text{ for some }y\in P_2(Y_j)\]
and
\[ \Psi_3(Z_k)\cap B_{r_{i+1}/4}(f_{W_{i+1}Y_jZ}(w_{i+1}y))\neq \emptyset\text{ some }y\in P_2(Y_j).\]
 \end{enumerate}
 Since $Y_j\in J_Y(i+1)$, we know $\rho'_{Y_j}\leq \min\{\rho'_{W_u},\rho'_{W_{i+1}}\}$.  By construction,  the centers of $B_{W_u}$ and $B_{W_{i+1}}$ have distance at least $10\rho'_{W_u}$ from each other.  Combining this with the fact that for all $y,y'\in P_2(Y_j)$, $d(y,y')\leq \rho'_{Y_j}$, we have by Claim \ref{cl:splitting} that
 \begin{align}\label{align:lg}
 d(f_{Y_jW_uZ}(yw_u),f_{W_{i+1}Y_jZ}(w_{i+1}y'))\geq 5\rho'_{W_u}\geq 5\rho'_{Z_k},
 \end{align}
But (\ref{align:lg}) implies that it is not possible to have both $B_{r_u/4}(f_{W_uY_jZ}(w_uy)))$ and $B_{r_{i+1}/4}(f_{W_{i+1}Y_jZ}(w_{i+1}y'))$  intersect $B_{Z_k}$, a contradiction.   This finishes our verification that $\calS_{W_u}\cap \calS_{W_{i+1}}=\emptyset$, and consequently the  $(i+1)$st step of our construction.

Clearly this will end at some stage $t^*$, where we will have that for some $X\in \{A,B,C\}$, $|\bigcup_{X_v\in L_X(t^*)}B'_{X_v}|\geq \e^{1/16}N$.  Say $1\leq i_1<\ldots<i_{s}\leq t^*$ are such that 
$$
W_{i_1},\ldots, W_{i_s}\subseteq\{X_i: i\in [t]\},
$$
 and $|\sum_{j=1}^sB'_{W_{i_j}}|\geq \e^{1/16}N$.  By the size axiom, this implies $\sum_{j=1}^s\rho'_{W_{i_j}}\geq \e^{1/16}/4p$.  Therefore, we have that 
\begin{align*}
|\calS|&\geq \sum_{j=1}^s|\calS_{W_{i_j}}|\geq (1-2\sqrt{\e})^2\e_3^{1/4}t^2(\sum_{j=1}^s\rho'_{W_{i_s}})/4p^3>\e_3^{3/4} t^2/20p^4>13\e_3t^2,
 \end{align*} 
where the last inequality uses that $\e\ll p^{-1}$ and $t\geq t_0$.  But this contradicts (\ref{al:sigmacont}), which finishes the proof. 
\end{proof}

We obtain as an immediate corollary that $\calH_{\HP}$ and $\calH_{\GS_p}$ require non-binary $\vdisc_3$-error.

\begin{corollary}\label{cor:nonbin}
$\calH_{\HP}$, as well as $\calH_{\GS_p}$ for any prime $p\geq 3$ require non-binary $\vdisc_3$-error.
\end{corollary}
\begin{proof}
This follows immediately from Lemmas \ref{prop:gsspecial} and \ref{prop:hpspecial}, and Theorem \ref{thm:irreg}.
\end{proof}

We will now use this fact to prove Theorems \ref{thm:ternarytriads} and \ref{thm:ternarytriads1}, which state that if $\calH$ is a hereditary $3$-graph property and either $\calH_{\HP}\subseteq \trip(\calH)$ or $\calH_{\GS_p}\subseteq \trip(\calH)$, then $\calH$ requires non-binary $\disc_{2,3}$-error.

\vspace{2mm}

\begin{proofof}{Theorems \ref{thm:ternarytriads} and \ref{thm:ternarytriads1}}
Suppose towards a contradiction there exists a hereditary $3$-graph property $\calH_2$ such that one of (1) or (2) below holds, and such that $\calH_2$ admits binary $\disc_{2,3}$-error. 
\begin{enumerate}
\item $\calH_1=\calH_{\HP}\subseteq \trip(\calH_2)$, or 
\item $\calH_1=\calH_{\GS_p}\subseteq \trip(\calH_2)$.
\end{enumerate} 

By Corollary \ref{cor:nonbin},  there is some $\e_1>0$ such that for all $T, M$, there is $N_{T,M}$ such that for all $n\geq N$, neither $\HP(n)$ nor $\GS_p(n)$ admit an $\e_1$-homogeneous partition $G_{T,M}$ into  at most $T$ parts with binary error.  

Let $0<\e_1\ll \e_1$, and choose $\e_2:\mathbb{N}\rightarrow (0,1)$ as in Proposition \ref{prop:wnipdens} for $\e_1'$ and $k=8$. Define $\e_2':\mathbb{N}\rightarrow (0,1]$ by choosing $\e_2'(x)\ll \e_1'\e_2(x)$ for all $x$.  By assumption, there are $T_1,L_1, M_1$ such that for all $G'=(V',E')\in \calH_2$ with $|V'|\geq M_1$, there is are $1\leq t\leq T_1$, $1\leq \ell\leq L_1$ and a $(t,\ell,\e'_1,\e'_2(\ell))$-decomposition of $V'$ which is $(\e_1',\e_2'(\ell))$-regular with respect to $G'$ with binary $\disc_{2,3}$-error.   

Let $T_2=T_1(\e'_1)^{-1}$ and choose $M_2\gg 3M_1T_1(\e'_1)^{-1}$.  Now fix $n\geq M_2,N_{T_2,M_2}$ and let $G=\HP(n)=(A\cup B\cup C, E)$ if (1) holds, and $G=\GS_p(n)=(A\cup B\cup C, E)$ if (2) holds.  

By assumption, $G\in \trip(\calH_2)$.  By definition this means there is some $H=(V,F)\in \calH_2$ and sets
\begin{align*}
A'&=\{v_{a}: a\in A\}\subseteq V',\text{ }B'=\{v_{b}: b\in B\}\subseteq V',\text{ and }C'=\{v_{c}: c\in C\}\subseteq V',
\end{align*}
such that $v_{a_g}v_{b_{g'}}v_{c_{g''}}\in F$ if and only if $abc\in E$. Note this implies $|V'|\geq n\geq M_1$.  Let $H'=H[A'\cup B'\cup C']$.  Then $H'\in \calH_2$, and $n\leq |V(H')|\leq 3n$.  Since $n\geq M_1$, there are $t_0\leq t\leq T_1$, $\ell_0\leq \ell\leq L_1$ and a $(t,\ell,\e'_1,\e'_2(\ell))$-decomposition $\calP$ of $V(H')$ which is $(\e'_1,\e'_2(\ell))$-regular with respect to $H'$, and which has binary $\disc_{2,3}$-error.  Say $\calP_{vert}=V'_1\cup \ldots \cup V'_t$ and $\calP_{edge}=\{P_{ij}^\alpha:ij\in {[t]\choose 2}, \alpha\leq \ell\}$.  Note that for each $1\leq i\leq t$, $n/t\leq |V_i|\leq 3 n/t$.  Define
\begin{align*}
\Sigma=\{ij\in {[t]\choose 2}: \text{ some }P_{ijk}^{\alpha,\beta,\gamma}\in \triads(\calP)\text{ fails }&\disc_{2,3}(\e_2(\ell),1/\ell)\\&\text{ with respect to }H'\}.
\end{align*}
 By assumption, $|\Sigma|\leq \e'_1 t^2$.  For each $1\leq i\leq t$, let 
$$
A_i=\{a\in A:v_a\in V_i'\},\text{ }B_i=\{b\in B:v_b\in V_i'\}\text{ and }C_i=\{c\in C:v_c\in V_i'\}.
$$
Then $\{A_i, B_i ,C_i: i\in [t]\}$ forms a partition of $A\cup B\cup C$.  Set $t_0=3\lceil \e_1^{-2}\rceil$ and $K=\lfloor 3n/t_0t\rfloor$.  For each $i\in [t]$ and $X\in \{A,B,C\}$, let $X_i=X_{i0}\cup X_{i1}\cup \ldots\cup X_{ik^X_i}$ be a partition with the property that $|X_{i0}|<K$ and so that for each $1\leq s\leq k^X_i$, $|X_{is}|=K$.  Let $k_X=\sum_{i=1}^tk^X_i$ and set $X_0=\bigcup_{i=1}^t X_{i0}$.  Note that by definition of $t_0$ and $K$, $|X_0|\leq tK\leq 3\e^2_1 n$.  Choose $\{X_0^{is}: 1\leq i\leq t,1\leq s\leq k^X_i\}$ any equipartition of  $X_0$ into $k$ parts, and for each $1\leq i\leq t$ and $1\leq s\leq k^X_i$, set $X_{is}'=X_{is}\cup X_{is}^0$.  Clearly $|X_{is}'\setminus X_{is}|\leq 3\e_2^2 |X_{is}|$.  

We now have an partition $\calQ_{vert}$  of $V(G)$ given by
\begin{align*}
\{A_{is}':1\leq i\leq t,1\leq s\leq k^A_i\}\cup \{B_{is}':1\leq i\leq t,&1\leq s\leq k^B_i\}\\
&\cup \{C_{is}':1\leq i\leq t,1\leq s\leq k^C_i\}.
\end{align*}
For each $1\leq \alpha\leq \ell$, $X\neq Y$ in $\{A,B,C\}$, $ij\in {[t]\choose 2}$, $1\leq s\leq k^X_i$, and $1\leq r\leq k^Y_{j}$, set $Q_{ij,sr,XY}^{\alpha}=P_{ij}^\alpha\cap K_2[X_{is}',Y_{jr}']$.  Note that by Lemma \ref{lem:sl}, if $P_{ij}^{\alpha}$ satisfies $\disc_2(\e_2'(\ell);1/\ell)$, then each $Q_{ij,sr,XY}^{\alpha}$ satisfies $\disc_2(\e_2(\ell);1/\ell)$.  

Then for each $X\in \{A,B,C\}$, $1\leq i\leq t$, and $sr\in {[k^X_i]\choose 2}$, choose any partition $K_2[W_{is}',W_{ir}']=\bigcup_{\alpha\leq \ell}Q_{is,ir',XX}^{\alpha}$ so that each $Q_{is,ir',XX}^{\alpha}$ has $\disc_2(\e_2(\ell);1/\ell)$ (such a partition exists by Lemma \ref{lem:3.8}).  We now set 
\begin{align*}
\calQ_{edge}=\{Q_{ij,sr, AB}^{\alpha}: 1\leq i,j\leq t, 1\leq & s\leq k^A_i, 1\leq r\leq k^B_j\} \\
\cup\{Q_{ij,sr, AC}^{\alpha}: 1\leq i,j\leq t, &1\leq s\leq k^A_i, 1\leq r\leq k^C_j\} \\
\cup\{Q_{ij,sr, BC}^{\alpha}&: 1\leq i,j\leq t, 1\leq s\leq k^B_i, 1\leq r\leq k^C_j\} \\
& \cup\{Q_{ij,sr, X}^{\alpha}: X\in \{A,B,C\}, 1\leq s,r\leq k^X_i\}.
\end{align*}
 It is straightforward to check that $\calQ$ is a $(tk,\ell,\e_1,\e_2(\ell))$-decomposition of $A\cup B\cup C$.  We then define 
\begin{align*}
\Sigma'=\{X_{is}Y_{jr}: & X\neq Y \in \{A,B,C\}, ij\in \Sigma, s\in [k^X_i], r\in [k^Y_j]\}\\
&\cup \{X_{is}Y_{ir}: X,Y\in \{A,B,C\}, s\in [k_i^X], r\in [k_i^Y]\}.
\end{align*}

It is straightforward to check that since $\e_1'\ll \e_1$ and by our choice of $t_0$, we have $|\Sigma'|\leq \e_1 |\calQ_{vert}|^2$.  We claim that any triple of sets from $\calQ_{vert}$ avoiding $\Sigma'$ is $\e_1$-homogeneous with respect to $G$.  Clearly this holds of any triple of the form $X_{is}X_{js}Y_{pq}$ for $X,Y\in \{A,B,C\}$ or $X_{is}X_{js}X_{pq}$ for some $X\in \{A,B,C\}$.  Indeed, any such triple will have edge density $0$ in $G$, since $G$ is $3$-partite with parts $A, B, C$.  

So fix some triple $A_{is}B_{jr}C_{pq}\in {\calQ_{vert}\choose 3}$ such that $A_{is}B_{jr}, B_{jr}C_{pq}, A_{is}C_{pq}\notin \Sigma'$, and a triad 
$$
\Gamma_{is,jr,pq}^{\alpha,\beta,\gamma}=(A_{is}\cup B_{jr}\cup C_{pq}, Q_{ij,sr, AB}^{\alpha}\cup Q_{jp,rq, BC}^{\beta}\cup Q_{ip,sq, AC}^{\gamma}).
$$
By construction, $ij, ip, jp\notin \Sigma$, so the triad $P_{ijp}^{\alpha,\beta,\gamma}=(V'_i\cup V'_j\cup V'_p, P_{ij}^{\alpha}\cup P_{jp}^{\beta}\cup P_{ip}^{\gamma})$ satisfies $\disc_{2,3}(\e_1',\e_2'(\ell))$ with respect to $H'$.  Let $d\in [0,1]$ be such that $|E(H')\cap K_3^{(2)}(P_{ijp}^{\alpha,\beta,\gamma})|=d|K_3^{(2)}(P_{ijp}^{\alpha,\beta,\gamma})|$. Now suppose $\Gamma'$ is a subgraph of $\Gamma_{is,jr,pq}^{\alpha,\beta,\gamma}$.  

Consider the graph $R_{ijp}^{\alpha,\beta,\gamma}$ with vertex set $V_i'\cup V_j'\cup V_p'$ and edge set 
$$
\{xyz\in K_3[V_i',V_j',V_p']: x=v_a, y=v_b, \text{ and }z=v_c\text{ for some }abc\in E(\Gamma_{is,jr,pq}^{\alpha,\beta,\gamma})\},
$$
and consider the subgraph $R'$ with $V_i'\cup V_j'\cup V_p'$ and edge set 
$$
\{xyz\in K_3[V_i',V_j',V_p']: x=v_a, y=v_b, \text{ and }z=v_c\text{ for some }abc\in K_3^{(2)}(\Gamma')\}.
$$
Since $R_{is,jr,pq}^{\alpha,\beta,\gamma}$ and $R'$ are subgraphs of $P_{ijp}^{\alpha,\beta,\gamma}$, by definition of $\disc_{2,3}(\e_1',\e_2'(\ell))$,
\[\big||E(H')\cap K_3^{(2)}(R_{is,jr,pq}^{\alpha,\beta,\gamma})|-d|K_3^{(2)}(R_{is,jr,pq}^{\alpha,\beta,\gamma})|\big|\leq \e_1'\frac{1}{\ell^3}|V'_i||V'_j||V'_p|\]
and
\[\big||E(H')\cap K_3^{(2)}(R')|-d|K_3^{(2)}(R')|\big|\leq \e_1'\frac{1}{\ell^3}|V'_i||V'_j||V'_p|.\]
Since $i,j,p$ are distinct, we have that 
$$
|K_3^{(2)}(R_{is,jr,pq}^{\alpha,\beta,\gamma}))|=|K_3^{(2)}(\Gamma_{is,jr,pq}^{\alpha,\beta,\gamma})|\text{ and }|K_3^{(2)}(R'))|=|K_3^{(2)}(\Gamma')|,
$$
as well as
\[|E(H')\cap K_3^{(2)}(R_{is,jr,pq}^{\alpha,\beta,\gamma})|=|E(G)\cap K_3^{(2)}(\Gamma_{is,jr,pq}^{\alpha,\beta,\gamma})|\]
 and
 \[|E(H')\cap K_3^{(2)}(R')|=|E(G)\cap K_3^{(2)}(\Gamma')|.\]
Consequently, if we let $d'$ be such that $|E(G)\cap K_3^{(2)}(\Gamma_{ijp}^{\alpha,\beta,\gamma})|=d'|K_3^{(2)}(\Gamma_{ijp}^{\alpha,\beta,\gamma})|$, then 
\begin{align*}
|d'-d|&\leq \e_1'\frac{1}{\ell^3}|V'_i||V'_j||V'_p|(|K_3^{(2)}(\Gamma_{is,jr,pq}^{\alpha,\beta,\gamma})|)^{-1}\\
&\leq \e_1'\frac{1}{\ell^3}|V_i'||V_j'||V_p'|3\e_1^{-2}\\
&\leq \e_1^8\frac{1}{\ell^3}|A_{is}||B_{jr}||C_{pq}|,
\end{align*}
where the second inequality is from Corollary \ref{cor:counting} and the definition of $\calQ_{vert}$, and the last inequality is since $\e_1'\ll \e_1$.  We also have from the above that 
$$
\big||E(G)\cap K_3^{(2)}(\Gamma')|-d|K_3^{(2)}(\Gamma')|\big|\leq \e_1'\frac{1}{\ell^3}|V'_i||V'_j||V'_p|\leq \e_1^8\frac{1}{\ell^3}|A_{is}||B_{jr}||C_{pq}|,
$$
where the last inequality is since $\e_1'\ll \e_1$ and the definition of $\calQ_{vert}$.  Combining these, we have that 
$$
\big||E(G)\cap K_3^{(2)}(\Gamma')|-d'|K_3^{(2)}(\Gamma')|\big|\leq 2 \e_1^8\frac{1}{\ell^3}|A_{is}||B_{jr}||C_{pq}|,
$$
which shows that $\Gamma_{is,jr,pq}^{\alpha,\beta,\gamma}$ satisfies $\disc_{2,3}(2\e_1^8, \e_2(\ell))$ with respect to $G$.  We can conclude by Proposition \ref{prop:wnipdens} (since $G$ has $\VC$-dimension at most $8$) that 
$$
|E(G)\cap K_3[A_{is},B_{jr},C_{pq}]|/|C_{is}||B_{jr}||C_{pq}|\in [0,\e_1^2)\cup (1-\e_1^2,1].
$$
This shows that any triple form $\calQ_{vert}$ avoiding $\Sigma'$ is $\e^2_1$-homogeneous with respect to $G$.  As defined $\calQ_{vert}$ may not be itself an equipartition, although its parts will differ in size by at most $\e_1^2n/|\calQ_{vert}|$.  It is now not too difficult to show that one can redistribute some vertices to make $\calQ_{vert}$ into an equipartition without changing the densities within any of the triads by very much.  The result will be an $\e_1$-homogeneous equipartition of $G$ with binary error and at most $T_2$ parts, a contradiction. 
\end{proofof}

\vspace{2mm}

We now prove as a Corollary \ref{cor:disc23hom}, which says that a property admitting binary, zero, or linear $\disc_{2,3}$-error is $\disc_{2,3}$-homogeneous.

\vspace{2mm}
\begin{proofof}{Corollary \ref{cor:disc23hom}}
Suppose $\calH$ admits zero $\disc_{2,3}$-error.  By Corollary \ref{cor:dischbark}, there is $k$ such that $\Hbar(k)\notin \trip(\calH)$.  By Fact \ref{fact:vc2universal}, $\calH$ has finite $\VC_2$-dimension, so by Theorem \ref{thm:dischom}, it is $\disc_{2,3}$-homogeneous.  Suppose now $\calH$ admits binary $\disc_{2,3}$-error.  By Theorem \ref{thm:ternarytriads}, there is $k\geq 1$ such that $\HP(k)\notin \trip(\calH)$.  By Fact \ref{fact:vc2universal}, $\calH$ has finite $\VC_2$-dimension, so by Theorem \ref{thm:dischom}, it is $\disc_{2,3}$-homogeneous.  Suppose now $\calH$ admits linear $\disc_{2,3}$-error.  By Theorem \ref{thm:FOP}, there is $k\geq 1$ such that $F(k)\notin \trip(\calH)$. By Fact \ref{fact:vc2universal}, $\calH$ has finite $\VC_2$-dimension, so by Theorem \ref{thm:dischom}, it is $\disc_{2,3}$-homogeneous.
\end{proofof}
\vspace{2mm}

A similar argument will be important in our proof of Theorem \ref{thm:disc3}.  In particular, we will use the following lemma.

\begin{lemma}\label{lem:disc3lem}
There are $\e_1>0$, $\e_2:\mathbb{N}\rightarrow (0,1]$, and $t_0,\ell_0\geq 1$ such that for all $T$, $L$ and $N$, there is $n\geq N$ such that for all $t_0\leq t\leq T_0$, and $\ell_0\leq \ell\leq T_0$, any $(t,\ell,\e_1,\e_2(\ell))$-decomposition of $V(\HP(n))$ which is $(\e_1,\e_2(\ell))$-regular, has a $\disc_{3}$-irregular triad.
\end{lemma}
\begin{proof}
Since $\calH_{\HP}$ has finite $\VC$-dimension it is $\disc_{2,3}$-homogeneous by Theorem \ref{thm:dischom}.  Suppose towards a contradiction that for all $\e_1>0$, $\e_2:\mathbb{N}\rightarrow (0,1]$, and $t_0,\ell_0\geq 1$ there are $T$, $L$ and $N$, such that for all $n\geq N$ there is $t_0\leq t\leq T_0$, and $\ell_0\ell\leq L_0$, any $(t,\ell,\e_2,\e_1)$-decomposition of $V(\HP(n))$ which is $(\e_1,\e_2)$-regular and has no $\disc_{3}$-irregular triad.  By Theorem \ref{thm:ternarytriads}, there are $\mu_1,\mu_2, t_0,\ell_0$ such that for all $T$, $L$ and $N$, there is $n\geq N$ such that for all $t_0\leq t\leq T_0$, and $\ell\leq T_0$, any $(t,\ell, \mu_1,\mu_2(\ell))$-decomposition of $V(\HP(n))$ which is $(\mu_1,\mu_2(\ell))$-regular, has non-binary $\disc_{2,3}$-error. 

By our assumption, and because $\calH_{\HP}$ is $\disc_{2,3}$-homogeneous, there are $T$, $L$ and $N$, so that for all $n\geq N$ there is $t_0\leq t\leq T_0$, and $\ell_0\leq \ell\leq L_0$, and a $(t,\ell,\mu_1,\mu_2(\ell))$-decomposition of $V(\HP(n))$ which is $(\mu_1,\mu_2(\ell))$-regular, with no $\disc_3$-irregular triads, and which is $\disc_{2,3}$-homogeneous.   Let $\Sigma$ be the set of pairs $V_iV_j$ where at least $(\e'_1)^{1/2} |V_i||V_j|$ of the $xy\in K_2[V_iV_j]$ are in a $\disc_2$-irregular $P\in \calP_{edge}$.  By assumption, $|\Sigma|\leq (\e_1')^{1/2}t^2$.  If $V_iV_j, V_iV_k, V_jV_k\notin \Sigma$, then by Proposition \ref{prop:wnipdens}, $|E\cap K_3[V_i,V_j,V_k]|/K_3[V_i,V_j,V_k]|\in [0,\e_1)\cup (1-\e_1,1]$.  But now $\calP$ is $\vdisc_3$-homogeneous with binary error, a contradiction.
\end{proof}

We now prove Theorem \ref{thm:disc3}.

\vspace{2mm}

\begin{proofof}{Theorem \ref{thm:disc3}}
Fix a hereditary $3$-graph property $\calH$, and assume it admits binary $\disc_{2,3}$-error.  Fix $\e_1>0$, $\e_2:\mathbb{N}\rightarrow (0,1]$ nonincreasing, and $t_0,\ell_0\geq 1$.  Let $\e_1'=\e_1^2$.  By assumption, there is $T$, $L$ and $N$ so that if $H\in \calH$ has at least $N$ vertices, then there exist $\ell_0\leq \ell\leq L$, $t_0\leq t\leq T$, and $\calP$ a $(t,\ell,\e'_1,\e_2(\ell))$-decomposition of $V(H)$ which is $(\e'_1,\e_2(\ell))$-regular, with binary $\disc_{2,3}$-error with respect to $H$.  Say $\calP$ consists of $\{V_i: i\in [t]\}$ and $\{P_{ij}^{\alpha}: ij\in {[t]\choose 2}, \alpha\leq \ell\}$.  Let $\Sigma\subseteq {[t]\choose 2}$ be such that every $\disc_{2,3}$-irregular triad of $\calP$ uses a pair $V_iV_j$ with $ij\in \Sigma$, and where $|\Sigma|\leq \e_1' t^2$.  

Let $\ell_1=\lfloor \ell/2\rfloor$ and $\ell_2=\lceil \ell/2 \rceil$.  For each $i\in [t]$, let $V_i=V_{i1}\cup V_{i2}$ be an equipartition.  For each $i<j$ with $ij\in \Sigma$, use Lemma \ref{lem:3.8} to choose partitions $V_{i1}\times V_j=\bigcup_{\alpha\leq \ell_1}Q_{ij}^{\alpha}$ and $V_{i2}\times V_j=\bigcup_{\ell_1+1\leq \alpha\leq \ell}Q_{ij}^{\alpha}$ so that for each $1\leq \alpha \leq \ell_1$, $Q_{ij}^{\alpha}$ has $\disc_2(\e_1(\ell_1);1/\ell_1)$ in $V_{i1}\times V_j$, and for each $\ell_1+1\leq \alpha \leq \ell$, $Q_{ij}^{\alpha}$ has $\disc_2(\e_1(\ell_2); 1/\ell_2)$ in $V_{i1}\times V_j$.  Clearly each $Q_{ij}^{\alpha}$ fails $\disc_2(\e_2(\ell);1/\ell)$.  Now define $\calQ$ to have vertex partition  $\{V_i: i\in [t]\}$ and edge partition 
$$
\Big\{P_{ij}^{\alpha}: ij\in {[t]\choose 2}\setminus \Sigma, \alpha\leq \ell\Big\}\cup \{Q_{ij}^{\alpha}: ij\in \Sigma, \alpha\leq \ell\}.
$$
By construction, any triad from $\calQ$ involving a pair $V_iV_j\in \Sigma$ will be $\disc_2$-irregular, and thus cannot be $\disc_3$-irregular.  Any triad from $\calQ$ not involving a pair $V_iV_j$ from $\Sigma$ will be $\disc_{2,3}$-regular, and thus cannot be $\disc_3$-irregular.  Thus $\calQ$ has no $\disc_3$-irregular triads.

Conversely, suppose $\calH$ has the property that for all for all $\ell_0,t_0\geq 1$, $\e_1>0$, all $\e_2:\mathbb{N}\rightarrow (0,1]$, there are $n_0$ and $T,L$ such that for all $n\geq n_0$, every $G\in \calH$ with at least $n_0$ vertices has an $(\e_1,\e_2)$-regular, $(t,\ell,\e_2,\e_1)$-decomposition for some $t_0\leq t\leq T$, and $\ell_0\leq \ell\leq L$, such that no triad of $\calP$ is $\disc_{3}$-irregular.  By Lemma \ref{lem:disc3lem}, there is $k$ such that $\HP(k)\notin \trip(\calH)$, so by Fact \ref{fact:vc2universal}, $\calH$ has finite $\VC_2$-dimension, and so by Theorem \ref{thm:dischom}, $\calH$ is $\disc_{2,3}$-homogeneous.  Fix $\e_1>0$, $\e_2:\mathbb{N}\rightarrow (0,1]$, and $t_0,\ell_0\geq 1$.  Without loss of generality, let us assume $\e_2$ is non-increasing (e.g. by replacing $\e_2$ with a non-increasing function bounded above by the original $\e_2$).

By the above, there exist $n_0$, $L_0$ and $T_0$ so that for all $n\geq n_0$, every $G\in \calH$ with at least $n_0$ vertices, there are $t_0\leq t\leq T_0$ and $\ell_0\leq \ell \leq L_0$ and a $(t,\ell,\e_1,\e_2(\ell))$-decomposition of $V(G)$ which is $(\e_1,\e_2(\ell))$-regular and $\disc_{2,3}$-homogeneous with respect to $G$, and such that, moreover, no triad from $\calP$ is $\disc_{3}$-irregular with respect to $G$.

Assume $n\geq n_0$ and $H=(V,E)\in  \calH_n$.  By assumption, there exists $t_0\leq t\leq T_0$, $\ell_0\leq \ell\leq L_0$ and $\calP$, a $(t,\ell,\e'_1,\e_2'(\ell))$-decomposition of $V$ which is $(\e_1',\e_2'(\ell))$-regular and $\disc_{2,3}$-homogeneous, and such that no triad of $\calP$ is $\disc_{3}$-irregular with respect to $H$.  Let $\Sigma$ be the set of pairs $V_iV_j$ where at least $(\e'_1)^{1/2} |V_i||V_j|$ of the $xy\in K_2[V_iV_j]$ are in a $\disc_2$-irregular $P\in \calP_{edge}$.  By assumption, $|\Sigma|\leq (\e_1')^{1/2}t^2$. 

Set $s=\lceil (\e_1')^{1/2}\ell \rceil$.  Fix $ij\notin \Sigma$.  By reindexing if necessary, we may assume that $P_{ij}^{s+1},\ldots, P_{ij}^\ell$ all satisfy $\disc_2(\e_2'(\ell);1/\ell)$.  Set $E'_0=\bigcup_{i=1}^sP_{ij}^{\alpha}$, and for each $s+1\leq u\leq \ell$, set $E'_{u-s}=P_{ij}^{u}$.  Note the density of $E_0'=\frac{s}{\ell}(1\pm s\e'_2(\ell))$.

Apply Lemma \ref{lem:3.8} to $E'_0$ with $\e=\delta=\e_2(\ell)$ and $p=1/(\ell-s)$ to obtain $R_{ij}^1,\ldots, R_{ij}^{\ell-s}$, where each $R_{ij}^u$ has $\disc_2(2s\e'_2(\ell)^{1/2}; s/(\ell(\ell-s)))$.  For each $1\leq m\leq \ell-s$, let $Q_{ij}^m=S_{ij}^m\cup R_{ij}^m$.  By Fact \ref{fact:adding}, each $Q_{ij}^m$ has $\disc_2(2s\e'_2(\ell)^{1/2}; 1/\ell+s/(\ell (\ell-s))$.  It is not hard to check that $1/\ell+s/(\ell (\ell-s)))=1/(\ell-s)$. Thus, each $Q_{ij}^m$ has $\disc_2(\e_2(k'); 1/k')$ (using choice of $\e_2'$), where $k'=\ell-s$.  Given $ij\notin \Sigma$, choose any partition $V_i\times V_j=\bigcup_{\alpha\leq k'}Q_{ij}^{\alpha}$.

Let $\calQ$ consist of $\{V_1,\ldots, V_t\}$, and $\{Q_{ij}^{\alpha}: ij\in {[t]\choose 2}\}$.  By construction, this is a $(t,k',\e_1,\e_2(k'))$-decomposition of $V$.  Suppose $ij, ik, jk\notin \Sigma$, and consider a triad of the form $Q_{ijk}^{\alpha,\beta,\gamma}=(V_i\cup V_j\cup V_k, Q_{ij}^{\alpha}\cup Q_{ik}^{\beta}\cup Q_{jk}^{\gamma})$, for some $1\leq\alpha,\beta,\gamma\leq k'$.  By construction, there was some triad of the form $P_{ijk}^{\alpha',\beta',\gamma'}$  which had $\disc_2(\e_2'(\ell);1/\ell)$ such that 
$$
|Q_{ij}^{\alpha}\setminus P_{ij}^{\alpha}|<\frac{s}{\ell k'}(1+2s\e_2'(\ell))|V_i||V_j|<\e^2_1|Q_{ij}^{\alpha}|,
$$
where the inequality is because $s\leq (\e_1')^{1/2}\ell+1$ and $\e_1'\ll \e_1$.   By assumption, $P_{ijk}^{\alpha',\beta',\gamma'}$ was a $\disc_{2,3}$-regular triad of $\calP$.  We show $Q_{ijk}^{\alpha,\beta,\gamma}$ is a $\disc_{2,3}$-regular triad of $\calQ$.  By Proposition \ref{prop:homimpliesrandome}, it suffices to show the density on $Q_{ijk}^{\alpha,\beta,\gamma}$ is in $[0,\e_1^2)\cup (1-\e_1^2,1]$. Since $P_{ijk}^{\alpha',\beta',\gamma'}$  has $\disc_2(\e_2'(\ell);1/\ell)$, it is by assumption, $\disc_{2,3}$-regular, and thus also $\disc_{2,3}$-homogeneous.  Thus there is $u\in \{0,1\}$ such that 
$$
|E^u\cap K_3^{(2)}(P_{ijk}^{\alpha',\beta',\gamma'})|\geq (1-\e_1')|K_3^{(2)}(P_{ijk}^{\alpha',\beta',\gamma'})|.
$$
Then, using Corollary \ref{cor:counting}, 
\begin{align*}
|E^u\cap K_3^{(2)}(Q_{ijk}^{\alpha,\beta,\gamma})|&\geq |K_3^{(2)}(Q_{ijk}^{\alpha,\beta,\gamma})|-\e_1'|K_3^{(2)}(P_{ijk}^{\alpha',\beta',\gamma'})|-\frac{9\e_1^2}{(k')^3}|V_i||V_j||V_k|\\
&\geq \frac{|V_i||V_j||V_k|}{(k')^3}\Big( 1-\e_1^2-\e_1'\frac{(k')^3}{\ell^3}-9\e_1^2\Big)\\
&\geq (1-\e_1)| K_3^{(2)}(Q_{ijk}^{\alpha,\beta,\gamma})|,
\end{align*}
where the inequality is by our choice of $\e_1'\ll \e_1$ and since $k'\geq (1-\e_1)\ell$. Thus $Q_{ijk}^{\alpha,\beta,\gamma}$ has $\disc_{2,3}(\e_1,\e_2(k'))$.  Thus in $\calQ$ all irregular triads use a pair from $\Sigma$, so it has binary $\disc_{2,3}$-error, as desired.
\end{proofof}

\vspace{2mm}

%% file: appendices.tex
\chapter{Facts about $\sim$-classes}\label{app:simclasses}

This appendix shows various facts about $\sim$-classes.   We begin with Proposition \ref{prop:simclasses2}, which says that all the decomposition properties of interest in this paper are closed under $\sim$-classes.  This boils down to the following two lemmas.

\begin{lemma}\label{lem:appsim1}
For all $\e>0$ and $T\geq 1$, there is $\delta>0$ so that the following holds.  Suppose $H=(V,E)$, and $H'=(V,E')$ are two graphs on the same vertex set $V$, and assume $H$ and $H'$ are $\delta$-close.  Suppose $V_1,\ldots, V_t$ is an equipartition of $V$ with $t\leq T$.  If $ijk\in {[t]\choose 3}$ and $V_iV_jV_k$ has $\vdisc_3(\e)$ in $H$ with density $r$, then it has $\vdisc_3(3\e)$ in $H'$ with density $r\pm \e$. 
\end{lemma}
\begin{proof}
Fix $\e>0$ and $T\geq 1$.  Define $\delta=\e^2/T^3$ and $N=\e^{-1}T$.  Suppose $H=(V,E)$, and $H'=(V,E')$ are two graphs on the same vertex set $V$ of size $n\geq N$, and assume $H$ and $H'$ are $\delta$-close.  Suppose $V_1,\ldots, V_t$ is an equipartition of $V$ with $t\leq T$.  Suppose $ijk\in {[t]\choose 3}$ and $V_iV_jV_k$ has $\vdisc_3(\e)$ in $H$ with density $r$. Then
\[|E'\cap K_3[V_i,V_j,V_k]|=|E\cap K_3[V_i,V_j,V_k]|\pm \delta |V|^3,\]
which in turn equals
\[r|V_i||V_j||V_k|\pm \delta n^3=(r\pm \delta t^3)|V_i||V_j||V_k|.\]
By definition of $\delta$ and since $t\leq T$, this shows the density of $V_iV_jV_k$ in $H'$ is within $\e$ of $r$.  Suppose now that $V_i'\subseteq V_i$, $V_j'\subseteq V_j$, and $V_k'\subseteq V_k$.  Then, since $H$ and $H'$ are $\delta$ close, and $V_iV_jV_k$ has $\vdisc_3(\e)$ in $H$, 
\begin{align*}
|E'\cap K_3[V'_i,V'_j,V'_k]|&=|E\cap K_3[V'_i,V'_j,V'_k]|\pm \delta |V|^3\\
&=(r\pm \e)|V'_i||V'_j||V'_k|\pm \delta n^3\\
&=(r\pm 2\e)|V_i||V_j||V_k|.
\end{align*}
Since $V_iV_jV_k$ has density $r\pm \e$, this shows $V_iV_jV_k$ satisfies $\vdisc_3(3\e)$ in $H'$.
\end{proof}

\begin{lemma}\label{lem:appsim2}
For all $\e>0$, there is $\e^*_2:\mathbb{N}\rightarrow (0,1]$ such that for all $\e_2:\mathbb{N}\rightarrow (0,1]$ satisfying $\e_2(x)\leq \e_2^*(x)$ for all $x\geq 1$, and all $L,T\geq 1$, there is $\delta>0$ so that the following holds.  Suppose $H=(V,E)$, and $H'=(V,E')$ are two graphs on the same vertex set $V$ of size $n\geq N$, such that $H$ and $H'$ are $\delta$-close.  Assume $\calP$ is a $(t,\ell,\e_1,\e_2(\ell))$-decomposition of $V$.  Then for any $G_{ijk}^{\alpha,\beta,\gamma}\in \triads(\calP)$ satisfying $\disc_{2,3}(\e_1,\e_2(\ell))$ with respect to $H$ with density $r$, then $G_{ijk}^{\alpha,\beta,\gamma}$  satisfies $\disc_{2,3}(3\e_1,\e_2(\ell))$ with respect to $H'$ with density $r\pm 2\e_1$.
\end{lemma}

\begin{proof}
Fix $\e>0$.  Define $\e_2^*:\mathbb{N}\rightarrow (0,1]$ so that $\e_2^*(x)\ll \e/x^2$ for all $x\geq 1$.  Now suppose $\e_2:\mathbb{N}\rightarrow (0,1]$ satisfies $\e_2(x)\leq \e_2^*(x)$ for all $x\geq 1$, and $L,T\geq 1$. Choose $\delta\ll \e_1^4\e_2(L)^2L^{-3}T^{-3}$.

Suppose $H=(V,E)$, and $H'=(V,E')$ are two graphs on the same vertex set $V$, and assume $H$ and $H'$ are $\delta$-close.  Suppose $\calP$ is a $(t,\ell,\e_1,\e_2(\ell))$-decomposition of $V$.  Suppose $G\in \triads(\calP)$ has $\disc_{2,3}(\e_1,\e_2(\ell))$ in $H$ with density $r$.  Then using Corollary \ref{cor:counting},
$$
|(E\Delta E')\cap K_3(G)|\leq \delta n^3\leq \e^2_1|K_3(G)|.
$$
Consequently, if $r'$ is the density of $G$ with respect to $H'$, then 
\[|r'-r|\leq \frac{|(E\Delta E')\cap K_3(G)|}{|K_3(G)|}\leq \e_1^2.\]  
Now suppose $G'\subseteq G$ is a subgraph.  Then since $G$ has $\disc_{2,3}(\e_1,\e_2(\ell))$ with respect to $H$ with density $r$, 
\begin{align*}
\big||E\cap K_3^{(2)}(G')|-r|K_3^{(2)}(G')|\big|\leq \e_1 \frac{n^3}{\ell^3t^3}.
\end{align*}
Thus
\begin{align*}
&\hspace{-20pt}\big||E'\cap K_3^{(2)}(G')|-r'|K_3^{(2)}(G')|\big|\\
&\leq \big||E'\cap K_2^{(2)}(G')|-|E\cap K_3^{(2)}(G')|\big|+\big||E\cap K_3^{(2)}(G')|-r'|K_3^{(2)}(G')|\big|\\
&\leq \delta n^3+\big||E\cap K_2^{(2)}(G')|-r|K_3^{(2)}(G')|\big|+|r-r'||K_3^{(2)}(G')|\\
&\leq \delta n^3 +\e_1 \frac{n^3}{\ell^3t^3}+\e_1^2|K_3^{(2)}(G')|\\
&\leq 3\e_1 \frac{n^3}{\ell^3t^3}.
\end{align*}
This shows that $G$ has $\disc_{2,3}(3\e_1,\e_2(\ell))$ with respect to $H'$ with density $r'$.
\end{proof}

\vspace{2mm}

\begin{proofof}{Proof of Proposition \ref{prop:simclasses2}} Assume $\calH$ and $\calH'$ are hereditary $3$-graph properties and $\calH$ is close to $\calH'$.  Suppose $\calH'$ satisfies one of the following.
\begin{enumerate}
\item[(a)] $\calH'$ is $\vdisc_3$-homogeneous.
\item[(b)] $\calH'$ admits binary $\vdisc_3$-error.
\item[(c)] $\calH'$ admits zero $\vdisc_3$-error.
\item[(d)] $\calH'$ is $\disc_{2,3}$-homogeneous.
\item[(e)] $\calH'$ admits linear $\disc_{2,3}$-error.
\item[(f)] $\calH'$ admits binary $\disc_{2,3}$-error.
\item[(g)] $\calH'$ admits zero $\disc_{2,3}$-error.
\end{enumerate}
Then we show $\calH$ satisfies each of these, respectively.  Assume first that (a), (b), or (c) hold.  Fix $\e>0$ and $t_0\geq 1$.  By assumption, there exist $T=T(\e/2,t_0)$, $N=N(\e/2,t_0)$ be such that for every $n\geq N$, the following holds.  For every $G\in\calH'_n$, there is a decomposition $V(G)=V_1\cup \ldots \cup V_t$ for some $t_0\leq t\leq T$, such that the following holds (respectively).
\begin{enumerate}
\item There is $\Sigma\subseteq {[t]\choose 3}$ such that $|\Sigma|\leq \e t^3/5$ and for all $ijk\notin \Sigma$, $V_iV_jV_k$ is $\e/5$-homogeneous with respect to $G$.
\item There is $\Gamma\subseteq {[t]\choose 2}$ such that $|\Gamma|\leq \e t^2/5$ and for all $ijk\in {[t]\choose 3}$ with $ij, jk, ik\notin \Gamma$, $(V_i,V_j,V_k)$ has $\vdisc_3(\e/5)$ with respect to $G$.  
\item For all $ijk\in {[t]\choose 3}$,   $V_iV_jV_k$ satisfies $\vdisc_3(\e/5)$ with respect to $G$. 
\end{enumerate}
Let $\delta$ be as in Lemma \ref{lem:appsim1} for $\e/5$ and $T$.  Choose $N'$ sufficiently large so that for every $n\geq N'$ and $G\in \calH_n$, there is $G'\in \calH'_n$ such that $G$ and $G'$ are $\delta$-close.  Now assume  $n\geq N'$ and $G=(V,E)\in \calH$ with $|V|=n$.  Then there is $G'=(V,E')\in \calH'$ which is $\delta$-close to $G$.  Further, there are $t_0\leq t\leq T$ and an equipartition  $V=V_1\cup \ldots \cup V_t$ such that the following holds (respectively).
\begin{enumerate}
\item There is $\Sigma\subseteq {[t]\choose 3}$ such that $|\Sigma|\leq \e t^3/5$ and for all $ijk\notin \Sigma$, $V_iV_jV_k$ is $\e/5$-homogeneous with respect to $G$.
\item There is $\Gamma\subseteq {[t]\choose 2}$ such that $|\Gamma|\leq \e t^2/5$ and for all $ijk\in {[t]\choose 3}$ with $ij, jk, ik\notin \Gamma$, $(V_i,V_j,V_k)$ has $\vdisc_3(\e/5)$ with respect to $G$.  
\item For all $ijk\in {[t]\choose 3}$,   $V_iV_jV_k$ satisfies $\vdisc_3(\e/5)$ with respect to $G$. 
\end{enumerate}
By Lemma \ref{lem:appsim1}, the following holds (respectively).
\begin{enumerate}
\item For all $ijk\notin \Sigma$, $V_iV_jV_k$ is $\e$-homogeneous with respect to $G'$.
\item For all $ijk\in {[t]\choose 3}$ with $ij, jk, ik\notin \Gamma$, $(V_i,V_j,V_k)$ has $\vdisc_3(\e)$ with respect to $G'$.  
\item For all $ijk\in {[t]\choose 3}$,   $V_iV_jV_k$ satisfies $\vdisc_3(\e)$ with respect to $G'$. 
\end{enumerate}
Thus we have shown that $\calH$ satisfies (a), (b), or (c), respectively. 

Assume now that (d), (e), (f), or (g) hold, respectively.  Fix $\e_1>0$, $\e_2:\mathbb{N}\rightarrow (0,1]$, and $t_0,\ell_0\geq 1$.  Let $\e_2^*$ be as in Lemma \ref{lem:appsim2} for $\e_1/5$.  Let $\e_2'(x)=\min\{\e_2(x), \e_2^*(x)\}$.  Given $\ell_0,t_0\geq 1$, choose $N, T,L$ such that for every $G=(V,E)\in \calH'$ on at least $N$ vertices,  there are $t_0\leq t\leq T$ and $\ell_0\leq \ell\leq L$ and an $(\e_1/5,\e'_2(\ell),t,\ell)$-decomposition $\calP$ of $V$  such that the following hold (respectively).
\begin{enumerate}
\item There is $\Sigma\subseteq {[t]\choose 3}$ such that $|\Sigma|\leq \e_1 t^3/5$ and for all $ijk\notin \Sigma$, every triad $P\in \triads(\calP)$ is $\e_1/5$-homogeneous with respect to $G$.
\item There is $\Sigma \subseteq {[t]\choose 3}$ such that $|\Sigma|\leq \e_1 t^3/5$ and for all $ijk\in {[t]\choose 3}\setminus \Sigma$, every triad $P\in \triads(\calP)$ satisfies $\disc_{2,3}(\e_1/5, \e_2(\ell))$ with respect to $G$.  
\item There is $\Gamma\subseteq {[t]\choose 2}$ such that $|\Gamma|\leq \e_1 t^2/5$ and for all $ijk\in {[t]\choose 3}$ with $ij, jk, ik\notin \Gamma$, every triad $P\in \triads(\calP)$ satisfies $\disc_{2,3}(\e_1/5, \e_2(\ell))$ with respect to $G$.  
\item Every triad $P\in \triads(\calP)$ satisfies $\disc_{2,3}(\e_1/5, \e_2(\ell))$ with respect to $G$. 
\end{enumerate}
Let $\delta$ be as in Lemma \ref{lem:appsim2} for $e_1/5$, $\e'_2$, $T$ and $L$.   Let $N$ be such that any $G\in \calH$ on at least $N$ vertices is $\delta$-close to some element of $\calH'$.  Now assume $n\geq N$ and $G=(V,E)\in \calH$ with $|V|=E$. Then there is $G'=(V,E')\in \calH'$ which is $\delta$-close to $G$.  By assumption there are $t_0\leq t\leq T$ and $\ell_0\leq \ell\leq L$ and a $(t,\ell,\e_1/5,\e'_2(\ell))$-decomposition $\calP$ of $V$ such that the following hold (respectively).
\begin{enumerate}
\item There is $\Sigma\subseteq {[t]\choose 3}$ such that $|\Sigma|\leq \e_1 t^3/5$ and for all $ijk\notin \Sigma$, every triad $P_{ijk}^{\alpha,\beta,\gamma}\in \triads(\calP)$ is $\e_1/5$-homogeneous with respect to $G'$.
\item There is $\Sigma \subseteq {[t]\choose 3}$ such that $|\Sigma|\leq \e_1 t^3/5$ and for all $ijk\in {[t]\choose 3}\setminus \Sigma$, every triad $P\in \triads(\calP)$ satisfies $\disc_{2,3}(\e_1/5, \e'_2(\ell))$ with respect to $G'$.  
\item There is $\Gamma\subseteq {[t]\choose 2}$ such that $|\Gamma|\leq \e_1 t^2/5$ and for all $ijk\in {[t]\choose 3}$ with $ij, jk, ik\notin \Gamma$, every triad $P\in \triads(\calP)$ satisfies $\disc_{2,3}(\e_1/5, \e'_2(\ell))$ with respect to $G'$.  
\item Every triad $P\in \triads(\calP)$ satisfies $\disc_{2,3}(\e_1/5, \e'_2(\ell))$ with respect to $G'$. 
\end{enumerate}
By Lemma \ref{lem:appsim2}, the following hold (respectively).
\begin{enumerate}
\item For all $ijk\notin \Sigma$, every $P\in \triads(\calP)$ is $\e_1$-homogeneous with respect to $G'$.
\item For all $ijk\in {[t]\choose 3}\setminus \Sigma$, every $P\in \triads(\calP)$ satisfies $\disc_{2,3}(\e_1, \e'_2(\ell))$ with respect to $G'$.  
\item For all $ijk\in {[t]\choose 3}$ with $ij, jk, ik\notin \Gamma$, every $P\in \triads(\calP)$ satisfies $\disc_{2,3}(\e_1, \e'_2(\ell))$ with respect to $G'$.  
\item Every $P\in \triads(\calP)$ satisfies $\disc_{2,3}(\e_1/5, \e'_2(\ell))$ with respect to $G'$. 
\end{enumerate}
This shows that $\calH$ has (d), (e), (f), or (g), respectively.
\end{proofof}

\vspace{2mm}

We now prove Theorem \ref{thm:simclasswnip}, whose proof is very similar to that of Theorem \ref{thm:simclassws}.

\vspace{2mm}

\begin{proofof}{Theorem \ref{thm:simclasswnip}} Suppose $\trip(\calH)$ contains $\Ubar(n)$ for all $n$.  Fix a SNIP property $\calH'$.  We show $\calH$ is far from $\calH'$.  By assumption, there is some $k$ so that $U^*(k)\notin \trip(\calH')$.  Let $\e=(1/2k)^{2^k+k+1}/|\calB(U^*(k))|$ and let $N\gg 2^k\e^{-1}$.  For all $n\geq N$, there is $G=(V,E)\in \calH$ such that $\trip(G)$ contains $\Ubar(n)$ and an induced sub-$3$-graph.  This means there are $\{a_i,c_i:i\in [n]\}\subseteq V$ and $\{b_S: S\subseteq \{c_1,\ldots, c_n\}\}\subseteq V$ such that $a_ib_Sc_t\in E$ if and only if $c_t\in S$.  Set $n_1=\lfloor n/(k+1)\rfloor$ and choose $C_1,\ldots, C_k$ to be disjoint subsets of $\{c_1,\ldots, c_{kn_1}\}$ of size $n_1$.  For each $X\subseteq [k]$, let 
$$
B_X=\{b_S: S\subseteq [n]\text{ and for each }i\in X, C_i\subseteq S\text{ and for each }i\notin X, C_i\cap S=\emptyset\}.
$$
 Note $|B_X|\geq 2^k\cdot 2^{n-kn_1}\geq 2^{n_1}$, so we can take $B_X'\subseteq B_X$ of size $n_1$.  We construct many induced subgraphs of $G$ as follows.
 \begin{enumerate}
 \item Choose $a_i$ for some $i\in [n]$.  There are $n$ choices.
 \item Choose some $(c_{i_1},\ldots, c_{i_k})\in C_1\times \ldots \times C_k$.  There are $n_1^k$ choices. 
 \item For each $X\subseteq [k]$, choose some $d_X\in B'_X$.  There are $n_1^{2^k}$ choices.
 \item Put $G[\{a_i\}\cup \{c_{i_1},\ldots, c_{i_k}\}\cup \{d_X: X\subseteq [k]\}]$ in $S$.
 \end{enumerate}
By construction, every element of $S$ is isomorphic to some element in $\calB(U^*(k))$ and 
$$
|S|\geq n_1^{k+2^k+1}>\e |\calB(U^*(k))|n^{k+2^k+1},
$$
where the inequality is by definition of $n_1$ and $\e$.  By the pigeonhole principle, there is some $G(n)\in \calB(U^*(k))$ so that $S$ contains at least $\e n^{k+2^k+1}$ elements isomorphic to $G(n)$.  By the pigeonhole principle, there is $G'\in \calB(U^*(k))$ so that $G'=G(n)$ for arbitrarily large $n$.  By definition, this implies $G'\in \Delta(\calH)$.  Clearly $G'\notin \Delta(\calH')$ (since $U^*(k)\notin \trip(\calH')$), so $\calH$ is far from $\calH'$.
 
 Conversely, suppose $\calH$ is far from every SNIP property $\calH'$.  By Lemma \ref{lem:clean}, $\calH'$ is SNIP if and only if, for some $k\geq 1$, $\calH'$ contains no clean copies of $U^*(k)$.  Thus, since $\calH$ is far from every SNIP property, we know that for all $k\geq 1$, there is $\delta>0$ so that for arbitrarily large $n$, there is $G(n)\in \calH$ and $\Gamma(n)\in\calB(U^*(k))$ so that $G(n)$ is   not $\delta/|\calB(U^*(k))|$-close to being $\Gamma(n)$-free.  For each $\Gamma\in \calB(U^*(k))$, let $c_k(\Gamma)$, $N_k(\Gamma)$ be as in Theorem \ref{thm:indrem} for $\Gamma$ and $\delta/|\calB(U^*(k))|$-close.  Then set $c_k=\min\{c_k(\Gamma):\Gamma\in \calB(U^*(k))\}$ and $N_k=\max\{N_k(\Gamma):\Gamma\in \calB(U^*(k))\}$.  Choose any $N>N_k$.  Then there is some $n\geq N$, so that that $G(n)$ is not $\delta$-close to being $\Gamma(n)$-free.  Say $V(\Gamma(n))=\{a,b_i, c_S: i\in [k], S\subseteq [k]\}$ and $ab_ic_S\in E(\Gamma(n))$ if and only if $i\in S$.  By Theorem \ref{thm:indrem},  $G$ contains at least $c_k n^{k+1+2^k}$ induced sub-$3$-graphs isomorphic to $\Gamma(n)$. Let $S$ be the set of these induced sub-$3$-graphs.  For each $H\in S$, there is a vertex $a_H$ and sets $D_H=\{d^H_i: i\in [k]\}$, $E_H=\{e_S^H: S\subseteq [k]\}$, such that $a_Hd_H^ie_S^H\in E$ if and only if $i\in S$.  Since $|S|\geq c_k n^{k+1+2^k}$, there are sets $D$ and $E$ and $X\subseteq S$ of size at least $c_k n/k!(2^k)!$ such that $(D,E)=(D_H,E_H)$ for all $H\in X$.  Then we must have that $\{a_H:H\in X\}|=|X|$. Let $Y\subseteq \{a_H: H\in X\}$ any set of size $k$.  Then $G[Y\cup D\cup E]\in \calH$ and $\trip(G[Y\cup D\cup E])\cong \Ubar(k)$.  Thus $\Ubar(k)\in \trip(\calH)$.
\end{proofof}

\chapter{Lemmas about combining and refining decompositions}\label{app:slicing}

In this appendix we prove the results stated in Section \ref{subsec:intersecting}.  We will also prove Theorem \ref{thm:counting} here, as it uses some of these lemmas.  

We begin with the proof of Corollary \ref{cor:counting}.

\vspace{2mm}

\begin{proofof}{Corllary \ref{cor:counting}}
Fix $t\geq 2$, $\e>0$ and $r\in (0,1]$. Without loss of generality, assume $\e\leq 1/2$.  Choose $\mu\ll_t \gamma \ll_t \e,r$.  

 Now assume $V=U_1\cup \ldots \cup U_t$ and $G=(V,E)$ is a $t$-partite graph with vertex partition $U_1,\ldots, U_t$, such that for each $1\leq i<j\leq t$, $G[U_i,U_j]$ has $\disc_2(\mu;r)$.  For each $1\leq i<j\leq t$, let $d_{ij}$ be such that $|E\cap K_2[U_i,U_j]|=d_{ij}|U_i||U_j|$.  By assumption, each $d_{ij}=r\pm \mu$.  By Proposition \ref{prop:counting}, 
$$
\Big|K_t(G)|-\prod_{ij\in {[t]\choose 2}}d_{ij}\prod_{i=1}^t|U_i|\Big|\leq \gamma \prod_{i=1}^t|U_i|.
$$
Let $d=\prod_{ij\in {[t]\choose 2}}d_{ij}$.  Note $d=(1\pm \mu)^{t\choose 2}r^{t\choose 2}=(1\pm \e/2)r^{t\choose 2}$, where the last equality is since $\mu\ll \e$.    Thus
$$
|K_t(G)|=(d \pm \gamma)\prod_{i=1}^t|U_i|=d\Big(1 \pm \gamma d^{-1} \Big)\prod_{i=1}^t|U_i|
$$
Since $d=(1\pm \e/2)r^{t\choose 2}$ and $\gamma \ll r,\e$, $\gamma d^{-1}\in (-2\e^2,2\e^2)$.  Thus 
$$
|K_t(G)|=(1\pm 2\e^2)d\prod_{i=1}^t|U_i|=(1\pm 2\e^2)(1 \pm \e/2)r^{t\choose 2}\prod_{i=1}^t|U_i|.
$$
Since $\e\leq 1/2$ this is equal to $(1\pm \e)r^{t\choose 2}\prod_{i=1}^t|U_i|$, as desired.
\end{proofof}

\vspace{2mm}

We now prove Lemma \ref{lem:subpairs}.

\vspace{2mm}
\begin{proofof}{Lemma \ref{lem:subpairs}}
Fix $k\geq 1$ and $\delta_3>0$.  Let $\delta_3'=\delta_3/32k^3$.  Suppose now $d_2,d_3\in (0,1]$ and $0<\delta_2\leq d_2^3/8$.  Set $\gamma=\delta_2/2k^3$ and let $\mu=\mu(3,\gamma,d_2)$ and $N_1=N_1(3,\gamma,d_2)$ as in Proposition \ref{prop:counting}.  Let $\delta'_2\ll \mu/4k^2, \delta_2/4k^2, d_2\delta_2/(6k^4)$, and choose $N\gg k, N_1$.

Suppose $V=V_1\cup V_2\cup V_3$ is a vertex set of size $n\geq N$, and $G=(V_1\cup V_2\cup V_3, E)$ is a $3$-partite graph.  Assume $H=(V,R)$ is a $3$-graph, and let $H'=(V, R\cap K_3^{(2)}(G))$.  Assume that for each $1\leq i< j\leq 3$, $G[V_i,V_j]$ has $\disc_2(\delta'_2;d_2)$, and $(H',G)$ has $\disc_{2,3}(\delta'_2,\delta'_3)$ with density $d_3$.  Suppose that for each $i\in [3]$, $V_i'\subseteq V_i$ satisfies $|V_i'|=m'_i\geq m_i/k$, where $m_i:=|V_i|$.  Set $G'=G[V'_1,V'_2,V'_3]$ and $H''=(V_1'\cup V_2'\cup V_3', R\cap K_3^{(2)}(G'))$.   

For each $1\leq i< j\leq 3$, let $d_{ij}$ be the density of $G[V_i,V_j]$.  By assumption $d_{ij}=d_2\pm \delta_2'$.  By Lemma \ref{lem:sl}, for each $ij\in {[3]\choose 2}$, $G'[V_i,V_j]$ has $\disc_2(2k^2\delta'_2;d'_{ij})$ for some $d'_{ij}=d_{ij}\pm k^2\delta'_2$.  Since $\delta_2'\ll  \delta_2/4k^2, \mu$, this means $G'[V_i',V_j']$ has $\disc_2(\mu;d_2)$.  Thus, by Proposition \ref{prop:counting}, we have
\begin{align*}
|K_3^{(2)}(G')|=d_2^3(1\pm \gamma)m'_1m'_2m'_3.
\end{align*}
Since $(H',G)$ has $\disc_{2,3}(\delta'_2,\delta'_3)$, 
$$
||R\cap K_3^{(2)}(G')|-d_3|K_3^{(2)}(G')||\leq \delta'_3 d_2^3 m_1m_2m_3.
$$

Therefore if $d_3'$ is such that $|R\cap K_3^{(2)}(G')|=d_3'|K_3^{(2)}(G')|$, then
$$
\Big|d_3'-d_3\Big|\leq \frac{\delta'_3 d_2^3 m_1m_2m_3}{|K_3^{(2)}(G')|}\leq \frac{\delta_3'd_2^3m_1m_2m_3}{d_2^3(1-\gamma)m'_1m'_2m'_3}\leq 2\delta_3',
$$
where the last inequality is by assumption, because $\gamma\leq 1/2$, and because each $m_i'\geq m_i/k$.  Suppose $G''\subseteq G'$.  Since $(H',G)$ has $\disc_{2,3}(\delta_2',\delta_3')$,
$$
\big||R\cap K_3^{(2)}(G'')|-d_3|K_3^{(2)}(G'')|\big|\leq \delta'_3 d_2^3 m_1m_2m_3\leq \delta_3d_2^3m'_1m'_2m'_3/2,
$$  
where the equality is by definition of $\delta_3'$ and since each $m_i\geq km_0$. Thus
\begin{align*}
&\hspace{-20pt}\big||R\cap K_3^{(2)}(G'')|-d_3'|K_3^{(2)}(G'')|\big|\\
&\leq\delta_3d_2^3m'_1m'_2m'_3/2+\big|d_3|K_3^{(2)}(G'')|-d_3'|K_3^{(2)}(G'')|\big|\\
&\leq \delta_3d_2^3m'_1m'_2m'_3/2+2\delta_3'|K_3^{(2)}(G'')|\\
&\leq \delta_3d_2^3m'_1m'_2m'_3/2+2\delta_3'd_2^3(1+ \gamma)m'_1m'_2m'_3\\
&\leq \delta_3d_2^8m'_1m'_2m'_3,
\end{align*}
where the last inequality is by definition of $\delta_2'$ and $\gamma$.  This finishes our verification that $(H'',G')$ has $\disc_{2,3}(\delta_3,\delta_2)$.
\end{proofof}
\vspace{2mm}

\begin{proofof}{Fact \ref{fact:adding}}
Assume $E_1,E_2$ are disjoint subsets of $K_2[U,V]$.  We begin with part (a).  Assume $(U\cup V, E_1)$ has $\disc_2(\e_1;d_1)$, $(U\cup V,E_2)$ has $\disc_2(\e_2;d_2)$.  We show $(U\cup V,E_1\cup E_2)$ has $\disc_2(\e_1+\e_2;d_1+d_2)$.  Suppose $U'\subseteq U$ and $V'\subseteq V$.  By assumption, $||E_1\cap K_2[U',V']|-d_1|U'||V'||\leq \e_1 |U|V|$, and $||E_2\cap K_2[U',V']|-d_2|U'||V'||\leq \e |U||V|$.  Then since $E_1\cap E_2=\emptyset$,
\begin{align*}
&\hspace{-20pt}|(E_1\cup E_2)\cap K_2[U',V']|\\
&=|E_1\cap K_2[U',V']|+|E_2\cap K_2[U',V']|-|(E_1\cap E_2)\cap K_2[U',V']|\\
&=|E_1\cap K_2[U',V']|+|E_2\cap K_2[U',V']|.
\end{align*}
By above, this is at most $(d_1+d_2)|U'||V'|+\e_1|U||V|+\e_2|U||V|$ and at least $(d_1+d_2)|U'||V'|-\e_1|U||V|-\e_2|U||V|$ since  $(U\cup V, E_1)$ has $\disc_2(\e_1;d_1)$, $(U\cup V, E_2)$ has $\disc_2(\e_2;d_2)$.  Thus $||(E_1\cup E_2)\cap K_2[U',V']|-(d_1+d_2)|U'||V'||\leq (\e_1+\e_2)|U||V|$, as desired.

We now show (b).  Assume $(U\cup V, E_1\cup E_2)$ has $\disc_2(\e;d)$, $(U\cup V,E_1)$ has $\disc_2(\e_1;d_1)$.  We show $(U\cup V,E_2)$ has $\disc_2(\e+\e_2;d-d_1)$.  Suppose $U'\subseteq U$ and $V'\subseteq V$.  Then
\begin{align*}
&\hspace{-10pt}\big||E_2\cap K_2[U',V']|-(d-d_1)|U'||V'|\big|\\
&=\big||(E_1\cup E_2)\cap K_2[U',V']|-|E_1\cap K_2[U',V']|-(d-d_1)|U'||V'|\big|\\
&\leq \big||(E_1\cup E_2)\cap K_2[U',V']|-d|U'||V'|\big|+\big||E_1\cap K_2[U',V']|-d_1|U'||V'|\big|.
\end{align*}
By assumption, this is at most $\e |U||V|+\e_1|U||V|=(\e+\e_1)|U||V|$.  This shows $(U\cup V,E_2)$ has $\disc_2(\e+\e_2;d-d_1)$.
\end{proofof}

\vspace{3mm}

Our next goal is to prove Lemma \ref{lem:refinement}. This will require two additional tools.  The first is a multi-colored version of Szemer\'{e}di's regularity lemma (see e.g.  \cite{Komlos.1996, Frankl.2002}).

\begin{theorem}[Multi-colored regularity lemma]\label{thm:mcreg}
For all $r\geq 1$, $\e>0$, and $m\geq 1$, there is $M=M(\e,m)$ such that the following holds.  Suppose $G=(V,E_1,\ldots, E_r)$ where $|V|\geq M$ and $E_1,\ldots, E_r$ are disjoint subsets of ${V\choose 2}$.  Then for any equipartition $\calP=\{ V_1, \ldots, V_m\}$ of $V$, there is are $m'\leq M$ and an equipartition, $\calP'=\{V_1', \ldots, V_{m'}'\}$ which refines $\calP$, such that all but at most $\e(m')^2$ pairs $(V_i',V_j')$, $(V_i'\cup V_j', E_u)$ satisfies $\disc_2(\e)$ for all $1\leq u\leq r$.  
\end{theorem}

Note that $M=M(\e,m)$ does not depend on the number of colors.  In fact, it is the same bound as in the usual regularity lemma.  Indeed, it is straightforward to prove Theorem \ref{thm:mcreg} by following the usual proof of the regularity lemma, with a slightly altered mean square density.  The mean square density  increases at each step in an identical way to the usual proof, so the bound ends up identical.  We now give a sketch of this proof, as we could not find this explicitly in the literature.

 We will follow closely the proof of the usual regularity lemma (see, for example \cite{Alon.2000}).  First, suppose that $r\geq 1$, for each $u\in [r]$,  $(V,E_u)$ is a graph on the vertex set $V$, where $|V|=n$, and $u\neq v$ implies $E_u\cap E_v=\emptyset$.  Given $X,Y\subseteq V$, define 
$$
d_u(X,Y)=\frac{|E_u\cap K_2[X,Y]|}{|X||Y|}.
$$ 
Note that $d_1(X,Y)+\ldots+d_r(X,Y)\leq 1$. We then define, for each $u\in [r]$,
$$
\msd_u(X,Y)=\frac{|X||Y|}{n^2}d^2_u(X,Y).
$$
Here, $\msd$ stands for ``mean square density.''  Note this is just the usual mean square density for the   graph $(V,E_u)$.  We then define
$$
\msd(X,Y)=\sum_{u=1}^r \msd_u(X,Y)=\frac{|X||Y|}{n^2}\sum_{u=1}^r d^2_u(X,Y).
$$
Note that since $0\leq d_1(X,Y)+\ldots +d_r(X,Y)\leq 1$, $\sum_{u=1}^rd_u^2(X,Y)\leq 1$, so $\msd(X,Y)\leq \frac{|X||Y|}{n^2}$.  If $\calP_X$ and $\calP_Y$ are partitions of $X$ and $Y$, define, for each $u\in [r]$,
$$
\msd_u(\calP_X,\calP_Y)=\sum_{U\in \calP_X}\sum_{W\in \calP_Y}\msd_u(U,W),
$$
 and set 
 $$
\msd(\calP_X,\calP_Y)=\sum_{u=1}^r \msd_u(\calP_X,\calP_Y).
$$  
For a partition $\calP=\{V_1,\ldots, V_k\}$ of $V$, let $\msd_u(\calP)=\msd_u(\calP,\calP)$ and $\msd(\calP)=\msd(\calP,\calP)$.  Note that
$$
\msd(\calP)=\sum_{i=1}^k\sum_{j=1}^k\frac{|V_i||V_j|}{n^2}(\sum_{u=1}^r d^2_u(V_i,V_j))\leq \sum_{i=1}^k\sum_{j=1}^k\frac{|V_i||V_j|}{n^2}=1
$$

\begin{lemma}
For any partitions $\calP_X$ of $X$ and $\calP_Y$ of $Y$, $\msd(\calP_X,\calP_Y)\geq \msd(X,Y)$.
\end{lemma}
\begin{proof}
It is well known that for each $1\leq u\leq r$, $\msd_u(\calP_X,\calP_Y)\geq \msd_u(X,Y)$  (see e.g. \cite[Corollary 3.7]{Gowers.20063gk}).  By definition, this implies $\msd(\calP_X,\calP_Y)\geq \msd(X,Y)$.
\end{proof}

\begin{lemma}\label{lem:refinereg}
If $(X\cup Y, E_u\cap K_2[X,Y])$ does not satisfy $\disc_2(\e)$, then there are $X_1\subseteq X$ and $Y_1\subseteq Y$, such that
$$
\msd(\{X_1,X\setminus X_1\},\{Y_1,Y\setminus Y_1\})> \msd(X,Y)+\e^4\frac{|X||Y|}{n^2}.
$$
\end{lemma}
\begin{proof}
Suppose $(X\cup Y, E_u\cap K_2[X,Y])$ fails $\disc_2(\e)$.  Then it is easy to check that $(X,Y)$ is not $\e$-regular, say this is witnessed by $X_1\subseteq X$ and $Y_1\subseteq Y$.  It is well known this implies that $\msd_u(\{X_1,X\setminus X_1\},\{Y_1,Y\setminus Y_1\})> \msd_u(X,Y)+\e^4\frac{|X||Y|}{n^2}$ (see e.g. \cite{Gowers.20063gk}).  Therefore
\begin{align*}
&\hspace{-20pt}\msd(\{X_1,X\setminus X_1\},\{Y_1,Y\setminus Y_1\})\\
&=\sum_{s=1}^r \msd_s((\{X_1,X\setminus X_1\},\{Y_1,Y\setminus Y_1\})\\
& > \msd_u(X,Y)+\e^4\frac{|X||Y|}{n^2}+\sum_{s\in [r]\setminus \{u\}}\msd_s(\{X_1,X\setminus X_1\},\{Y_1,Y\setminus Y_1\})\\
&\geq \msd_u(X,Y)+\e^4\frac{|X||Y|}{n^2}+\sum_{s\in [r]\setminus \{u\}}\msd_s(X,Y)\\
&=\msd(X,Y)+\e^4\frac{|X||Y|}{n^2},
\end{align*}
where the last inequality is due to Lemma \ref{lem:refinereg}.
\end{proof}

One can now deduce Theorem \ref{thm:mcreg} exactly as in usual proofs of the regularity lemma (see e.g. \cite{Rodl.2010}, \cite{Alon.2000}).  The second tool used in the proof of Lemma \ref{lem:refinement}  is  Lemma \ref{lem:even} below, which is a corollary of Lemma \ref{lem:3.8}.

\begin{lemma}\label{lem:even}
Suppose $\e_1>0$, $\ell\geq 1$, and $0<\e_2<\e_1/(64 \ell')$, where $\ell':=[32 \ell^2\e_1^{-2}]$.  There is $m_0=m_0(\e_1,\e_2, \ell, \delta)$ such that the following holds.  Suppose $|U|=|V|=m\geq m_0$, $s\leq \ell$, and $K_2(U,V)=E_0\cup \ldots \cup E_s$ is a partition such that for each $i\in [s]$,  $(U\cup V, E_i)$ has $\disc(\e_2)$, and $|E_0|\leq \e_1 m^2$.  Then there is a partition 
$$
K_2(U,V)=E_{0,1}'\cup \ldots E_{0,k}'\cup E'_1\cup \ldots \cup E'_{k'},
$$
where $k+k'=\ell'$, such that the following hold.
 \begin{enumerate}[label=\normalfont(\roman*)]
 \item For each $1\leq i\leq k$, and $1\leq \alpha\leq k'$, $(U\cup V, E_{0,i}')$ and $(U\cup V, E_\alpha)$ have $\disc_2(\e_2^{1/3};1/\ell')$.
\item For each $1\leq \alpha\leq k'$, there is $1\leq j\leq s$ such that $|E_j|\geq (\sqrt{2/32}\e_1/s) m^2$ and $|E_\alpha'\setminus E_j|<(2\e^{1/2}_2/\ell') m^2$.  Moreover, if $k'<\ell'$, then $E_\alpha'\setminus E_j=\emptyset$.
 \item If $E_0'=\bigcup_{j=1}^kE_{0,j}'$, then $|E_0'|\leq 2\e_1 m^2$.
 \end{enumerate}
\end{lemma}
\begin{proof}
 Fix $\e_1>0$, $\ell\geq 1$, and $0<\e_2<\e_1/(64 \ell')$, where $\ell':=[32 \ell^2\e_1^{-2}]$. Set $\delta=\e_2/2(\ell')^2$, and let $m_0=m_0(\e_2/2\ell',\ell',\delta/2\ell')$ be as in Lemma \ref{lem:3.8}. Suppose $s\leq \ell$, $|U|=|V|=m\geq m_0$, and $K_2(U,V)=E_0\cup \ldots \cup E_s$ is a partition such that for each $1\leq i \leq s$,  $(U\cup V, E_i)$ has $\disc_2(\e_2)$, and such that $|E_0|\leq \e_1 m^2$.  For each $i\in \{0,\ldots, s\}$, let $\rho_i=|E_i|/m^2$, and define
$$
\Sigma=\Big\{i\in [s]: \rho_i\leq \sqrt{\frac{2}{\ell'}}\Big\}.
$$
Given $i\in [s]\setminus \Sigma$, set $p(i)=1/\rho_i\ell'$ and $u(i)=[ 1/p(i)]$.  Observe that for all $i\in [s]\setminus \Sigma$, $\rho_i\geq \sqrt{2/\ell'}$ implies $p(i)<\rho_i/2$, and further, $\e_2\leq \e_1$ implies $\rho_i\geq 2\e_2$.  Therefore, for each $i\in [s]\setminus \Sigma$, Lemma \ref{lem:3.8} implies  there is a partition $E_i=R_{i,0}\cup \ldots \cup R_{i,u(i)}$ such that $|R_{i,0}|<p(i)\rho_i (1+\delta)m^2$ and such that for each $1\leq j\leq u(i)$, $R_{i,j}$ has $\disc_2(\e_2)$ and density $\rho_ip(i)(1\pm \delta)=\frac{1}{\ell'}(1\pm \delta)$. 

Let $k'=\sum_{i\in [s]\setminus \Sigma}u(i)$. Since $\delta<1/(\ell')^2$ and 
$$
m^2\geq |\bigcup_{i=1}^s \bigcup_{j=1}^{u(i)}R_{i,j}|\geq (1-k'\delta)\frac{k'}{\ell'}m^2,
$$
we must have that $k'\leq \ell'$.  Fix an enumeration 
$$
\{R_{i,j}: i\in [s]\setminus \Sigma,1\leq j\leq u(i)\}=\{\Gamma_1,\ldots, \Gamma_{k'}\},
$$
and set $\Gamma_0= E_0\cup \bigcup_{i\in [s]\setminus \Sigma}R_{i,0}\cup \bigcup_{i\in  \Sigma}E_i$.  Observe that
\begin{align*}
|\Gamma_0|&\leq \e_1m^2+(s-|\Sigma|)(\frac{1}{\ell'}(1+\delta)m^2)+|\Sigma|\sqrt{\frac{2}{\ell'}}m^2\\
&\leq m^2(\e_1+\frac{s}{\ell'}(1+\delta)+|\Sigma|(\sqrt{\frac{2}{\ell'}}- \frac{1}{\ell'}(1+\delta))\\
&\leq 4\e_1 m^2,
\end{align*}
where the inequality is because $s\leq \ell\leq \e_1\ell'$, and $\sqrt{2}|\Sigma|/\sqrt{\ell'}\leq \sqrt{2}s\e_1/\ell\leq 3\e_1$.   On other hand, since for each $i\in [k']$, $\Gamma_i$ has density $1/\ell' (1\pm \delta)$, we have that 
\begin{align}\label{gamma}
|\Gamma_0|=\frac{\ell'-k'}{\ell'}(1\pm k'\delta)m^2.
\end{align}
If $k'=\ell'$, then (\ref{gamma}) implies that $|\Gamma_0|\leq k'\delta m^2<\e_2m^2/2\ell'$.   In this case, choose any equipartition of $\Gamma_0$ into $\ell'$ pieces, say $\Gamma_0=E_1''\cup \ldots \cup E_{\ell'}''$, and for each $\alpha\in [\ell']$, set $E'_\alpha=\Gamma_\alpha\cup E_\alpha''$.  Then set $k=0$, so  $k+k'=\ell'$.  By construction, for each $\alpha\in [\ell']$, there is $j\in [s]$ such that $|E_{\alpha}'\setminus E_j|=|E_\alpha''|\leq \e_2 m^2/2\ell'$, so (ii) holds.  Further, since each $E_\alpha$ has $\disc_2(\e_2;1/\ell')$ and $|E_{\alpha}'\Delta E_{\alpha}|\leq \e_2 m^2/2\ell'$, it follows that $E_\alpha'$ has $\disc_2(2\e^{1/2}_2)$ and density $1/\ell'\pm 2\sqrt{\e_2}$, so (i) is satisfied.   Since (iii) holds trivially, we are done.

If $k'=\ell'-1$, set $k=1$, so $k+k'=\ell'$.  Define $E_{0,1}=\Gamma_0$, and for each $1\leq i\leq \ell'$, set $E_i'=\Gamma_i$.  Note that by Fact \ref{fact:adding}, $\Gamma_1\cup \ldots \cup \Gamma_{k'}$ has $\disc_2(k'\e_2;k'/\ell')$.  Since $E_{0,1}$ is the compliment of this union, it has $\disc_2(k'\e_2;1-k'/\ell')$.  Consequently, $E_{0,1}$ has $\disc_2(\e_2^{1/2};1/\ell')$ (since $k'=\ell'-1$).  Thus (i)-(iii) hold for $E'_{0,1}\cup E'_1\cup \ldots \cup E'_{k'}$. 

Finally, assume $k'\leq \ell'-2$.  In this case, set $k=\ell'-k'$ and define $\rho_0'=|\Gamma_0|/m^2$.  Since $|\Gamma_0|\leq 2\e_1 m^2$, $\rho_0'\leq 2\e_1$. By  (\ref{gamma}), this implies $(1-2\e_1)\ell'\leq k'$.  Let $\mu=\e_2^2/2\ell'$, $u=k$, and $p=1/u$.  Using that $(1-2\e_1)\ell'\leq k'$ and (\ref{gamma}), we have that $p<\rho_0'/2$.  By definition, $2\mu\leq p^{-1}$.  Thus we may apply Lemma \ref{lem:3.8} to $\Gamma_0$ with $\mu$, $p$, $u$, and $\delta'=\delta/(\ell')^2$.  From this we obtain a partition $\Gamma_0=E_{0,1}\cup \ldots \cup E_{0,u}$ such that each $E_{0,i}$ has $\disc_2(\mu)$ and density $\rho_0'\frac{1}{u}(1\pm \delta')=\frac{1}{\ell'}(1\pm \delta)$.  Thus each $E_{0,i}$ has $\disc_2(\e_2^{1/2};1/\ell')$, as desired.  Thus (i)-(iii) are satisfied with this $k+k'=\ell'$ and $E_{0,1}\cup \ldots \cup E_{0,k}\cup \Gamma_1\cup \ldots \cup \Gamma_{k'}$. 
\end{proof}

We now prove Lemma \ref{lem:refinement}.

\vspace{2mm}

\begin{proofof}{Lemma \ref{lem:refinement}} Fix $\e_1>0$ and $\e_2:\mathbb{N}\rightarrow (0,1]$.  Without loss of generality, let us assume $\e_2$ is non-increasing (e.g. by replacing $\e_2$ with a non-increasing function bounded above by the original $\e_2$). Choose $\mu_1\ll \e_1$.  Define $\mu_2:\mathbb{N}\rightarrow (0,1]$ by setting $\mu_2(x)\ll \mu_1\e_2(\mu_1^{-1}x^2)/x^2$ for all $x\geq 1$.  Let $M_{Sz}(x,y)$ denote the bound from Theorem \ref{thm:mcreg} with parameters $x$ and $y$.  Then set $f(x_1,y_1,x_2,y_2)=32(y_1y_2)^2\mu_1^{-1}$, and set $g(x_1,y_1,x_2,y_2)=\mu_1^{-1}M_{Sz}(\mu_1^{-1}x_1x_2, \mu_2f(x_1,y_1,x_2,y_2))$. 

Now fix $\ell, t, r, s\geq 1$.  Set $L=g(t,\ell, s,r)$, $T= f(t,\ell, s, r)$, and choose $N\gg T, L, \mu_2(L)^{-1}$.  

Suppose $|V|=n\geq N$, $\calP$ is a $(t,\ell)$-decomposition of $V$ consisting of $\calP_{vert}=\{V_i: i\in [t]\}$, $\calP_{edge}=\{P_{ij}^{\alpha}: \alpha\leq \ell\}$, and assume $\calQ$ is an $(s,r)$-decomposition of $V$ consisting of $\calQ_{vert}=\{Q_i: i\in [s]\}$, $\calQ_{edge}=\{Q_{ij}^\alpha: \alpha\in [r]\}$.  For each $\alpha\in [\ell]$ and $\beta\in [r]$, define $P^{\alpha}=\bigcup_{ij\in {[t]\choose 2}}P_{ij}^{\alpha}$ and $Q^{\beta}=\bigcup_{ij\in {[s]\choose 2}}Q_{ij}^{\alpha}$, and set $E^{\alpha,\beta}=P^{\alpha}\cap Q^{\beta}$.

For each $(i,j)\in [t]\times [s]$, define  $X_{ij}=V_i\cap Q_j$, and then let $\calX=\{X_{ij}: ij\in [t]\times [s]\}$.  Set $m=(\mu_1)^2n/ts$, and for each $X_{ij}$ choose a partition $X_{ij}=X_{ij}^0\cup \ldots \cup X_{ij}^{s_{ij}}$ with the property that for each $1\leq u\leq s_{ij}$, $|X_{ij}^u|=m$ and $|X_{ij}^0|<m$.  Setting $X^0=\bigcup_{ij\in [t]\times [s]}X_{ij}^0$, it is easy to see that $|X^0|\leq (\mu_1)^2n$.  Fix an enumeration $U_1,\ldots, U_{\tau}$ of the set $\{X_{ij}^u: 1\leq u\leq s_{ij}\}$, and note we must have $\tau\leq n/m\leq (\mu_1)^{-2}ts$.  Let $\ell'=[32 (\ell r)^2 (\mu_1)^{-2}]$, and define $\rho_1=\mu_1$ and $\rho_2=\mu_2(\ell')$.  Note that $\ell'\leq L$ and $\rho_2<\rho_1/(64\ell')$.

Consider the edge-colored graph with vertex set $V'=\bigcup_{i=1}^{\tau}U_i$, edge colors $F^{\alpha,\beta}:=E^{\alpha,\beta}\cap{V'\choose 2}$, for $\alpha\in [\ell]$ and $\beta\in [r]$.  By construction $\calU=\{U_1,\ldots, U_{\tau}\}$ is an equipartition of $V'$.  By Theorem \ref{thm:mcreg}, there is $K\leq T$, and $\calS=\{S_1,\ldots, S_K\}$ an equipartition refining $\calU$ such that all but at most $\rho_2 K^2$ pairs from $\calS$ satisfy $\disc_2(\rho_2)$ with respect to each of $F^{\alpha,\beta}$.  Let $\Sigma_1$ be the set of pairs $(S_i,S_j)$ from $\calS$ which do not satisfy $\disc_2(\rho_2)$ with respect to all the $F^{\alpha,\beta}$. By assumption, $|\Sigma_1|\leq\rho_2K^2$.  

Our next step is to partition each $K_2[S_i,S_j]$ into equally sized quasirandom graphs.  

For each $(S_i,S_j)\notin \Sigma_1$, define $I_{ij}$ to be the set of $F^{\alpha,\beta}$ with density at most $\rho_1/\ell s$ in $(S_i,S_j)$, and set $E_{ij}^0=(\bigcup_{I_{ij}}E^{\alpha,\beta})\cap K_2[S_i,S_j]$.  Note $|E^0_{ij}|\leq \rho_1 |S_i||S_j|$.       By Lemma \ref{lem:even}, there is some $k_{ij}+k'_{ij}=\ell'$ and a partition 
$$
S_i\times S_j=W^{0,1}_{ij}\cup \ldots \cup W^{0,k_{ij}}_{ij}\cup W^1_{ij}\cup \ldots \cup W_{ij}^{k_{ij}'}
$$
such that each $W^u_{ij}$ and $W^{0,v}_{ij}$ have $\disc_2(\rho_2^{1/3};1/\ell')$, such that for each $1\leq u\leq k_{ij}'$, there is $(\alpha,\beta)\in [\ell]\times [r]$ with $|W^u_{ij}\setminus E^{\alpha,\beta}|<(2\rho_2^{1/3}/\ell')|S_i||S_j|$, and such that $W^0_{ij}:=\bigcup_{v=1}^{k_{ij}} W_{ij}^{0,v}$ has size at most $\rho_1 |S_i||S_j|$.  Since $|W_{ij}^0|=(1\pm k_{ij}\rho^{1/3}_2)m^2(k_{ij}/\ell')$ and $\rho_2\ll \rho_1$, this implies $k_{ij}\leq 2\rho_1 \ell'$.

Now let $q=\min\{k_{ij}': (S_i,S_j)\notin \Sigma_1\}$.  Note that by above, $q\geq (1-2\rho_1)\ell'$.  For each $(S_i,S_j)\notin \Sigma_1$, consider $K_2[S_i,S_j]\setminus (\bigcup_{u=1}^qW_{ij}^u)$.  Let 
$$
\Gamma_{ij}:=K_2[S_i,S_j]\setminus (\bigcup_{u=1}^qW_{ij}^u),
$$
and note that by construction, $\Gamma_{ij}$ is a union of $r_{ij}:=k_{ij}+(k_{ij}'-q)$ many graphs, each with $\disc_2(\rho_2^{1/3}; 1/\ell')$.  If $r_{ij}=0$, then let $E_{ij}^u=W_{ij}^u$ for each $u\in [q]$.  Otherwise, apply Lemma \ref{lem:3.8} to obtain a partition 
$$
\Gamma_{ij}=X_{ij}^1\cup \ldots \cup X_{ij}^q,
$$
each with $\disc_2(\rho_2)$ and density $\frac{r_{ij}}{\ell'q}(1\pm \rho_2)$.  Now for each $1\leq u\leq q$, define $E_{ij}^u=W_{ij}^u\cup X_{ij}^u$.  By Fact \ref{fact:adding}, each $E_{ij}^u$ satisfies $\disc_2(2\rho^{1/3};1/q)$.   Note $K_2[S_i,S_j]=\bigcup_{u\leq q}E_{ij}^u$.

Now, for each $(S_i,S_j)\in \Sigma_1$, choose any partition $S_i\times S_j=\bigcup_{u=1}^{q}E_{ij}^u$ such that $E_{ij}^u$ has $\disc_2(\rho_2; 1/q)$ (such a partition exists by Lemma \ref{lem:3.8}).  

We now deal with the left over vertices, i.e. those in $X^0$. Choose an equipartition $X^0=X^0_1\cup \ldots\cup X^0_K$.   Suppose first there is some $X_i^0$ with $|X_i^0|<\rho_2^{1/3}n/K$.  Then for each $1\leq i\leq K$, $|X^0_i|\leq \rho_2^{1/3}n/K$.  In this case, set $S_i'=S_i\cup X_i^0$ for each $i\in [K]$, and for each $1\leq i<j\leq K$, choose any partition $K_2[X_i^0, S_j']\cup K_2[S_i', X_j^0]=\bigcup_{\alpha=1}^{q}B_{ij}^{\alpha}$.  Note that for each $\alpha$, $|B_{ij}^{\alpha}|\leq \rho_2^{1/3}m(n/K)/q\leq \rho_1|S_i'||S_j'|/q$.  For each $ij\in {[K]\choose 2}$, set $R_{ij}^u=E_{ij}^u\cup B_{ij}^u$.  Note that for each $S_iS_j\notin \Sigma_1$, and each $1\leq u\leq q$, there is some $E^{\alpha,\beta}$ so that
$$
|R_{ij}^u\setminus E^{\alpha,\beta}|\leq |B_{ij}^{u}|+3\rho_1|E_{ij}^u| \leq 4\rho_1|R_{ij}^u|.
$$
In fact, $R_{ij}^u\setminus E^{\alpha,\beta}$ has the form $Y_{ij}^u\cup (W_{ij}^{u'}\setminus E^{\alpha,\beta}) \cup B_{ij}^u$, for some $1\leq u'\leq k_{ij}'$, and $\alpha,\beta$.  Thus since $Y_{ij}^u$ has $\disc_2(\rho_2^{1/3})$ and since $|(W_{ij}^{u'}\setminus E^{\alpha,\beta}) \cup B_{ij}^u|<\rho_2^{1/3}|W_{ij}^{u'}|$, it is easy to check that $R_{ij}^u\setminus E^{\alpha,\beta}$ has $\disc_2(\rho_2^{1/6})$. 

Suppose now that each $X_i^0$ has size at least $\rho_2^{1/6}n/K$.  Let $b=\lceil \rho_2^{1/6}n/K\rceil$.  For each $i\in [K]$, choose any partitions $X_i^0=X_{i,0}^0\cup X_{i,1}^0\cup \ldots \cup X_{i,a}^0$ and $S'_i=S_{i,0}\cup S_{i,1}\cup \ldots S_{i,c}$ such that for each $1\leq j\leq a$ and $1\leq j'\leq c$, $|X_{i,j}^0|=|S_{i,j'}|=b$, and $|X_{i,0}|,|S_{i,0}|<b$.  Note $a,c\leq 2\rho_2^{-1/6}$.

Given $1\leq i<j\leq K$, for each $1\leq i'\leq c$ and $1\leq j'\leq a$, use Lemma \ref{lem:3.8} with parameters $q$ and $\delta=\e=\rho_2$ to choose  partitions $K_2[X_{i,i'}^0,S_{j,j'}]=D_{ii',jj'}^1\cup \ldots \cup D_{ii',jj'}^{q}$ so that for each $\alpha\leq q$, $(X_{i,i'}^0\cup S_{j,j'}, D_{ii',jj'}^u)$ satisfies $\disc_2(\rho_2;1/q)$.  Note this means that each $|D_{ii',jj'}^u|\leq \frac{ac}{q}(1+\rho_2)<2\rho_{1/3}|S_i'||S_j'|/q$.  Then choose an arbitrary partition 
$$
K_2[X_{i,0}^0,S_j']=\bigcup_{\alpha\leq \ell'}B_{ij}^{\alpha}.
$$
Note that each $|B_{ij}^{\alpha}|<\rho_2^{1/6}|S_i'||S_j'|$.  Now define, for each $1\leq u\leq q$, 
$$
R_{ij}^u=E_{ij}^u\cup (\bigcup_{i'=1}^{c}\bigcup_{j'=1}^{a}D_{ii',jj'}^u)\cup B_{ij}^u.
$$
By Fact \ref{fact:adding}, $R_{ij}^u\setminus B_{ij}^u$ has $\disc_2(2\rho_2^{-1/6}\rho_2^{1/3})=\disc_2(\rho_2^{1/6})$, and size
$$
|E_{ij}^u|+ac \frac{b^2}{q}(1\pm \rho_2^{1/6})=\frac{1}{q}(1\pm 5\rho_2^{1/6}).
$$
Thus since each $B_{ij}^u$ is very small, we have that $|R_{ij}^u|=\frac{1}{q}(1\pm 5\rho_2^{1/6})$, and $R_{ij}^u$ has $\disc_2(4\rho_2^{1/6}(q))$.  Observe that by construction, for each $u$, there is $E^{\alpha,\beta}$ so that
$$
|R_{ij}^u\setminus E^{\alpha,\beta}|\leq |E_{ij}^u\setminus E^{\alpha,\beta}|+ac \frac{b^2}{q}(1\pm \rho_2^{1/6})<4\rho_1|R_{ij}^u|.
$$
Setting $\calR_{vert}=\{S_i': i\in [K]\}$ and $\calR_{edge}=\{E_{ij}^u: ij\in {[K]\choose 2}, u\in [q]\}$, we have that $\calR$ is an approximate $(\e_1,\e_2(q))$-refinement of both $\calP$ and $\calQ$.
\end{proofof}

\vspace{2mm}

We now prove Lemma \ref{lem:intersecting}.

\vspace{2mm}

\begin{proofof}{Lemma \ref{lem:intersecting}}
Fix $\e_1>0$ and $\e_2:\mathbb{N}\rightarrow (0,1]$.  As in the proof of Lemma \ref{lem:refinement}, we may assume $\e_2$ is non-increasing.  Choose $\e_1'\ll \e_1$.  Define $\rho,\e_2',\e_2'':\mathbb{N}\rightarrow (0,1]$ by taking $\e''_2(x)\ll \e_1\e_1'/x$,  $\rho_2(x)\ll \min\{\e_2(x), \e_2''(x)\}$, and  $\e_2'(x)\ll \rho_2(x)/x$, for all $x\in \mathbb{N}$. Fix $t,\ell\geq 1$ and $t',\ell'\geq 1$, and choose $N\gg \ell, t(\e_1')^{-1}\e_2'(\ell)^{-1}, \ell',t'$.

Suppose $H=(V,E)$ is a $3$-graph with $|V|=n\geq N$ and $\calQ$ is a $(t,\ell,\e_1',\e_2'(\ell))$-decomposition of $V$ which is $\e_1'$-homogeneous with respect to $H$.  Suppose $\calP$ is a $(t',\ell',\e_1',\e_2'(\ell'))$-decomposition of $V$ which is a $(\e'_1,\e'_2(\ell'))$-approximate refinement of $\calQ$.  Say $\calQ_{vert}=\{V_1,\ldots, V_{t}\}$, and $\calQ_{edge}=\{Q_{ij}^{\alpha}:ij\in {[t]\choose 2}, \alpha\leq \ell\}$ and $\calP_{vert}=\{W_1,\ldots, W_{t'}\}$, and $\calP_{edge}=\{P_{ij}^{\alpha}:ij\in {[t']\choose 2}, \alpha\leq \ell'\}$.  To ease notation, given $i,j,k\in [t]$ and $\alpha,\beta,\gamma\in [\ell]$, we let $Q_{ijk}^{\alpha,\beta,\gamma}$ denote the triad $(V_i\cup V_j\cup V_k, Q_{ij}^{\alpha} \cup Q_{ik}^{\beta}\cup Q_{jk}^{\gamma})$.  Similarly, given $i,j,k\in [t']$ and $\alpha,\beta,\gamma\in [\ell']$, we let $P_{ijk}^{\alpha,\beta,\gamma}$ denote the triad $(W_i\cup W_j\cup W_k, P_{ij}^{\alpha}\cup P_{ik}^{\beta}\cup P_{jk}^{\gamma})$. Define 
\begin{align*}
\Gamma_1=\{Q_{ijk}^{\alpha,\beta,\gamma}: Q_{ijk}^{\alpha,\beta,\gamma}\text{ is not }&\disc_{2,3}(\e_1',\e_2'(\ell))\text{-regular  and} \\ & \e_1'\text{-homogeneous with respect to }H\}
\end{align*}
and
\begin{align*}
\Gamma_2=\{P_{ijk}^{\alpha,\beta,\gamma}: P_{ijk}^{\alpha,\beta,\gamma}\text{ is not }&\disc_{2,3}(\e_1,\e_2(\ell'))\text{-regular and} \\ & \e_1\text{-homogeneous with respect to }H\}.
\end{align*}
Setting $\Omega_1=\bigcup_{Q_{ijk}^{\alpha,\beta,\gamma}\in \Gamma_1}K_2^{(2)}(Q_{ijk}^{\alpha,\beta,\gamma})$, we have by assumption that $|\Omega_1|\leq \e_1' n^3$.  Our goal is to show that $|\Omega_2|\leq \e_1 n^3$, where $\Omega_2:=\bigcup_{P_{ijk}^{\alpha,\beta,\gamma}\in \Gamma_2}K_2^{(2)}(P_{ijk}^{\alpha,\beta,\gamma})$.

Let $\Sigma\subseteq {\calP_{vert}\choose 2}$ witness that $\calQ$ is an $(\e'_1,\e'_2(\ell'))$-approximate refinement of $\calQ$.   By the definition of an approximate refinement, for each $V_iV_j\notin \Sigma$, and each $P_{ij}^{\alpha}\in \calP_{edge}$, there is some $f(P_{ij}^{\alpha})\in \calQ_{edge}$ such that $|P_{ij}^{\alpha}\setminus f(P_{ij}^{\alpha})|<\e'_1|P_{ij}^{\alpha}|$ and such that $P_{ij}^{\alpha}\setminus f(P_{ij}^{\alpha})$ has $\disc_2(\e'_2(\ell'))$.  Similarly, for all $P_{ijk}^{\alpha, \beta,\gamma}\in \triads(\calP)$ with $V_iV_j, V_iV_k,V_jV_k\notin  \Sigma$, we let $f(P_{ijk}^{\alpha, \beta,\gamma})$ denote the triad $Q_{uvw}^{\alpha',\beta',\gamma'}$ with the property that $f(P_{ij}^{\alpha})=Q_{uv}^{\alpha'}$, $f(P_{ik}^{\beta})=Q_{uw}^{\beta'}$, and $f(P_{jk}^{\gamma})=Q_{vw}^{\gamma'}$.  We now define 
$$
\Sigma'=\Sigma\cup \{P_{ij}^{\alpha}\in \calP_{edge}: P_{ij}^{\alpha} \text{ fails }\disc_2(\e_1',\e_2'(\ell'))\}.
$$
Note that we must have $|\Sigma'|\leq \e'_1(t')^2+\sqrt{\e_1'}(t')^2\leq 2\sqrt{\e_1'}(t')^2$, where the first inequality is by definition of a $(\e'_1,\e'_2(\ell'))$-approximate refinement, and the second is because $\calP$ is a $(t',\ell, \e_1',\e_2'(\ell'))$-decomposition of $V$.  Now let
$$
E^*=\bigcup_{ij\notin \Sigma'}P_{ij}^{\alpha}\setminus f(P_{ij}^{\alpha})\;\text{ and }\;V^*=\bigcup_{W_i\in \calP_{vert}}W_i\setminus V_i.
$$
By assumption, $|E^*|\leq \e_2'(\ell') |V|^2$ and $|V^*|\leq \e_1'|V|$.  Define
$$
\Sigma_{\calQ}=\{Q_{ij}^{\alpha}\in \calQ_{edge}: |Q_{ij}^{\alpha}\cap E^*|\geq \sqrt{\e_1'}|Q_{ij}^{\alpha}|\text{ or } \min\{|V_i\cap V^*|,|V_j\cap V^*|\}\geq \sqrt{\e_1'}|V_j|\}.
$$
Clearly, the bounds on $E^*$ and $V^*$ implies $|\Sigma_{\calQ}|\leq 3\sqrt{\e'_1}\ell t^2$. Now set
$$
\Gamma_1'=\{Q_{ijs}^{\alpha,\beta,\gamma}: Q_{ij}^{\alpha}, Q_{is}^{\beta},\text{ or }Q_{js}^{\gamma}\in \Sigma_{\calQ}\}.
$$
For each $ij\in {[t]\choose 2}$, set $I_{ij}^{\alpha}=\{P_{uv}^{\beta}\in \calP_{edge}: f(P_{uv}^{\beta})=Q_{ij}^{\alpha}\}$, and define, for each $ijs\in {[t]\choose 3}$ and $\alpha,\beta,\gamma\leq \ell$,
$$
\calI_{ijs}^{\alpha,\beta,\gamma }=\{P_{uvw}^{\rho,\gamma,\tau}\in \triads(\calP): P_{uv}^{\rho}\in I_{ij}^{\alpha}, P_{uw}^{\mu}\in I_{is}^{\beta}, P_{vw}^{\tau}\in I_{js}^{\gamma}\}.
$$

Suppose $Q_{ijs}^{\alpha,\beta,\gamma}\notin \Gamma_1\cup \Gamma_1'$.  We claim that almost all triples contained in $K_3^{(2)}(Q_{ijs}^{\alpha,\beta,\gamma})$ are in an $\e_1$-homogeneous triad of $\calP$.  First, we observe that by definition of $\Gamma_1$, there is a $\delta=\delta(ijs,\alpha\beta\gamma)\in \{0,1\}$ so that 
$$
|E^{\delta}\cap K_3^{(2)}(Q_{ijs}^{\alpha,\beta,\gamma})|\geq (1-\e_1')|K_3^{(2)}(Q_{ijs}^{\alpha,\beta,\gamma})|.
$$
Next, note that
$$
\bigcup_{P_{uvw}^{\rho,\gamma,\tau}\in \calI_{ijs}^{\alpha,\beta,\gamma }}K_3^{(2)}(P_{uvw}^{\rho,\gamma,\tau}\cap Q_{ijs}^{\alpha,\beta,\gamma})\subseteq K_3^{(2)}(Q_{ijs}^{\alpha,\beta,\gamma}).
$$
This implies that 
$$
|\calI_{ijs}^{\alpha,\beta,\gamma}|\geq (1-3\sqrt{\e_1'})\frac{n^3}{t^3\ell^3}\Big/(1-\e_2(\ell))\frac{n^3}{(t')^3(\ell')^3}\geq (1-4\sqrt{\e_1'})(t')^3(\ell')^3/t^3\ell^3.
$$
For all $P_{uvw}^{\rho,\gamma,\tau}\in \calI_{ijs}^{\alpha,\beta,\gamma }$, each of $ P_{uv}^{\rho}\setminus f(P_{uv}^{\rho})$, $P_{uw}^{\mu}\setminus f(P_{uw}^{\mu})$ and $P_{vw}^{\tau}\setminus f(P_{vw}^{\tau})$ satisfy $\disc_2(\e'_2(\ell))$ and have size at most $\e_1'\frac{n^2}{(t')^2}$.  Combining this with the observation above and Corollary \ref{cor:counting}, we have 
$$
(1-3\e_1')^3\frac{n^3}{(t')^3(\ell')^3}|\calI_{ijs}^{\alpha,\beta,\gamma }|\leq |K_2^{(2)}(Q_{ijs}^{\alpha,\beta,\gamma})|\leq (1+\e'_1)\frac{n^3}{\ell^3t^3}.
$$
Consequently, using that $\e_1'$ is small, we have that $|\calI_{ijs}^{\alpha,\beta,\gamma }|\leq \frac{(t')^3(\ell')^3}{t^3\ell^3}(1-\sqrt{\e_1'})$.  On the other hand, since $Q_{ijs}^{\alpha,\beta,\gamma}\notin  \Gamma_1'$,
$$
|K_2^{(2)}(Q_{ijs}^{\alpha,\beta,\gamma})\setminus \bigcup_{P_{uvw}^{\rho,\gamma,\tau}\in \calI_{ijs}^{\alpha,\beta,\gamma }}K_3^{(2)}(P_{uvw}^{\rho,\gamma,\tau}\cap Q_{ijs}^{\alpha,\beta,\gamma})|\leq 3\sqrt{\e_1'}\frac{n^3}{t^3\ell^3}.
$$
Using the inequalities above and Corollary \ref{cor:counting}, we have 
\begin{align*}
\sum_{\triads(I_{ijs}^{\alpha,\beta,\gamma})}|E^{\delta}\cap K_3^{(2)}(P_{uvw}^{\rho,\mu,\tau}))|&\geq |E^{\delta}\cap\bigcup_{P_{uvw}^{\rho,\gamma,\tau}\in \calI_{ijs}^{\alpha,\beta,\gamma }} K_3^{(2)}(P_{uvw}^{\rho,\gamma,\tau}\cap Q_{ijs}^{\alpha,\beta,\gamma})|\\
&\geq (1-3\sqrt{\e_1'})|E^{\delta}\cap K_3^{(2)}(Q_{ijs}^{\alpha,\beta,\gamma})|\\
&\geq (1-3\sqrt{\e_1'})(1-\e_1')|K_3^{(2)}(Q_{ijs}^{\alpha,\beta,\gamma})|\\
&\geq (1-3\sqrt{\e_1'})(1-\e_1')^2\frac{n^3}{t^3\ell^3}.
\end{align*}
This implies
$$
\sum_{\calI_{ijs}^{\alpha,\beta,\gamma}}\frac{|E^{\delta}\cap K_3^{(2)}(P_{uvw}^{\rho,\mu,\tau}))|}{|K_3^{(2)}(P_{uvw}^{\rho,\mu,\tau}))|}\geq \frac{t^3\ell^3}{(t')^3(\ell')^3}(1-4(\e_1')^{1/2})\geq (1-(\e_1')^{1/4})|\calI_{ijs}^{\alpha,\beta,\gamma}|.
$$
Since each term in the summand above is in $[0,1]$, and $\e_1'\ll \e_1$, this implies that there is a set $\calJ_{ijs}^{\alpha,\beta,\gamma}\subseteq \calI_{ijs}^{\alpha,\beta,\gamma}$ such that $|\calJ_{ijs}^{\alpha,\beta,\gamma}|\geq (1-\e_1^2)|\calI_{ijs}^{\alpha,\beta,\gamma}|$, and for all $P_{uvw}^{\rho,\mu,\tau}\in \calJ_{ijs}^{\alpha,\beta,\gamma}$
$$
\frac{|E^{\delta}\cap K_3^{(2)}(P_{uvw}^{\rho,\mu,\tau}))|}{|K_3^{(2)}(P_{uvw}^{\rho,\mu,\tau})|}\geq (1-\e_1^2).
$$
Thus, the number of triples from $Q_{ijs}^{\alpha,\beta,\gamma}$ which are in an $\e_1^2$-homogeneous triad of $\calP$ is at least 
\begin{align*}
|\calJ_{ijs}^{\alpha,\beta,\gamma}|(1-\e_1^2)\frac{n^3}{(t')^3(\ell')^3}&\geq (1-\e_1^2)|\calI_{ijs}^{\alpha,\beta,\gamma}|(1-\e_1^2)\frac{n^3}{(t')^3(\ell')^3}\\
&\geq (1-\e_1^2)(1-4\sqrt{\e_1'})((t')^3(\ell')^3/t^3\ell^3)(1-\e_1^2)\frac{n^3}{(t')^3(\ell')^3}\\
&=(1-\e_1^2)^2(1-4\sqrt{\e_1'})^2\frac{n^3}{t^3\ell^3}.
\end{align*}
We now have that 
\begin{align*}
|\Omega_2|&\geq |\triads(\calQ)\setminus (\Gamma_1\cup \Gamma_1')|(1-\e_1^2)^2(1-4\sqrt{\e_1'})^2\frac{n^3}{t^3\ell^3}\\
&\geq (1-10\sqrt{\e_1'})t^3\ell^3(1-\e_1^2)^2(1-4\sqrt{\e_1'})^2\frac{n^3}{t^3\ell^3}\\
&\geq (1-\e_1)n^3,
\end{align*}
which finishes the proof.
\end{proofof}

\vspace{3mm}

Finally, we prove Theorem \ref{thm:counting}.

\vspace{3mm}

\begin{proofof}{Theorem \ref{thm:counting}}
Given $\xi, d_3>0$ and $t\geq 1$, let $\xi_0>0$ be sufficiently small so that $(1+\xi_0)^{t^3+t^1}\leq 1+\xi$.  Choose $\delta_3=\delta_3(t,\xi_0,d_3/2)$ as in Theorem \ref{thm:countinginducedHG}.  Fix $d_2>0$, and choose $\delta_2=\delta_2(t,\xi_0,d_3/2,\delta_3,d_2)$ and $n_0(t,\xi_0,d_3/2,\delta_3,d_2)$ as in Theorem \ref{thm:countinginducedHG}.  Let $m_0$ be as in Fact \ref{fact:rg} for $\delta_2$ and $d_2$, and let $\delta_2'=d_2^8\delta_2/4$.  

Let $\delta_3'=\delta_3'(\min\{\delta_3,d_3\xi_0\},2)$, $\delta_2''=\delta_2''(\delta_2',\min\{\delta_3,d_3\xi_0\},d_3,2))$, and  $M_0=M_0(\delta_2',\min\{\delta_3,d_3\xi_0\},d_3,2)$ be from Lemma \ref{lem:subpairs}.  Now set $N=tm_0n_0M_0(\delta''_2)^{-4}$.

Fix $t_a,t_b,t_c\geq 1$ such that $t_a+t_b+t_c=t$, and $F=(U\cup W\cup Z,R_F)$ a $3$-partite $3$-graph with $U=\{u_1,\ldots, u_{t_a}\}$, $W=\{w_1,\ldots, w_{t_b}\}$, and $Z=\{z_1,\ldots, z_{t_c}\}$.  Let $V$ be a set of size $n\geq N$, equipped with a partition 
$$
V=A_1\cup \ldots \cup A_{t_a}\cup B_1\cup \ldots \cup B_{t_b}\cup C_1\cup \ldots \cup C_{t_c},
$$
where each part has size $\frac{n}{t}(1\pm \delta_2'')$.  Let $A=\bigcup_{i=1}^{t_a}A_i$,  $B=\bigcup_{i=1}^{t_b}B_i$, and  $C=\bigcup_{i=1}^{t_c}C_i$.  Suppose   $H=(V,R)$ is a $3$-graph and $G=(V,E)$ is a $t$-partite graph with the vertex partition displayed above.  Assume $H$ and $G$ satisfy the hypotheses of Theorem \ref{thm:counting} with $\delta''_2, d_2, d_3, \delta'_3$.  

Let $H'=(V, E\setminus ({A\choose 2}\cup {B\choose 2}\cup {C\choose 2}))$. We define a new graph $G'$ on $V$ as follows.  For each $X\in \{A,B,C\}$ and $ij\in {[t_x]\choose 2}$, define $G'[X_i,X_j]$ to be any bipartite graph on $X_i\cup X_j$ with $\disc_2(\delta'_2;d_2)$ (such a graph exists by Fact \ref{fact:rg}).  For $X\neq Y\in \{A,B,C\}$, $i\in [t_x]$ and $j\in [t_y]$, define $G'[X_i,Y_j]=G[X_i,Y_j]$.  For each $XYZ\in \{A,B,C\}^3$ and $(i,j,s)\in [t_x]\times [t_y]\times [t_z]$,  set $G_{XYZ,ijs}:=G'[X_i,Y_j,Z_s]$ and $H_{XYZ, ijs}=(X_i\cup Y_j\cup Z_s, R\cap K_3^{(2)}(G_{XYZ, ijs}))$, and let $d_{XYZ,ijs}$ be such that $|R\cap K_3^{(2)}(G_{XYZ,ijs})|=d_{XYZ,ijs}|K_3^{(2)}(G_{XYZ,ijs})|$.  If $XYZ=ABC$, then by assumption, $(H_{XYZ,ijs},G_{XYZ,ijs})$ satisfies $\disc_{2,3}(\e_1,\e_2(\ell))$.  If $XYZ\neq ABC$, then by definition $(H_{XYZ,ijs},G_{XYZ,ijs})$ has density $0$, so by Proposition \ref{prop:homimpliesrandome}, it satisfies $\disc_{2,3}(\delta_3',\delta_2'')$.  Set
$$
d=(\prod_{u_iv_jw_s\in R_F}d_{ABC,ijs})(\prod_{u_iv_jw_s\in R_F}(1-d_{ABC,ijs})).
$$

Let $m=\min\{|A_i|, |B_j|,|C_k|:i\in [t_a], j\in [t_b], k\in [t_c]\}$, and for each $(i,j,k)\in [t_a]\times [t_b]\times [t_c]$, choose any $m$-element subsets $A'_i$, $B_j'$, and $C_k'$ of $A_i$, $B_j$, and $C_k$ respectively.  Note $m\geq (1-\delta_2')\frac{n}{t}\geq \max\{n_0,m_0\}$.  Define $V'=\bigcup_{i=1}^{t_a}A_i'\cup \bigcup_{i=1}^{t_b}B_i'\cup \bigcup_{i=1}^{t_c}C_i'$, and set  $G''=G'[V']$ and $H''=H'[V']$.  

Then for each $XYZ\in \{A,B,C\}^3$ and $(i,j,s)\in [t_x]\times [t_y]\times [t_z]$, let $G'_{XYZ,ijs}=G_{XYZ,ijs}[V']$ and $H'_{XYZ,ijs}=H_{XYZ,ijs}[V']$.  Note that by Lemma \ref{lem:subpairs}, each $(G'_{XYZ,ijs},H'_{XYZ,ijs})$ has $\disc_{2,3}(\delta_3,\delta_2')$ with density $d_{XYZ,ijs}(1\pm \xi_0)$.

By Theorem \ref{thm:countinginducedHG}, there is a set $\Sigma\subseteq  \prod_{i=1}^{t_1}A'_i\times \prod_{i=1}^{t_2}B'_i\times \prod_{i=1}^{t_3}C'_i$ such that for all $\abar\bbar\cbar\in \Sigma$, $x_iy_jc_s\in {\abar\bbar\cbar\choose3}$ if and only if $x_iy_jz_s\in E(H')$, and such that
$$
|\Sigma|=d_2^{t\choose 2}d' (1\pm \xi_0)m^t,
$$
where $d'=(\prod_{(u_i,v_j,w_s)\in R} d'_{ABC, ijs})(\prod_{(u_i,v_j,w_s)\notin R} (1-d'_{ABC, ijs}))$.  Recall that each $d_{ABC,ijs}'=d_{ABC,ijs}(1\pm \xi_0)$.  Consequently, $d'=d(1\pm \xi_0)^{t^3}$.  Thus
\begin{align*}
|\Sigma|&=d_2^{t\choose 2}d (1\pm \xi_0)m^t\\
&=d_2^{t\choose 2}d (1\pm \xi_0)^{t^3+t}\prod_{i=1}^{t_1}|A_i| \prod_{i=1}^{t_2}|B_i| \prod_{i=1}^{t_3}|C_i|\\
&=d_2^{t\choose 2}d (1\pm \xi)\prod_{i=1}^{t_1}|A_i| \prod_{i=1}^{t_2}|B_i| \prod_{i=1}^{t_3}|C_i|.
\end{align*}
By construction, for all $\abar\bbar\cbar\in \Sigma$, $a_ib_jc_s\in R$ if and only if $u_iw_jz_s\in R_F$, so this finishes the proof.
\end{proofof}

\chapter{Proofs of Ramsey facts}\label{app:ramsey}

In this appendix we prove the Ramsey facts from Section \ref{subsec:ramseyfacts}. For convenience, we will use model theoretic language along with the compactness theorem.  This is essentially only to cut down on the notation required, and all the proofs can be easily finitized by replacing infinite sets with sufficiently large finite ones.

Let $\calL=\{R(x,y,z)\}$ consist of a single ternary relation symbol.  Recall that for any hereditary $3$-graph property, there is a universal $\calL$-theory $T_{\calH}$ such that $\calH$ is the class of finite models of $T_{\calH}$.  

\vspace{2mm}

\begin{proofof}{Fact \ref{fact:vc2universal2}}
Suppose $\calH$ has $\IP_2$.  By compactness, there is a an infinite $3$-graph $\calM=(V,F)\models  T_{\calH}$ with sets $A=\{a_i: i\in \mathbb{N}\}\subseteq V$, $B=\{b_i: i\in \mathbb{N}\}\subseteq V$, and $C=\{c_S: S\subseteq A\times B, S\text{ finite}\}\subseteq V$ so that $a_ib_jc_S\in E$ if and only if $(i,j)\in S$.   Clearly we may choose $A'\subseteq A$ and $B'\subseteq B$ such that each of $A'$ and $B'$ are infinite and $A'\cap B'=\emptyset$.  Suppose $H=(V_1\cup V_2\cup V_3, E)$ is a finite $3$-partite $3$-graph.

Let $f:V_1\rightarrow A'$, $g:V_2\rightarrow B'$ be any injections.  By construction $f(V_1)\cap g(V_2)=\emptyset$.  For each $w\in V_3$, set $S_w=\{c_S\in C: S\cap (f(V_1)\times g(V_2))=\{f(u)f(v): uvw\in E\}\}$.  There are infinitely many distinct $S\subseteq A\times B$ with the property that $S\cap (f(V_1)\times g(V_2))=\{f(u)f(v): uvw\in E\}$, and thus $S_w$ is infinite for all $w\in V_3$.  Thus, for all $w\in V_3$, we have that $S_w'=S_w\setminus (f(V_1)\cup g(V_2))$ is also infinite.  Since $V_3$ is finite, we can clearly choose elements $h(w)\in S_w'$ so that $h(V_3)=\{h(w):w\in V_3\}$ has size $V_3$.  By construction, $H':=\calM[f(V_1)\cup g(V_2)\cup h(V_3)]$ is a clean copy of $H$.  Since $\calM\models T_{\calH}$ and universal sentences are closed under substructures, $H'\models T_{\calH}$, and consequently, $H'\in  \calH$ is as desired.
\end{proofof}

\vspace{2mm}

\begin{proofof}{Lemma \ref{lem:clean}}
We begin with (1).  Suppose $H^*(k)\in \trip(\calH)$ for all $k\geq 1$.  By compactness there is a $3$-graph $\calM=(V,E)\models T_{\calH}$, a vertex $c\in M$, and subsets $A=\{a_i: i\in \mathbb{N}\}$, $B=\{b_i: i\in \mathbb{N}\}\subseteq V$ such that $ca_ib_j\in E$ if and only if $i\leq j$.  It is easy to see that all elements in $A$ and $B$ must be pairwise distinct, respectively.  Given $a_i\in A$, since $\calM$ is a $3$-graph, for any $j\geq i$, $ca_ib_j\in E$ implies $a_i\neq b_j$.  Thus, there can be at most one element in $B$ equal to $a_i$, and it would have to be in the set $\{b_1,\ldots, b_i\}$.  Similarly, since $c$ is in an edge with every element of $A$ and $B$, it is distinct from all of $A\cup B$.

Fix $k\geq 1$.  We define sequences $(m_k,\ldots, m_1)$ and $(I_k,\ldots, I_1)$ as follows.

Start with $m_k=k^3$, and $I_k=\{j<k^3: b_j\neq a_{k^3}\}$.  

Suppose now that $0\leq i\leq k-1$, and assume we have defined $m_k,\ldots, m_{k-i}$ and $I_k,\ldots, I_{k-i}$.  Let $m_{k-i-1}=\min I_{k-i}$ and let  $I_{k-i-1}=\{j<m_{k-i-1}: b_j\neq a_{k-i-1}\}$.  At the end we will have that $ca_{m_u}b_{m_v}$ holds if and only if $u\leq v$, and $X=\{c\}\cup \{a_{m_u},b_{m_u}: u\in [k]\}$ has size $2k+1$.  Thus $\calM[X]$ is a clean copy of $H^*(k)$.  As above, $\calM[X]\in \calH$, so we are done.

For part (2), suppose $\Hbar(k)\in \trip(\calH)$ for all $k\geq 1$.  By compactness there is a $3$-graph $\calM=(V,E)\models T_{\calH}$ with subsets $C=\{c_i: i\in \mathbb{N}\}$, $A=\{a_i: i\in \mathbb{N}\}$, and $B=\{b_i: i\in \mathbb{N}\}\subseteq M$ such that $c_sa_ib_j\in E$ if and only if $i\leq j$, and $c_i\neq c_j$ for each $i\neq j$.  It is easy to see that all elements in $A$ and $B$ must be pairwise distinct, respectively.  Given $a_i\in A$, since $\calM$ is a $3$-graph, for any $j\geq i$, $c_1a_ib_j\in E$ implies $a_i\neq b_j$.  Thus, there can be at most one element in $B$ equal to $a_i$, and it would have to be in the set $\{b_1,\ldots, b_i\}$.  Similarly, since each $c_s$ is in an edge with every element of $A$ and also some edge with every element of $B$, we automatically know $C\cap A=C\cap B=\emptyset$. 

Fix $k\geq 1$.  Construct a sequence $m_k,\ldots, m_{1}$ and $I_k,\ldots, I_{1}$ exactly as in part (1).  Then for all $s\in [k]$, $c_sa_{m_u}b_{m_v}$ holds if and only if $u\leq v$, and $X=\{c_1,\ldots c_k\}\cup \{a_{m_u},b_{m_u}: u\in [k]\}$ has size $3k$.  Thus $\calM[X]$ is a clean copy of $\Hbar(k)$.  As above, $\calM[X]\in \calH$, so we are done.

We now do part (3).  Suppose $F(k)\in \trip(\calH)$ for all $k\geq 1$.  By compactness there is a $3$-graph $\calM=(V,E)\models T_{\calH}$ with subsets $C=\{c_i: i\in \mathbb{N}\}$, $A=\{a_i: i\in \mathbb{N}\}$, and $B=\{b^f_i: f:\mathbb{N}^2\rightarrow \mathbb{N}, i\in \mathbb{N}\}\subseteq M$ such that $a_ib^f_jc_s\in E$ if and only if $s\leq f(i,j)$. It is easy to see this implies that the elements of $C$ and $A$ are pairwise distinct, respectively.  Further, since every $a_i$ is in an edge with every $c_s$ (if $g$ is the constant function $g(x,y)=s$, then $a_ib_1^gc_s\in E$), we must have $A\cap C=\emptyset$.  

Fix $k\geq 1$ and let $N\gg k$.  Let $I$ be the set of functions from $[k]^2\rightarrow [k]$.  Fix any enumeration $I=\{f_u: u\in k^{2k}\}$.  For each $f_u\neq f_v\in I$, let $(\alpha_{uv},\beta_{uv})$ be the least element of $[k]^2$ (in the lexicographic order) such that $f_u(\alpha_{uv},\beta_{uv})\neq f_v(\alpha_{uv},\beta_{uv})$.  Given $f_u:[k]^2\rightarrow [k]$, and $w\in [k^2]$ define $g^w_{f_u}:[N^2]\rightarrow [N]$ as follows.  For all $i,j\in [k]$, set $g^2_{f_u}(i,j)=f_u(i,j)$, and for each $u\neq v\in k^{2k}$, and $j\in [k]$, define $g^w_{f_u}(k+3^u4^v,j)=g^w_{f_u}(k+3^v4^u,j)=f_u(\alpha_{uv},\beta_{uv})$, and for each $j\neq j'\in [k]$, let $g^w_{f_u}(k+5^j7^{j'},j)=1$ and  $g^w_{f_u}(k+5^j7^{j'},j')=2$, and then for each $w< w'\leq N$, define $g^w_{f_u}(k+11^w13^{w'},j)=w$ and $g^{w'}_{f_u}(k+11^w13^{w'},j)=w'$.  The point of all this is that $|\{e_i^{g^w_{f_u}}: w\in [N],f_u\in I, i\in [k]\} |=|I|k^2N$.  For each $f_u\in I$, and each $j\in [j]$, let $X_{u,j}=\{e_j^{g^w_{f_u}}: w\in [N]\}$.  Then $X_{u,j}$ has size $N$, so it contains a subset $X_{u,j}'$ of size at least $N-2k$ which is disjoint from $\{c_i,a_i': i\in [k]\}$. 

Choose $(b_1^{f_1},\ldots, b_k^{f_k})\in X_{1,1}\times \ldots \times X_{1,k}$ so that the coordinates are pairwise distinct.  This is possible by the size of each $X_{1,j}$.

At step $u+1$, choose $(b_1^{f_{u+1}},\ldots, b_k^{f_{u+1}})\in X_{u+1,1}\times \ldots \times X_{u+1,k}$ which are pairwise distinct, and also distinct from all $b_j^v$ for all $1\leq v\leq u$ and $j\in [k]$, which is possible by how large the $X_{u+1,j}$ are. 

In the end, we have a set $V'' =\{c_i,a_i',b_i^{f_u}: i\in [k], f_u\in I\}$ such that $a_ib_j^{f_u}c_s\in E$ if and only if $s\leq f_u(i,j)$, and such that the elements of $V''$ are all pairwise distinct.  In other words, $\calM[V'']$ is a clean copy of $F(k)$.  As before, this tells us $\calM[V'']\in \calH$, as desired.
\end{proofof}

\vspace{2mm}

\begin{proofof}{Lemma \ref{lem:ubarkuniversal}}
Suppose $\Ubar(k)\in \trip(\calH)$ for all $k\geq 1$.  By compactness there is a $3$-graph $\calM=(V,E)\models T_{\calH}$ with subsets $C=\{c_i: i\in \mathbb{N}\}$, $A=\{a_i: i\in \mathbb{N}\}$, and $B=\{b_S:S\subseteq \mathbb{N}\}$ such that $c_ja_ib_S\in E$ if and only if $i\in S$.  It is easy to see that all elements in $A$ and $B$ must be pairwise distinct, respectively, and the $c_i$ are pairwise distinct by assumption.  Given $a_i$, $c_j$, and $b_S$, with $i\in S$, since $c_ja_ib_S\in E$, we know $a_i\neq c_j$, $b_S\neq c_j$.  Thus, $(B\setminus \{b_{\emptyset}\})\cap A=(B\setminus \{b_{\emptyset}\})\cap C=\emptyset$.  

Fix a finite bipartite graph $G=(U\cup V, E')$ and $n\geq 1$.  Choose $\{A_v: v\in V\}$ so that each $A_v\subseteq A$ has size $|V|^n$ and so that $v\neq v'$ implies $A_v\cap A_{v'}=\emptyset$.   For each $u\in U$, let 
$$
B_{u}=\Big\{b_S: \bigcup_{v\in N_G(u)}A_v\subseteq S\text{ and }S\cap (\bigcup_{v\in V\setminus N_G(u)}A_v)=\emptyset\Big\}.
$$
Note that $B_{u}$ is infinite for every $u\in U$.  Thus, we may clearly choose a tuple $(b_u)_{u\in U}$ and $(x_v)_{v\in V}$ so that $x_v\in A_v$ for each $v\in V$, $b_u\in B_u$ for each $u\in U$, and so that all the elements of $\{b_u: u\in U\}\cup \{x_v: v\in V\}$ are pairwise distinct.  Now choose any $c_1,\ldots, c_n$ also disjoint from this set (this is possible as $C$ is inifinite).  We now have that $c_jb_ux_v\in E$ if and only if $v\in N_G(u)$, i.e. if and only if $uv\in E'$.  It is now clear that $\calM[\{b_u: u\in U\}\cup \{x_v: v\in V\}\cup \{c_1,\ldots, c_n\}]\in \calH$ is a clean copy of $n\otimes G$, as desired.
\end{proofof}

\chapter{General properties of $\NFOP_2$-formulas}\label{app:fopgen}

This appendix contains some general properties of $\NFOP_2$ formulas. We will use model theoretic notation.  First, we observe that formally, all three variables play a distinct role in the definition of $\FOP_2$.  For this reason, when we give the definition of $\FOP_2$ in general, we must partition the variables into three parts (rather than the more typical two).

\begin{definition}\label{def:fopgeneral}
Suppose $\calL$ is a first-order language, $\calM$ is an $\calL$-structure, and $\phi(x,y,z)$ is an $\calL$-formula.  We say the tripartitioned formula $\phi(x;y;z)$ \emph{has $k$-$\FOP_2$ in $\calM$} if there are $b_1,\ldots, b_k,c_1,\ldots, c_k\in M$, and for each $f:[k]^2\rightarrow [k]$, $a_1^f,\ldots, a_k^f\in M$ such that $\calM\models \phi(a^f_i,b_j,c_s)$ if and only if $s\leq f(i,j)$.
\end{definition}

It is easy to see that if a formula $\phi(x,y,z)$ has $\VC_2$-dimension at least $\ell$, then it has $\ell$-$\FOP_2$ (see e.g. Fact \ref{fact:vc2universal}).  Let $\calM$ be the Fra\"{i}ss\'{e} limit of all $3$-graphs.  It is well known that $Th(\calM)$ is simple (this can be shown via an argument similar to Corollary 7.3.14 in \cite{Tent.2012}).  Since the edge relation in this theory has unbounded $\VC_2$-dimension, this shows that simple $\nRightarrow \text{NFOP}_2$.  On the other hand, using the machinery of this paper (namely Theorem \ref{thm:FOP} and Section \ref{subsec:rbinary}), it is possible to show that any quantifier-free formula $\phi(x;y;z)$ in the graph language is $\NFOP_2$ in any graph. This implies, for example, that any quantifier-free formula in three variables in the Henson graph is $\NFOP_2$.  Consequently, $\NFOP_2\nRightarrow$ simple.  We also show in \cite[Lemma 5.11]{Terry.2021a} that if a formula has $\ell$-$\FOP_2$ then it has $\VC$-dimension at least $\ell$.

It is well known that stable formulas are closed under boolean combinations. Our next results address the question of whether $\NFOP_2$-formulas are similarly closed under boolean combinations.  In our next proposition, we give an equivalent formulation of $\ell$-$\FOP_2$, and also show that $\NFOP_2$ is closed under negations.  

\begin{proposition}\label{prop:basicprop}
Suppose $\calL$ is a first-order language, $\calM$ is an $\calL$-structure, and $\phi(x,y,z)$ is an $\calL$-formula.   
\begin{enumerate}[label=\normalfont(\roman*)]
\item If $\phi(x_1;x_2;x_3)$ has $\ell$-$\FOP_2$ in $\calM$, then $\neg \phi(x_1;x_2;x_3)$ has $(\ell-1)$-$\FOP_2$.
\item The following are equivalent.
\begin{enumerate}[label=\normalfont(\alph*)]
\item $\phi(x_1;x_2;x_3)$ has $\ell$-$\FOP_2$ for all $\ell$. 
\item For all $\ell$, there are $b_1,\ldots,b_{\ell},c_1,\ldots,c_{\ell}$ and $a_f$ for $f\in [\ell]^{[\ell]}$ such that $\phi(a^f,b_i,c_k)$ holds if and only if $k\leq f(i)$. 
\end{enumerate}
\end{enumerate}
\end{proposition}
\begin{proof}
We begin with (i).  Let $y_1,\ldots, y_\ell,z_1,\ldots, z_\ell$ and $\{x_i^f: i\in [\ell], f:[\ell]^2\rightarrow [\ell]\}$ be such that $\phi(x^f_i,y_j, z_k)$ holds in $\calM$ if and only if $k\leq f(i,j)$.  For each $f:[\ell-1]^2\rightarrow [\ell-1]$, let $g_f:[\ell]^2\rightarrow [\ell]$ be the function defined by $g_f(x,y)=\ell-f(x,y)+1$ for $x,y\in[\ell-1]$, and $g_f(x,y)=1$ otherwise.  For each $i\in [\ell-1]$, let $w_i=z_{\ell-i+1}$, $v_i=x_{i}$, and for each $f:[\ell-1]^2\rightarrow [\ell-1]$, let $u_i^f=x_i^{g_f}$.  Then for each $i,j,k\in [\ell]$, $\phi(u_i^f, v_j,w_k)=\phi(x^{g_f}_{i}, y_j, z_{\ell-k+2})$ holds if and only if $\ell-k+2< g_f(i,j)$, if and only if $\ell-k+1< \ell-f(i,j)+1$, if and only if $f(i,j)<  k$. Thus $\neg \phi(u^f_i,v_j,w_k)$ holds if and only if $k\leq f(i,j)$.  This shows $\neg \phi(x;y;z)$ has $(\ell-1)$-$\FOP_2$. 

We now prove (ii).  We being with the direction $(a)\Rightarrow (b)$. Enumerate $[\ell]^{[\ell]}=\{f_1,\ldots, f_m\}$. Define $f\colon [m]^2\to [m]$ such that $f(i,j)=f_i(j)$ if $j\leq \ell$ and $f(i,j)=\ell$ if $j>\ell$. Choose $a^f_1,\ldots,a^f_m,b_1,\ldots,b_m,c_1,\ldots,c_m$ so that $\phi(a^f_i,b_j,c_k)$ holds if and only if $k\leq f(i,j)$, i.e. if and only if $k\leq f_i(j)$.   
For each $t\in [m]$, set  $a_{f_t}=a^f_t$. Then for all $j,k\leq \ell$, we have $\phi(a_{f_t},b_{j},c_k)$  if and only  if and only if $k\leq f_t(i)$. 

We now show $(b)\Rightarrow (a)$. Fix $\ell\geq 1$. Let $b_1,\ldots,b_\ell,c_1,\ldots,c_\ell$ and $a_f$ for $f\in [\ell]^{[\ell]}$ be as in $(b)$. Fix $f\colon[\ell]^2\to[\ell]$ and for each $j\leq\ell$,  define $f_j\colon[\ell]\to[\ell]$ so that $f_j(i)=f(i,j)$. For each $i\in [\ell]$, set $a^f_i=a_{f_i}$. Then $\phi(a^f_i,b_j,c_k)$ holds if and only if $k\leq f_j(i)$, i.e., $k\leq f(i,j)$. 
\end{proof}

Proposition \ref{prop:basicprop} implies that being $\NFOP_2$ is closed under negations.  We now show that it is also closed under conjunctions.   We note that while the proof below is finitary, it passes through regularity and counting lemmas, which means that the quantitative relationship between the parameters is likely to be very poor.  It would be of interest to find a simpler finitary proof.  On the other hand, a more recent proof in the infinitary setting has appeared in \cite{Aldaim.2023}. This relies on the closure of stability under boolean combinations along with several Ramsey theoretic results.   

\begin{theorem}\label{thm:and}
For all $k\geq 1$ there is $t\geq 1$ so that the following holds.  Suppose $\calL$ is a first-order language, $\calM$ is an $\calL$-structure, and $\phi_1(x,y,z),\phi_2(x,y,z)$ are $\calL$-formulas.  Suppose $\phi_1(x;y;z)$ and $\phi_2(x;y;z)$ have no $k$-$\FOP_2$ in $\calM$.  Then $(\phi_1\wedge \phi_2)(x;y;z)$ has no $t$-$\FOP_2$. 
\end{theorem}
\begin{proof}
Suppose towards a contradiction that there exists $k\geq 1$ so that for all $t\geq 1$, there exists an $\calL$-structure $\calM(t)$ in which neither $\phi_1(x;y;z)$ nor $\phi_2(x;y;z)$ have $k$-$\FOP_2$, but in which $(\phi_1\wedge \phi_2)(x;z;y)$ has  $t$-$\FOP_2$.

By compactness, there exists an $\calL$-structure $\calM$ in which neither $\phi_1(x;y;z)$ nor $\phi_2(x;y;z)$ have $k$-$\FOP_2$, but where $\phi_1\wedge \phi_2(x;z;y)$ has $m$-$\FOP_2$ for all $m\geq 1$. 

Let $s=s(k)$ be such that the conjunction of two $k$-stable formulas is $s$-stable. Let $\e_1,\e_2$ be as in Proposition \ref{prop:sufffop} for $k$ and choose $S\gg \e_1^{-1}s$.  Choose $\e_1''\ll \e_1'\ll \e_1/S$, and define $\e_2'',\e'_2:\mathbb{N}\rightarrow (0,1]$ by defining, for each $x$, $\e_2''(x)\ll \e_2'(x)\ll \e_2(x)/kx$.   

Let $\mu'_1,\mu_1,\mu_2,\mu_2'$ be as in Lemma \ref{lem:intersecting} for $\e_1'',\e_2''$.  Choose $f,g:\mathbb{N}^4\rightarrow (0,1]$ as in Lemma \ref{lem:refinement} for $\mu_1'',\mu_2''$, where $\mu_1''=\min\{\mu_1,\mu_1'\}$ and $\mu_2''(x)=\min\{\mu_2(x), \mu_2'(x)\}$.  Define $\e_2''''(x)\ll \e_2'''(x)\ll \mu_2''(x)\e_2''(x)$ for all $x$.

Let $T_1, L_1, N_1$ be as in Proposition \ref{thm:FOPfinite} for $k$, $\mu''_1$, and $\mu_2''$, and set $T=f_1(T_1,L_1,T_1,L_1)$ and $L=g_1(T_1,L_1,T_1,L_1)$.  Choose $S\gg L$ and choose $N\gg S LT\mu_2''(LT)^{-1}(\mu_1'')^{-1} N_1N_2$.  Define $n=n(k,N)$ from Lemma \ref{lem:otherway}.  

Now let $U=\{u_1,\ldots, u_N\}$, $W=\{w_1,\ldots, w_N\}$ be two disjoint sets, each of size $N$, and let $\bigcup_{\alpha\leq S}P_{UW}^{\alpha}$ be a partition of $K_2[U,W]$ such that each $P_{UW}^{\alpha}$ has $\disc_2(\e''''_2(LT); 1/S)$ (such a partition exists by Lemma \ref{lem:3.8}).  

By assumption and Lemma \ref{lem:clean}, there are disjoint sets $A,B,C\subseteq M$, with 
$$
A=\{a_i^f: i\in [n], f:[n]^2\rightarrow [n]\}, B=\{b_i:i\in [n]\},\text{ and }C=\{c_i: i\in [n]\},
$$
so that $\calM\models (\phi_1\wedge \phi_2)(a^f_i,b_j, c_k)$ if and only if $k\leq f(i,j)$.  Let 
$$
E=\{abc\in K_3[A,B,C]: \calM\models (\phi_1\wedge \phi_2)(a,b,c)\},
$$
and for each $\alpha\in \{1,2\}$, let $E_{\alpha}=\{abc\in K_3[A,B,C]: \calM\models \phi_{\alpha}(a,b,c)\}$.  Note that $(A\cup B\cup C, E)$, $(A\cup B\cup C, E_0)$ and $(A\cup B\cup C, E_1)$ are each $3$-partite $3$-graphs.

By Lemma \ref{lem:otherway}, there are subsets $X=\{x_i: i\in [N]\}\subseteq A$, $Y=\{y_i: i\in [N]\}\subseteq B$, and $Z_1,\ldots, Z_{S_1}\subseteq C$ each of size $N$, so that if we define $Q_{XY}^{\alpha}:=\{x_iy_j: u_iw_j\in P_{UW}^{\alpha}\}$ for each $\alpha\in [S]$, then for all $x_iy_j\in Q^{\alpha}_{XY}$ and $z\in Z_{\beta}$, we have $x_iy_jz\in E$ if and only if $\beta\leq \alpha$.  

Let $V'=X\cup Y\cup Z_1\cup \ldots \cup Z_S$, and let $\calQ$ be the $(S+2, S)$-decomposition of $V'$ consisting of $\{X,Y,Z_1,\ldots, Z_S\}$ and $\{Q_{XY}^{\alpha}: \alpha,m\leq S\}$.  Define  $H_1=(V',E_1[V'])$, $H_2=(V',E_2[V'])$, and $H_{12}=(V',(E_1\wedge E_2)[V'])$.  Note that by assumption, $H_1$ and $H_2$ have no $k$-$\FOP_2$, and thus both have $\VC_2$-dimension at most $k$. 

By Theorem \ref{thm:FOPfinite}, there are integers $t_1\leq T_1$, $\ell_1\leq L_1$ and a  $(t_1,\ell_1,\mu''_1,\mu''_2(\ell_1))$-decomposition $\calP_1$ of $V'$ which is $\mu''_1$-homogeneous with respect to $H_1$ and which has linear error, say witnessed by $\Sigma_1$.  Similarly, by Theorem \ref{thm:FOPfinite}, there are $t_2\leq T_1$, $\ell_2\leq L_1$ and an $(\mu''_1,\mu''_2(\ell_2),t_2,\ell_2)$-decomposition  $\calP_2$ of $V'$ which is $\mu''_1$-homogeneous and which has linear error with respect $H_2$, say witnessed by $\Sigma_2$.  Let $\calP_4=\calP_1\wedge \calP_2$ be the coarsest decomposition of $V$ refining both $\calP_1$ and $\calP_2$.  Let $\calP_3$ be the decomposition of $V$ with $(\calP_3)_{vert}=\calQ_{vert}$ and $(\calP_3)_{edge}=\{K_2[F,F']: FF'\in {\calQ_{vert}\choose 2}\}$. 

By Lemma \ref{lem:refinement}, there exist integers $t_3\leq f_1(t_1t_2,\ell_1\ell_2,S+2,1)$ and $\ell_3\leq g_1(t_1t_2,\ell_1\ell_2,S+2,1)$, as well as a $(t_3,\ell_3,\mu''_1,\mu''_2(\ell_3))$-decomposition $\calR$ of $V'$ which is an $(\mu''_1,\mu_2''(\ell_3))$-approximate refinement of $\calP_3$ and $\calP_4$.  Say this is witnessed by $\Sigma\subseteq {\calR_{vert}\choose 2}$, as in Definition \ref{def:refinement}.

By Lemma \ref{lem:intersecting}, $\calR$ is $(\e_1'',\e_2''(\ell_3))$-regular and $\e''_1$-homogeneous with respect to both $E_1$ and $E_2$.  Since $\calR_{vert}$ refines $\calQ_{vert}$, it has the form 
$$
\calR_{vert}=\{X_i: i\in [t]\}\cup \{Y_i: i\in [t]\}\cup \{Z_{i,j}: i\in [S], j\in [t]\},
$$
for some $t\leq t_3$, and where each $X_i\subseteq X$, $Y_i\subseteq Y$ and $Z_{i,j}\subseteq Z_i$.   Let $\calR_X=\{X_i: i\in [t]\}$, $\calR_Y= \{Y_i: i\in [t]\}$, and $\calR_Z= \{Z_{i,j}: i\in [S], j\in [t]\}$.  For each $i,j\in [t]$ and $i'\in [S]$, we also have the partitions $K_2[X_i,Y_j]=\bigcup_{\alpha\leq \ell_3}R_{X_iY_j}^{\alpha} $, $K_2[X_i,Z_{i',j}]=\bigcup_{\alpha\leq \ell_3}R_{X_iZ_{i',j}}^{\alpha} $ and $K_2[Y_i,Z_{i',j}]=\bigcup_{\alpha\leq \ell_3}R_{Y_iZ_{i',j}}^{\alpha} $ coming from $\calR_{edge}$.  For each $1\leq i,j\leq t$, set $Q_{X_iY_j}^{\alpha}=Q_{XY}^{\alpha}[X_i,Y_j]$.  Note that by Lemma \ref{lem:sl}, each $Q_{X_iY_j}^{\alpha}$ satisfies $\disc_2(\e_2'''(S);1/S)$. 

For each $X_iY_j\notin \Sigma$ and $\tau\in [\ell_3]$, let $f_1(R_{X_iY_j}^{\tau})\in (\calP_1)_{edge}$ and $f_2(R_{X_iY_j}^{\tau})\in (\calP_2)_{edge}$ be such that $|R_{X_iY_j}^{\tau}\setminus (f_1(R_{X_iY_j}^{\tau})\cap f_2(R_{X_iY_j}^{\tau}))|\leq \mu_1''|R_{X_iY_j}^{\tau}|$, and so that $R_{X_iY_j}^{\tau}\setminus (f_1(R_{X_iY_j}^{\tau})\cap f_2(R_{X_iY_j}^{\tau}))$ has $\disc_2(\mu_2''(\ell_3))$.  Setting 
$$
C(X_i,Y_j)=\bigcup_{\alpha\in \ell_4} R_{X_iY_j}^{\tau}\setminus (f_1(R_{X_iY_j}^{\tau})\cap f_2(R_{X_iY_j}^{\tau})),
$$
we have, by definition of a $(\mu_1'',\mu_2''(\ell_4))$-approximate refinement, that for all $X_iY_j\notin \Sigma$, $|C(X_i,Y_j)|\leq \mu_1''|X_i||Y_j|$.  Define 
$$
J(X_i,Y_j)=\{Q_{XY}^{\alpha}\in \calQ_{edge}: |Q_{XY}^{\alpha}\cap C(X_i,Y_j)|\geq \sqrt{\mu_1''}|Q_{XY}^{\alpha}|\}.
$$
Then $|\{xy\in K_2[X_i,Y_j]: xy\in Q_{XY}^{\alpha}\text{ for some }Q_{XY}^{\alpha}\in J(X_i,Y_j)\}|\leq \sqrt{\mu_1''}|X_i||Y_j|$, and consequently, $|J(X_i,Y_j)|\leq 2\sqrt{\mu_1''}S$.  We now define
\begin{align*}
\mathbf{\Gamma}=\Big\{G\in \triads(\calR): \text{ for each }& u\in \{0,1\}, (H_u|G,G)\text{ has }\disc_{2,3}(\e_1'',\e_2''(\ell_3))\\
&\;\; \text{ and } \frac{|E_u\cap K_3^{(2)}(G)|}{|K_3^{(2)}(G)|}\in [0,\e_1'')\cup (1-\e_1'',1]\Big\},
\end{align*}
and set $\Gamma=\{xyz\in K_3[X,Y,Z]:xyz\in K_3^{(2)}(G) \text{ for some }G\in \mathbf{\Gamma}\}$.  By assumption, $|\Gamma|\geq (1-\e_1'')|V|^3$.  Therefore, setting 
$$
I=\Big\{xy\in K_2[X,Y]: |\{z\in Z: xyz\in \Gamma\}|\geq (1-\e_1')|Z|, xy\notin \bigcup_{X_iY_j\in \Sigma}C(X_i,Y_j)\Big\},
$$
we must have that $|I|\geq (1-\e_1')|V|^2$ (here we are using that $\mu_1''\ll \e_1''\ll \e_1'$).  Therefore, there must be some $X_i\in \calR_X$ and $Y_j\in\calR_Y$ so that $|I\cap K_2[X_i,Y_j]|\geq (1-\e_1')|X_i||Y_j|$.  Given $\alpha,\beta,\gamma\leq \ell_3$, $b\in [t]$, and $u\in [S]$, let 
$$
R_{X_iY_jZ_{u,b}}^{\alpha,\beta,\gamma}=(X_i\cup Y_j\cup Z_{u,b}, R_{X_iY_j}^{\alpha}\cup R_{X_iZ_{u,b}}^{\beta}\cup R_{Y_jZ_{u,b}}^{\gamma}).
$$
Let $\mathbf{\Gamma}'=\mathbf{\Gamma}\cap \{R_{X_iY_jZ_{u,b}}^{\alpha,\beta,\gamma}: \alpha,\beta,\gamma\leq \ell_3, u\in [S], b\in [t]\}$.  We now define a bipartite graph $\mathbf{G}=(U\cup W, F)$, where
\[U=\{R_{X_iY_j}^{\alpha}:\alpha\leq \ell_3\},\text{ }W=\{R_{X_iZ_{u,b}}^{\beta}R_{Y_jZ_{u,b}}^{\gamma}:\gamma,\beta\leq \ell_3b\in [t], u\in [S]\}\]
and
\[F= \{ R_{X_iY_j}^{\alpha}(R_{X_iZ_{u,b}}^{\beta}R_{Y_jZ_{u,b}}^{\gamma})\in K_2[U,W]: R_{X_iY_jZ_{u,b}}^{\alpha,\beta,\gamma}\in \mathbf{\Gamma}'\}.\]
By our choice of $X_iY_j$, 
$$
\Big|\Big(\bigcup_{G\in \mathbf{\Gamma}'}K_3^{(2)}(G)\Big)\cap K_3[X_i,Y_j,Z]\Big|\geq (1-\e_1')^2|X_i||Y_j||Z|.
$$
Consequently, using Corollary \ref{cor:counting},
$$
|F|\geq \frac{(1-\e_1')^2|X_i||Y_j||Z|}{(1+\e_1')|X_i||Y_j||Z|/St\ell_3^3}\geq (1-\e_1')^3 St\ell_3^3\geq (1-3\e_1')St\ell_3^3.
$$
Set
$$
\Omega=\{Z_\alpha: \alpha\in [S] \text{ and for some }c\in W,  |N_{\mathbf{G}}(c)|\geq (1-\sqrt{3\e'_1})\ell_3\}.
$$
Then $|F|\geq (1-3\e_1')St\ell_3^3$ implies that $|\Omega|\geq (1-\sqrt{3\e_1'})S$.  Now set
$$
\Omega'=\{Z_\alpha\in \Omega: Q_{X_iY_j}^\alpha\notin J(X_i,Y_j)\}.
$$
Then by above, $|\Omega'|\geq (1-\sqrt{3\e_1'}-2\sqrt{\mu_1''})\geq (1-10\sqrt{\e_1'})S$. Choose $Z_{i_1},\ldots, Z_{i_{s}}\in \Omega'$ with $i_1<\ldots<i_s$, and for each $1\leq v\leq s$, choose any $c_v=R_{X_iZ_{i_v,b_v}}^{\beta_v}R_{Y_jZ_{i_v,b_v}}^{\gamma_v}\in Z_{i_v}$ so that $|N_{\mathbf{G}}(c)|\geq (1-\sqrt{3\e_1'})\ell_3$. Consider
$$
O=\bigcap_{v=1}^sN_{\mathbf{G}}(c_v).
$$
By construction, $|O|\geq (1-s\sqrt{3\e'_1})\ell_3$.  Since $\e_1'\ll \e/S$ and $S\gg s$, $|O|\geq (1-\e_1)\ell_3$.  For each $v\in [s]$, consider 
$$
L_v=\{xy: \text{there is } P\in O\text{ with } xy\in P\}\cap Q_{XY}^{i_v}[X_i,Y_j].
$$
By definition of $\Omega'$, and because $\e_1'\ll \e_1/S$, 
\begin{align*}
|L_v|\geq |Q_{XY}^{i_v}[X_i,Y_j]|-|Q_{XY}^{i_v}[X_i,Y_j]\cap C(X_i,Y_j)|&\geq  |Q_{XY}^{i_v}[X_i,Y_j]|-\sqrt{\e_1'}|X_i||Y_j|\\
&\geq \Big(\frac{1}{S}-2\sqrt{\e_1'}\Big)|X_i||Y_j|\\
&\geq (1-\e_1)|Q_{XY}^{i_v}[X_i,Y_j]|.
\end{align*}
Clearly this implies there is some $R_{X_iY_j}^{\tau_v}\in O$ with $|R_{X_iY_j}^{\tau_v}\setminus Q_{XY}^{i_v}|\leq \e_1 |R_{X_iY_j}^{\tau_v}|$.  It is not difficult to see now that for each $u\leq v$, 
$$
|(E_1\cap E_2)\cap K^{(2)}_3(R_{X_iY_jZ_{i_u}}^{\tau_v\alpha_u\beta_u})|\geq (1-\e_1)|K^{(2)}_3(R_{X_iY_jZ_{i_u}}^{\tau_v\alpha_u\beta_u})|,
$$
while if $v<u$, then 
$$
|(E_1\cap E_2)\cap K^{(2)}_3(R_{X_iY_jZ_{i_u}}^{\tau_v\alpha_u\beta_u})|< \e_1|K^{(2)}_3(R_{X_iY_jZ_{i_u}}^{\tau_v\alpha_u\beta_u})|.
$$
Since each of these triads $R_{X_iY_jZ_{i_u}}^{\tau_v\alpha_u\beta_u}$ are in $\mathbf{\Gamma}$, we have that for all $u\leq v$, 
$$
\min\{|E_1\cap K^{(2)}_3(R_{X_iY_jZ_{i_u}}^{\tau_v\alpha_u\beta_u})|, |E_2\cap K^{(2)}_3(R_{X_iY_jZ_{i_u}}^{\tau_v\alpha_u\beta_u})|\}\geq (1-\mu_1'')|K^{(2)}_3(R_{X_iY_jZ_{i_u}}^{\tau_v\alpha_u\beta_u})|,
$$
and for all $v<u$, there $w_{u,v}\in \{0,1\}$ so that 
$$
|E_{w_{u,v}}\cap K^{(2)}_3(R_{X_iY_jZ_{i_u}}^{\tau_v\alpha_u\beta_u})|\leq \mu_1''|K^{(2)}_3(R_{X_iY_jZ_{i_u}}^{\tau_v\alpha_u\beta_u})|.
$$
By our choice of $s$, there is a subsequence $1\leq i_1'<\ldots<i_k'\leq s$ and $w\in \{0,1\}$ so that so $E_{w_{i_u',i_v'}}=E_w$ for all $v<u$.  Let $f:\{a_1,\ldots, a_k\}\rightarrow \calR_{edge}$ and $g:\{b_1,\ldots, b_k\}\rightarrow \calR_{cnr}$ be defined by $f(a_v)=R_{X_{i}Y_{i}}^{\tau_{i_v'}}$ and $g(b_u)=R_{X_iZ_{i_u'}}^{\alpha_{i_u}'}R_{Y_jZ_{i_u'}}^{\beta_{i_u'}}$, respectively. Then $(f,g)$ is an encoding of $H(k)$ in $(H_w,\calP)$ with corresponding partition $(X_i,Y_j,\bigcup_{u=1}^kZ_{i_u'})$.  By Proposition \ref{prop:sufffop}, there are $(c_1,\ldots, c_k)\in Z_{i_1'}\times \ldots \times Z_{i_k'}$, $\{a_1^h,\ldots, a_k^h | h:[k]^2\rightarrow [k]\}\subseteq X_i$, and $\{b_1,\ldots, b_k\}\subseteq Y_j$ such that for all $h:[k]^2\rightarrow [k]$, $a_i^hb_jc_k\in E_w$ if and only if $k\leq h(i,j)$.  By definition of $E_w$, this implies that $\phi_w(x;y;z)$ has $k$-$\FOP_2$ in $\calM$, a contradiction.
\end{proof}

Combining  Proposition \ref{prop:basicprop}, Theorem \ref{thm:and}, and the fact that $\wedge$, $\neg$ form a complete system of connectives, we have thus shown that being $\NFOP_2$ is closed under finite boolean combinations. 

\begin{corollary}
In a first order structure $\calM$, the set of $\NFOP_2$ formulas is closed under boolean combinations.
\end{corollary}

We end by constructing for all $k\geq 1$, a $3$-graph which has $k$-$\FOP_2$  and $\VC_2$-dimension $1$.

\begin{example}\label{exa:vcfop}
Choose $\e_1>0$ and $\e_2:\mathbb{N}\rightarrow (0,1]$ as in Proposition \ref{prop:sufffop} for $k$, and let $N$ be as in Proposition \ref{prop:sufffop} for $\e_1,\e_2,\ell=k$ and $t=k+2$, and let $n\gg N$.  Let $U=\{u_1,\ldots, u_n\}$, $W=\{w_1,\ldots, w_n\}$ and let $K_2[U, W]=\bigcup_{\alpha\leq N}P_{UW}^{\alpha}$ be a partition so that each $P_{UW}^{\alpha}$ has $\disc_2(\e_2(k); 1/k)$ (such a partition exists by Lemma \ref{lem:3.8}).  Let $Z_1,\ldots, Z_k$ be disjoint sets, each of size $n$, and set $Z=\bigcup_{i=1}^kZ_i$.   We now let $V=U\cup W\cup Z$ and define $H=(V,E)$ to be the $3$-partite $3$-graph where $E=\{uwz: z\in Z_i, uw\in P_{UW}^{\alpha}\text{ for some }\alpha\leq i\}$.  By Proposition \ref{prop:sufffop}, $H$ has $k$-$\FOP_2$.  We show it has $\VC_2$-dimension at most $1$.

Fix any $x_1,x_2,y_1,y_2\in V$, and let $X=\{x_1,x_2\}$ and $Y=\{y_1,y_2\}$ be two $2$-element subsets.  We claim that the edge relation in $H$ cannot shatter $X\times Y$.  Suppose towards a contradiction this was not the case.  Then there are, for each $S\subseteq X\times Y$, some $z_S\in V$ so that $z_Sx_iy_j\in E$ if and only if $(i,j)\in S$.  Clearly this implies $x_1\neq x_2$ and $y_1\neq y_2$.  Further, since $z_{X\times Y}x_iy_j\in E$ for each $i,j\in [2]$, and since $H$ is $3$-partite, we must have that there is a relabeling $\{A,B,C\}=\{U,W,Z\}$ so that $x_1,x_2\in A$,  $y_1,y_2\in B$, and $z_S\in C$ for all $S\neq \emptyset$.  

Suppose first that $A=Z$.  Note that for all $i\in [k]$ and $u\in U$ and $w\in W$, either $Z_i\subseteq N_H(uw)$ or $Z_i\cap N_H(uw)=\emptyset$.  Thus we clearly cannot have $x_1,x_2$ in the same $Z_i$.  Therefore, without loss of generality, we may assume that for some $1\leq i<j\leq k$, $x_1\in Z_i$ and $x_2\in Z_j$.  But now there exists no $uw\in K_2[U,W]$ with $x_2uw\in E$ and $x_1uw\notin E$, a contradiction.  Thus, $A\neq Z$.  A symmetric argument shows $B\neq Z$. 

Thus $\{A,B\}=\{U,W\}$.  Without loss of generality, assume $A=U$ and $B=W$.  For each $(i,j)\in [2]^2$, let $\alpha_{ij}$ be such that $x_iy_j\in P_{UW}^{\alpha_{ij}}$.  Clearly we must have $\alpha_{11}\neq \alpha_{12}$, since $z_{\{(1,1)\}}x_1y_1\in E$ but $z_{\{(1,2)\}}x_1y_1\notin E$.  Assume $\alpha_{11}<\alpha_{12}$ (the other case is similar).  But now there exists no $z\in Z$ with $zx_1y_2\in E$, and $ax_1y_1\notin E$.  Again this is a contradiction, and finishes the proof.
\end{example}